\documentclass[a4paper,12pt]{article}
\usepackage{amssymb}
\usepackage{amsthm}
\usepackage{amsmath}
\usepackage[mathscr]{eucal}
\newtheorem{thm}{Theorem}[section]
\newtheorem{lem}[thm]{Lemma}
\newtheorem{cor}[thm]{Corollary}
\newtheorem{prop}[thm]{Proposition}
\newtheorem{rem}[thm]{Remark}
\newtheorem{exmp}[thm]{Example}
\newtheorem{defn}[thm]{Definition}

\title{Centralizers of reflections and reflection-independence of Coxeter groups}
\author{Koji Nuida}
\date{}
\begin{document}
\maketitle
\begin{abstract}
A Coxeter group $W$ is called reflection independent if its reflections are uniquely determined by $W$ only, independently on the choice of the generating set.
We give a new sufficient condition for the reflection independence, and examine this condition for Coxeter groups in certain classes, possibly of infinite ranks.
We also determine the finite irreducible components of another Coxeter group, that is a subgroup of $W$ generated by the reflections centralizing a given generator of $W$.
Determining such a subgroup makes our criterion efficient.
\end{abstract}

\section{Introduction}

Let $(W,S)$ be a Coxeter system, namely $W$ is a Coxeter group with corresponding generating set $S$.
In this paper, we do not place any restriction on $W$ or $S$ unless otherwise noticed; for example, the cardinality of $W$ or of $S$ may be infinite.
An element of $W$ conjugate to some element of $S$ is called a \emph{reflection} in $W$.
This definition of reflections seems to depend on the choice of the generating set $S$; and it actually does, as shown by a well-known example provided by the symmetric group $S_6$ of degree $6$.
Following Patrick Bahls \cite{Bah}, a Coxeter group $W$ is called \emph{reflection independent} if the notion of reflections in $W$ is independent of the generating set $S$.
The aim of the present paper is to investigate the notion of reflections in $W$, obtaining some sufficient conditions for reflection independence and revealing Coxeter groups in certain classes to be reflection independent.

As the name suggests, reflections in Coxeter groups are of not only group-theoretical, but also geometrical importance.
Any Coxeter group $W$ is realized as a reflection group in a vector space $V$ equipped with a (not necessarily nondegenerate) symmetric bilinear form, where an element of $W$ is a reflection in the above sense if and only if it is geometrically a reflection in $V$.
However, this realization does depend on the set of reflections in $W$ in the above sense.
This phenomenon would exhibit the importance of reflection independence of Coxeter groups.
Another aspect exhibiting the importance is a relationship with the isomorphism problem of Coxeter groups.
This is the problem of deciding which Coxeter groups are isomorphic as abstract groups; the Weyl group of type $G_2$, which is isomorphic to another Weyl group $S_2 \times S_3$, implies that the problem is not so obvious.
It is conjectured in {\cite[Conjecture 8.1]{Bra-McC-Muh-Neu}} that any two generating sets of a Coxeter group would be related via ``diagram twistings'' if these sets define the same set of reflections.
Assuming this conjecture, the variation of generating sets of a Coxeter group will be apparent if it is reflection independent.

Here we introduce some preceding results on this topic.
Although some stronger results on certain Coxeter groups yielding their reflection independence had been given (e.g.\ {\cite[Main theorem]{Cha-Dav}}), the notion of reflection independence itself had not attracted attention until recently.
The main strategy of preceding studies is to examine maximal finite subgroups of given Coxeter groups; the key property is given in {\cite[Lemma 1.6]{Fra}}.
Reflection independence of Coxeter groups in some classes have been investigated along this line; see e.g.\ \cite{Bah,Fra-How-Muh,Muh-Wei}.
However, the strategy based on maximal finite subgroups would not work in a study of general Coxeter groups of infinite ranks, since it may happen that the group possesses no maximal finite subgroups.
On the other hand, the strategy of the present paper is different and basically applicable to an arbitrary case, including the case of infinite ranks.
For example, we prove in Theorem \ref{thm:refindep_2spherical} that, if a Coxeter group $W$ is infinite, irreducible, and the product of any pair of two generators has finite order (i.e.\ ``$2$-spherical''), then $W$ is reflection independent even if it has infinite rank.

A key ingredient of our strategy established below is the subgroup $W^{\perp x}$ of $W$ associated to each $x \in S$.
The subgroup $W^{\perp x}$ is generated by the reflections in $W$ other than $x$ itself which commute with $x$.
This is a Coxeter group, due to a result of Vinay V.\ Deodhar \cite{Deo_refsub} or Matthew Dyer \cite{Dye}.
Let $(W^{\perp x})_{\mathrm{fin}}$ denote the product of all finite irreducible components of $W$ (Definition \ref{defn:finitepart}).
Roughly speaking, our sufficient condition says that $W$ is reflection independent whenever the group $(W^{\perp x})_{\mathrm{fin}}$ is sufficiently small for every $x \in S$ (Theorem \ref{thm:Nui_ext_scforrefindep}).
To make the strategy efficient, we determine the structure of the subgroups $(W^{\perp x})_{\mathrm{fin}}$ completely, as summarized in Section \ref{sec:finitepart}.
Note that the group $(W^{\perp x})_{\mathrm{fin}}$ is a part of the centralizer $Z_W(x)$ of $x$, so our description of $(W^{\perp x})_{\mathrm{fin}}$ is based on the results on the centralizers given by Brigitte Brink \cite{Bri} and by the author \cite{Nui_centra}.

The present paper is organized as follows.
In Section \ref{sec:preliminaries}, we summarize basics of Coxeter groups and some further definitions and properties, and fix some notations and terminology for graphs, posets and groupoids.
Section \ref{sec:Zofreflections} introduces the preceding results on the centralizers mentioned above, and prepares some further observations.
In Section \ref{sec:mainresult_specialcase}, we study the subgroups $(W^{\perp x})_{\mathrm{fin}}$ for some special cases.
These results are not only parts of the results for general cases, but also crucial ingredients of the following arguments proceeded inductively.
Section \ref{sec:finitepart} determines the structure of $(W^{\perp x})_{\mathrm{fin}}$ completely; this description is divided into six special cases (Theorems \ref{thm:finpart_O_typeB}, \ref{thm:finpart_O_trunk_acyclic}, \ref{thm:finpart_O_trunk_cyclic}, \ref{thm:finpart_nonO_Otrunk_F4}, \ref{thm:finpart_nonO_Otrunk_B}, \ref{thm:finpart_nonO_Onontrunk}) and a `generic' case (Theorem \ref{thm:finpart_generalcase}).
Finally, in Section \ref{sec:resultsonrefindep}, we give a new sufficient condition for reflection independence mentioned above, and verify the condition for Coxeter groups in certain classes (including ``$2$-spherical'' cases), owing to the description of $(W^{\perp x})_{\mathrm{fin}}$ given in the preceding sections.

\noindent
\textbf{Acknowledgement.}
The author would like to express his deep gratitude to Professor Itaru Terada and Professor Kazuhiko Koike for their precious advice and encouragement.
Moreover, the author was supported by JSPS Research Fellowship throughout this research.

\section{Preliminaries}
\label{sec:preliminaries}

In this paper, we often denote a set $\{x\}$ with one element also by $x$ for simplicity, unless some ambiguity occurs.
\subsection{On graphs, posets and groupoids}
\label{sec:graphsandgroupoids}

In this subsection, we briefly summarize the notion of groupoids, and fix some notations and convention for graphs, groupoids and posets.

In this paper, a graph $G$ signifies an edge-unoriented graph with no loops, where it may possess edge-labels or parallel edges.
The vertex set of $G$ is denoted by $V(G)$.
Let $\vec{E}(G)$ denote the set of the edges of $G$ with some orientation, and for $x,y \in V(G)$, let $\vec{E}_{x,y}(G)$ be the set of the edges of $G$ \emph{from $y$ to $x$} (this convention, chosen here for technical reason, is reverse what is usual).
For $x \in V(G)$, let $G_{\sim x}$ denote the connected component of $G$ containing $x$.
If $G$ has no parallel edges, let $(x_n,\dots,x_1,x_0)$ denote the unique path in $G$ which starts from a vertex $x_0$, visits vertices $x_1,\dots,x_{n-1}$ in this order and ends at a vertex $x_n$.
For example, the (closed) path of length $0$ starting from $x$ is denoted by $(x)$.

A \emph{groupoid} is a small category whose morphisms are all invertible; namely a family of sets $\mathcal{G}=\{\mathcal{G}_{x,y}\}_{x,y \in V(\mathcal{G})}$ with some index set (or \emph{vertex set}) $V(\mathcal{G})$ endowed with (1) multiplications $\mathcal{G}_{x,y} \times \mathcal{G}_{y,z} \to \mathcal{G}_{x,z}$ satisfying the associativity law, (2) an identity element $1=1_x$ in each $\mathcal{G}_{x,x}$, and (3) an inverse $g^{-1} \in \mathcal{G}_{y,x}$ for every $g \in \mathcal{G}_{x,y}$.
Each $\mathcal{G}_{x,x}$ is a group, called a \emph{vertex group} of $\mathcal{G}$ and denoted by $\mathcal{G}_x$.
We write $g \in \mathcal{G}$ to signify $g \in \bigsqcup_{x,y \in V(\mathcal{G})}\mathcal{G}_{x,y}$, and $X \subseteq \mathcal{G}$ for $X \subseteq \bigsqcup_{x,y \in V(\mathcal{G})}\mathcal{G}_{x,y}$.
A \emph{homomorphism} between groupoids is a covariant functor between them regarded as categories.
Notions such as \emph{isomorphisms} and \emph{subgroupoids} are defined as usual.

The \emph{fundamental groupoid} of a graph $G$, denoted here by $\pi_1(G;*,*)$, is a groupoid with vertex set $V(G)$.
The elements of $\vec{E}_{x,y}(G)$ define elements of $\pi_1(G;x,y)=\pi_1(G;*,*)_{x,y}$, and $\pi_1(G;*,*)$ is freely generated by these elements.
Note that the same edge endowed with two opposite orientations give inverse elements to each other.
A vertex group $\pi_1(G;x)=\pi_1(G;x,x)$ is the \emph{fundamental group} of $G$ at $x$, which is a free group.
For technical reasons, \emph{a path in a graph is written from right to left in this paper}, hence a path $e_n \cdots e_2e_1$, where $e_i \in \vec{E}_{x_i,x_{i-1}}$ for $1 \leq i \leq n$, determines an element of $\pi_1(G;x_n,x_0)$.

Let $P$ be a poset (partially ordered set), with partial order denoted by $\preceq$.
Here the notation $x \preceq y$ includes the case $x=y$; we write $x \prec y$ to signify ``$x \preceq y$ and $x \neq y$''.
Two element $x,y \in P$ are called \emph{comparable} if $x \preceq y$ or $y \preceq x$, and \emph{incomparable} otherwise.
An element $x \in P$ \emph{covers} $y \in P$ if $y \prec x$ and there is no element $z \in P$ with $y \prec z \prec x$.
If $P$ has a unique minimal element $0_P$, then the elements of $P$ which cover $0_P$ are called \emph{atoms} of $P$.
A chain $y_0 \prec y_1 \prec \cdots \prec y_n$ in $P$ from $y_0 \in P$ to $y_n \in P$ is called \emph{saturated} if $y_i$ covers $y_{i-1}$ for $1 \leq i \leq n$.
A subset $Q \subseteq P$ is an \emph{order ideal} of $P$ if $y \preceq x \in Q$ yields $y \in Q$.
For $x \in P$, let $\wedge^x=\{y \in P \mid y \preceq x\}$ denote the principal order ideal generated by $x$.
Dually, a subset $Q \subseteq P$ is a \emph{filter} of $P$ if $Q \ni x \preceq y$ yields $y \in Q$; let $\vee_x=\{y \in P \mid x \preceq y\}$ denote the principal filter generated by $x \in P$.
A unique maximal element of the intersection $\wedge^x \cap \wedge^y$ is called (if exists) the \emph{meet} of $x,y$ and denoted by $x \wedge y$.
If $x \wedge y$ exists for all $x,y \in P$, then $P$ is called a \emph{meet-semilattice}.

An example of posets is the \emph{tree order} on a tree $T$.
This order depends on a specified vertex (\emph{root vertex}) $x_0$ of $T$.
We define $x \preceq y$ in $T$ if the unique non-backtracking path in $T$ from $x_0$ to $y$ contains $x$.
Since $T$ is a tree, this is a partial order and makes $T$ a meet-semilattice with unique minimal element $x_0$.
For $y \in T$, the subset $\wedge^y$ is a finite saturated chain from $x_0$ to $y$, which is the non-backtracking path from $x_0$ to $y$.
Thus each $\wedge^x \cap \wedge^y$ is also a finite chain, yielding the existence of $x \wedge y$.
Two distinct vertices of the tree $T$ are adjacent if and only if one of them covers the other in the poset $T$; in particular, the atoms of the poset $T$ are the neighbors of $x_0$ in the tree $T$.
\begin{rem}
\label{rem:decompintotrees}
For any connected graph $G$, let $K$ be a connected, nonempty subset of $V(G)$ which contains the vertices of all simple closed paths in $G$.
Then it is straightforward to show that the graph $G$ is the union of its full subgraph $G|_K$ with vertex set $K$ and mutually disjoint trees $T_x$ indexed by $x \in K$, with $T_x \cap K=\{x\}$.
(A vertex $y \in V(G) \smallsetminus K$ belongs to $T_x$ where $x$ is the unique vertex in $K$ nearest from $y$.)
In this case, each $T_x$ is equipped with the tree order with root vertex $x$, denoted by $\preceq_x$.
Decompositions of graphs of this type will play an important role in later sections.
\end{rem}
The following graph-theoretic lemma will be needed later.
\begin{lem}
\label{lem:intersectionoftree}
Let $T$ be a tree, containing finitely many nonempty subtree $T_i$ with $1 \leq i \leq n$.
If any two subtrees $T_i$ and $T_j$ intersect with each other, then the intersection of all the subtrees $T_i$ is also a nonempty subtree of $T$.
\end{lem}
\begin{proof}
First, note that the intersection $T'=\bigcap_{i=1}^{n}T_i$ is also connected.
Indeed, if $x,y \in T'$, then every tree $T_i$ contains $x$ and $y$, so it also contains the unique path $P$ between $x$ and $y$ in $T$, therefore $P \subseteq T'$.

Thus it suffices to show that $T' \neq \emptyset$.
We proceed the proof by induction on $n$, the case $n=1$ being obvious.
The result in the previous paragraph proves the case $n=2$.
Moreover, once the case $n=3$ is proved, the claim for the case $n \geq 4$ will follow inductively.
Indeed, the induction assumption implies that the three subtrees $T_1$, $T_2$ and $\bigcap_{i=3}^{n}T_i$ satisfy the hypothesis, so the claim follows.

So suppose that $n=3$.
We denote by $\preceq$ the tree order on $T$ with fixed root vertex $x_0$ (see above for terminology).
Let $y_{ij}$ be a minimal element of $T'_{ij}=T_i \cap T_j$ with respect to $\preceq$.
Now observe that, for any $z,z' \in T$, the unique path $P(z,z')$ in $T$ between $z$ and $z'$ is the union of the two paths $P(z,z \wedge z')$ and $P(z',z \wedge z')$.

We show that either $y_{23} \preceq y_{12}$ or $y_{23} \preceq y_{13}$.
If $y_{23} \not\preceq y_{12}$ and $y_{23} \not\preceq y_{13}$, then we have $y_{23} \wedge y_{1i} \prec y_{23}$ for any $i \in \{2,3\}$, so $z \in P(y_{23},y_{1i})$ where $z$ is the unique element of the poset $T$ which is covered by $y_{23}$ (note that $y_{23} \wedge y_{1i} \preceq z$).
Now $P(y_{23},y_{1i}) \subseteq T_i$ since $T_i$ is connected and $y_{23},y_{1i} \in T_i$, so $z \in T_2 \cap T_3=T'_{23}$, contradicting the minimality of $y_{23}$.

Thus we have shown that $y_{23} \preceq y_{12}$ or $y_{23} \preceq y_{13}$.
Assume by symmetry that $y_{23} \preceq y_{12}$.
Then the claim follows if $y_{23}=y_{12}$ (namely $y_{23}=y_{12} \in T'$); so suppose that $y_{23} \prec y_{12}$.
Let $z$ be the unique element of $T$ which is covered by $y_{12}$.
Then we have $y_{23} \preceq z \prec y_{12}$, so $z \in P(y_{12},y_{23}) \subseteq T_2$.
Thus the minimality of $y_{12}$ implies that $z \not\in T_1$.
Now if $y_{12} \wedge y_{13} \prec y_{12}$, then $y_{12} \wedge y_{13} \preceq z \prec y_{12}$ and so $z \in P(y_{12},y_{13}) \subseteq T_1$, a contradiction.
This implies that $y_{12} \wedge y_{13}=y_{12}$, so $y_{23} \prec y_{12} \preceq y_{13}$, therefore $y_{12} \in P(y_{23},y_{13}) \subseteq T_3$.
Hence we have $y_{12} \in T'$, concluding the proof.
\end{proof}
\subsection{Coxeter groups}
\label{sec:Coxetergroups}

The basic definitions and facts summarized in this subsection are found in \cite{Hum} unless otherwise specified.

A pair $(W,S)$ of a group $W$ and its generating set $S$, \emph{possibly of infinite cardinality}, is called a \emph{Coxeter system} if $W$ admits the following presentation
\[
W=\langle S \mid (st)^{m_{s,t}}=1 \textrm{ for all } s,t \in S \textrm{ such that } m_{s,t}<\infty \rangle,
\]
where the $m_{s,t} \in \{1,2,\dots\} \sqcup \{\infty\}$ are symmetric in $s,t \in S$, and $m_{s,t}=1$ if and only if $s=t$.
A group $W$ is a \emph{Coxeter group} if some $S \subseteq W$ makes $(W,S)$ a Coxeter system.
The cardinality of $S$ is called the \emph{rank} of $(W,S)$, or of $W$ if the set $S$ is obvious from the context.
An \emph{isomorphism of Coxeter systems} $f:(W,S) \to (W',S')$ is a group isomorphism $f:W \to W'$ which maps $S$ onto $S'$.
It is well known that $m(s,t)$ is precisely the order of $st$ in $W$; so every generator $s \in S$ has order $2$.
Thus the pair $(W,S)$ determines uniquely the data $m(s,t)$, so the Coxeter systems $(W,S)$ are in one-to-one correspondence with the \emph{Coxeter graphs} $\Gamma$ (up to isomorphism); $\Gamma$ is a graph with vertex set $S$, and $s,t \in S$ are joined by a unique edge with label $m_{s,t}$ if and only if $m_{s,t} \geq 3$ (by convention, the label `3' is omitted when drawing a picture).
On the other hand, it is \emph{not} true that two Coxeter systems $(W,S)$ and $(W',S')$ have isomorphic Coxeter graphs whenever $W \simeq W'$ as abstract groups.
An aim of this paper is to study some aspects of the correspondence between Coxeter groups and Coxeter graphs.

For $w \in W$, let $\ell(w)=\ell_S(w)$ denote the \emph{length} of $w$ with respect to $S$, that is the minimal number $r$ such that $w=s_1s_2 \cdots s_r$ for some $s_i \in S$.
We have $\ell(1)=0$, $\ell(w^{-1})=\ell(w)$ and $\ell(ws)=\ell(w) \pm 1$ for $s \in S$.
Such an expression $w=s_1s_2 \cdots s_r$ with $r=\ell(w)$ is called a \emph{reduced expression} of $w$ (with respect to $S$).
\begin{defn}
\label{defn:supportofw}
Define the \emph{support} $\mathrm{supp}(w)$ of $w \in W$ by
\[
\mathrm{supp}(w)=\{s_1,\dots,s_r\} \subseteq S \textrm{ if } w=s_1 \cdots s_r \textrm{ is a reduced expression.}
\]
\end{defn}
It is well known that this definition is independent of the reduced expression.
This yields the following fact:
\begin{lem}
\label{lem:supportofproduct}
Let $w_1,\dots,w_n \in W$ such that $\ell(w_1 \cdots w_n)=\sum_{i=1}^{n}\ell(w_i)$.
Then we have $\mathrm{supp}(w_1 \cdots w_n)=\bigcup_{i=1}^{n}\mathrm{supp}(w_i)$.
\end{lem}

For a subset $I \subseteq S$, let $W_I=\langle I \rangle$ denote the \emph{standard parabolic subgroup} of $W$ generated by $I$.
A subgroup of $W$ conjugate to some $W_I$ is called a \emph{parabolic subgroup}.
(In some context, the term ``parabolic subgroup'' signifies a subgroup of the form $W_I$ only.)
Then $(W_I,I)$ is also a Coxeter system; its length function $\ell_I$ is the restriction of $\ell$ to $W_I$, and its Coxeter graph $\Gamma_I$ is the full subgraph of $\Gamma$ with vertex set $I$.
It is well known that $W_I \cap W_J=W_{I \cap J}$ for $I,J \subseteq S$; for $w \in W$, the subgroup $W_{\mathrm{supp}(w)}$ is the unique minimal standard parabolic subgroup of $W$ containing $w$.

The following theorem of Jacques Tits is useful.
Its proof is found in several books and papers; see e.g.\ {\cite[Lemma 1.2]{Bah_book}}.
Note that the assumption on the finiteness of ranks in the original version is not necessary, since the finite subgroup in the statement is always contained in a standard parabolic subgroup of finite rank.
\begin{thm}
[Tits]
\label{thm:Tits_finitesubgroup}
Any finite subgroup of a Coxeter group $W$ is contained in a finite parabolic subgroup.
\end{thm}

If a subset $I$ of $S$ is (the vertex set of) a connected component of the Coxeter graph $\Gamma$, then $W_I$ is called an \emph{irreducible component} of $W$ (with respect to $S$) or of $(W,S)$; in this case, the subset $I$ is referred to an \emph{irreducible component} of $S$.
Then $W$ is the (restricted) direct product of all irreducible components.
If $\Gamma$ is connected, then $W$, $S$ or $(W,S)$ is called \emph{irreducible}.
\begin{defn}
\label{defn:finitepart}
Let $W_{\mathrm{fin}}$ be the product of all finite irreducible components of $W$, called the \emph{finite part} of $W$.
\end{defn}
Although the decomposition of $W$ into irreducible components depends on the set $S$, the \emph{subset} $W_{\mathrm{fin}}$ of $W$ (not only its group structure) is uniquely determined by $W$ only, independent of $S$ (see {\cite[Theorem 3.4]{Nui_isom}} or {\cite[Example 3.3]{Nui_ext}}).
The finite part of a Coxeter group will play a central role in our argument below.

For a Coxeter system $(W,S)$, let $V$ denote the \emph{geometric representation space} (over $\mathbb{R}$) endowed with basis $\Pi=\{\alpha_s \mid s \in S\}$ and symmetric bilinear form $\langle \,,\, \rangle$ determined by
\[
\langle \alpha_s, \alpha_t \rangle=
\begin{cases}
-\cos(\pi/m(s,t)) & \textrm{if } m(s,t)<\infty;\\
-1 & \textrm{if } m(s,t)=\infty.
\end{cases}
\]
Then $\langle \,,\, \rangle$ is $W$-invariant, where $W$ acts on $V$ by $s \cdot v=v-2\langle \alpha_s, v \rangle \alpha_s$ for $s \in S$ and $v \in V$.
Let $\Phi=W \cdot \Pi$, $\Phi^+=\Phi \cap \mathbb{R}_{\geq 0}\Pi$ and $\Phi^-=-\Phi^+$ denote the \emph{root system}, the set of \emph{positive roots} and the set of \emph{negative roots}, respectively.
Each element of $\Phi$ (a \emph{root}) is a unit vector with respect to $\langle \,,\, \rangle$, and $\Phi=\Phi^+ \sqcup \Phi^-$.
The set $\Pi$ is called a \emph{simple system}, and its elements are called \emph{simple roots}.
For any subset $\Psi \subseteq \Phi$ and $w \in W$, write
\[
\Psi^+=\Psi \cap \Phi^+,\ \Psi^-=\Psi \cap \Phi^- \textrm{ and } \Psi\left[w\right]=\{\gamma \in \Psi^+ \mid w \cdot \gamma \in \Phi^-\}.
\]
Then the length $\ell(w)$ of $w \in W$ coincides with the cardinality of $\Phi\left[w\right]$.

For a vector $v=\sum_{s \in S}c_s\alpha_s \in V$, put
\begin{equation}
\label{eq:notationforcoefficient}
\left[\alpha_s\right]v=c_s \textrm{ and } \mathrm{supp}(v)=\{s \in S \mid c_s \neq 0\} \textrm{ (the \emph{support} of $v$)}.
\end{equation}
For a subset $I \subseteq S$, let $V_I$ denote the subspace of $V$ spanned by the set $\Pi_I=\{\alpha_s \mid s \in I\}$, and put $\Phi_I=\Phi \cap V_I$.
Then it is known that
\[
\Phi_I=W_I \cdot \Pi_I,
\]
the root system of $(W_I,I)$ (see e.g.\ {\cite[Lemma 4]{Fra-How_nearly}}).
This implies that the support $\mathrm{supp}(\gamma)$ of any root $\gamma$ is irreducible (see e.g.\ {\cite[Section 2.3]{Nui_isom}}).

For a root $\gamma=w \cdot \alpha_s \in \Phi$, let $s_\gamma$ denote the \emph{reflection} $wsw^{-1} \in W$ along the root $\gamma$ acting on $V$ by $s_\gamma \cdot v=v-2\langle \gamma, v \rangle \gamma$ for $v \in V$.
Note that $s_{\alpha_s}=s$ for all $s \in S$.
It is well known that, for $\beta,\gamma \in \Phi$, we have $s_\beta=s_\gamma$ if and only if $\beta=\pm \gamma$, and $s_\beta s_\gamma=s_\gamma s_\beta$ if and only if either $\beta=\pm \gamma$ or $\beta$ is orthogonal to $\gamma$.
Write
\[
S^W=\bigcup_{w \in W}wSw^{-1},
\]
the set of reflections in $W$ with respect to $S$.
It is known that
\[
S^W \cap W_I=I^{W_I} \textrm{ for } I \subseteq S
\]
(see {\cite[Section 5.8, Exercise 4]{Hum}}), so since $\Phi_I=W_I \cdot \Pi_I$, we have
\[
\mathrm{supp}(\gamma)=\mathrm{supp}(s_\gamma) \textrm{ for } \gamma \in \Phi.
\]

The following lemma will be used later.
\begin{lem}
[{\cite[Lemma 2.7]{Nui_ext}}]
\label{lem:wdoesnotfixgamma}
Let $1 \neq w \in W$, $I=\mathrm{supp}(w)$, $\gamma \in \Phi^+$ and $J=\mathrm{supp}(\gamma)$.
Suppose that $I \cap J=\emptyset$ and $J$ is adjacent to $I$ in the Coxeter graph $\Gamma$.
Then $w \cdot \gamma \in \Phi_{I \cup J}^+ \smallsetminus \Phi_J$, hence $w \cdot \gamma \neq \gamma$.
\end{lem}

The following subgraph of a Coxeter graph is introduced.
\begin{defn}
\label{defn:oddCoxetergraph}
For a subset $I \subseteq S$, we define the \emph{odd Coxeter graph} $\Gamma_I^{\mathrm{odd}}$ to be the subgraph of the Coxeter graph $\Gamma_I$ of $(W_I,I)$ obtained by removing all edges with labels either even or $\infty$.
Moreover, for $s \in I$, let $I_{\sim s}^{\mathrm{odd}}$ denote the vertex set of the connected component $(\Gamma_I^{\mathrm{odd}})_{\sim s}$ of $\Gamma_I^{\mathrm{odd}}$ containing $s$.
\end{defn}
Then we have the following well-known property:
\begin{prop}
[See {\cite[Section 5.3, Exercise]{Hum}}]
\label{prop:oddCoxetergraph}
For $I \subseteq S$ and $s,t \in I$, the following three conditions are equivalent:
(1) $s$ is conjugate to $t$ in $W_I$; (2) $\alpha_s \in W_I \cdot \alpha_t$; (3) $s \in I_{\sim t}^{\mathrm{odd}}$.
\end{prop}

Let $W_I$ be a standard parabolic subgroup of $W$.
We say that $(W_I,I)$, $W_I$, $I$ or $\Gamma_I$ is of \emph{finite type} if $W_I$ is a finite group.
As is summarized in {\cite[Chapter 2]{Hum}}, the finite irreducible Coxeter groups are completely classified, into the following seven families; a standard parabolic subgroup $W_I$ of finite rank (or $(W_I,I)$, $I$ or $\Gamma_I$) is \emph{of type $X$} if, under a labelling $I=\{s_1,s_2,\dots,s_n\}$ of elements of $I$, the data $m_{i,j}=m_{s_i,s_j}$ satisfy the following condition:
\begin{equation}
\label{eq:ruleforfinitepart}
\begin{cases}
m_{i,i+1}=3\quad (1 \leq i \leq n-1) & \textrm{if } X=A_n,\ n \geq 1\\
m_{i,i+1}=3\quad (1 \leq i \leq n-2),\ m_{n-1,n}=4 & \textrm{if } X=B_n,\ n \geq 2\\
m_{i,i+1}=m_{n-2,n}=3\quad (1 \leq i \leq n-2) & \textrm{if } X=D_n,\ n \geq 4\\
m_{1,3}=m_{2,4}=m_{i,i+1}=3\quad (3 \leq i \leq n-1) & \textrm{if } X=E_n,\ 6 \leq n \leq 8\\
m_{1,2}=m_{3,4}=3,\ m_{2,3}=4 & \textrm{if } X=F_4,\ n=4\\
m_{i,i+1}=3\quad (1 \leq i \leq n-2),\ m_{n-1,n}=5 & \textrm{if } X=H_n,\ n=3,4\\
m_{1,2}=m<\infty & \textrm{if } X=I_2(m),\ n=2
\end{cases}
\end{equation}
and $m_{i,j}=2$ if $i \neq j$ and $m_{i,j}$ is not determined by the above list.
When we emphasize the chosen labelling of elements of $I$ as well as the type of $I$, we say that \emph{a sequence $(s_1,s_2,\dots,s_n)$ is of type $X$}.
Let $W(X)$ denote the unique Coxeter group (up to isomorphism) of type $X$.

\begin{figure}
\centering
\begin{picture}(150,50)(0,-50)
\put(70,-10){$A_\infty$}
\put(10,-30){\circle{6}}\put(40,-30){\circle{6}}\put(70,-30){\circle{6}}\put(100,-30){\circle{6}}
\put(10,-45){\hbox to0pt{\hss$1$\hss}}\put(40,-45){\hbox to0pt{\hss$2$\hss}}\put(70,-45){\hbox to0pt{\hss$3$\hss}}\put(100,-45){\hbox to0pt{\hss$4$\hss}}
\put(13,-30){\line(1,0){24}}\put(43,-30){\line(1,0){24}}\put(73,-30){\line(1,0){24}}\put(103,-30){\line(1,0){20}}\put(128,-33){$\cdots$}
\end{picture}
\qquad
\begin{picture}(170,50)(0,-50)
\put(70,-10){$A_{\pm \infty}$}
\put(40,-30){\circle{6}}\put(70,-30){\circle{6}}\put(100,-30){\circle{6}}\put(130,-30){\circle{6}}
\put(40,-45){\hbox to0pt{\hss$-1$\hss}}\put(70,-45){\hbox to0pt{\hss$0$\hss}}\put(100,-45){\hbox to0pt{\hss$1$\hss}}\put(130,-45){\hbox to0pt{\hss$2$\hss}}
\put(0,-33){$\cdots$}\put(37,-30){\line(-1,0){20}}\put(43,-30){\line(1,0){24}}\put(73,-30){\line(1,0){24}}\put(103,-30){\line(1,0){24}}\put(133,-30){\line(1,0){20}}\put(158,-33){$\cdots$}
\end{picture}\\
\ \\
\begin{picture}(150,50)(0,-50)
\put(70,-10){$B_\infty$}
\put(10,-30){\circle{6}}\put(40,-30){\circle{6}}\put(70,-30){\circle{6}}\put(100,-30){\circle{6}}
\put(10,-45){\hbox to0pt{\hss$0$\hss}}\put(40,-45){\hbox to0pt{\hss$1$\hss}}\put(70,-45){\hbox to0pt{\hss$2$\hss}}\put(100,-45){\hbox to0pt{\hss$3$\hss}}
\put(13,-30){\line(1,0){24}}\put(43,-30){\line(1,0){24}}\put(73,-30){\line(1,0){24}}\put(103,-30){\line(1,0){20}}\put(128,-33){$\cdots$}
\put(25,-26){\hbox to0pt{\hss$4$\hss}}
\end{picture}
\qquad
\begin{picture}(160,50)(0,-50)
\put(70,-10){$D_\infty$}
\put(20,-15){\circle{6}}\put(20,-45){\circle{6}}\put(50,-30){\circle{6}}\put(80,-30){\circle{6}}\put(110,-30){\circle{6}}
\put(8,-17){\hbox to0pt{\hss$0'$\hss}}\put(8,-47){\hbox to0pt{\hss$0''$\hss}}\put(50,-45){\hbox to0pt{\hss$1$\hss}}\put(80,-45){\hbox to0pt{\hss$2$\hss}}\put(110,-45){\hbox to0pt{\hss$3$\hss}}
\put(47,-28){\line(-2,1){24}}\put(47,-32){\line(-2,-1){24}}\put(53,-30){\line(1,0){24}}\put(83,-30){\line(1,0){24}}\put(113,-30){\line(1,0){20}}\put(138,-33){$\cdots$}
\end{picture}
\caption{Four Coxeter groups of infinite rank}
\label{fig:nonf.g.Coxetergroup}
\end{figure}
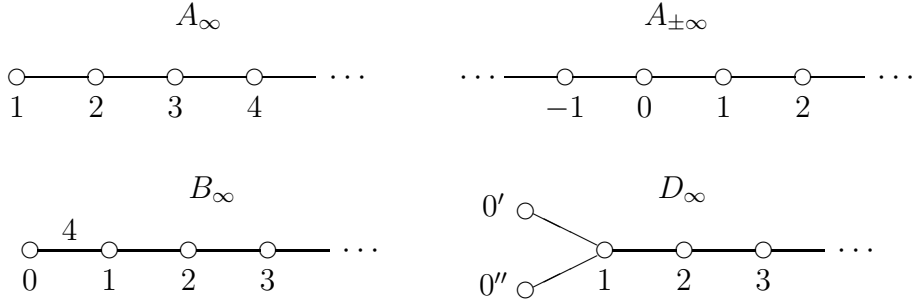

On the other hand, the four Coxeter groups of infinite (countable) ranks in Figure \ref{fig:nonf.g.Coxetergroup} will play an important role in our argument below.
We also adopt the above terminology; for example, an infinite sequence $(s'_0,s''_0,s_1,s_2,\dots)$ of elements of $S$ is of type $D_\infty$ if $m_{s'_0,s_1}=m_{s''_0,s_1}=m_{s_i,s_{i+1}}=3$ for all $i \geq 1$, and the value of the $m$ is $2$ for any other pair of distinct elements of $S$.

If $(W_I,I)$ is of finite type, let $w_0(I)$ denote the \emph{longest element} of $W_I$ (with respect to $S$), which is an involution and maps $\Pi_I$ onto $-\Pi_I$; this induces a graph automorphism $\sigma_I \in \mathrm{Aut}\,\Gamma_I$ such that
\[
w_0(I)sw_0(I)=\sigma_I(s) \textrm{ and } w_0(I) \cdot \alpha_s=-\alpha_{\sigma_I(s)} \textrm{ for all } s \in I.
\]
The following fact is well known.
\begin{prop}
\label{prop:centerofW}
Let $(W,S)$ be an irreducible Coxeter system.
\begin{enumerate}
\item The center $Z(W)$ of $W$ is nontrivial if and only if $(W,S)$ is of type $A_1$, $B_n$ (with $2 \leq n<\infty$), $D_n$ (with $4 \leq n<\infty$ and $n$ even), $E_7$, $E_8$, $F_4$, $H_3$, $H_4$ or $I_2(m)$ (with $m$ even).
In particular, $(W,S)$ is of finite type if $Z(W) \neq 1$.
\item If $Z(W) \neq 1$, then $Z(W)$ has order two and is generated by $w_0(S)$, and $\sigma_S$ is the identity map on $\Gamma$.
\item If $Z(W)=1$ and $|W|<\infty$, namely $(W,S)$ is of type $A_n$ (with $2 \leq n<\infty$), $D_n$ (with $4 \leq n<\infty$ and $n$ odd), $E_6$ or $I_2(m)$ (with $m$ odd), then $\sigma_S$ is the unique non-identity automorphism of $\Gamma$.
\end{enumerate}
\end{prop}
Moreover, we will use the following result later.
\begin{prop}
[See {\cite[Corollary 2.22 (3)]{Nui_isom}}]
\label{prop:rootformaximalparabolic}
Suppose that $(W,S)$ is irreducible and not of finite type.
Then for $s \in S$, the only elements of $W$ which normalizes $W_I$ (where $I=S \smallsetminus s$) are the elements of $W_I$.
Hence there does not exist a root $\gamma \in \Phi$ which is orthogonal to $\Pi_I$.
\end{prop}
\begin{proof}
The former part is {\cite[Corollary 2.22 (3)]{Nui_isom}} (note that this is also an easy consequence of the main result of \cite{Bri-How_normal}).
For the latter part, observe that such a root $\gamma$ cannot belong to $\Phi_I$, so $s_\gamma \not\in I^{W_I}$; while $s_\gamma$ normalizes (actually, centralizes) $W_I$, so $s_\gamma \in S^W \cap W_I=I^{W_I}$ by the former part.
This contradiction proves the claim.
\end{proof}
\section{Notes on centralizers of reflections}
\label{sec:Zofreflections}

In this section, we summarize some preceding results on the structure of the centralizer $Z_W(x)$ of a generator $x \in S$ in an arbitrary Coxeter group $W$, which play a central role of this paper, together with some additional remarks.
Descriptions of the centralizer $Z_W(x)$ are found in several papers; namely, in \cite{How} by Robert B.\ Howlett for finite $W$; in \cite{Bri} by Brink; in \cite{Bor} by Richard E.\ Borcherds; in \cite{Bri-How_normal} by Brink and Howlett; and in \cite{Nui_centra} by the author.
Note that the centralizer of a reflection in $W$ is conjugate to some $Z_W(x)$ since the reflection is (by definition) conjugate to some $x \in S$.
\subsection{Structure of the centralizers}
\label{sec:structureofZ}

The first uniform description of $Z_W(x)$ for general $W$ and $x \in S$ was given by Brink \cite{Bri} (note that the general result in \cite{Nui_centra} agrees with hers when specializing to this case; see {\cite[Example 4.12]{Nui_centra}}).
Her result shows that, if $x \in I \subseteq S$, then the centralizer $Z_{W_I}(x)$ of $x$ in $W_I$ admits the following decomposition:
\[
Z_{W_I}(x)=(\langle x \rangle \times W_I^{\perp x}) \rtimes Y_x^{(I)}=\langle x \rangle \times (W_I^{\perp x} \rtimes Y_x^{(I)}).
\]
(In what follows, we omit the subscripts or superscripts `$S$' when $I=S$.)
Here we denote by $W_I^{\perp x}$ the subgroup of $W$ generated by the reflections along roots in the set $\Phi_I^{\perp x}$ defined by
\[
\Phi_I^{\perp x}=\{\gamma \in \Phi_I \mid \langle \gamma, \alpha_x \rangle=0\}.
\]
On the other hand, we put
\[
Y^{(I)}_{y,z}=\{w \in W_I \mid \alpha_y=w \cdot \alpha_z \textrm{ and } \Phi_I^{\perp z}\left[w\right]=\emptyset\} \textrm{ for } y,z \in I.
\]
Since $(\Phi_I^{\perp y})^+=w \cdot (\Phi_I^{\perp z})^+$ for any $w \in Y^{(I)}_{y,z}$, the family $Y^{(I)}=\{Y^{(I)}_{y,z}\}_{y,z \in I}$ forms a groupoid (see Section \ref{sec:graphsandgroupoids} for notion of groupoids) with vertex set $I$, and the third factor $Y_x^{(I)}$ of $Z_{W_I}(x)$ is the vertex group of $Y^{(I)}$ at $x \in I$.

This group $Y_x^{(I)}$ is described in terms of the odd Coxeter graph $\Gamma_I^{\mathrm{odd}}$ of $W_I$ (see Definition \ref{defn:oddCoxetergraph}).
Brink gave an isomorphism $\pi^{(I)}:\mathscr{F}^{(I)} \to Y^{(I)}$ from the fundamental groupoid
\[
\mathscr{F}^{(I)}=\pi_1(\Gamma_I^{\mathrm{odd}};*,*)
\]
of the graph $\Gamma_I^{\mathrm{odd}}$ (see Section \ref{sec:graphsandgroupoids} for definition) to $Y^{(I)}$ in the following manner:
\begin{thm}
[\cite{Bri}]
\label{thm:Brink_pipreserveslength}
The groupoid isomorphism $\pi^{(I)}:\mathscr{F}^{(I)} \to Y^{(I)}$ is defined by $\pi^{(I)}(y)=y$ for vertices $y \in I$,
\[
\pi^{(I)}(y,z)=(zy)^{(m-1)/2} \in Y^{(I)}_{y,z} \textrm{ for distinct } y,z \in I \textrm{ with } m=m_{y,z} \textrm{ odd}
\]
and extending it multiplicatively to the whole $\mathscr{F}^{(I)}$.
Moreover, if a path $(y_0,y_1,\dots,y_n) \in \mathscr{F}^{(I)}_{y_0,y_n}$ is non-backtracking, i.e.\ $y_{i+2} \neq y_i$ for $0 \leq i \leq n-2$, then $\ell(\pi^{(I)}(y_0,y_1,\dots,y_n))=\sum_{i=0}^{n-1}\ell(\pi^{(I)}(y_i,y_{i+1}))$.
\end{thm}
\begin{cor}
\label{cor:supportofelementsofY}
If a path $(y_0,y_1,\dots,y_n) \in \mathscr{F}^{(I)}_{y_0,y_n}$ is nontrivial (i.e.\ $n \geq 1$) and non-backtracking, then $\mathrm{supp}(\pi^{(I)}(y_0,y_1,\dots,y_n))=\{y_0,y_1,\dots,y_n\}$.
\end{cor}
\begin{proof}
By Theorem \ref{thm:Brink_pipreserveslength} and Lemma \ref{lem:supportofproduct}, we have
\[
\mathrm{supp}(\pi^{(I)}(y_0,y_1,\dots,y_n))=\bigcup_{i=0}^{n-1}\mathrm{supp}(\pi^{(I)}(y_i,y_{i+1})),
\]
while $\mathrm{supp}(\pi^{(I)}(y_i,y_{i+1}))=\{y_i,y_{i+1}\}$ by definition, proving the claim.
\end{proof}
This $\pi^{(I)}$ induces a group isomorphism from the fundamental group $\mathscr{F}^{(I)}_x=\pi_1(\Gamma_I^{\mathrm{odd}};x)$ of $\Gamma_I^{\mathrm{odd}}$ at $x \in I$ to the group $Y^{(I)}_x$.
Note that, for distinct $y,z \in I$ with $m_{y,z}$ odd, we have
\[
\pi^{(I)}(y,z)=w_0(\{y,z\})z=s_{\widetilde{\alpha}_{y,z}}z
\]
where $\widetilde{\alpha}_{y,z}$ denotes the highest root of a Coxeter group $W_{\{y,z\}}$; that is
\begin{equation}
\label{eq:highestroot_misodd}
\widetilde{\alpha}_{y,z}=\frac{1}{2\sin(\pi/(2m))}(\alpha_y+\alpha_z) \textrm{ with } m=m_{y,z} \textrm{ odd}.
\end{equation}

We consider the second factor $W_I^{\perp x}$ of $Z_{W_I}(x)$.
As is mentioned in \cite{Bri}, a general theorem of Deodhar \cite{Deo_refsub} or of Dyer \cite{Dye} implies that $W_I^{\perp x}$ is a Coxeter group with the set $\Phi_I^{\perp x}$ playing a role of the root system.
More precisely, let $\Pi_x^{(I)}$ be the set of the roots in $(\Phi_I^{\perp x})^+$ which cannot be a positive linear combination of other roots in $(\Phi_I^{\perp x})^+$; in other words, $\Pi_x^{(I)}$ is the ``simple system'' in $\Phi_I^{\perp x}$.
Define
\[
R_x^{(I)}=\{s_\gamma \mid \gamma \in \Pi_x^{(I)}\}.
\]
\begin{thm}
\label{thm:propertyofWperpx}
Let $I \subseteq J \subseteq S$ and $x,y \in I$.
\begin{enumerate}
\item \label{thm_item:WperpxisCoxetergroup}
\textbf{(\cite{Deo_refsub}; see also {\cite[Theorem 2.3]{Nui_centra}})}
The pair $(W_I^{\perp x},R_x^{(I)})$ is a Coxeter system.
For $w \in W_I^{\perp x}$, the length $\widetilde{\ell}_{I,x}(w)$ of $w$ with respect to $R_x^{(I)}$ coincides with the cardinality of $\Phi_I^{\perp x}\left[w\right]$.
Moreover, for $\gamma \in (\Phi_I^{\perp x})^+$, we have $\widetilde{\ell}_{I,x}(ws_\gamma)<\widetilde{\ell}_{I,x}(w)$ if and only if $w \cdot \gamma \in \Phi^-$.
\item \label{thm_item:orderandbilinearform}
\textbf{({\cite[Theorem 4.4]{Dye}})}
Let $\beta,\gamma \in \Pi_x^{(I)}$, and $m$ the order of $s_\beta s_\gamma \in W$.
Then we have
\[
\begin{cases}
\langle \beta, \gamma \rangle=-\cos(\pi/m) & \textrm{if } m<\infty;\\
\langle \beta, \gamma \rangle \leq -1 & \textrm{if } m=\infty.
\end{cases}
\]
In particular, we have $\langle \beta, \gamma \rangle \leq 0$ for distinct $\beta,\gamma \in \Pi_x^{(I)}$.
\item \label{thm_item:conjugateinduceisomorphism}
\textbf{({\cite[Theorem 4.6 (4)]{Nui_centra}})}
If $w \in Y_{x,y}^{(I)}$, then $\Pi_x^{(I)}=w \cdot \Pi_y^{(I)}$, hence the map $u \mapsto wuw^{-1}$ is an isomorphism of Coxeter systems from $(W_I^{\perp y},R_y^{(I)})$ to $(W_I^{\perp x},R_x^{(I)})$.
\item \label{thm_item:inclusionofWperp}
We have $W_I^{\perp x}=W_J^{\perp x} \cap W_I$, $\Pi_x^{(I)}=\Pi_x^{(J)} \cap \Phi_I$ and $R_x^{(I)}=R_x^{(J)} \cap W_I$.
\end{enumerate}
\end{thm}
\begin{proof}
We prove Claim \ref{thm_item:inclusionofWperp}.
For the second equality, if $\gamma \in (\Phi_I^{\perp x})^+$ is expressed as a positive linear combination of some $\beta \in (\Phi_J^{\perp x})^+$ different from $\gamma$, then $\mathrm{supp}(\gamma)$ is the union of all the sets $\mathrm{supp}(\beta)$, deducing that $\beta \in \Phi_I^{\perp x}$.
Thus $\Pi_x^{(I)} \subseteq \Pi_x^{(J)} \cap \Phi_I$, while $\Pi_x^{(J)} \cap \Phi_I \subseteq \Pi_x^{(I)}$ by definition, therefore $\Pi_x^{(I)}=\Pi_x^{(J)} \cap \Phi_I$.
Now the third equality also follows, since $J^{W_J} \cap W_I=I^{W_I}$.

Finally, we prove the first equality.
It suffices to show that $w \in W_I^{\perp x}$ for all $w \in W_J^{\perp x} \cap W_I$; we proceed the proof by induction on $n=\widetilde{\ell}_{J,x}(w)$, the case $n=0$ being obvious.
So suppose that $n \geq 1$, and take a reduced expression $w=s_{\gamma_1}s_{\gamma_2} \cdots s_{\gamma_n}$ with respect to $R_x^{(J)}$, where $\gamma_i \in \Pi_x^{(J)}$.
Then we have $\widetilde{\ell}_{J,x}(ws_{\gamma_n})<\widetilde{\ell}_{J,x}(w)$, so $w \cdot \gamma_n \in \Phi^-$ by Claim \ref{thm_item:WperpxisCoxetergroup}, therefore $\gamma_n \in \Phi_I$ since $w \in W_I$.
Thus $s_{\gamma_n} \in R_x^{(J)} \cap W_I=R_x^{(I)}$, so $ws_{\gamma_n} \in W_J^{\perp x} \cap W_I$.
Now the induction assumption implies that $ws_{\gamma_n} \in W_I^{\perp x}$, while $s_{\gamma_n} \in R_x^{(I)}$ as above, hence $w \in W_I^{\perp x}$ as desired.
\end{proof}
The result in the author's papar \cite{Nui_centra} describes the Coxeter group $W_I^{\perp x}$ more precisely.
We prepare some further notations.
Choose and fix a maximal tree $\mathcal{T}$ in the connected component $(\Gamma_I^{\mathrm{odd}})_{\sim x}$ of $\Gamma_I^{\mathrm{odd}}$ containing $x$, and for $y,z \in I_{\sim x}^{\mathrm{odd}}$ (see Definition \ref{defn:oddCoxetergraph} for notation), let $p_{y,z}^{(I)} \in \mathscr{F}_{y,z}^{(I)}$ denote the unique non-backtracking path in $\mathcal{T}$ from $z$ to $y$.
Define
\[
\mathcal{E}_x^{(I)}=\{(y,s) \mid y \in I_{\sim x}^{\mathrm{odd}},\ s \in I \textrm{ and } m_{y,s} \textrm{ is even}\},
\]
and put
\begin{equation}
\label{eq:highestroot_miseven}
\widetilde{\alpha}_{y,s}=\frac{1}{\sin(\pi/m_{y,s})}(\cos(\pi/m_{y,s})\alpha_y+\alpha_s) \in \Phi_{\{y,s\}}^+ \textrm{ for } (y,s) \in \mathcal{E}_x^{(I)}.
\end{equation}
Note that this $\widetilde{\alpha}_{y,s}$ is the unique root in $(\Phi_{\{y,s\}}^{\perp y})^+$ (see {\cite[Lemma 4.1]{Nui_centra}}), so we have $\widetilde{\alpha}_{y,s} \in \Pi_y^{(I)}$.
Note also that $\widetilde{\alpha}_{y,s}=\alpha_s$ if $m_{y,s}=2$, while $\mathrm{supp}(\widetilde{\alpha}_{y,s})=\{y,s\}$ if $m_{y,s} \neq 2$.
Now define
\[
\gamma_x^{(I)}(c;y,s)=\pi^{(I)}(cp_{x,y}^{(I)}) \cdot \widetilde{\alpha}_{y,s} \textrm{ for } c \in \mathscr{F}_x^{(I)} \textrm{ and } (y,s) \in \mathcal{E}_x^{(I)}.
\]
Since $cp_{x,y}^{(I)} \in \mathscr{F}_{x,y}^{(I)}$ and $\widetilde{\alpha}_{y,s} \in \Pi_y^{(I)}$, it follows from Theorems \ref{thm:Brink_pipreserveslength} and \ref{thm:propertyofWperpx} (\ref{thm_item:conjugateinduceisomorphism}) that $\gamma_x^{(I)}(c;y,s) \in \Pi_x^{(I)}$.
Let $r_x^{(I)}(c;y,s)$ denote the reflection along the root $\gamma_x^{(I)}(c;y,s)$, so $r_x^{(I)}(c;y,s) \in R_x^{(I)}$.
For simplicity, we abbreviate $\gamma_x^{(I)}(1;y,s)$ to $\gamma_x^{(I)}(y,s)$ and $r_x^{(I)}(1;y,s)$ to $r_x^{(I)}(y,s)$.
\begin{exmp}
\label{ex:supportofgammaxys_length2}
Suppose that $I=\{x,y,s\}$, $m_{x,y}$ is odd, $m_{y,s}$ is even and $m_{x,s} \neq 2$.
Put $\beta=\widetilde{\alpha}_{x,y}$ and $\beta'=\widetilde{\alpha}_{y,s}$.
Then we have $(y,s) \in \mathcal{E}_x^{(I)}$ and
\[
\gamma_x^{(I)}(y,s)=\pi^{(I)}(x,y) \cdot \beta'=s_\beta y \cdot \beta'=s_\beta \cdot \beta'=\beta'-2\langle \beta, \beta' \rangle \beta
\]
since $\langle \beta', \alpha_y \rangle=0$.
Moreover, we have
\[
\langle \beta, \beta' \rangle=(\left[\alpha_x\right]\beta) \langle \alpha_x, \beta' \rangle=(\left[\alpha_x\right]\beta)\bigl((\left[\alpha_y\right]\beta') \langle \alpha_x, \alpha_y \rangle+(\left[\alpha_s\right]\beta') \langle \alpha_x, \alpha_s \rangle\bigr)<0
\]
since $\left[\alpha_x\right]\beta>0$, $\left[\alpha_y\right]\beta' \geq 0$, $\left[\alpha_s\right]\beta'>0$, $\langle \alpha_x, \alpha_y \rangle<0$ and $\langle \alpha_x, \alpha_s \rangle<0$.
Thus we have $\mathrm{supp}(\gamma_x^{(I)}(y,s))=\mathrm{supp}(\beta) \cup \mathrm{supp}(\beta')=I$.
\end{exmp}
\begin{exmp}
\label{ex:supportofgammaxys_length3}
Suppose that $I=\{x,y,z,s\}$, $m_{x,y}$ and $m_{y,z}$ are odd, $m_{x,z}$, $m_{x,s}$ and $m_{z,s}$ are even, and $m_{y,s} \neq 2$.
Then we have $(z,s) \in \mathcal{E}_x^{(I)}$ and $\gamma_x^{(I)}(z,s)=\pi^{(I)}(x,y,z) \cdot \widetilde{\alpha}_{z,s}$.
Now we have $\beta=\pi^{(I)}(y,z) \cdot \widetilde{\alpha}_{z,s} \in \Pi_y^{(I)}$ and $\mathrm{supp}(\beta)=I \smallsetminus x$ by Example \ref{ex:supportofgammaxys_length2}; so $\langle \widetilde{\alpha}_{x,y}, \beta \rangle<0$ by a similar argument to Example \ref{ex:supportofgammaxys_length2} (since $\langle \alpha_x, \alpha_y \rangle<0$).
Hence we have
\[
\gamma_x^{(I)}(z,s)=\pi^{(I)}(x,y) \cdot \beta=s_{\widetilde{\alpha}_{x,y}} \cdot \beta=\beta-2\langle \widetilde{\alpha}_{x,y}, \beta \rangle \widetilde{\alpha}_{x,y}
\]
and so $\mathrm{supp}(\gamma_x^{(I)}(z,s))=\mathrm{supp}(\beta) \cup \mathrm{supp}(\widetilde{\alpha}_{x,y})=I$.
\end{exmp}
By generalizing the arguments in Examples \ref{ex:supportofgammaxys_length2} and \ref{ex:supportofgammaxys_length3}, we obtain the following observation.
\begin{lem}
\label{lem:conditionfornonorthogonal}
Let $\beta,\gamma \in \Pi_x^{(I)}$, and suppose that $(\mathrm{supp}(\beta) \smallsetminus x) \cap \mathrm{supp}(\gamma)=\emptyset$, and $\mathrm{supp}(\beta) \smallsetminus x$ is adjacent to $\mathrm{supp}(\gamma)$ in the Coxeter graph $\Gamma_I$.
Then $\langle \beta, \gamma \rangle<0$; hence (by Theorem \ref{thm:propertyofWperpx} (\ref{thm_item:orderandbilinearform})) two generators $s_\beta$ and $s_\gamma$ of $W_I^{\perp x}$ do not commute.
\end{lem}

\begin{figure}
\centering
\begin{picture}(360,570)(0,-570)
\put(50,-10){\hbox to0pt{\hss{\bf (I)}\hss}}
\put(20,-30){\circle{6}}\put(20,-20){\hbox to0pt{\hss$y=z$\hss}}
\put(20,-90){\circle{6}}\put(20,-105){\hbox to0pt{\hss$s$\hss}}
\put(80,-60){\circle{6}}\put(90,-63){$t$}
\put(23,-88){\line(2,1){54}}\put(45,-88){$m$}
\put(50,-120){\hbox to0pt{\hss$q=(z)$,\ $k=m$\hss}}
\put(50,-135){\hbox to0pt{\hss$2 \leq m<\infty$\hss}}
\put(170,-10){\hbox to0pt{\hss{\bf (II)}\hss}}
\put(140,-30){\circle{6}}\put(140,-20){\hbox to0pt{\hss$y=z$\hss}}
\put(140,-90){\circle{6}}\put(140,-105){\hbox to0pt{\hss$s$\hss}}
\put(200,-60){\circle{6}}\put(210,-63){$t$}
\put(143,-32){\line(2,-1){54}}\put(165,-40){$m$}
\put(170,-120){\hbox to0pt{\hss$q=(z)$,\ $k=2$\hss}}
\put(170,-135){\hbox to0pt{\hss$m$ even\hss}}
\put(290,-10){\hbox to0pt{\hss{\bf (III)}\hss}}
\put(260,-30){\circle{6}}\put(260,-20){\hbox to0pt{\hss$y=z$\hss}}
\put(260,-90){\circle{6}}\put(260,-105){\hbox to0pt{\hss$s$\hss}}
\put(320,-60){\circle{6}}\put(330,-63){$t$}
\put(263,-32){\line(2,-1){54}}\put(285,-40){$4$}
\put(263,-88){\line(2,1){54}}
\put(290,-120){\hbox to0pt{\hss$q=(z)$,\ $k=4$\hss}}
\put(50,-160){\hbox to0pt{\hss{\bf (IV)}\hss}}
\put(20,-180){\circle{6}}\put(20,-170){\hbox to0pt{\hss$y$\hss}}
\put(20,-240){\circle{6}}\put(20,-255){\hbox to0pt{\hss$s=t$\hss}}
\put(80,-210){\circle{6}}\put(90,-213){$z$}
\put(23,-182){\line(2,-1){54}}\put(45,-190){$m$}
\put(50,-270){\hbox to0pt{\hss$q=(y,z)$,\ $k=1$\hss}}
\put(50,-285){\hbox to0pt{\hss$m$ odd\hss}}
\put(170,-160){\hbox to0pt{\hss{\bf (V)}\hss}}
\put(140,-180){\circle{6}}\put(140,-170){\hbox to0pt{\hss$y=t$\hss}}
\put(140,-240){\circle{6}}\put(140,-255){\hbox to0pt{\hss$s=z$\hss}}
\put(200,-210){\circle{6}}\put(210,-213){$y'$}
\put(143,-182){\line(2,-1){54}}
\put(143,-238){\line(2,1){54}}
\put(170,-270){\hbox to0pt{\hss$q=(y,y',z)$,\ $k=1$\hss}}
\put(290,-160){\hbox to0pt{\hss{\bf (VI)}\hss}}
\put(260,-180){\circle{6}}\put(260,-170){\hbox to0pt{\hss$y=t$\hss}}
\put(260,-240){\circle{6}}\put(260,-255){\hbox to0pt{\hss$s=z$\hss}}
\put(320,-210){\circle{6}}\put(330,-213){$y'$}
\put(263,-182){\line(2,-1){54}}
\put(263,-238){\line(2,1){54}}\put(285,-238){$5$}
\put(290,-270){\hbox to0pt{\hss$q=(y,y',z)$,\ $k=2$\hss}}
\put(50,-310){\hbox to0pt{\hss{\bf (VII)}\hss}}
\put(20,-330){\circle{6}}\put(20,-320){\hbox to0pt{\hss$y=t$\hss}}
\put(20,-390){\circle{6}}\put(20,-405){\hbox to0pt{\hss$s=z$\hss}}
\put(80,-360){\circle{6}}\put(90,-363){$y'$}
\put(23,-332){\line(2,-1){54}}\put(45,-340){$5$}
\put(23,-388){\line(2,1){54}}
\put(50,-420){\hbox to0pt{\hss$q=(y,y',z)$,\ $k=2$\hss}}
\put(170,-310){\hbox to0pt{\hss{\bf (VIII)}\hss}}
\put(140,-330){\circle{6}}\put(140,-320){\hbox to0pt{\hss$y$\hss}}
\put(140,-390){\circle{6}}\put(140,-405){\hbox to0pt{\hss$s=t$\hss}}
\put(200,-360){\circle{6}}\put(210,-363){$z$}
\put(143,-332){\line(2,-1){54}}
\put(143,-388){\line(2,1){54}}\put(165,-388){$4$}
\put(170,-420){\hbox to0pt{\hss$q=(y,z)$,\ $k=2$\hss}}
\put(290,-310){\hbox to0pt{\hss{\bf (IX)}\hss}}
\put(260,-330){\circle{6}}\put(260,-320){\hbox to0pt{\hss$y=z$\hss}}
\put(260,-390){\circle{6}}\put(260,-405){\hbox to0pt{\hss$s$\hss}}
\put(320,-360){\circle{6}}\put(330,-363){$t$}
\put(260,-333){\line(0,-1){54}}\put(255,-363){\hbox to0pt{\hss$m$}}
\put(290,-420){\hbox to0pt{\hss$q=(z)$,\ $k=2$\hss}}
\put(290,-435){\hbox to0pt{\hss$m$ even\hss}}
\put(50,-460){\hbox to0pt{\hss{\bf (X)}\hss}}
\put(20,-480){\circle{6}}\put(20,-470){\hbox to0pt{\hss$y=z$\hss}}
\put(20,-540){\circle{6}}\put(20,-555){\hbox to0pt{\hss$s$\hss}}
\put(80,-510){\circle{6}}\put(90,-513){$t$}
\put(20,-483){\line(0,-1){54}}\put(15,-513){\hbox to0pt{\hss$4$}}
\put(23,-538){\line(2,1){54}}
\put(50,-570){\hbox to0pt{\hss$q=(z)$,\ $k=4$\hss}}
\put(170,-460){\hbox to0pt{\hss{\bf (XI)}\hss}}
\put(140,-480){\circle{6}}\put(140,-470){\hbox to0pt{\hss$y$\hss}}
\put(140,-540){\circle{6}}\put(140,-555){\hbox to0pt{\hss$s=t$\hss}}
\put(200,-510){\circle{6}}\put(210,-513){$z$}
\put(140,-483){\line(0,-1){54}}\put(135,-513){\hbox to0pt{\hss$4$}}
\put(143,-482){\line(2,-1){54}}
\put(170,-570){\hbox to0pt{\hss$q=(y,z)$,\ $k=2$\hss}}
\end{picture}
\caption{The eleven subgraphs of $(\Gamma_I^{\mathrm{odd}})_{\sim x}$}
\label{fig:localsubgraph}
\end{figure}
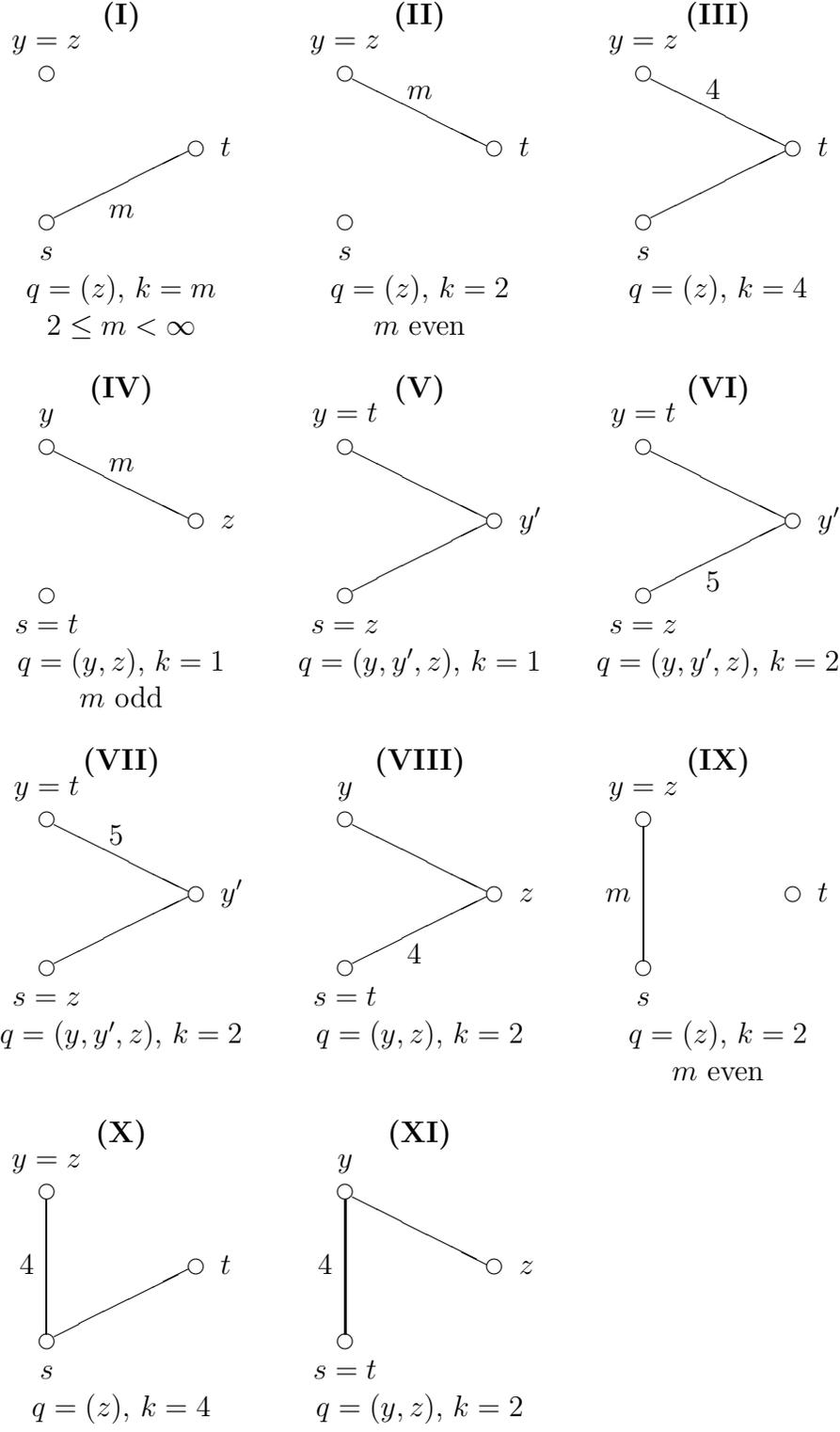

For an integer $k \geq 1$, we define relations $\overset{k}{\sim}_I$ on $\mathcal{E}_x^{(I)}$ and on $\mathscr{F}_x^{(I)} \times \mathcal{E}_x^{(I)}$, as follows.
For a path $q \in \mathscr{F}_{y,z}^{(I)}$ in $(\Gamma_I^{\mathrm{odd}})_{\sim x}$, put $q_{(x)}=p_{x,y}^{(I)}qp_{z,x}^{(I)} \in \mathscr{F}_x^{(I)}$.
Then for $(y,s),(z,t) \in \mathcal{E}_x^{(I)}$ and $c \in \mathscr{F}_x^{(I)}$, we define $(c;y,s) \overset{k}{\sim}_I (cq_{(x)};z,t)$ if $(\Gamma_I^{\mathrm{odd}})_{\sim x}$ possesses a subgraph (with three vertices) of one of the eleven forms in Figure \ref{fig:localsubgraph} (where the path $q$ and the integer $k$ are determined in the figure); in this case, we also define $(y,s) \overset{k}{\sim}_I (z,t)$.
Observe that all of them are symmetric relations.
Moreover, let $\sim_I$ denote the transitive closure of $\overset{1}{\sim}_I$, which is an equivalence relation.
The next theorem describes the structure of $W_I^{\perp x}$ in terms of these relations.
\begin{thm}
[{\cite[Theorems 4.13 and 4.14]{Nui_centra}}]
\label{thm:presentationofWperpx}
We have
\[
\Pi_x^{(I)}=\{\gamma_x^{(I)}(c;y,s) \mid c \in \mathscr{F}_x^{(I)} \textrm{ and } (y,s) \in \mathcal{E}_x^{(I)}\}.
\]
For $c,c' \in \mathscr{F}_x^{(I)}$ and $(y,s),(y',s') \in \mathcal{E}_x^{(I)}$, we have $r_x^{(I)}(c;y,s)=r_x^{(I)}(c';y',s')$ if and only if $(c;y,s) \sim_I (c';y',s')$.
Moreover, for $2 \leq k<\infty$, the product $r_x^{(I)}(c;y,s)r_x^{(I)}(c';y',s')$ has order $k$ if and only if there exist $d,d' \in \mathscr{F}_x^{(I)}$ and $(z,t),(z',t') \in \mathcal{E}_x^{(I)}$ such that $(c;y,s) \sim_I (d;z,t) \overset{k}{\sim}_I (d';z',t') \sim_I (c';y',s')$.
\end{thm}

For $(y,s) \in \mathcal{E}_x^{(I)}$, let $\left[y,s\right]_I$ denote the $\sim_I$-equivalence class of $(y,s)$ in $\mathcal{E}_x^{(I)}$.
For subsets $X,Y \subseteq \mathcal{E}_x^{(I)}$ and $1 \leq k<\infty$, write $X \overset{k}{\sim}_I Y$ to signify that $\xi \overset{k}{\sim}_I \zeta$ for some $\xi \in X$ and $\zeta \in Y$.
Moreover, write $X \overset{\infty}{\sim}_I Y$ to signify that $X \not\overset{k}{\sim}_I Y$ for any $k<\infty$.
Then Theorem \ref{thm:presentationofWperpx} implies that
\[
r_x^{(I)}(c;y,s)r_x^{(I)}(d;z,t) \textrm{ has infinite order if } (z,t) \not\in \left[y,s\right]_I \textrm{ and } \left[y,s\right]_I \overset{\infty}{\sim}_I \left[z,t\right]_I.
\]

We define an auxiliary graph $\mathcal{I}_x^{(I)}$ with vertex set $\mathcal{E}_x^{(I)}$, introduced in {\cite[Section 4.4]{Nui_centra}}.
Let $(y,s),(z,t) \in \mathcal{E}_x^{(I)}$.
Then for each subgraph of $(\Gamma_I^{\mathrm{odd}})_{\sim x}$ of type IV or V in Figure \ref{fig:localsubgraph} (which are the ones with $k=1$), we draw an edge $\widetilde{e} \in \vec{E}_{(y,s),(z,t)}(\mathcal{I}_x^{(I)})$ of $\mathcal{I}_x^{(I)}$; we denote the specified path $q \in \mathscr{F}_{y,z}^{(I)}$ by $\iota(\widetilde{e})$.
Now the same edge with opposite orientation (denoted by $\widetilde{e}^{-1}$) is also drawn.
We identify every $\widetilde{e}^{-1}$ with $\widetilde{e}$; this identification turn $\mathcal{I}_x^{(I)}$ into an edge-unoriented graph.
Note that the vertex set of the connected component $(\mathcal{I}_x^{(I)})_{\sim (y,s)}$ is the $\sim_I$-equivalence class $\left[y,s\right]_I$.
Now {\cite[Remark 4.15]{Nui_centra}} says that $\iota$ extends to an injective groupoid homomorphism from the fundamental groupoid $\pi_1(\mathcal{I}_x^{(I)};*,*)$ of $\mathcal{I}_x^{(I)}$ to $\mathscr{F}^{(I)}$, with $\iota(y,s)=y$ for $(y,s) \in \mathcal{E}_x^{(I)}$.

Now a part of Theorem \ref{thm:presentationofWperpx} can be restated as follows:
\begin{cor}
\label{cor:Wperpxandiota}
Let $c,d \in \mathscr{F}_x^{(I)}$ and $(y,s),(z,t) \in \mathcal{E}_x^{(I)}$.
Then $r_x^{(I)}(c;y,s)=r_x^{(I)}(d;z,t)$ if and only if there is a path $\widetilde{q} \in \pi_1(\mathcal{I}_x^{(I)};(y,s),(z,t))$ such that $c^{-1}d=\iota(\widetilde{q})_{(x)}$.
Hence, specializing to the case $(y,s)=(z,t)$, we have $r_x^{(I)}(c;y,s)=r_x^{(I)}(d;y,s)$ if and only if $(c^{-1}d)_{(y)} \in \iota(\pi_1(\mathcal{I}_x^{(I)};(y,s)))$.
\end{cor}
\begin{exmp}
\label{ex:misevenneq2}
Let $(y,s) \in \mathcal{E}_x^{(I)}$ such that $m_{y,s} \geq 4$.
Then $(y,s)$ has no neighbor in $\mathcal{I}_x^{(I)}$ (see Figure \ref{fig:localsubgraph}), so $\left[y,s\right]_I=\{(y,s)\}$ and $\iota(\pi_1(\mathcal{I}_x^{(I)};(y,s)))=1$.
Hence, if in addition $s \in I_{\sim x}^{\mathrm{odd}}$, then we have $(s,y) \in \mathcal{E}_x^{(I)}$ and $\left[s,y\right]_I=\{(s,y)\}$ as well, so $\left[y,s\right]_I \overset{\infty}{\sim}_I \left[s,y\right]_I$ (see Figure \ref{fig:localsubgraph}).
\end{exmp}
\begin{exmp}
\label{ex:mis2}
Let $(y,s) \in \mathcal{E}_x^{(I)}$ such that $m_{y,s}=2$ and $s \not\in O=I_{\sim x}^{\mathrm{odd}}$.
Then by Figure \ref{fig:localsubgraph}, the neighbors of $(y,s)$ in the graph $\mathcal{I}_x^{(I)}$ are the vertices $(z,s)$ with $z \in O \smallsetminus y$, $m_{y,z}$ odd and $m_{z,s}=2$, where the edge joining them is of type IV in Figure \ref{fig:localsubgraph}.
Thus it follows inductively that all edges of $(\mathcal{I}_x^{(I)})_{\sim (y,s)}$ are of type IV, and we have
\[
\left[y,s\right]_I=\{(z,s) \mid z \in J_{\sim y}^{\mathrm{odd}}\} \textrm{ and } \iota(\pi_1(\mathcal{I}_x^{(I)};(y,s)))=\mathscr{F}_y^{(K)} \subseteq \mathscr{F}_y^{(I)},
\]
where $J=\{z \in O \mid m_{z,s}=2\}$ and $K=J_{\sim y}^{\mathrm{odd}}$.
\end{exmp}
\subsection{Finite irreducible components of $W_I^{\perp x}$}
\label{sec:finitecomponentofWperpx}

\begin{figure}
\centering
\begin{picture}(140,70)(0,-70)
\put(10,-50){\circle{6}}\put(50,-50){\circle{6}}\put(90,-50){\circle{6}}\put(130,-50){\circle{6}}
\put(10,-65){\hbox to0pt{\hss$1$\hss}}\put(50,-65){\hbox to0pt{\hss$2$\hss}}\put(90,-65){\hbox to0pt{\hss$3$\hss}}\put(130,-65){\hbox to0pt{\hss$4$\hss}}
\put(13,-50){\line(1,0){34}}\put(53,-50){\line(1,0){34}}\put(93,-50){\line(1,0){34}}
\put(30,-45){\hbox to0pt{\hss$m$\hss}}
\put(10,-47){\line(1,3){10}}\put(20,-17){\line(1,0){60}}\put(90,-47){\line(-1,3){10}}
\put(50,-12){\hbox to0pt{\hss$\infty$\hss}}
\end{picture}
\caption{Coxeter group of type $P(m)$}
\label{fig:typeP}
\end{figure}
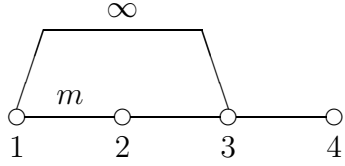

To determine the structure of the finite part $(W_I^{\perp x})_{\mathrm{fin}}$ (see Definition \ref{defn:finitepart}) in later sections, here we give some preliminary observations.
First, let $(\Pi_x^{(I)})_{\mathrm{fin}}$ be the set of all $\gamma \in \Pi_x^{(I)}$ such that the corresponding generator $s_\gamma \in R_x^{(I)}$ belongs to $(W_I^{\perp x})_{\mathrm{fin}}$.
The following property is shown in \cite{Nui_centra}.
\begin{thm}
[{\cite[Theorem 7.1]{Nui_centra}}]
\label{thm:Ycentralizefinitepart}
Any element of $Y_x^{(I)}$ fixes the set $(\Pi_x^{(I)})_{\mathrm{fin}}$ pointwise.
\end{thm}
Secondly, by definition, we have
\[
\gamma_x^{(I)}(d;y,s)=\pi^{(I)}(dc^{-1}) \cdot \gamma_x^{(I)}(c;y,s) \textrm{ for } c,d \in \mathscr{F}_x^{(I)} \textrm{ and } (y,s) \in \mathcal{E}_x^{(I)},
\]
so $r_x^{(I)}(d;y,s)$ is the conjugate of $r_x^{(I)}(c;y,s)$ by $\pi^{(I)}(dc^{-1}) \in Y_x^{(I)}$.
Thus by Theorem \ref{thm:propertyofWperpx} (\ref{thm_item:conjugateinduceisomorphism}), we have $\gamma_x^{(I)}(c;y,s) \in (\Pi_x^{(I)})_{\mathrm{fin}}$ if and only if $\gamma_x^{(I)}(d;y,s) \in (\Pi_x^{(I)})_{\mathrm{fin}}$.
Let $(\mathcal{E}_x^{(I)})_{\mathrm{fin}}$ be the set of all $(y,s) \in \mathcal{E}_x^{(I)}$ such that $\gamma_x^{(I)}(c;y,s) \in (\Pi_x^{(I)})_{\mathrm{fin}}$ for some (or equivalently, all) $c \in \mathscr{F}_x^{(I)}$.
\begin{rem}
\label{rem:Efinisindependentoftree}
Although the definition of $\gamma_x^{(I)}(c;y,s)$ depends on the choice of the maximal tree $\mathcal{T}$, the set $(\mathcal{E}_x^{(I)})_{\mathrm{fin}}$ does not.
Indeed, in the expression $\gamma_x^{(I)}(c;y,s)$ of an element of $\Pi_x^{(I)}$, any change of the tree $\mathcal{T}$ gives rise to a change of $c \in \mathscr{F}_x^{(I)}$ only, leaving $(y,s) \in \mathcal{E}_x^{(I)}$ unchanged.
\end{rem}
\begin{rem}
\label{rem:Efinisindependentofx}
By definition, if $y \in I_{\sim x}^{\mathrm{odd}}$, then we have $\mathcal{E}_y^{(I)}=\mathcal{E}_x^{(I)}$; and the isomorphism $W_I^{\perp x} \to W_I^{\perp y}$, induced by the conjugation by an element $\pi^{(I)}(p_{y,x}^{(I)})$ of $Y_{y,x}^{(I)}$, maps $r_x^{(I)}(z,s) \in R_x^{(I)}$ to $r_y^{(I)}(z,s) \in R_y^{(I)}$ (see Theorem \ref{thm:propertyofWperpx} (\ref{thm_item:conjugateinduceisomorphism})).
Thus we have $r_x^{(I)}(z,s) \in (W_I^{\perp x})_{\mathrm{fin}}$ if and only if $r_y^{(I)}(z,s) \in (W_I^{\perp y})_{\mathrm{fin}}$, hence $(\mathcal{E}_y^{(I)})_{\mathrm{fin}}=(\mathcal{E}_x^{(I)})_{\mathrm{fin}}$.
\end{rem}
Since $\pi^{(I)}(c) \in Y_x^{(I)}$ for all $c \in \mathscr{F}_x^{(I)}$, Theorem \ref{thm:Ycentralizefinitepart} yields the following:
\begin{prop}
\label{prop:gammaisindependentofc}
If $(y,s) \in (\mathcal{E}_x^{(I)})_{\mathrm{fin}}$, then $\gamma_x^{(I)}(c;y,s)=\gamma_x^{(I)}(y,s)$ for all $c \in \mathscr{F}_x^{(I)}$.
Moreover, by Corollary \ref{cor:Wperpxandiota}, the conclusion is equivalent to the following:
\begin{equation}
\label{eq:gammaisindependentofc}
\textrm{the homomorphism } \iota:\pi_1(\mathcal{I}_x^{(I)};(y,s)) \to \mathscr{F}_y^{(I)} \textrm{ is surjective}.
\end{equation}
\end{prop}
For example, (\ref{eq:gammaisindependentofc}) is satisfied if $(\Gamma_I^{\mathrm{odd}})_{\sim x}$ is acyclic (since now $\mathscr{F}_x^{(I)}=1$).
In the case of (\ref{eq:gammaisindependentofc}), we have the following lemma.
\begin{lem}
\label{lem:ifgammaisindependentofc}
Suppose that $(y,s) \in \mathcal{E}_x^{(I)}$ satisfies the condition (\ref{eq:gammaisindependentofc}).
Let $c' \in \mathscr{F}_x^{(I)}$ and $(y',s') \in \mathcal{E}_x^{(I)}$.
\begin{enumerate}
\item \label{lem_item:whengammaisindependentofc_generator}
We have $r_x^{(I)}(c';y',s')=r_x^{(I)}(y,s)$ if and only if $(y',s') \in \left[y,s\right]_I$.
Hence (by Proposition \ref{prop:gammaisindependentofc}) $(y',s')$ also satisfies (\ref{eq:gammaisindependentofc}) if $(y',s') \in \left[y,s\right]_I$.
\item \label{lem_item:whengammaisindependentofc_relation}
For $2 \leq k<\infty$, the product $r_x^{(I)}(y,s)r_x^{(I)}(c';y',s')$ has order $k$ if and only if $\left[y',s'\right]_I \overset{k}{\sim}_I \left[y,s\right]_I$.
\end{enumerate}
\end{lem}
\begin{proof}
\textbf{(\ref{lem_item:whengammaisindependentofc_generator})} The ``only if'' part follows immediately from Theorem \ref{thm:presentationofWperpx}.
For the other part, we can take a path $\widetilde{q} \in \vec{E}_{(y',s'),(y,s)}(\mathcal{I}_x^{(I)})$.
Then Corollary \ref{cor:Wperpxandiota} shows that $r_x^{(I)}(c';y',s')=r_x^{(I)}(c'\iota(\widetilde{q})_{(x)};y,s)$, while (\ref{eq:gammaisindependentofc}) implies that $r_x^{(I)}(c'\iota(\widetilde{q})_{(x)};y,s)=r_x^{(I)}(y,s)$ (see Proposition \ref{prop:gammaisindependentofc}), proving the claim.

\noindent
\textbf{(\ref{lem_item:whengammaisindependentofc_relation})} The ``only if'' part also follows from Theorem \ref{thm:presentationofWperpx}.
For the other part, take $(z,t) \in \left[y,s\right]_I$ and $(z',t') \in \left[y',s'\right]_I$ such that $(z,t) \overset{k}{\sim}_I (z',t')$.
Then by an argument similar to (\ref{lem_item:whengammaisindependentofc_generator}), we have $r_x^{(I)}(c';y',s')=r_x^{(I)}(d_1;z',t')$ for some $d_1 \in \mathscr{F}_x^{(I)}$, the product $r_x^{(I)}(d_1;z',t')r_x^{(I)}(d_2;z,t)$ has order $k$ for some $d_2 \in \mathscr{F}_x^{(I)}$, and $r_x^{(I)}(d_2;z,t)=r_x^{(I)}(d_3;y,s)$ for some $d_3 \in \mathscr{F}_x^{(I)}$ (see Theorem \ref{thm:presentationofWperpx}).
This yields the conclusion, since now $r_x^{(I)}(y,s)=r_x^{(I)}(d_3;y,s)$.
\end{proof}
\begin{cor}
\label{cor:necessaryconditionforfinitecomponent}
Let $(y,s) \in (\mathcal{E}_x^{(I)})_{\mathrm{fin}}$.
Then $(y,s)$ satisfies (\ref{eq:gammaisindependentofc}), and for any $(z,t) \in \mathcal{E}_x^{(I)}$, we have either $(z,t) \in \left[y,s\right]_I$ or $\left[z,t\right]_I \overset{k}{\sim}_I \left[y,s\right]_I$ for some $2 \leq k<\infty$.
\end{cor}
\begin{proof}
The former claim is shown in Proposition \ref{prop:gammaisindependentofc}.
For the latter, since $r_x^{(I)}(y,s)$ belongs to a finite irreducible component of $W_I^{\perp x}$, we have either $r_x^{(I)}(z,t)=r_x^{(I)}(y,s)$, or $r_x^{(I)}(z,t)r_x^{(I)}(y,s)$ has finite order.
Thus the claim follows from Lemma \ref{lem:ifgammaisindependentofc}.
\end{proof}
\begin{cor}
\label{cor:inthesamecomponentofWperpx}
Let $(y,s),(z,t) \in \mathcal{E}_x^{(I)}$, and suppose that either $\left[y,s\right]_I=\left[z,t\right]_I$ or $\left[y,s\right]_I \overset{k}{\sim}_I \left[z,t\right]_I$ for some $3 \leq k<\infty$.
Then $(y,s) \in (\mathcal{E}_x^{(I)})_{\mathrm{fin}}$ if and only if $(z,t) \in (\mathcal{E}_x^{(I)})_{\mathrm{fin}}$.
\end{cor}
\begin{proof}
Suppose that $(y,s) \in (\mathcal{E}_x^{(I)})_{\mathrm{fin}}$, so $(y,s)$ satisfies (\ref{eq:gammaisindependentofc}) by Proposition \ref{prop:gammaisindependentofc}.
Then by the hypothesis, Lemma \ref{lem:ifgammaisindependentofc} implies that either $r_x^{(I)}(z,t)=r_x^{(I)}(y,s)$, or $r_x^{(I)}(z,t)r_x^{(I)}(y,s)$ has order $k$.
In any case, $r_x^{(I)}(z,t)$ and $r_x^{(I)}(y,s)$ belong to the same irreducible component of $W_I^{\perp x}$, showing that $(z,t) \in (\mathcal{E}_x^{(I)})_{\mathrm{fin}}$.
Thus the claim follows by symmetry.
\end{proof}
\begin{cor}
\label{cor:miseven}
Let $(y,s) \in (\mathcal{E}_x^{(I)})_{\mathrm{fin}}$, and put $O=I_{\sim x}^{\mathrm{odd}}$.
If $m_{y,s} \geq 4$, then $s \not\in O$ and $\Gamma_O^{\mathrm{odd}}$ is acyclic.
\end{cor}
\begin{proof}
Note that $(y,s)$ satisfies (\ref{eq:gammaisindependentofc}) by Proposition \ref{prop:gammaisindependentofc}.
Thus the claim follows from Example \ref{ex:misevenneq2} and Corollary \ref{cor:necessaryconditionforfinitecomponent}.
\end{proof}
\begin{exmp}
\label{ex:WperpxforAn}
Let $I=(x_1,x_2,\dots)$ be a (finite or infinite) sequence of type $A_n$, with $3 \leq n \leq \infty$ (see (\ref{eq:ruleforfinitepart}) and Figure \ref{fig:nonf.g.Coxetergroup} for terminology).
Put $\xi_{i,j}=(x_i,x_j)$ for simplicity.
Then for every $x_i$, all elements of the set $\mathcal{E}_{x_i}^{(I)}$ satisfy the condition (\ref{eq:gammaisindependentofc}), and $\mathcal{E}_{x_i}^{(I)}$ consists of distinct $\sim_I$-equivalence classes $\left[\xi_{1,j}\right]_I$ with $3 \leq j \leq n$ (or $3 \leq j<\infty$ if $n=\infty$); indeed, the relations $\overset{1}{\sim}_I$ in $\mathcal{E}_{x_i}^{(I)}$ are precisely the followings:
\begin{itemize}
\item $\xi_{j,k} \overset{1}{\sim}_I \xi_{j+1,k}$, with $k \leq j-2$ or $k \geq j+3$; and
\item $\xi_{j,j+2} \overset{1}{\sim}_I \xi_{j+2,j}$
\end{itemize}
(see Figure \ref{fig:localsubgraph}), so we have
\[
\xi_{1,j} \overset{1}{\sim}_I \xi_{2,j} \overset{1}{\sim}_I \cdots \overset{1}{\sim}_I \xi_{j-2,j} \overset{1}{\sim}_I \xi_{j,j-2} \overset{1}{\sim}_I \xi_{j+1,j-2} \overset{1}{\sim}_I \cdots.
\]
Moreover, we have $\left[\xi_{1,j}\right]_I \overset{3}{\sim}_I \left[\xi_{1,j+1}\right]_I$ for all $j \geq 3$, since $\xi_{1,j} \overset{3}{\sim}_I \xi_{1,j+1}$; while $\left[\xi_{1,j}\right]_I \overset{2}{\sim}_I \left[\xi_{1,k}\right]_I$ if $|j-k| \geq 2$, since $\xi_{1,j} \overset{2}{\sim}_I \xi_{1,k}$ (see also Figure \ref{fig:localsubgraph}).
Thus by Lemma \ref{lem:ifgammaisindependentofc}, the Coxeter group $W_I^{\perp x_i}$ is generated by the distinct generators $r_{x_i}^{(I)}(\xi_{1,j})$ with $3 \leq j \leq n$, and $W_I^{\perp x_i}$ is of type $A_{n-2}$ (or $A_\infty$ if $n=\infty$).
\end{exmp}
\begin{exmp}
\label{ex:WperpxforH4}
Let $I=(x_1,x_2,x_3,x_4)$ be a sequence of type $H_4$ (see (\ref{eq:ruleforfinitepart})).
Put $\xi_{i,j}=(x_i,x_j)$.
Then for every $x_i$, we have
\begin{eqnarray*}
\mathcal{E}_{x_i}^{(I)}=\left[\xi_{4,2}\right]_I \sqcup \left[\xi_{1,3}\right]_I \sqcup \left[\xi_{1,4}\right]_I,\ \left[\xi_{4,2}\right]_I=\{\xi_{4,2}\},\ \left[\xi_{1,3}\right]_I=\{\xi_{1,3},\xi_{3,1},\xi_{4,1}\},\\
\left[\xi_{1,4}\right]_I=\{\xi_{1,4},\xi_{2,4}\},\ \left[\xi_{4,2}\right]_I \overset{3}{\sim}_I \left[\xi_{1,3}\right]_I \overset{5}{\sim}_I \left[\xi_{1,4}\right]_I \textrm{ and } \left[\xi_{4,2}\right]_I \overset{2}{\sim}_I \left[\xi_{1,4}\right]_I
\end{eqnarray*}
(note that $\xi_{4,2} \overset{3}{\sim}_I \xi_{4,1}$, $\xi_{1,3} \overset{5}{\sim}_I \xi_{1,4}$ and $\xi_{4,2} \overset{2}{\sim}_I \xi_{2,4}$ by Figure \ref{fig:localsubgraph}).
Thus, by a similar argument to Example \ref{ex:WperpxforAn}, the Coxeter group $W_I^{\perp x_i}$ is generated by three generators $r_{x_i}^{(I)}(\xi_{4,2})$, $r_{x_i}^{(I)}(\xi_{1,3})$ and $r_{x_i}^{(I)}(\xi_{1,4})$, and $W_I^{\perp x_i}$ is of type $H_3$.
\end{exmp}
\begin{exmp}
\label{ex:WperpxforP}
We say that a sequence $(x_1,x_2,x_3,x_4)$ of distinct elements of $S$ is of \emph{type $P(m)$} (with $m$ odd) if $m_{1,2}=m$, $m_{2,3}=m_{3,4}=3$, $m_{1,3}=\infty$ and $m_{1,4}=m_{2,4}=2$, where $m_{i,j}=m_{x_i,x_j}$ (see Figure \ref{fig:typeP}).
In this case, we also say that $I=\{x_1,\dots,x_4\}$, $W_I$, $(W_I,I)$ or $\Gamma_I$ is of type $P(m)$.

Now let $I=(x_1,x_2,x_3,x_4)$ be a sequence of type $P(m)$, and put $\xi_{i,j}=(x_i,x_j)$.
Then for every $x_i$, all elements of $\mathcal{E}_{x_i}^{(I)}$ satisfy (\ref{eq:gammaisindependentofc}), and we have
\begin{eqnarray*}
\mathcal{E}_{x_i}^{(I)}=\{\xi_{1,4},\xi_{2,4},\xi_{4,1},\xi_{4,2}\}=\left[\xi_{1,4}\right]_I \sqcup \left[\xi_{4,1}\right]_I,\\
\left[\xi_{1,4}\right]_I=\{\xi_{1,4},\xi_{2,4},\xi_{4,1}\},\ \left[\xi_{4,2}\right]_I=\{\xi_{4,2}\} \textrm{ and } \left[\xi_{1,4}\right]_I \overset{m}{\sim}_I \left[\xi_{4,1}\right]_I
\end{eqnarray*}
(note that $\xi_{4,2} \overset{m}{\sim}_I \xi_{4,1}$ by Figure \ref{fig:localsubgraph}).
Thus $W_I^{\perp x_i}$ is generated by two generators $r_{x_i}^{(I)}(\xi_{1,4})$ and $r_{x_i}^{(I)}(\xi_{4,1})$, and $W_I^{\perp x_i}$ is of type $I_2(m)$ (so it is finite).
\end{exmp}

In our argument below, we often derive some information on the set $(\mathcal{E}_x^{(I)})_{\mathrm{fin}}$ from that on another set $(\mathcal{E}_x^{(J)})_{\mathrm{fin}}$ with $x \in J \subseteq I$.
This is based on the following observation; note that $\mathcal{E}_x^{(J)} \subseteq \mathcal{E}_x^{(I)}$ whenever $x \in J \subseteq I$.
\begin{lem}
\label{lem:restrictionoffinitepart}
If $x \in J \subseteq I \subseteq S$, then $(\mathcal{E}_x^{(I)})_{\mathrm{fin}} \cap \mathcal{E}_x^{(J)} \subseteq (\mathcal{E}_x^{(J)})_{\mathrm{fin}}$; hence $\mathcal{E}_x^{(J)} \smallsetminus (\mathcal{E}_x^{(J)})_{\mathrm{fin}} \subseteq \mathcal{E}_x^{(I)} \smallsetminus (\mathcal{E}_x^{(I)})_{\mathrm{fin}}$.
\end{lem}
\begin{proof}
By Remark \ref{rem:Efinisindependentoftree}, we may choose the maximal tree $\mathcal{T}$ in $(\Gamma_I^{\mathrm{odd}})_{\sim x}$ such that its full subgraph with vertex set $J_{\sim x}^{\mathrm{odd}}$ is also a maximal tree in $(\Gamma_J^{\mathrm{odd}})_{\sim x}$.
Now if $(y,s) \in (\mathcal{E}_x^{(I)})_{\mathrm{fin}} \cap \mathcal{E}_x^{(J)}$, then we have $\gamma_x^{(I)}(y,s)=\gamma_x^{(J)}(y,s)$ (since $p_{y,x}^{(I)}=p_{y,x}^{(J)}$ by the choice of $\mathcal{T}$), so $\gamma_x^{(J)}(y,s) \in (\Pi_x^{(I)})_{\mathrm{fin}} \cap \Pi_x^{(J)}$.
Moreover, since $R_x^{(J)} \subseteq R_x^{(I)}$ (by Theorem \ref{thm:propertyofWperpx} (\ref{thm_item:inclusionofWperp})), the irreducible component of $W_J^{\perp x}$ containing the generator $r_x^{(J)}(y,s)$ is contained in that of $W_I^{\perp x}$ containing the same $r_x^{(J)}(y,s)$.
Now the latter irreducible component is finite (since $\gamma_x^{(J)}(y,s) \in (\Pi_x^{(I)})_{\mathrm{fin}}$), so the former one is also finite, showing that $r_x^{(J)}(y,s) \in (W_J^{\perp x})_{\mathrm{fin}}$ and $(y,s) \in (\mathcal{E}_x^{(J)})_{\mathrm{fin}}$.
Hence the claim holds.
\end{proof}
Finally, we prepare the following lemma.
\begin{lem}
\label{lem:deniedtriangle}
Suppose that $(y,s) \in (\mathcal{E}_x^{(I)})_{\mathrm{fin}}$.
Then there does not exist an element $t \in I$ such that $m_{y,t}$ is even and $m_{s,t}=\infty$.
\end{lem}
\begin{proof}
We may assume that $y=x$ (see Remark \ref{rem:Efinisindependentofx}).
Assume contrary that such $t$ exists, and put $J=\{y,s,t\}$.
Then $\mathcal{E}_y^{(J)}=\{(y,s),(y,t)\}$; and by Figure \ref{fig:localsubgraph}, we have $(y,s) \overset{\infty}{\sim}_J (y,t)$, so $\left[y,s\right]_J=\{(y,s)\}$, $\left[y,t\right]_J=\{(y,t)\}$ and $\left[y,s\right]_J \overset{\infty}{\sim}_J \left[y,t\right]_J$.
This contradicts Corollary \ref{cor:necessaryconditionforfinitecomponent}, so the claim follows.
\end{proof}
\section{The structure of $(W^{\perp x})_{\mathrm{fin}}$ for special cases}
\label{sec:mainresult_specialcase}

In this section, we determine the structure of the finite part $(W_I^{\perp x})_{\mathrm{fin}}$ of the Coxeter group $W_I^{\perp x}$ for certain special cases.
The results in this section are of independent interest, and also important steps for the general cases considered in Section \ref{sec:finitepart}.
\subsection{The $2$-spherical case}
\label{sec:specialcase_2spherical}

We start with the following lemma.
We say that a subset $I$ of $S$ is \emph{$2$-spherical} if $m_{s,t}<\infty$ for all $s,t \in I$.
Recall the graph-theoretic fact that, if $\Gamma_I$ is connected and its subgraph $\Gamma_I^{\mathrm{odd}}$ is acyclic, then $\Gamma_I$ possesses a maximal tree which contains $\Gamma_I^{\mathrm{odd}}$.
Recall also the definition of the tree order on a tree, introduced in Section \ref{sec:graphsandgroupoids}.
\begin{lem}
\label{lem:2spherical_conditionfortypeA}
Let $I$ be an irreducible, $2$-spherical subset of $S$ with at least three elements, such that $\Gamma_I^{\mathrm{odd}}$ is acyclic (so (\ref{eq:gammaisindependentofc}) is satisfied by any $(y,s) \in \mathcal{E}_x^{(I)}$).
Let $T$ be a maximal tree of $\Gamma_I$ containing $\Gamma_I^{\mathrm{odd}}$, $x$ a vertex of $T$ with degree $1$, and $x'$ the unique neighbor of $x$ in $T$.
\begin{enumerate}
\item \label{lem_item:2sp_cforA_root}
If $(y,s) \in \mathcal{E}_x^{(I)}$ and $x' \not\in \mathrm{supp}(\gamma_x^{(I)}(y,s))$, then $(y,s) \in \left[x,t\right]_I$ for some $t \in I \smallsetminus \{x,x'\}$.
Hence if $\Pi_x^{(I)} \subseteq \Phi_{I \smallsetminus x'}$, then $\mathcal{E}_x^{(I)}=\bigcup_{t \in I \smallsetminus \{x,x'\}}\left[x,t\right]_I$.
\item \label{lem_item:2sp_cforA_wholeset}
If $\mathcal{E}_x^{(I)}=\bigcup_{t \in I \smallsetminus \{x,x'\}}\left[x,t\right]_I$, then $I$ is of type $A_n$ with $3 \leq n \leq \infty$ (see Figure \ref{fig:nonf.g.Coxetergroup} for definition of type $A_\infty$).
\end{enumerate}
\end{lem}
\begin{proof}
\textbf{(\ref{lem_item:2sp_cforA_root})}
The second claim follows from the first one.
For the first claim, the hypothesis means that $\gamma_x^{(I)}(y,s) \in \Pi_x^{(I)} \cap \Phi_J=\Pi_x^{(J)}$ where $J=I \smallsetminus x'$ (see Theorem \ref{thm:propertyofWperpx} (\ref{thm_item:inclusionofWperp})).
Now we have $J_{\sim x}^{\mathrm{odd}}=\{x\}$ by the choice of $x'$, so $\mathcal{E}_x^{(J)}=\{(x,t) \mid t \in J \smallsetminus x\}$ (since $I$ is $2$-spherical) and $\gamma_x^{(I)}(y,s)=\gamma_x^{(J)}(x,t)$ for some $t \in J \smallsetminus x$.
Since $\gamma_x^{(J)}(x,t)=\widetilde{\alpha}_{x,t}=\gamma_x^{(I)}(x,t)$, Theorem \ref{thm:presentationofWperpx} implies that $(y,s) \sim_I (x,t)$ and $(y,s) \in \left[x,t\right]_I$, as desired.

\noindent
\textbf{(\ref{lem_item:2sp_cforA_wholeset})}
First, we show that $m_{y,z}$ is odd if $y$ and $z$ are adjacent vertices of the tree $T$, yielding $T=\Gamma_I^{\mathrm{odd}}$.
Assume contrary that $m_{y,z}$ is even (recall that $I$ is $2$-spherical); by symmetry, we may assume further that $y \prec z$ where $\preceq$ denote the tree order on $T$ with root vertex $x$.
Moreover, by taking such a pair $y,z$ with $y$ as minimal (with respect to $\preceq$) as possible, we may assume in addition that the path in $T$ from $x$ to $y$ consists of edges with odd labels, namely $y \in I_{\sim x}^{\mathrm{odd}}$.
This implies that $(y,z) \in \mathcal{E}_x^{(I)}$, while $m_{y,z} \neq 2$, so we have $\left[y,z\right]_I=\{(y,z)\}$ (see Example \ref{ex:misevenneq2}).
However, now $(y,z)$ cannot belong to any $\left[x,t\right]_I$ with $t \in I \smallsetminus \{x,x'\}$, since $y=x$ implies that $z=x'$.
This contradiction proves the claim.

Secondly, we show that $m_{y,z}$ is odd or $2$ for any $y,z \in I$, yielding $\Gamma_I=\Gamma_I^{\mathrm{odd}}$.
If $m_{y,z}$ is even and $m_{y,z} \neq 2$, and we assume by symmetry that $y \neq x$, then we have $(y,z) \in \mathcal{E}_x^{(I)}$ (since $y \in I_{\sim x}^{\mathrm{odd}}$ as above) and $\left[(y,z)\right]_I=\{(y,z)\}$ (see Example \ref{ex:misevenneq2}), contradicting the hypothesis.
Thus this claim follows.

We show that any pair of elements of the poset $T$ are comparable, therefore $T$ is a (finite or infinite) chain $x_0 \prec x_1 \prec x_2 \prec \cdots$, with $x_0=x$ and $x_1=x'$.
If this claim fails, then $T$ possesses two incomparable elements $y,z$ both of which cover their meet $y \wedge z$.
Note that $(y,z) \in \mathcal{E}_x^{(I)}$.
Now by Figure \ref{fig:localsubgraph}, it is straightforward to verify that the set
\[
\{(y',z) \mid y' \in \vee_y\} \cup \{(z',y) \mid z' \in \vee_z\}
\]
is a nonempty subset of $\mathcal{E}_x^{(I)}$ which is closed under the relation $\overset{1}{\sim}_I$.
This means that the subset is the union of at least one $\sim_I$-equivalence classes; however, this contradicts the hypothesis, since the subset contains no element of the form $(x,t)$.
Thus this claim holds.

Finally, we show that $m_{x_{i-1},x_i}=3$ for any $i$ with $x_{i-1},x_i \in T$, concluding the proof.
To see this, put $j=2$ if $i=1$, and $j=i$ if $i \geq 2$.
Note that $(x_j,x_{j-2}) \in \mathcal{E}_x^{(I)}$ (since $|I| \geq 3$).
Now if $m_{x_{i-1},x_i} \neq 3$, then it is also straightforward to verify that
\[
\left[x_j,x_{j-2}\right]_I=\{(x_k,x_{j-2}) \mid k \geq j \textrm{ and } x_k \in T\},
\]
which contains no element of the form $(x,t)$, contradicting the hypothesis.
Hence we have $m_{x_{i-1},x_i}=3$ as desired.
\end{proof}
The main result in this subsection is the following.
\begin{thm}
\label{thm:specialcase_2spherical_acyclic}
Let $x \in I \subseteq S$, and suppose that $I$ is irreducible, $2$-spherical and not of finite type, with $\Gamma_I^{\mathrm{odd}}$ acyclic.
Then $(W_I^{\perp x})_{\mathrm{fin}} \neq 1$ if and only if one of the following two conditions is satisfied (up to labelling of elements):
\begin{enumerate}
\item \label{thm_item:specialcase_2sp_ac_B}
$I=(x_0,x_1,\dots)$ is of type $B_\infty$ and $x \neq x_0$ (see Figure \ref{fig:nonf.g.Coxetergroup}); in this case, we have $W_I^{\perp x} \simeq W(A_1) \times W(B_\infty)$, $(\mathcal{E}_x^{(I)})_{\mathrm{fin}}=\left[x_1,x_0\right]_I=\{(x_1,x_0)\}$ and $(W_I^{\perp x})_{\mathrm{fin}}$ is generated by a single generator
\[
r_x^{(I)}(x_1,x_0)=wr_{x_1}^{(I)}(x_1,x_0)w^{-1}=wx_0x_1x_0w^{-1} \textrm{ where } w=\pi^{(I)}(p_{x,x_1}^{(I)}).
\]
\item \label{thm_item:specialcase_2sp_ac_D}
$I=(x'_0,x''_0,x_1,\dots)$ is of type $D_\infty$ (see Figure \ref{fig:nonf.g.Coxetergroup}); in this case, we have $W_I^{\perp x} \simeq W(A_1) \times W(D_\infty)$, $(\mathcal{E}_x^{(I)})_{\mathrm{fin}}=\left[x'_0,x''_0\right]_I=\{(x'_0,x''_0),(x''_0,x'_0)\}$ and $(W_I^{\perp x})_{\mathrm{fin}}$ is generated by a single generator
\[
r_x^{(I)}(x'_0,x''_0)=wr_{x'_0}^{(I)}(x'_0,x''_0)w^{-1}=wx''_0w^{-1} \textrm{ where } w=\pi^{(I)}(p_{x,x'_0}^{(I)}).
\]
\end{enumerate}
\end{thm}
\begin{proof}
The `if' part is verified by direct computation.
In the case \ref{thm_item:specialcase_2sp_ac_B}, we have
\[
\mathcal{E}_x^{(I)}=\left[x_1,x_0\right]_I \sqcup \left[x_2,x_0\right]_I \sqcup \bigsqcup_{i=3}^{\infty}\left[x_1,x_i\right]_I \textrm{ and } \left[x_1,x_0\right]_I=\{(x_1,x_0)\},
\]
while $\left[x_2,x_0\right]_I \overset{4}{\sim}_I \left[x_1,x_3\right]_I$, $\left[x_1,x_i\right]_I \overset{3}{\sim}_I \left[x_1,x_{i+1}\right]_I$ (for $i \geq 3$), and any other pair is related via $\overset{2}{\sim}_I$.
This yields the description of $W_I^{\perp x}$.
On the other hand, in the case \ref{thm_item:specialcase_2sp_ac_D}, we have
\[
\mathcal{E}_x^{(I)}=\left[x'_0,x''_0\right]_I \sqcup \left[x'_0,x_2\right]_I \sqcup \left[x''_0,x_2\right]_I \sqcup \bigsqcup_{i=3}^{\infty}\left[x'_0,x_i\right]_I
\]
and $\left[x'_0,x''_0\right]_I=\{(x'_0,x''_0),(x''_0,x'_0)\}$, while $\left[x'_0,x_i\right]_I \overset{3}{\sim}_I \left[x'_0,x_{i+1}\right]_I$ (for $i \geq 2$), $\left[x''_0,x_2\right]_I \overset{3}{\sim}_I \left[x'_0,x_3\right]_I$, and any other pair is related via $\overset{2}{\sim}_I$.
This yields the description of $W_I^{\perp x}$.
Hence the `if' part is proved.

From now, we show the ``only if'' part.
Fix a finite irreducible component $G_x$ of $W_I^{\perp x}$.
For each $y \in I_{\sim x}^{\mathrm{odd}}$, put $w_y=\pi^{(I)}(p_{y,x}^{(I)}) \in Y_{y,x}^{(I)}$, where $p_{y,x}^{(I)} \in \mathscr{F}_{y,x}^{(I)}$ denotes the unique non-backtracking path between $x$ and $y$ in the acyclic graph $\Gamma_I^{\mathrm{odd}}$.
Write $G_y=w_yG_xw_y^{-1}$ and $K_y=\mathrm{supp}(\langle y \rangle \times G_y)$, so $|K_y|<\infty$ and (by Theorem \ref{thm:propertyofWperpx} (\ref{thm_item:conjugateinduceisomorphism})) $G_y$ is a finite irreducible component of $W_I^{\perp y}$.

First we show the following lemma:
\begin{lem}
\label{thm_lem:special_2sp_ac_minimalK}
In this setting, there is an element $y \in I_{\sim x}^{\mathrm{odd}}$ such that $K_y$ is of finite type and $K_z=K_y$ for all $z \in (K_y)_{\sim y}^{\mathrm{odd}}$.
\end{lem}
\begin{proof}
[Proof of Lemma \ref{thm_lem:special_2sp_ac_minimalK}.]
Put $G'_y=\langle y \rangle \times G_y$.
Since $G'_x$ is a finite subgroup of $W_I$, Theorem \ref{thm:Tits_finitesubgroup} implies that there are some $u \in W_I$ and a subset $J \subseteq I$ of finite type such that $uG'_x u^{-1} \subseteq W_J$.
Now $uxu^{-1} \in W_J \cap S^W=J^{W_J}$, so $u'ux(u'u)^{-1} \in J$ for some $u' \in W_J$, and $u'uG'_x(u'u)^{-1} \subseteq W_J$.
By replacing $u$ with $u'u$, we may assume that $uG'_x u^{-1} \subseteq W_J$ and $uxu^{-1} \in J$.

Put $y=uxu^{-1}$, so $y \in I_{\sim x}^{\mathrm{odd}}$ (by Proposition \ref{prop:oddCoxetergraph}) and $w_y^{-1}u \in Z_{W_I}(x)$.
Now $\mathscr{F}_x^{(I)}=1$ since $\Gamma_I^{\mathrm{odd}}$ is acyclic, so we have $Y_x^{(I)}=1$ (by Theorem \ref{thm:Brink_pipreserveslength}) and $w_y^{-1}u \in \langle x \rangle \times W_I^{\perp x}$, therefore $w_y^{-1}uG'_x(w_y^{-1}u)^{-1}=G'_x$ since $G'_x$ is a direct factor of $\langle x \rangle \times W_I^{\perp x}$.
Thus $W_J \supseteq uG'_x u^{-1}=w_yG'_xw_y^{-1}=G'_y$, so $K_y \subseteq J$ by definition, therefore $K_y$ is of finite type.

Moreover, we have $G'_z=w_zw_y^{-1}G'_y(w_zw_y^{-1})^{-1} \subseteq W_{K_y}$ for $z \in (K_y)_{\sim y}^{\mathrm{odd}}$ (indeed, $w_zw_y^{-1}=\pi^{(I)}(p_{z,y}^{(I)}) \in W_{K_y}$ since the path $p_{z,y}^{(I)}$ is contained in $\Gamma_{K_y}^{\mathrm{odd}}$), so $K_z \subseteq K_y$ by definition.
Thus if $K_z \neq K_y$ for some $z$, then we can replace $y$ with $z$, making $K_z$ smaller.
Since $|K_y|<\infty$, this process must terminate, with resulting element of $(K_y)_{\sim y}^{\mathrm{odd}} \subseteq I_{\sim x}^{\mathrm{odd}}$ having the desired property.
\end{proof}
So we may assume that $K=K_x$ is of finite type and $K_y=K$ for all $y \in K_{\sim x}^{\mathrm{odd}}$.
Note that $K \neq I$ since $I$ is not of finite type.
Now put
\[
\Psi=\{\gamma \in \Pi_x^{(I)} \mid s_\gamma \in G_x\} \textrm{ and } \Psi'=\Psi \cup \{\alpha_x\},
\]
so by definition of $K$, we have $\Psi' \subseteq \Phi_K^+$ and $\Psi'$ cannot be contained in $\Phi_J$ for a proper subset $J \subset K$.
Moreover, since $G_x$ is an irreducible component of $W_I^{\perp x}$, we have
\begin{equation}
\label{thm_eq:specialcase_2sp_ac_Psiisirrcomp}
\gamma \in \Psi \textrm{ whenever } \beta \in \Psi,\ \gamma \in \Pi_x^{(I)} \textrm{ and } \langle \beta, \gamma \rangle<0
\end{equation}
(the last inequality means that $s_\beta$ and $s_\gamma$ do not commute).
Now we have the following:
\begin{lem}
\label{thm_lem:special_2sp_ac_expandK}
In this setting, if $y \in K \smallsetminus x$ is not contained in but adjacent (in $\Gamma_I$) to the support $\mathrm{supp}(\gamma)$ of $\gamma \in \Pi_x^{(I)}$, then $(K \smallsetminus x) \cap \mathrm{supp}(\gamma) \neq \emptyset$.
\end{lem}
\begin{proof}
[Proof of Lemma \ref{thm_lem:special_2sp_ac_expandK}.]
Since $\Psi' \not\subseteq \Phi_{K \smallsetminus y}$, there is a root $\beta \in \Psi$ with $y \in \mathrm{supp}(\beta)$.
Now if we assume contrary that $(K \smallsetminus x) \cap \mathrm{supp}(\gamma)=\emptyset$, then the choice of $\beta$ implies that $(\mathrm{supp}(\beta) \smallsetminus x) \cap \mathrm{supp}(\gamma)=\emptyset$ (since $\Psi \subseteq \Phi_K$) and $\mathrm{supp}(\beta) \smallsetminus x$ is adjecent to $\mathrm{supp}(\gamma)$ in $\Gamma_I$.
Thus Lemma \ref{lem:conditionfornonorthogonal} shows that $\langle \beta, \gamma \rangle<0$, so $\gamma \in \Psi \subseteq \Phi_K$ by (\ref{thm_eq:specialcase_2sp_ac_Psiisirrcomp}), while $\mathrm{supp}(\gamma) \not\subseteq K$ by our assumption (note that $\mathrm{supp}(\gamma)$ must contain an element other than $x$).
This contradiction proves the claim.
\end{proof}
Take a maximal tree $T$ in the connected graph $\Gamma_I$ containing $\Gamma_I^{\mathrm{odd}}$.
Let $\preceq$ denote the tree order on $T$ with root vertex $x$, $\Lambda$ the set of atoms of the poset $T$ and $\Omega_y$ (for each $y \in \Lambda$) the set of the elements of $T$ covering $y$.
Note that $\Lambda \neq \emptyset$ since $I \neq \{x\}$.
Now we have the following observations:
\begin{lem}
\label{thm_lem:special_2sp_ac_KandT}
In this setting, we have the followings.
\begin{enumerate}
\item \label{thm_lem_item:special_2sp_ac_adjacenttoK}
If $s \in T \smallsetminus (x \cup \Lambda)$ is adjacent to some $y \in K$ in $T$, then $s \in K$.
\item \label{thm_lem_item:special_2sp_ac_height2}
For $y \in \Lambda$ and $z \in \Omega_y$, we have $\vee_z \subseteq K$ or $\vee_z \cap K=\emptyset$.
\item \label{thm_lem_item:special_2sp_ac_ifatomisin}
If $y \in \Lambda \cap K$, then $\vee_y \subseteq K$.
\item \label{thm_lem_item:special_2sp_ac_atomeven}
If $y \in \Lambda$ and $m_{x,y}$ is even, then $\vee_y \subseteq K$ or $\vee_y \cap K=\emptyset$.
\item \label{thm_lem_item:special_2sp_ac_atomodd}
If $y \in \Lambda$ and $m_{x,y}$ is odd, then $y \not\in K$.
\item \label{thm_lem_item:special_2sp_ac_simpleroot}
We have $(x,s) \in (\mathcal{E}_x^{(I)})_{\mathrm{fin}}$ and $\gamma_x^{(I)}(x,s) \in \Psi$ for all $s \in K \smallsetminus x$.
\end{enumerate}
\end{lem}
\begin{proof}
[Proof of Lemma \ref{thm_lem:special_2sp_ac_KandT}.]
\textbf{(\ref{thm_lem_item:special_2sp_ac_adjacenttoK})}
Since the tree $T$ contains $\Gamma_I^{\mathrm{odd}}$ and $I$ is $2$-spherical, we have $(x,s) \in \mathcal{E}_x^{(I)}$, $y \neq x$ and $s \in \mathrm{supp}(\gamma_x^{(I)}(x,s)) \subseteq \{x,s\}$ (recall that $\widetilde{\alpha}_{x,s}$).
Thus $(K \smallsetminus x) \cap \{x,s\} \neq \emptyset$ and $s \in K$ by Lemma \ref{thm_lem:special_2sp_ac_expandK}.

\noindent
\textbf{(\ref{thm_lem_item:special_2sp_ac_height2})}
If $\vee_z \cap K \neq \emptyset$, then repeated use of Claim \ref{thm_lem_item:special_2sp_ac_adjacenttoK} implies that $\vee_z \subseteq K$, since $\vee_z \cap (x \cup \Lambda)=\emptyset$ and $\vee_z$ is connected in $T$.

\noindent
\textbf{(\ref{thm_lem_item:special_2sp_ac_ifatomisin})}
Now Claim \ref{thm_lem_item:special_2sp_ac_adjacenttoK} implies that $\Omega_y \subseteq K$, so Claim \ref{thm_lem_item:special_2sp_ac_height2} proves the claim.

\noindent
\textbf{(\ref{thm_lem_item:special_2sp_ac_atomeven})}
Suppose that $\vee_y \cap K \neq \emptyset$.
If $y \in K$, then $\vee_y \subseteq K$ by Claim \ref{thm_lem_item:special_2sp_ac_ifatomisin}.
On the other hand, if $s \in \vee_y \cap K$ and $s \neq y$, then $z \in K$ by Claim \ref{thm_lem_item:special_2sp_ac_height2}, where $z \in \Omega_y$ and $z \preceq s$.
Now $(x,y) \in \mathcal{E}_x^{(I)}$ and $z \in K$ is adjacent to $\mathrm{supp}(\gamma_x^{(I)}(x,y))=\{x,y\}$, so $y \in K$ by Lemma \ref{thm_lem:special_2sp_ac_expandK}, therefore $\vee_y \subseteq K$ as above.

\noindent
\textbf{(\ref{thm_lem_item:special_2sp_ac_atomodd})}
Assume contrary that $y \in K$, and deduce a contradiction $K=I$ as follows.
By Claim \ref{thm_lem_item:special_2sp_ac_ifatomisin}, it suffices to show that $y' \in K$ for all $y' \in \Lambda \smallsetminus y$.
Now since $K_y=K$, we have $x \in K_y$, while $y'$ is adjacent to $x$ but not adjacent to $y$ in $T$.
Thus we have $y' \in K_y=K$ as desired, by applying Claim \ref{thm_lem_item:special_2sp_ac_ifatomisin} to the new root vertex $y$ of $T$ instead of $x$.

\noindent
\textbf{(\ref{thm_lem_item:special_2sp_ac_simpleroot})}
Since $T$ contains $\Gamma_I^{\mathrm{odd}}$, Claim \ref{thm_lem_item:special_2sp_ac_atomodd} implies that $K_{\sim x}^{\mathrm{odd}}=\{x\}$, so (since $I$ is $2$-spherical) we have $\mathcal{E}_x^{(K)}=\{(x,s) \mid s \in K \smallsetminus x\}$.
Thus for $s \in K \smallsetminus x$, the root $\gamma_x^{(I)}(x,s)=\gamma_x^{(K)}(x,s)$ is the unique element of $\Pi_x^{(I)} \cap \Phi_K=\Pi_x^{(K)}$ (see Theorem \ref{thm:propertyofWperpx} (\ref{thm_item:inclusionofWperp})) whose support contains $s$.
Since $\Psi \subseteq \Phi_K$ and $\Psi \not\subseteq \Phi_{K \smallsetminus s}$, the set $\Psi$ must contain $\gamma_x^{(I)}(x,s)$, so $r_x^{(I)}(x,s) \in G_x \subseteq (W_I^{\perp x})_{\mathrm{fin}}$, concluding the proof.
\end{proof}
We divide the remaining proof of Theorem \ref{thm:specialcase_2spherical_acyclic} into two cases.

\noindent
\textbf{Case 1: $|\Lambda| \geq 2$.}

In this case, we verify the condition \ref{thm_item:specialcase_2sp_ac_B} of the statement.

\noindent
\textbf{Step 1-1: $m_{x,y}$ is even and $y \in K$ for some $y \in \Lambda$.}

Take some $s \in K \smallsetminus x$.
If $s \in \Lambda$, then Lemma \ref{thm_lem:special_2sp_ac_KandT} (\ref{thm_lem_item:special_2sp_ac_atomodd}) shows that $s$ is the desired element.
So suppose that $s \not\in \Lambda$, and take $y \in \Lambda$ and $z \in \Omega_y$ such that $z \preceq s$.
Lemma \ref{thm_lem:special_2sp_ac_KandT} (\ref{thm_lem_item:special_2sp_ac_atomeven}) implies that $y$ is the desired element if $m_{x,y}$ is even (since now $s \in \vee_y \subseteq K$); so assume contrary that $m_{x,y}$ is odd.
Then $z \in K$ and $y \not\in K$ by Lemma \ref{thm_lem:special_2sp_ac_KandT} (\ref{thm_lem_item:special_2sp_ac_height2}) and (\ref{thm_lem_item:special_2sp_ac_atomodd}), so $\gamma_x^{(I)}(x,z) \in \Psi$ by Lemma \ref{thm_lem:special_2sp_ac_KandT} (\ref{thm_lem_item:special_2sp_ac_simpleroot}).
However, by taking $y' \in \Lambda \smallsetminus y$, we have $(y,y') \in \mathcal{E}_x^{(I)}$ and $\mathrm{supp}(\gamma_x^{(I)}(y,y'))=\{x,y,y'\}$ (see Example \ref{ex:supportofgammaxys_length2}), so $\gamma_x^{(I)}(y,y') \not\in \Psi \subseteq \Phi_K$ and $\langle \gamma_x^{(I)}(x,z), \gamma_x^{(I)}(y,y') \rangle<0$ (by Lemma \ref{lem:conditionfornonorthogonal}).
This contradicts (\ref{thm_eq:specialcase_2sp_ac_Psiisirrcomp}), so Step 1-1 is concluded.

\noindent
\textbf{Step 1-2: $m_{x,y'}$ is odd for all $y' \in \Lambda \smallsetminus y$.
Hence $\Lambda \cap K=\{y\}$ by Lemma \ref{thm_lem:special_2sp_ac_KandT} (\ref{thm_lem_item:special_2sp_ac_atomodd}) and Step 1-1.}

If $m_{x,y'}$ is even, then Example \ref{ex:misevenneq2} and Figure \ref{fig:localsubgraph} show that
\[
\left[x,y\right]_I=\{(x,y)\} \overset{\infty}{\sim}_I \{(x,y')\}=\left[x,y'\right]_I
\]
(since $m_{x,y} \neq 2$ and $m_{x,y'} \neq 2$), while $(x,y) \in (\mathcal{E}_x^{(I)})_{\mathrm{fin}}$ by Lemma \ref{thm_lem:special_2sp_ac_KandT} (\ref{thm_lem_item:special_2sp_ac_simpleroot}).
This contradicts Corollary \ref{cor:necessaryconditionforfinitecomponent}, so Step 1-2 is concluded.

\noindent
\textbf{Step 1-3: We have $\Omega_y=\emptyset$, so $\vee_y=\{y\}$.}

Take $y' \in \Lambda \smallsetminus y$, so $m_{x,y'}$ is odd and $y' \not\in K$ by Step 1-2.
Assume contrary that some $z \in \Omega_y$ exists.
Then, since $y \in K$, Lemma \ref{thm_lem:special_2sp_ac_KandT} (\ref{thm_lem_item:special_2sp_ac_atomeven}) shows that $z \in K$, so $\gamma_x^{(I)}(x,z) \in \Psi$ by Lemma \ref{thm_lem:special_2sp_ac_KandT} (\ref{thm_lem_item:special_2sp_ac_simpleroot}).
Now we have $(y',y) \in \mathcal{E}_x^{(I)}$, $\mathrm{supp}(\gamma_x^{(I)}(y',y))=\{x,y,y'\}$ (see Example \ref{ex:supportofgammaxys_length2}) and $\langle \gamma_x^{(I)}(x,z), \gamma_x^{(I)}(y',y) \rangle<0$ (by Lemma \ref{lem:conditionfornonorthogonal}), while $\gamma_x^{(I)}(y',y) \not\in \Psi$ since $y' \not\in K$.
This contradicts (\ref{thm_eq:specialcase_2sp_ac_Psiisirrcomp}), so Step 1-3 is concluded.

\noindent
\textbf{Step 1-4: If $\gamma \in \Pi_x^{(I)}$ and $y \not\in \mathrm{supp}(\gamma)$, then $\Lambda \cap \mathrm{supp}(\gamma)=\emptyset$.
Hence we have $\Pi_x^{(I \smallsetminus y)} \subseteq \Phi_{I \smallsetminus \Lambda}$.}

Assume contrary that $y' \in \Lambda \cap \mathrm{supp}(\gamma)$ (so $y' \neq y$).
Then $y' \not\in K$ (by Step 1-2) and $\gamma \not\in \Psi$; while $\gamma_x^{(I)}(x,y) \in \Psi$ by Lemma \ref{thm_lem:special_2sp_ac_KandT} (\ref{thm_lem_item:special_2sp_ac_simpleroot}), $y'$ is adjacent to $\mathrm{supp}(\gamma_x^{(I)}(x,y))=\{x,y\}$ and so $\langle \gamma_x^{(I)}, \gamma \rangle<0$ by Lemma \ref{lem:conditionfornonorthogonal}.
This contradicts (\ref{thm_eq:specialcase_2sp_ac_Psiisirrcomp}), so Step 1-4 is concluded.

\noindent
\textbf{Step 1-5: We have $|\Lambda|=2$.}

If $\Lambda \smallsetminus y$ contains two distinct elements $y',y''$, then $m_{x,y'}$ is odd by Step 1-2, so $(y',y'') \in \mathcal{E}_x^{(I)}$; however, by Example \ref{ex:supportofgammaxys_length2}, the root $\gamma_x^{(I)}(y',y'') \in \Pi_x^{(I)}$ contradicts Step 1-4.
Thus Step 1-5 is concluded.

Owing to Step 1-5, put $\Lambda=\{y,y'\}$.
Note that $m_{x,y'}$ is odd by Step 1-2.

\noindent
\textbf{Step 1-6: $I \smallsetminus y$ is of type $A_n$ for some $3 \leq n \leq \infty$.}

First, we show that $\Omega_{y'} \neq \emptyset$.
If $\Omega_{y'}=\emptyset$, then $I=\{x,y,y'\}$ and $\mathcal{E}_x^{(I)}=\{(x,y),(y',y)\}$, while $(x,y) \in (\mathcal{E}_x^{(I)})_{\mathrm{fin}}$ by Lemma \ref{thm_lem:special_2sp_ac_KandT} (\ref{thm_lem_item:special_2sp_ac_simpleroot}).
Now Corollary \ref{cor:necessaryconditionforfinitecomponent} implies that $(x,y) \overset{m}{\sim}_I (y',y)$ for some $m<\infty$; since $m_{x,y} \neq 2$, this is possible only if $m_{x,y}=4$, $m_{y',y}=2$ and $m_{x,y'}=3$ (see Figure \ref{fig:localsubgraph}), namely $I$ is of type $B_3$.
This contradicts the hypothesis, so we have $\Omega_{y'} \neq \emptyset$ and $|I \smallsetminus y| \geq 3$.
Thus by Step 1-4 and Lemma \ref{lem:2spherical_conditionfortypeA} (\ref{lem_item:2sp_cforA_root}), the set $I \smallsetminus y$, the full subgraph of $T$ with vertex set $I \smallsetminus y$, and the two elements $x$ and $y'$ satisfy the hypothesis of Lemma \ref{lem:2spherical_conditionfortypeA} (\ref{lem_item:2sp_cforA_wholeset}).
Hence the lemma proves the claim.

\noindent
\textbf{Step 1-7: $m_{x,y}=4$ and $m_{s,y}=2$ for all $s \in I \smallsetminus \{x,y\}$; hence $I$ is of type $B_{n+1}$ (if $n<\infty$) or $B_\infty$.}

First, since $(x,y) \in (\mathcal{E}_x^{(I)})_{\mathrm{fin}} \cap \mathcal{E}_x^{(J)} \subseteq (\mathcal{E}_x^{(J)})_{\mathrm{fin}}$ (where $J=\{x,y,y'\}$) and $\mathcal{E}_x^{(J)}=\{(x,y),(y',y)\}$, we have $m_{x,y}=4$ by a similar argument to Step 1-6.
Moreover, if we assume contrary that $m_{s,y} \neq 2$ for some $s \in I \smallsetminus \{x,y\}$, then we have $(s,y) \in \mathcal{E}_x^{(I)}$ and
\[
\left[x,y\right]_I=\{(x,y)\} \overset{\infty}{\sim}_I \{(s,y)\}=\left[s,y\right]_I.
\]
This contradicts Corollary \ref{cor:necessaryconditionforfinitecomponent}, since $(x,y) \in (\mathcal{E}_x^{(I)})_{\mathrm{fin}}$ by Lemma \ref{thm_lem:special_2sp_ac_KandT} (\ref{thm_lem_item:special_2sp_ac_simpleroot}).
Hence we have $m_{s,y}=2$ for all $s \in I \smallsetminus \{x,y\}$, concluding the proof.

Since $I$ is not of finite type, Step 1-7 actually shows that $I$ is of type $B_\infty$.
Hence the condition \ref{thm_item:specialcase_2sp_ac_B} in the statement is satisfied.

\noindent
\textbf{Case 2: $|\Lambda|=1$.}

Put $\Lambda=\{y\}$, so $I=x \cup \vee_y$.
In this case, we verify the condition \ref{thm_item:specialcase_2sp_ac_D} of the statement. 
Since $\emptyset \neq K \neq I$, we have $y \not\in K$ by Lemma \ref{thm_lem:special_2sp_ac_KandT} (\ref{thm_lem_item:special_2sp_ac_ifatomisin}), so $s \in K$ for some $s \in I \smallsetminus \{x,y\}$.
Take $z \in \Omega_y$ with $z \preceq s$, so $\vee_z \subseteq K$ by Lemma \ref{thm_lem:special_2sp_ac_KandT} (\ref{thm_lem_item:special_2sp_ac_height2}).

\noindent
\textbf{Step 2-1: We have $|\Omega_y| \geq 2$, so $|I \smallsetminus \vee_z| \geq 3$.}

Assume contrary that $\Omega_y=\{z\}$, so $I=\{x,y\} \cup \vee_z$ and $K=I \smallsetminus y$.
Note that $|I|=|K|+1<\infty$ since $K$ is of finite type.
First, we show that $y \not\in \mathrm{supp}(\gamma)$ for any $\gamma \in \Pi_x^{(I)}$.
If $y \in \mathrm{supp}(\gamma)$, then $\gamma \not\in \Psi \subseteq \Phi_K$, so $\langle \gamma_x^{(I)}(x,s), \gamma \rangle=0$ for all $s \in K \smallsetminus x$ (since $\gamma_x^{(I)}(x,s) \in \Psi$ by Lemma \ref{thm_lem:special_2sp_ac_KandT} (\ref{thm_lem_item:special_2sp_ac_simpleroot}); see also (\ref{thm_eq:specialcase_2sp_ac_Psiisirrcomp}) and Theorem \ref{thm:propertyofWperpx} (\ref{thm_item:orderandbilinearform})).
Since $s \in \mathrm{supp}(\gamma_x^{(I)}(x,s)) \subseteq \{x,s\}$, this means that $\langle \alpha_s, \gamma \rangle=0$ for $s \in K \smallsetminus x$, while $\langle \alpha_x, \gamma \rangle=0$.
Thus $\langle \alpha_s, \gamma \rangle=0$ for all $s \in K=I \smallsetminus y$, contradicting Proposition \ref{prop:rootformaximalparabolic} since $I$ is irreducible and not of finite type.
Hence $y \not\in \mathrm{supp}(\gamma)$ for any $\gamma \in \Pi_x^{(I)}$, so $\Pi_x^{(I)} \subseteq \Phi_{I \smallsetminus y}$.

Now Lemma \ref{lem:2spherical_conditionfortypeA} (\ref{lem_item:2sp_cforA_root}) shows that $I$, $T$, $x$ and $y$ satisfy the hypothesis of Lemma \ref{lem:2spherical_conditionfortypeA} (\ref{lem_item:2sp_cforA_wholeset}), so it follows that $I$ is of type $A_n$ with $3 \leq n<\infty$ (recall that $|I|<\infty$ as above).
This is a contradiction, since $I$ is not of finite type.
Hence Step 2-1 is concluded.

\noindent
\textbf{Step 2-2: $I \smallsetminus \vee_z$ is of type $A_n$ for some $3 \leq n \leq \infty$.}

Put $J=I \smallsetminus \vee_z$, so $|J| \geq 3$ by Step 2-1.
As well as Step 2-1, owing to Lemma \ref{lem:2spherical_conditionfortypeA}, it suffices to show that $y \not\in \mathrm{supp}(\gamma)$ for all $\gamma \in \Pi_x^{(J)} \subseteq \Pi_x^{(I)}$.
Now if $\gamma \in \Pi_x^{(J)}$ and $y \in \mathrm{supp}(\gamma)$, then we have $\gamma \not\in \Psi$ (since $y \not\in K$); while $\gamma_x^{(I)}(x,z) \in \Psi$ (by Lemma \ref{thm_lem:special_2sp_ac_KandT} (\ref{thm_lem_item:special_2sp_ac_simpleroot})) and $\langle \gamma_x^{(I)}(x,z), \gamma \rangle<0$ (by Lemma \ref{lem:conditionfornonorthogonal}), contradicting (\ref{thm_eq:specialcase_2sp_ac_Psiisirrcomp}).
Hence Step 2-2 is concluded.

In particular, we have $|\Omega_y|=2$.
Put $\Omega_y=\{z,z'\}$.

\noindent
\textbf{Step 2-3: We have $\vee_z=\{z\}$.}

Assume contrary that $\vee_z \neq \{z\}$, and take an element $s$ of $\vee_z$ covering $z$ in $T$.
Now $s \in K$ and $\gamma_x^{(I)}(x,s) \in \Psi$ (by Lemma \ref{thm_lem:special_2sp_ac_KandT} (\ref{thm_lem_item:special_2sp_ac_simpleroot})); while $z' \in I_{\sim x}^{\mathrm{odd}}$ (by Step 2-2), $\mathrm{supp}(\gamma_x^{(I)}(z',z))=\{x,y,z,z'\}$ (see Example \ref{ex:supportofgammaxys_length3}), $\gamma_x^{(I)}(z',z) \not\in \Psi$ (since $y \not\in K$) and $\langle \gamma_x^{(I)}(x,s), \gamma_x^{(I)}(z',z) \rangle<0$ (by Lemma \ref{lem:conditionfornonorthogonal}).
This contradicts (\ref{thm_eq:specialcase_2sp_ac_Psiisirrcomp}), so Step 2-3 is concluded.

\noindent
\textbf{Step 2-4: $m_{x,z}=2$, $m_{y,z}=3$ and $m_{s,z}=2$ for $s \in \vee_{z'}$; hence $I$ is of type $D_{n+1}$ (if $n<\infty$) or $D_\infty$.}

Note that $(x,z) \in (\mathcal{E}_x^{(I)})_{\mathrm{fin}}$ by Lemma \ref{thm_lem:special_2sp_ac_KandT} (\ref{thm_lem_item:special_2sp_ac_simpleroot}).
First, we show that $m_{x,z}=2$ and $m_{y,z}=3$.
If this fails, then we have $\left[x,z\right]_I=\{(x,z)\}$, while $(z',z) \in \mathcal{E}_x^{(I)}$ and so the subset
\[
\mathcal{E}'=\{(z,z')\} \cup \{(s,z) \mid s \in \vee_{z'}\}
\]
of $\mathcal{E}_x^{(I)}$ contains $(z',z)$ and is closed under $\overset{1}{\sim}_I$, therefore $\left[z',z\right]_I \subseteq \mathcal{E}'$.
However, now we have $(x,z) \not\in \mathcal{E}'$ and $(x,z) \overset{\infty}{\sim}_I \mathcal{E}'$ (indeed, if $(x,z) \overset{m}{\sim}_I (t,t')$ for some $m<\infty$, then either $(t,t')=(z,x)$ or $(t,t')=(x,t'')$ with $t'' \in \vee_{z'}$); so $(x,z) \not\in \left[z',z\right]_I$ and $\left[x,z\right]_I \overset{\infty}{\sim}_I \left[z',z\right]_I$, contradicting Corollary \ref{cor:necessaryconditionforfinitecomponent}.
Thus we have $m_{x,z}=2$ and $m_{y,z}=3$, so it follows that $\left[x,z\right]_I=\{(x,z),(z,x)\}$.

Finally, we show that $m_{s,z}=2$ for all $s \in \vee_{z'}$, concluding Step 2-4.
Note that $(s,z) \in \mathcal{E}_x^{(I)}$.
Now if $m_{s,z} \neq 2$, then we have $\left[s,z\right]_I=\{(s,z)\}$ (see Example \ref{ex:misevenneq2}), so $\left[x,z\right]_I \overset{\infty}{\sim}_I \left[s,z\right]_I$ (see Figure \ref{fig:localsubgraph}), contradicting Corollary \ref{cor:necessaryconditionforfinitecomponent}.
Thus $m_{s,z}=2$ for all $s \in \vee_{z'}$, as desired.

Since $I$ is not of finite type, Step 2-4 actually shows that $I$ is of type $D_\infty$.
Hence the condition \ref{thm_item:specialcase_2sp_ac_D} in the statement is satisfied.

Now the proof of Theorem \ref{thm:specialcase_2spherical_acyclic} is concluded.
\end{proof}
\begin{cor}
\label{cor:no2sphericalacyclicsupset}
Suppose that $x \in J \subseteq S$ and $(y,s) \in (\mathcal{E}_x^{(J)})_{\mathrm{fin}}$.
Then there does not exist a finite, irreducible, $2$-spherical subset $I \subseteq J$ such that $y,s \in I$, the graph $\Gamma_I^{\mathrm{odd}}$ is acyclic and $I$ is not of finite type.
\end{cor}
\begin{proof}
If such a subset $I$ exists, then we have $(y,s) \in (\mathcal{E}_y^{(J)})_{\mathrm{fin}} \cap \mathcal{E}_y^{(I)}$ by Theorem \ref{thm:propertyofWperpx} (\ref{thm_item:conjugateinduceisomorphism}), so $(y,s) \in (\mathcal{E}_y^{(I)})_{\mathrm{fin}}$ by Lemma \ref{lem:restrictionoffinitepart}, while $(\mathcal{E}_y^{(I)})_{\mathrm{fin}}=\emptyset$ by Theorem \ref{thm:specialcase_2spherical_acyclic}.
This contradiction proves the claim.
\end{proof}
\subsection{The case with $\Gamma^{\mathrm{odd}}$ a path}
\label{sec:specialcase_oddpath}

This subsection treats the case that the odd Coxeter graph $\Gamma_I^{\mathrm{odd}}$ is a simple non-closed path.
Namely, we show the following result.
\begin{thm}
\label{thm:Wperpxforoddpath}
Let $I \subseteq S$ and $n=|I|$.
Suppose that $3 \leq n<\infty$ and $\Gamma_I^{\mathrm{odd}}$ is a simple non-closed path $(x_1,x_2,\dots,x_n)$ with $m_{x_1,x_n}<\infty$.
Let $1 \leq i \leq n$.
\begin{enumerate}
\item \label{thm_item:Wperpxforoddpath_H3}
If $\mathrm{type}\,I=H_3$ (with $n=3$), then $W_I^{\perp x_i}$ is finite and generated by two commuting generators (namely $W_I^{\perp x_i}$ is of type $A_1 \times A_1$).
\item \label{thm_item:Wperpxforoddpath_AHP}
If $\mathrm{type}\,I=A_n$, $H_4$ or $P(m)$, with $n=4$ in the latter two cases (see Example \ref{ex:WperpxforP} for definition of type $P(m)$), then $W_I^{\perp x_i}$ is finite, irreducible and of type $A_{n-2}$, $H_3$ or $I_2(m)$, respectively.
Hence for $I \subseteq J \subseteq S$, all elements of $R_{x_i}^{(I)}$ belong to the same irreducible component of $W_J^{\perp x_i}$ (see Theorem \ref{thm:propertyofWperpx} (\ref{thm_item:inclusionofWperp})); thus for $(y,s),(z,t) \in \mathcal{E}_{x_i}^{(I)}$, we have $(y,s) \in (\mathcal{E}_{x_i}^{(J)})_{\mathrm{fin}}$ if and only if $(z,t) \in (\mathcal{E}_{x_i}^{(J)})_{\mathrm{fin}}$.
\item \label{thm_item:Wperpxforoddpath_other}
Otherwise, we have $(\mathcal{E}_{x_i}^{(I)})_{\mathrm{fin}}=\emptyset$ and $(W_I^{\perp x_i})_{\mathrm{fin}}=1$.
Hence for $I \subseteq J \subseteq S$, we have $(\mathcal{E}_{x_i}^{(J)})_{\mathrm{fin}} \cap \mathcal{E}_{x_i}^{(I)}=\emptyset$ (see Lemma \ref{lem:restrictionoffinitepart}).
\end{enumerate}
\end{thm}
\begin{proof}
Put $m_{i,j}=m_{x_i,x_j}$ and $\xi_{i,j}=(x_i,x_j)$.
We proceed the proof by induction on $n$.
Note that, since $\Gamma_I^{\mathrm{odd}}$ is connected, Remark \ref{rem:Efinisindependentofx} shows that it suffices to prove the claim for only one $x_i$.
Moreover, since $\Gamma_I^{\mathrm{odd}}$ is acyclic, any element of $\mathcal{E}_{x_i}^{(I)}$ satisfies the condition (\ref{eq:gammaisindependentofc}).

First, suppose that $n=3$.
Then $\mathcal{E}_{x_1}^{(I)}=\{\xi_{1,3},\xi_{3,1}\}$ since $m_{1,3}<\infty$.
Now a direct computation shows that $\xi_{1,3} \sim_I \xi_{3,1}$ if $\mathrm{type}\,I=A_3$; $\xi_{1,3} \overset{2}{\sim}_I \xi_{3,1}$ if $\mathrm{type}\,I=H_3$; and $\xi_{1,3} \not\sim_I \xi_{3,1}$ and $\xi_{1,3} \overset{\infty}{\sim}_I \xi_{3,1}$ otherwise.
Thus the claim follows from Lemma \ref{lem:ifgammaisindependentofc} (note that $r_{x_i}^{(I)}(\xi_{1,3})r_{x_i}^{(I)}(\xi_{3,1})$ has infinite order in the third case).

So suppose that $n \geq 4$.
If $I$ satisfies the hypothesis of Claim \ref{thm_item:Wperpxforoddpath_AHP}, then the claim follows from Examples \ref{ex:WperpxforAn}, \ref{ex:WperpxforH4} and \ref{ex:WperpxforP}.
So suppose that the hypothesis of Claim \ref{thm_item:Wperpxforoddpath_AHP} is not satisfied.
This is the case of Claim \ref{thm_item:Wperpxforoddpath_other}, so our aim is to show that $(\mathcal{E}_{x_i}^{(I)})_{\mathrm{fin}}=\emptyset$.
We divide the proof into three cases.

\noindent
\textbf{Case 1: $m_{1,n-1} \neq 2$ and $m_{2,n} \neq 2$.}

Since $n \geq 4$, this assumption implies by a direct computation that $\left[\xi_{1,n}\right]_I=\{\xi_{1,n}\}$, $\left[\xi_{n,1}\right]_I=\{\xi_{n,1}\}$ and $\xi_{1,n} \overset{\infty}{\sim}_I \xi_{n,1}$.
Thus $\xi_{1,n},\xi_{n,1} \not\in (\mathcal{E}_{x_i}^{(I)})_{\mathrm{fin}}$ by Corollary \ref{cor:necessaryconditionforfinitecomponent}.
On the other hand, if $m_{1,n-1}<\infty$, then Corollary \ref{cor:miseven} implies that $\xi_{1,n-1},\xi_{n-1,1} \not\in (\mathcal{E}_{x_i}^{(I)})_{\mathrm{fin}}$ (since $m_{1,n-1} \neq 2$ and $I_{\sim {x_i}}^{\mathrm{odd}}=I$).
By symmetry, a similar result follows for $\xi_{2,n}$ and $\xi_{n,2}$.

By these and symmetry, it suffices to show that $\xi_{j,k} \not\in (\mathcal{E}_{x_i}^{(I)})_{\mathrm{fin}}$ for indices $j<k$ such that $\xi_{j,k} \in \mathcal{E}_{x_i}^{(I)} \smallsetminus \{\xi_{1,n},\xi_{1,n-1},\xi_{2,n}\}$.
Assume contrary that $\xi_{j,k} \in (\mathcal{E}_{x_i}^{(I)})_{\mathrm{fin}}$.
Then, since $\left[\xi_{n,1}\right]_I=\{\xi_{n,1}\}$, Corollary \ref{cor:necessaryconditionforfinitecomponent} implies that $\left[\xi_{j,k}\right]_I \overset{m}{\sim}_I \xi_{n,1}$ for some $2 \leq m<\infty$; so a direct computation shows that $\left[\xi_{j,k}\right]_I$ contains either $\xi_{n,\ell}$ (with $2 \leq \ell \leq n-2$) or $\xi_{n-1,1}$.
Since $\xi_{n-1,1} \not\in (\mathcal{E}_{x_i}^{(I)})_{\mathrm{fin}}$, Corollary \ref{cor:inthesamecomponentofWperpx} denies the latter possibility, so $\left[\xi_{j,k}\right]_I$ contains some $\xi_{n,\ell}$.
This is possible only if all $m_{p,k}$ (with $j \leq p \leq k-2$) and $m_{k-2,p}$ (with $k \leq p \leq n$) are $2$ and $m_{k-2,k-1}=m_{k-1,k}=3$.
(Indeed, by a similar argument to Example \ref{ex:WperpxforAn}, the connected component $(\mathcal{I}_{x_i}^{(I)})_{\sim \xi_{j,k}}$ of the graph $\mathcal{I}_{x_i}^{(I)}$ containing $\xi_{j,k}$ is a full subgraph of the path
\[
(\xi_{1,k},\xi_{2,k},\dots,\xi_{j,k},\dots,\xi_{k-2,k},\xi_{k,k-2},\xi_{k+1,k-2},\dots,\xi_{n,k-2});
\]
and if some of these conditions fails, then some of the above edges between $\xi_{j,k}$ and $\xi_{n,k-2}$ does not exist and so $(\mathcal{I}_{x_i}^{(I)})_{\sim \xi_{j,k}}$ cannot achieve to $\xi_{n,k-2}$.)

These conditions and the assumption of Case 1 imply that $k-2 \neq 1,2$, so $k \geq 5$.
Now by these conditions, a direct computation shows that
\[
\left[\xi_{j,k}\right]_I \ni \xi_{k-2,k} \overset{m_{k,k+1}}{\sim}_I \xi_{k-2,k+1} \overset{m_{k+1,k+2}}{\sim}_I \cdots \overset{m_{n-1,n}}{\sim}_I \xi_{k-2,n},
\]
so $\xi_{k-2,n} \in (\mathcal{E}_{x_i}^{(I)})_{\mathrm{fin}}$ by Corollary \ref{cor:inthesamecomponentofWperpx}.
Since $\left[\xi_{1,n}\right]_I=\{\xi_{1,n}\}$, Corollary \ref{cor:necessaryconditionforfinitecomponent} implies that $\left[\xi_{k-2,n}\right]_I \overset{m'}{\sim}_I \xi_{1,n}$ for some $2 \leq m'<\infty$, forcing $\left[\xi_{k-2,n}\right]_I$ to contain either $\xi_{1,\ell'}$ (with $3 \leq \ell' \leq n-1$) or $\xi_{2,n}$.
However, since $k-2 \geq 3$ and $m_{2,n} \neq 2$, this is impossible by a similar argument to the previous paragraph.
This contradiction proves the claim in Case 1.

\noindent
\textbf{Case 2: $m_{1,n-1} \neq 2$ and $m_{2,n}=2$. (By symmetry, the case that $m_{1,n-1}=2$ and $m_{2,n} \neq 2$ is similarly observed.)}

This case is divided into the following two subcases.

\noindent
\textbf{Case 2-1: $m_{1,k}<\infty$ for some $3 \leq k \leq n-2$.}

Now $n \geq 5$.
Let $k$ be the largest index with this property.
First, for any $3 \leq j \leq k$ with $m_{1,j}<\infty$ (e.g.\ $j=k$), a similar argument to the second paragraph of Case 1 shows that $\left[\xi_{n,1}\right]_I=\{\xi_{n,1}\} \overset{\infty}{\sim}_I \left[\xi_{j,1}\right]_I$ (since $m_{1,n-1} \neq 2$).
Thus $\xi_{n,1},\xi_{j,1} \not\in (\mathcal{E}_{x_i}^{(I)})_{\mathrm{fin}}$ by Corollary \ref{cor:necessaryconditionforfinitecomponent}.

Secondly, we show that $\xi_{1,j} \not\in (\mathcal{E}_{x_i}^{(I)})_{\mathrm{fin}}$ for such $j$.
By the first remark of the proof, we may assume that $i=1$.
If $J=\{x_1,\dots,x_k\}$ is of type $A_k$, $H_4$ or $P(m)$, then since $\xi_{j,1} \not\in (\mathcal{E}_{x_1}^{(I)})_{\mathrm{fin}}$ and $m_{1,k}<\infty$, the induction assumption (applied to $J$) proves the claim.
If $J$ is of type $H_3$ (so $j=k=3$), then we have $\left[\xi_{1,j}\right]_I=\{\xi_{1,j}\}$, $\left[\xi_{n,1}\right]_I=\{\xi_{n,1}\}$ and $\xi_{1,j} \overset{\infty}{\sim}_I \xi_{n,1}$ by a direct computation, proving the claim (by Corollary \ref{cor:necessaryconditionforfinitecomponent}).
In the other cases, the claim follows from the induction assumption (applied to $J$).

Moreover, since $m_{1,n-1} \neq 2$, Corollary \ref{cor:miseven} implies that $\xi_{1,n-1},\xi_{n-1,1} \not\in (\mathcal{E}_{x_i}^{(I)})_{\mathrm{fin}}$ whenever $m_{1,n-1}<\infty$.
Similarly, we have $\xi_{1,n},\xi_{n,1} \not\in (\mathcal{E}_{x_i}^{(I)})_{\mathrm{fin}}$ if $m_{1,n} \neq 2$, while $\xi_{1,n} \overset{1}{\sim}_I \xi_{2,n}$ if $m_{1,n}=2$ (note that $m_{2,n}=2$).
Thus by Corollary \ref{cor:inthesamecomponentofWperpx} and the first remark of the proof, our remaining task is to show that $(\mathcal{E}_{x_n}^{(I)})_{\mathrm{fin}} \cap \mathcal{E}_{x_n}^{(K)}=\emptyset$ where $K=I \smallsetminus x_1$.

Note that now $m_{2,n}<\infty$ and $|K| \geq 4$; so $K$ does not satisfy the hypothesis of Claim \ref{thm_item:Wperpxforoddpath_H3}.
On the other hand, if $K$ satisfies the hypothesis of Claim \ref{thm_item:Wperpxforoddpath_other}, then the claim follows from the induction assumption (applied to $K$).
Finally, suppose that $K$ satisfies the hypothesis of Claim \ref{thm_item:Wperpxforoddpath_AHP}.
Note that $\xi_{n,2} \in \mathcal{E}_{x_n}^{(I)}$.
Then by the induction assumption (applied to $K$), it suffices to show that $\xi_{n,2} \not\in (\mathcal{E}_{x_n}^{(I)})_{\mathrm{fin}}$.
Now Lemma \ref{lem:conditionfornonorthogonal} implies that $r_{x_n}^{(I)}(\xi_{n,2})$ does not commute with $r_{x_n}^{(I)}(\xi_{n,1})$ (since $x_1$ and $x_2$ are adjacent); while $\xi_{n,1} \not\in (\mathcal{E}_{x_n}^{(I)})_{\mathrm{fin}}$ as above, so $r_{x_n}^{(I)}(\xi_{n,1}) \not\in (W_I^{\perp x_n})_{\mathrm{fin}}$.
Thus we have $r_{x_n}^{(I)}(\xi_{n,2}) \not\in (W_I^{\perp x_n})_{\mathrm{fin}}$ and $\xi_{n,2} \not\in (\mathcal{E}_{x_n}^{(I)})_{\mathrm{fin}}$, as desired.

Hence the claim holds in Case 2-1.

\noindent
\textbf{Case 2-2: $m_{1,k}=\infty$ for all $3 \leq k \leq n-2$.}

In the first two paragraph, we show that $\xi_{1,n-1},\xi_{1,n},\xi_{n-1,1},\xi_{n,1} \not\in (\mathcal{E}_{x_n}^{(I)})_{\mathrm{fin}}$.
Note that $\left[\xi_{n,1}\right]_I=\{\xi_{n,1}\}$ (since $m_{1,n-1} \neq 2$), $n \geq 4$ and
\[
\left[\xi_{1,n}\right]_I \subseteq \mathcal{E}' \textrm{ where } \mathcal{E}'=\{\xi_{1,n},\xi_{2,n},\dots,\xi_{n-2,n},\xi_{n,n-2}\}.
\]
First, suppose that $m_{1,n-1}<\infty$.
Then, since $m_{1,n-1} \neq 2$, we have $\left[\xi_{1,n-1}\right]_I=\{\xi_{1,n-1}\}$, $\left[\xi_{n-1,1}\right]_I=\{\xi_{n-1,1}\}$ and $\xi_{1,n-1},\xi_{n-1,1} \not\in (\mathcal{E}_{x_n}^{(I)})_{\mathrm{fin}}$ (see Corollary \ref{cor:miseven}).
Similarly, we have $\left[\xi_{n,1}\right]_I \overset{\infty}{\sim}_I \left[\xi_{1,n-1}\right]_I$ (since $\xi_{n,1} \overset{\infty}{\sim}_I \xi_{1,n-1}$ by a direct computation) and so $\xi_{n,1} \not\in (\mathcal{E}_{x_n}^{(I)})_{\mathrm{fin}}$.
Moreover, we have $\left[\xi_{1,n}\right]_I \overset{\infty}{\sim}_I \left[\xi_{n-1,1}\right]_I$ (since $\xi_{n-1,1} \overset{\infty}{\sim}_I \mathcal{E}'$ by a direct computation) and so $\xi_{1,n} \not\in (\mathcal{E}_{x_n}^{(I)})_{\mathrm{fin}}$ as desired.

Secondly, suppose that $m_{1,n-1}=\infty$, so $\xi_{1,n-1},\xi_{n-1,1} \not\in \mathcal{E}_{x_n}^{(I)}$.
Now it suffices to show that $\left[\xi_{1,n}\right]_I \overset{\infty}{\sim}_I \xi_{n,1}$ (note that $\left[\xi_{n,1}\right]_I=\{\xi_{n,1}\}$).
If this fails, then we have $\xi_{n,n-2} \in \left[\xi_{1,n}\right]_I$ and $\xi_{n,n-2} \overset{m}{\sim}_I \xi_{n,1}$ for some $m<\infty$, since now we have $\xi_{n,1} \overset{\infty}{\sim}_I (\mathcal{E}' \smallsetminus \xi_{n,n-2})$.
The latter relation implies that $m_{1,n-2}<\infty$, so $n=4$ by the assumption of Case 2-2.
Moreover, by a similar argument to the second paragraph of Case 1, the property $\xi_{4,2} \in \left[\xi_{1,4}\right]_I$ implies that $\xi_{1,4} \overset{1}{\sim}_I \xi_{2,4} \overset{1}{\sim}_I \xi_{4,2}$, therefore $m_{1,4}=m_{2,4}=2$ and $m_{2,3}=m_{3,4}=3$; namely, $(x_1,\dots,x_4)$ is of type $P(m_{1,2})$.
However, this is a contradiction, since we are assuming that the hypothesis of Claim \ref{thm_item:Wperpxforoddpath_AHP} is not satisfied.

Thus we have shown that $\xi_{1,n-1},\xi_{1,n},\xi_{n-1,1},\xi_{n,1} \not\in (\mathcal{E}_{x_n}^{(I)})_{\mathrm{fin}}$.
By the assumption of Case 2-2, it now suffices to show that $(\mathcal{E}_{x_n}^{(I)})_{\mathrm{fin}} \cap \mathcal{E}_{x_n}^{(J)}=\emptyset$ where $J=I \smallsetminus \{x_1\}$.
Now by the induction assumption (applied to $J$), the claim follows immediately if $J$ satisfies the hypothesis of Claim \ref{thm_item:Wperpxforoddpath_other}; similarly, if $J$ satisfies the hypothesis of Claim \ref{thm_item:Wperpxforoddpath_AHP}, then we have $\xi_{n,2} \not\in (\mathcal{E}_{x_n}^{(I)})_{\mathrm{fin}}$ by a similar argument to Case 2-1 (since $\xi_{n,1} \not\in (\mathcal{E}_{x_n}^{(I)})_{\mathrm{fin}}$ as above), so the claim also follows.
Finally, if $J$ is of type $H_3$ (so $n=4$), then $\mathcal{E}_{x_4}^{(J)}=\{\xi_{2,4},\xi_{4,2}\}$ and we have $\xi_{4,2} \not\in (\mathcal{E}_{x_4}^{(I)})_{\mathrm{fin}}$ similarly; while we have $\left[\xi_{2,4}\right]_I \subseteq \{\xi_{1,4},\xi_{2,4}\} \overset{\infty}{\sim}_I \{\xi_{4,1}\}=\left[\xi_{4,1}\right]_I$, so $\left[\xi_{2,4}\right]_I \overset{\infty}{\sim}_I \left[\xi_{4,1}\right]_I$ and $\xi_{2,4} \not\in (\mathcal{E}_{x_4}^{(I)})_{\mathrm{fin}}$.

Hence the claim holds in Case 2-2.

\noindent
\textbf{Case 3: $m_{1,n-1}=m_{2,n}=2$.}

Put $J=I \smallsetminus x_n$ and $K=I \smallsetminus x_1$.
This case is divided into the following two subcases.

\noindent
\textbf{Case 3-1: $m_{1,n} \neq 2$.}

Note that $m_{1,n}<\infty$ by the hypothesis, so we have $\left[\xi_{1,n}\right]_I=\{\xi_{1,n}\}$, $\left[\xi_{n,1}\right]_I=\{\xi_{n,1}\}$ and $\xi_{1,n},\xi_{n,1} \not\in (\mathcal{E}_{x_i}^{(I)})_{\mathrm{fin}}$ (see Corollary \ref{cor:miseven}).
Thus by symmetry, it suffices to show that $(\mathcal{E}_{x_1}^{(I)})_{\mathrm{fin}} \cap \mathcal{E}_{x_1}^{(J)}=\emptyset$.
Now we have $\xi_{1,n-1} \not\in (\mathcal{E}_{x_1}^{(I)})_{\mathrm{fin}}$ by a similar argument to Case 2-1 (since $\xi_{1,n} \not\in (\mathcal{E}_{x_1}^{(I)})_{\mathrm{fin}}$); so the claim follows from the induction assumption (applied to $J$) whenever $J$ satisfies the hypothesis of Claims \ref{thm_item:Wperpxforoddpath_AHP} or \ref{thm_item:Wperpxforoddpath_other}.
On the other hand, if $\mathrm{type}\,J=H_3$ (so $n=4$), then a direct computation shows that $\left[\xi_{3,1}\right]_I=\{\xi_{3,1}\}$ (since $m_{1,4} \neq 2$) and $\xi_{3,1} \overset{\infty}{\sim}_I \xi_{1,4}$, therefore $\left[\xi_{3,1}\right]_I \overset{\infty}{\sim}_I \left[\xi_{1,4}\right]_I$ (since $\left[\xi_{1,4}\right]_I=\{\xi_{1,4}\}$) and so $\xi_{3,1} \not\in (\mathcal{E}_{x_1}^{(I)})_{\mathrm{fin}}$.
Hence the claim holds in Case 3-1.

\noindent
\textbf{Case 3-2: $m_{1,n}=2$.}

Observe that $\xi_{1,n-1} \overset{m_{n-1,n}}{\sim}_I \xi_{1,n} \overset{1}{\sim}_I \xi_{2,n}$ and $\xi_{n,2} \overset{m_{1,2}}{\sim}_I \xi_{n,1} \overset{1}{\sim}_I \xi_{n-1,1}$, so it suffices to show that $(\mathcal{E}_{x_2}^{(I)})_{\mathrm{fin}} \cap \mathcal{E}_{x_2}^{(J)}=\emptyset$ and $(\mathcal{E}_{x_2}^{(I)})_{\mathrm{fin}} \cap \mathcal{E}_{x_2}^{(K)}=\emptyset$ (since then it follows from Corollary \ref{cor:inthesamecomponentofWperpx} that $\xi_{1,n},\xi_{n,1} \not\in (\mathcal{E}_{x_2}^{(I)})_{\mathrm{fin}}$).

If $J$ satisfies the hypothesis of Claim \ref{thm_item:Wperpxforoddpath_other}, then the induction assumption (applied to $J$) yields the former claim, particularly $\xi_{1,n-1},\xi_{n-1,1} \not\in (\mathcal{E}_{x_2}^{(I)})_{\mathrm{fin}}$.
Now the above observation and Corollary \ref{cor:inthesamecomponentofWperpx} imply that $\xi_{2,n},\xi_{n,2} \not\in (\mathcal{E}_{x_2}^{(I)})_{\mathrm{fin}}$.
Thus the other claim follows immediately if $\mathrm{type}\,K=H_3$ (since now $n=4$ and $\mathcal{E}_{x_2}^{(K)}=\{\xi_{2,4},\xi_{4,2}\}$), while this claim follows from the induction assumption (applied to $K$) if $K$ satisfies the hypothesis of Claims \ref{thm_item:Wperpxforoddpath_AHP} or \ref{thm_item:Wperpxforoddpath_other}.
Thus the claim holds in this case; and it holds similarly if $K$ satisfies the hypothesis of Claim \ref{thm_item:Wperpxforoddpath_other}.

Finally, we consider the case that each of $J$ and $K$ satisfies the hypothesis of Claims \ref{thm_item:Wperpxforoddpath_H3} or \ref{thm_item:Wperpxforoddpath_AHP}.
Note that $|J|=|K|$.
Now we are assuming that $I$ satisfies the hypothesis of Claim \ref{thm_item:Wperpxforoddpath_other}, so $I$ is not of finite type.
Thus by Theorem \ref{thm:specialcase_2spherical_acyclic}, the claim holds whenever $I$ is $2$-spherical.
So assume contrary that $I$ is not $2$-spherical.
Now by the above assumption on $J$ and $K$, at least one of $J$ and $K$ is of type $P(m)$, so $n=5$ and there are the four possibilities (up to symmetry) listed in Figure \ref{fig:listforWperpxforoddpath}.
In the list, the vertices $x_1,\dots,x_5$ of the Coxeter graph $\Gamma_I$, depicted in the top rows, are ordered from left to right, and $m$, $m_1$, $m_2$ are odd numbers.
Graphs in the middle rows signify the elements of $\mathcal{E}_{x_i}^{(I)}$ (where we abbreviate $\xi_{p,q}$ to $pq$) and their relations $\overset{k}{\sim}_I$; the relations $\overset{1}{\sim}_I$ are denoted by duplicated (vertical) lines.
Thus by Lemma \ref{lem:ifgammaisindependentofc}, the Coxeter graph of the Coxeter system $(W_I^{\perp x_i},R_{x_i}^{(I)})$ is as depicted in the bottom rows; so it is connected and not of finite type, therefore $(\mathcal{E}_{x_i}^{(I)})_{\mathrm{fin}}=\emptyset$ as desired.

Hence the proof of Theorem \ref{thm:Wperpxforoddpath} is concluded.
\end{proof}
\begin{figure}
\centering
\begin{picture}(40,250)(0,-250)
\put(0,-20){$\Gamma_I=$}
\put(0,-110){$\mathcal{E}_x^{(I)}=$}
\put(0,-220){$W_I^{\perp x}=$}
\end{picture}
\begin{picture}(150,250)(0,-250)
\put(10,-30){\circle{6}}\put(40,-30){\circle{6}}\put(70,-30){\circle{6}}\put(100,-30){\circle{6}}\put(130,-30){\circle{6}}
\put(13,-30){\line(1,0){24}}\put(43,-30){\line(1,0){24}}\put(73,-30){\line(1,0){24}}\put(103,-30){\line(1,0){24}}
\put(70,-27){\line(0,1){15}}\put(70,-12){\line(1,0){60}}{\put(130,-27){\line(0,1){15}}}
\put(100,-8){\hbox to0pt{\hss$\infty$\hss}}\put(115,-40){\hbox to0pt{\hss$m$\hss}}
\put(10,-80){\hbox to0pt{\hss$13$\hss}}\put(10,-110){\hbox to0pt{\hss$31$\hss}}\put(10,-140){\hbox to0pt{\hss$41$\hss}}\put(10,-170){\hbox to0pt{\hss$51$\hss}}
\put(10,-83){\line(0,-1){15}\,\line(0,-1){15}}\put(10,-113){\line(0,-1){15}\,\line(0,-1){15}}\put(10,-143){\line(0,-1){15}\,\line(0,-1){15}}
\put(60,-80){\hbox to0pt{\hss$14$\hss}}\put(60,-110){\hbox to0pt{\hss$24$\hss}}\put(60,-140){\hbox to0pt{\hss$42$\hss}}\put(60,-170){\hbox to0pt{\hss$52$\hss}}
\put(60,-83){\line(0,-1){15}\,\line(0,-1){15}}\put(60,-113){\line(0,-1){15}\,\line(0,-1){15}}\put(60,-143){\line(0,-1){15}\,\line(0,-1){15}}
\put(110,-80){\hbox to0pt{\hss$15$\hss}}\put(110,-110){\hbox to0pt{\hss$25$\hss}}
\put(110,-83){\line(0,-1){15}\,\line(0,-1){15}}
\put(20,-76){\line(1,0){30}}\put(35,-72){\hbox to0pt{\hss$3$\hss}}
\put(20,-136){\line(1,0){30}}\put(35,-132){\hbox to0pt{\hss$3$\hss}}
\put(20,-166){\line(1,0){30}}\put(35,-162){\hbox to0pt{\hss$3$\hss}}
\put(70,-76){\line(1,0){30}}\put(85,-72){\hbox to0pt{\hss$m$\hss}}
\put(70,-106){\line(1,0){30}}\put(85,-102){\hbox to0pt{\hss$m$\hss}}
\put(10,-220){\circle{8}}\put(60,-220){\circle{8}}\put(110,-220){\circle{8}}
\put(14,-220){\line(1,0){42}}\put(64,-220){\line(1,0){42}}\put(85,-230){\hbox to0pt{\hss$m$\hss}}
\put(10,-216){\line(0,1){15}}\put(10,-201){\line(1,0){100}}{\put(110,-216){\line(0,1){15}}}\put(60,-197){\hbox to0pt{\hss$\infty$\hss}}
\end{picture}
\qquad\quad
\begin{picture}(150,250)(0,-250)
\put(10,-30){\circle{6}}\put(40,-30){\circle{6}}\put(70,-30){\circle{6}}\put(100,-30){\circle{6}}\put(130,-30){\circle{6}}
\put(13,-30){\line(1,0){24}}\put(43,-30){\line(1,0){24}}\put(73,-30){\line(1,0){24}}\put(103,-30){\line(1,0){24}}
\put(70,-27){\line(0,1){15}}\put(70,-12){\line(1,0){60}}{\put(130,-27){\line(0,1){15}}}
\put(25,-40){\hbox to0pt{\hss$5$\hss}}\put(100,-8){\hbox to0pt{\hss$\infty$\hss}}\put(115,-40){\hbox to0pt{\hss$m$\hss}}
\put(10,-110){\hbox to0pt{\hss$31$\hss}}\put(10,-140){\hbox to0pt{\hss$41$\hss}}\put(10,-170){\hbox to0pt{\hss$51$\hss}}
\put(10,-113){\line(0,-1){15}\,\line(0,-1){15}}\put(10,-143){\line(0,-1){15}\,\line(0,-1){15}}
\put(40,-80){\hbox to0pt{\hss$13$\hss}}
\put(80,-80){\hbox to0pt{\hss$14$\hss}}\put(80,-110){\hbox to0pt{\hss$24$\hss}}\put(80,-140){\hbox to0pt{\hss$42$\hss}}\put(80,-170){\hbox to0pt{\hss$52$\hss}}
\put(80,-83){\line(0,-1){15}\,\line(0,-1){15}}\put(80,-113){\line(0,-1){15}\,\line(0,-1){15}}\put(80,-143){\line(0,-1){15}\,\line(0,-1){15}}
\put(130,-80){\hbox to0pt{\hss$15$\hss}}\put(130,-110){\hbox to0pt{\hss$25$\hss}}
\put(130,-83){\line(0,-1){15}\,\line(0,-1){15}}
\put(15,-102){\line(1,1){20}}\put(20,-92){\hbox to0pt{\hss$2$\hss}}
\put(50,-76){\line(1,0){20}}\put(60,-72){\hbox to0pt{\hss$3$\hss}}
\put(20,-136){\line(1,0){50}}\put(45,-132){\hbox to0pt{\hss$5$\hss}}
\put(20,-166){\line(1,0){50}}\put(45,-162){\hbox to0pt{\hss$5$\hss}}
\put(90,-76){\line(1,0){30}}\put(105,-72){\hbox to0pt{\hss$m$\hss}}
\put(90,-106){\line(1,0){30}}\put(105,-102){\hbox to0pt{\hss$m$\hss}}
\put(10,-220){\circle{8}}\put(40,-220){\circle{8}}\put(80,-220){\circle{8}}\put(130,-220){\circle{8}}
\put(44,-220){\line(1,0){32}}\put(84,-220){\line(1,0){42}}\put(105,-218){\hbox to0pt{\hss$m$\hss}}
\put(10,-216){\line(0,1){25}}\put(10,-191){\line(1,0){120}}{\put(130,-216){\line(0,1){25}}}\put(70,-187){\hbox to0pt{\hss$\infty$\hss}}
\put(10,-216){\line(1,2){8}}\put(18,-200){\line(1,0){54}}\put(80,-216){\line(-1,2){8}}\put(45,-210){\hbox to0pt{\hss$5$\hss}}
\put(40,-224){\line(0,-1){15}}\put(40,-239){\line(1,0){90}}\put(130,-224){\line(0,-1){15}}\put(85,-236){\hbox to0pt{\hss$\infty$\hss}}
\end{picture}
\\
\underline{\hbox to\textwidth{\hss\,\hss}}
\\
\ \\
\begin{picture}(40,250)(0,-250)
\put(0,-20){$\Gamma_I=$}
\put(0,-130){$\mathcal{E}_x^{(I)}=$}
\put(0,-220){$W_I^{\perp x}=$}
\end{picture}
\begin{picture}(150,250)(0,-250)
\put(10,-30){\circle{6}}\put(40,-30){\circle{6}}\put(70,-30){\circle{6}}\put(100,-30){\circle{6}}\put(130,-30){\circle{6}}
\put(13,-30){\line(1,0){24}}\put(43,-30){\line(1,0){24}}\put(73,-30){\line(1,0){24}}\put(103,-30){\line(1,0){24}}
\put(10,-27){\line(1,2){8}}\put(18,-11){\line(1,0){44}}{\put(70,-27){\line(-1,2){8}}}
\put(70,-27){\line(1,2){8}}\put(78,-11){\line(1,0){44}}{\put(130,-27){\line(-1,2){8}}}
\put(25,-40){\hbox to0pt{\hss$m_1$\hss}}\put(40,-8){\hbox to0pt{\hss$\infty$\hss}}\put(100,-8){\hbox to0pt{\hss$\infty$\hss}}\put(115,-40){\hbox to0pt{\hss$m_2$\hss}}
\put(10,-140){\hbox to0pt{\hss$41$\hss}}\put(10,-170){\hbox to0pt{\hss$51$\hss}}
\put(10,-143){\line(0,-1){15}\,\line(0,-1){15}}
\put(60,-80){\hbox to0pt{\hss$14$\hss}}\put(60,-110){\hbox to0pt{\hss$24$\hss}}\put(60,-140){\hbox to0pt{\hss$42$\hss}}\put(60,-170){\hbox to0pt{\hss$52$\hss}}
\put(60,-83){\line(0,-1){15}\,\line(0,-1){15}}\put(60,-113){\line(0,-1){15}\,\line(0,-1){15}}\put(60,-143){\line(0,-1){15}\,\line(0,-1){15}}
\put(110,-80){\hbox to0pt{\hss$15$\hss}}\put(110,-110){\hbox to0pt{\hss$25$\hss}}
\put(110,-83){\line(0,-1){15}\,\line(0,-1){15}}
\put(20,-136){\line(1,0){30}}\put(35,-132){\hbox to0pt{\hss$m_1$\hss}}
\put(20,-166){\line(1,0){30}}\put(35,-162){\hbox to0pt{\hss$m_1$\hss}}
\put(70,-76){\line(1,0){30}}\put(85,-72){\hbox to0pt{\hss$m_2$\hss}}
\put(70,-106){\line(1,0){30}}\put(85,-102){\hbox to0pt{\hss$m_2$\hss}}
\put(10,-220){\circle{8}}\put(60,-220){\circle{8}}\put(110,-220){\circle{8}}
\put(14,-220){\line(1,0){42}}\put(64,-220){\line(1,0){42}}
\put(35,-230){\hbox to0pt{\hss$m_1$\hss}}\put(85,-230){\hbox to0pt{\hss$m_2$\hss}}
\put(10,-216){\line(0,1){15}}\put(10,-201){\line(1,0){100}}{\put(110,-216){\line(0,1){15}}}\put(60,-197){\hbox to0pt{\hss$\infty$\hss}}
\end{picture}
\qquad
\begin{picture}(150,250)(0,-250)
\put(10,-30){\circle{6}}\put(40,-30){\circle{6}}\put(70,-30){\circle{6}}\put(100,-30){\circle{6}}\put(130,-30){\circle{6}}
\put(13,-30){\line(1,0){24}}\put(43,-30){\line(1,0){24}}\put(73,-30){\line(1,0){24}}\put(103,-30){\line(1,0){24}}
\put(40,-27){\line(0,1){15}}\put(40,-12){\line(1,0){60}}{\put(100,-27){\line(0,1){15}}}
\put(70,-8){\hbox to0pt{\hss$\infty$\hss}}
\put(10,-70){\hbox to0pt{\hss$13$\hss}}\put(10,-100){\hbox to0pt{\hss$31$\hss}}\put(10,-130){\hbox to0pt{\hss$41$\hss}}\put(10,-160){\hbox to0pt{\hss$51$\hss}}
\put(10,-73){\line(0,-1){15}\,\line(0,-1){15}}\put(10,-103){\line(0,-1){15}\,\line(0,-1){15}}\put(10,-133){\line(0,-1){15}\,\line(0,-1){15}}
\put(50,-100){\hbox to0pt{\hss$14$\hss}}
\put(90,-130){\hbox to0pt{\hss$52$\hss}}
\put(130,-70){\hbox to0pt{\hss$15$\hss}}\put(130,-100){\hbox to0pt{\hss$25$\hss}}\put(130,-130){\hbox to0pt{\hss$35$\hss}}\put(130,-160){\hbox to0pt{\hss$53$\hss}}
\put(130,-73){\line(0,-1){15}\,\line(0,-1){15}}\put(130,-103){\line(0,-1){15}\,\line(0,-1){15}}\put(130,-133){\line(0,-1){15}\,\line(0,-1){15}}
\put(20,-66){\line(1,0){100}}\put(70,-62){\hbox to0pt{\hss$2$\hss}}
\put(20,-100){\line(4,-1){100}}\put(70,-110){\hbox to0pt{\hss$2$\hss}}
\put(20,-156){\line(1,0){100}}\put(70,-168){\hbox to0pt{\hss$2$\hss}}
\put(20,-70){\line(1,-1){24}}\put(40,-84){\hbox to0pt{\hss$3$\hss}}
\put(120,-72){\line(-3,-1){62}}\put(80,-82){\hbox to0pt{\hss$3$\hss}}
\put(20,-152){\line(3,1){62}}\put(60,-136){\hbox to0pt{\hss$3$\hss}}
\put(122,-152){\line(-1,1){24}}\put(100,-144){\hbox to0pt{\hss$3$\hss}}
\put(10,-210){\circle{8}}\put(50,-210){\circle{8}}\put(90,-210){\circle{8}}\put(130,-210){\circle{8}}
\put(14,-210){\line(1,0){32}}\put(54,-210){\line(1,0){32}}\put(94,-210){\line(1,0){32}}
\put(70,-206){\hbox to0pt{\hss$\infty$\hss}}
\put(10,-206){\line(0,1){15}}\put(10,-191){\line(1,0){80}}{\put(90,-206){\line(0,1){15}}}
\put(50,-214){\line(0,-1){15}}\put(50,-229){\line(1,0){80}}{\put(130,-214){\line(0,-1){15}}}
\end{picture}
\caption{List for the proof of Theorem \ref{thm:Wperpxforoddpath}}
\label{fig:listforWperpxforoddpath}
\end{figure}
\subsection{The case with $\Gamma^{\mathrm{odd}}$ connected}
\label{sec:specialcase_oddconnected}

In this and the next two cases, we determine the structure of the Coxeter group $(W_O^{\perp x})_{\mathrm{fin}}$ where $O=S_{\sim x}^{\mathrm{odd}}$.
Recall that $V(G)$ denotes the vertex set of a graph $G$.
First we prepare some preliminary observations.
\begin{lem}
\label{lem:closedpathandfinitepart}
Let $(y,s) \in \mathcal{E}_x^{(I)}$, and suppose that the graph $\Gamma_I^{\mathrm{odd}}$ contains a non-backtracking nontrivial closed path $P=(z_n,\dots,z_1,z_0)$ starting from $z_0=z_n=y$ such that $s$ is not contained in but adjacent to $V(P)$ in the Coxeter graph $\Gamma_I$.
Then $(y,s) \not\in (\mathcal{E}_x^{(I)})_{\mathrm{fin}}$.
\end{lem}
\begin{proof}
By Remark \ref{rem:Efinisindependentofx}, we may assume without loss of generality that $x=y$.
Put $w=\pi^{(I)}(P) \in Y_x^{(I)}$, so $w \cdot \alpha_x=\alpha_x$.
Then the hypothesis implies that $\mathrm{supp}(w)=V(P)$ (by Corollary \ref{cor:supportofelementsofY}) and $w \cdot \alpha_s \neq \alpha_s$ (by Lemma \ref{lem:wdoesnotfixgamma}), so we have $w \cdot \gamma_x^{(I)}(x,s) \neq \gamma_x^{(I)}(x,s)$ since $s \in \mathrm{supp}(\widetilde{\alpha}_{x,s}) \subseteq \{x,s\}$.
Thus $(x,s) \not\in (\mathcal{E}_x^{(I)})_{\mathrm{fin}}$ by Theorem \ref{thm:Ycentralizefinitepart}.
\end{proof}
\begin{lem}
\label{lem:cyclecontainsyonly}
Let $(y,s) \in \mathcal{E}_x^{(I)}$, and suppose that $s \in I_{\sim x}^{\mathrm{odd}}$ and there is a simple closed path $P=(y_n,\dots,y_1,y_0)$ (with $n \geq 3$ and $y_n=y_0$) in $\Gamma_I^{\mathrm{odd}}$ which contains $y$ but does not contain $s$.
Then $(y,s) \not\in (\mathcal{E}_x^{(I)})_{\mathrm{fin}}$.
\end{lem}
\begin{proof}
By Remark \ref{rem:Efinisindependentofx} and symmetry, we may assume that $x=y=y_0$.
Then the claim follows immediately from Lemma \ref{lem:closedpathandfinitepart} if $s$ is adjacent to $V(P)$ in $\Gamma_I$; so suppose not.
Let $P'=(z_r,\dots,z_1,z_0)$ be the shortest path in $\Gamma_I^{\mathrm{odd}}$ from $z_0=s$ to $V(P)$ (such a path indeed exists since $s \in I_{\sim x}^{\mathrm{odd}}$), so $P'$ is simple, non-closed, non-backtracking and satisfies that $P' \cap P=\{z_r\}$.
Since $s$ is not adjacent to $V(P)$ in $\Gamma_I$, we have $(y_i,s) \in \left[x,s\right]_I$ for $0 \leq i \leq n$ (see Figure \ref{fig:localsubgraph}), so Corollary \ref{cor:inthesamecomponentofWperpx} allows us to assume that $x$ is adjacent to $z_r$ in $P$, namely $z_r=y_1$.

Now if $m_{x,z_i}$ is odd for some $1 \leq i \leq r-1$, then the graph $\Gamma_J^{\mathrm{odd}}$ (where $J=\{x,z_r,\dots,z_1,s\}$) contains a simple closed path which contains $x$, does not contain $s$ and is closer (in $\Gamma_I^{\mathrm{odd}}$) to $s$ than $P$.
Thus the claim follows by induction on the distance $r$ in $\Gamma_I^{\mathrm{odd}}$ between $P$ and $s$.

On the other hand, suppose that $m_{x,z_i}$ is even or $\infty$ for $1 \leq i \leq r-1$.
Note that $m_{x,s}$ is even and $r \geq 2$ (since $s$ is not adjacent to $V(P)$).
Now the graph $\Gamma_J^{\mathrm{odd}}$ is a simple non-closed path $(x,z_r,\dots,z_1,s)$ with $|J|=r+2 \geq 4$.
Thus by Theorem \ref{thm:Wperpxforoddpath} (applied to $J$), our claim holds unless $J$ is of type $A_{r+2}$, $H_4$ or $P(m)$; so suppose that $J$ is of one of these types.
Now if $m_{x,z_{r-1}}<\infty$, then we have $(x,z_{r-1}) \in \mathcal{E}_x^{(I)}$ and $(x,z_{r-1}) \not\in (\mathcal{E}_x^{(I)})_{\mathrm{fin}}$ by the first argument of the proof, so $(x,s) \not\in (\mathcal{E}_x^{(I)})_{\mathrm{fin}}$ by Theorem \ref{thm:Wperpxforoddpath} (\ref{thm_item:Wperpxforoddpath_AHP}).
On the other hand, suppose that $m_{x,z_{r-1}}=\infty$, so $\mathrm{type}\,J=P(m)$ and $r=2$.
Then $z_1$ is the unique neighbor of $s$ in $\Gamma_K^{\mathrm{odd}}$ (where $K=V(P) \cup J$) and $m_{x,z_1}=\infty$, so we have $\left[s,x\right]_K=\{(s,x)\}$, while $\mathscr{F}_s^{(K)} \neq 1$.
Thus Proposition \ref{prop:gammaisindependentofc} implies that $(s,x) \not\in (\mathcal{E}_x^{(K)})_{\mathrm{fin}}$, so $(s,x) \not\in (\mathcal{E}_x^{(I)})_{\mathrm{fin}}$ (by Lemma \ref{lem:restrictionoffinitepart}) and $(x,s) \not\in (\mathcal{E}_x^{(I)})_{\mathrm{fin}}$ (by Theorem \ref{thm:Wperpxforoddpath} (\ref{thm_item:Wperpxforoddpath_AHP}) applied to $J$).
Hence the proof is concluded.
\end{proof}
\begin{lem}
\label{lem:cyclenotcontainy}
Let $(y,s) \in \mathcal{E}_x^{(I)}$, and suppose that $s \in I_{\sim x}^{\mathrm{odd}}$ and there is a simple closed path $P=(y_n,\dots,y_1,y_0)$ (with $n \geq 3$ and $y_0=y_n$) in $(\Gamma_I^{\mathrm{odd}})_{\sim x}$ which does not contain $y$.
Then $(y,s) \not\in (\mathcal{E}_x^{(I)})_{\mathrm{fin}}$.
\end{lem}
\begin{proof}
By Remark \ref{rem:Efinisindependentofx}, we may assume that $x=y$.
Moreover, we may assume without loss of generality that the graph $\Gamma_{V(P)}^{\mathrm{odd}}$ itself is a simple closed path, i.e.\ $\Gamma_{V(P)}^{\mathrm{odd}}$ possesses no extra edges other than those of $P$.
Take one of the shortest paths $P'=(z_r,\dots,z_1,z_0)$ in $\Gamma_I^{\mathrm{odd}}$ from $z_0=x$ to $V(P)$, so the path $P'$ is simple, $r \geq 1$ and $P' \cap P=\{z_r\}$.
By symmetry, we may assume that $z_r=y_0=y_r$.

First, we prove the claim in the case that $s \not\in P \cup P'$.
If $s$ is adjacent to $V(P')$ in $\Gamma_I$, then the claim follows from Lemma \ref{lem:closedpathandfinitepart} (applied to the closed path $P'^{-1}PP'$).
Otherwise, we have $(z_i,s) \in \left[x,s\right]_I$ for all $i$, while $(z_r,s) \not\in (\mathcal{E}_x^{(I)})_{\mathrm{fin}}$ by Lemma \ref{lem:cyclecontainsyonly} (since $s \not\in P$), so the claim follows from Corollary \ref{cor:inthesamecomponentofWperpx}.

So suppose that $s \in P \cup P'$.
Note that $m_{y_i,z_j}$ is even or $\infty$ for any $i$ and $j \leq r-2$, by the choice of $P'$.
We divide the proof into two cases.

\noindent
\textbf{Case 1: $r=1$.}

Note that $s \in P$ since $s \in P \cup P'=P \cup x$.
Let $y_{i_0}$, $y_{i_1},\dots,y_{i_{k-1}}$ (where $0=i_0<i_1< \cdots <i_{k-1} \leq n-1$) be the vertices $y_j$ of $P$ such that $m_{x,y_j}$ is odd, and put $i_k=n$.
Note that $s \neq y_{i_j}$ for $0 \leq j \leq k$, since $m_{x,s}$ is even.
Now if $k \geq 2$, then there is some $1 \leq j \leq k$ such that the simple closed path $(x,y_{i_j},y_{i_j-1},\dots,y_{i_{j-1}+1},y_{i_{j-1}},x)$ in $\Gamma_I^{\mathrm{odd}}$ does not contain $s$.
Thus the claim follows from Lemma \ref{lem:cyclecontainsyonly}.

So suppose that $k=1$, namely $m_{x,y_j}$ is even or $\infty$ for $1 \leq j \leq n-1$.
By symmetry, we may assume that $s=y_i$ with $1 \leq i \leq n-2$.
Put $J=\{x,y_0,y_1,\dots,y_i\}$, so the graph $\Gamma_J^{\mathrm{odd}}$ is a simple non-closed path $Q=(s,y_{i-1},\dots,y_0,x)$ with $i+2$ vertices.
Note that $m_{x,s}<\infty$.
Now the claim follows from Theorem \ref{thm:Wperpxforoddpath} if $Q$ satisfies the hypothesis of Claim \ref{thm_item:Wperpxforoddpath_other} of the theorem.
If $Q$ satisfies the hypothesis of Claim \ref{thm_item:Wperpxforoddpath_AHP} of the theorem, then we have $(s,x) \not\in (\mathcal{E}_x^{(I)})_{\mathrm{fin}}$ by Lemma \ref{lem:cyclecontainsyonly}, so the claim follows from the theorem.

Finally, suppose that $Q$ is of type $H_3$, so $i=1$ and $s=y_1$.
Then, since $y_0$ is the unique neighbor of $x$ in $\Gamma_K^{\mathrm{odd}}$ (where $K=V(P) \cup x$) and $\mathscr{F}_x^{(K)} \neq 1$, we have $\left[x,s\right]_K=\{(x,s)\}$ (see Figure \ref{fig:localsubgraph}) and $(x,s) \not\in (\mathcal{E}_x^{(K)})_{\mathrm{fin}}$ by Proposition \ref{prop:gammaisindependentofc}, so the claim follows from Lemma \ref{lem:restrictionoffinitepart}.

\noindent
\textbf{Case 2: $r \geq 2$.}

First, if $m_{z_{r-1},y_i}$ is odd for some $1 \leq i \leq n-1$, then the graph $\Gamma_J^{\mathrm{odd}}$ (where $J=V(P) \cup V(P')$) contains a simple closed path (with at least three vertices) which does not contain $x$ and is closer to $x$ than $P$ in $\Gamma_I^{\mathrm{odd}}$.
Thus the claim is deduced by induction on the distance $r$ from $x$ to $P$ in $\Gamma_I^{\mathrm{odd}}$.

So suppose that $m_{z_{r-1},y_i}$ is even or $\infty$ for $1 \leq i \leq n-1$, namely $\Gamma_J^{\mathrm{odd}}$ is precisely the union of the path $P'$ and the closed path $P$.
Now if the shortest path $Q$ in $\Gamma_J^{\mathrm{odd}}$ from $x$ to $s$ satisfies the hypothesis of Theorem \ref{thm:Wperpxforoddpath} (\ref{thm_item:Wperpxforoddpath_other}), then the claim follows from the theorem.
If $Q$ satisfies the hypothesis of Theorem \ref{thm:Wperpxforoddpath} (\ref{thm_item:Wperpxforoddpath_AHP}), then we have $(s,x) \in (\mathcal{E}_x^{(I)})_{\mathrm{fin}}$ by either Lemma \ref{lem:cyclecontainsyonly} (if $s \in P$) or induction on the distance $r$ from $x$ to $P$ in $\Gamma_I^{\mathrm{odd}}$ (if $s \in P' \smallsetminus z_r$); so the claim follows from the theorem.
Finally, if $Q$ is of type $H_3$, then $z_1$ is the unique neighbor of $x$ in $\Gamma_J^{\mathrm{odd}}$, and $s=z_2$; so we have $\left[x,s\right]_J=\{(x,s)\}$ (see Figure \ref{fig:localsubgraph}) and $(x,s) \not\in (\mathcal{E}_x^{(J)})_{\mathrm{fin}}$ by Proposition \ref{prop:gammaisindependentofc} (since $\mathscr{F}_x^{(J)} \neq 1$), therefore the claim follows from Lemma \ref{lem:restrictionoffinitepart}.
Hence the proof is concluded.
\end{proof}
\begin{lem}
\label{lem:cyclecontainsyands}
Let $(y,s) \in \mathcal{E}_x^{(I)}$ and $J \subseteq I_{\sim x}^{\mathrm{odd}}$, and suppose that $\Gamma_J^{\mathrm{odd}}$ is a simple closed path $(y_n,\dots,y_1,y_0)$ with $n \geq 4$, $y_0=y_n=y$ and $s \in J$.
Put $m_{i,j}=m_{y_i,y_j}$.
If $(y,s) \in (\mathcal{E}_x^{(I)})_{\mathrm{fin}}$, then we have $n=4$, $s=y_2$, $m_{0,1}=m_{1,2}=m_{2,3}=m_{3,0}=3$, $m_{0,2}=2$ and $m_{1,3}=\infty$ (see Figure \ref{fig:diamond}).
\end{lem}
\begin{proof}
By Remark \ref{rem:Efinisindependentofx}, we may assume that $x=y$.
Put $s=y_k$ (so $2 \leq k \leq n-2$) and $P_{j,i}=(y_j,y_{j-1},\dots,y_{i+1},y_i)$ for $0 \leq i \leq j \leq n$.
First, we show that $k=2$ and $V(P_{2,0})$ is of type $A_3$.
Since $(x,s) \in (\mathcal{E}_x^{(I)})_{\mathrm{fin}}$, Theorem \ref{thm:Wperpxforoddpath} implies that $V(P_{k,0})$ (containing $x$ and $s$) satisfies the hypothesis of Claim \ref{thm_item:Wperpxforoddpath_H3} or \ref{thm_item:Wperpxforoddpath_AHP} of the theorem.
If $V(P_{k,0})$ is of type $A_{k+1}$ (with $k+1 \geq 4$), $H_4$ or $P(m)$, then a direct computation shows that $y_{k-2} \in \mathrm{supp}(\gamma)$ where
\[
\gamma=\pi^{(I)}(P_{k-1,0}) \cdot \gamma_x^{(I)}(x,s)=\pi^{(I)}(P_{k-1,0}) \cdot \alpha_{y_k}.
\]
Since $y_{k-2} \not\in \mathrm{supp}(\pi^{(I)}(P_{n,k-1}))$, we have $y_{k-2} \in \mathrm{supp}(\beta)$ where
\[
\beta=\pi^{(I)}(P_{n,k-1}) \cdot \gamma=\pi^{(I)}(P_{n,0}) \cdot \gamma_x^{(I)}(x,s),
\]
so $\beta \neq \alpha_{y_k}=\gamma_x^{(I)}(x,s)$.
However, this contradicts Theorem \ref{thm:Ycentralizefinitepart}, since $\pi^{(I)}(P_{n,0}) \in Y_x^{(I)}$ and $(x,s) \in (\mathcal{E}_x^{(I)})_{\mathrm{fin}}$.
Moreover, if $V(P_{k,0})$ is of type $H_3$ (so $k=2$), then a direct computation also shows that
\[
y_1 \in \mathrm{supp}(\pi^{(I)}(P_{2,0}) \cdot \gamma_x^{(I)}(x,s)) \textrm{ and } y_1 \not\in \mathrm{supp}(\pi^{(I)}(P_{n,2})),
\]
so $\pi^{(I)}(P_{n,0}) \cdot \gamma_x^{(I)}(x,s) \neq \alpha_{y_2}=\gamma_x^{(I)}(x,s)$ similarly, contradicting Theorem \ref{thm:Ycentralizefinitepart} as well.
Thus we have $k=2$ and $V(P_{2,0})$ is of type $A_3$.
By symmetry, it also follows that $n=4$ and $V(P_{4,2})$ is also of type $A_3$.

Now our remaining task is to show that $m_{1,3}=\infty$.
Note that $(x,s) \in (\mathcal{E}_x^{(J)})_{\mathrm{fin}}$ by Lemma \ref{lem:restrictionoffinitepart}.
If $m_{1,3}<\infty$, then $(y_1,y_3) \in \mathcal{E}_x^{(J)}$, and a direct computation shows that $\left[x,s\right]_J=\{(y_0,y_2),(y_2,y_0)\}$ and $(y_1,y_3) \overset{\infty}{\sim}_J \left[x,s\right]_J$; this contradicts Corollary \ref{cor:necessaryconditionforfinitecomponent}.
Hence the claim follows.
\end{proof}
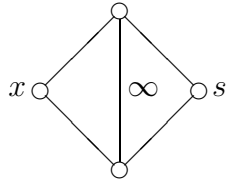
\begin{figure}
\centering
\begin{picture}(100,70)(0,-70)
\put(20,-35){\circle{6}}\put(50,-5){\circle{6}}\put(50,-65){\circle{6}}\put(80,-35){\circle{6}}
\put(22,-33){\line(1,1){26}}\put(22,-37){\line(1,-1){26}}\put(78,-33){\line(-1,1){26}}\put(78,-37){\line(-1,-1){26}}\put(50,-8){\line(0,-1){54}}
\put(8,-37){$x$}\put(85,-37){$s$}\put(53,-37){$\infty$}
\end{picture}
\caption{Coxeter graph in Lemma \ref{lem:cyclecontainsyands}}
\label{fig:diamond}
\end{figure}

Now we state the main results in this subsection, whose proofs are given in the next two subsections.
Recall notations and terminology for posets summarized in Section \ref{sec:graphsandgroupoids}.
\begin{thm}
\label{thm:finitepart_oddconnected_cyclic}
Put $O=S_{\sim x}^{\mathrm{odd}}$, and suppose that the graph $\Gamma_O^{\mathrm{odd}}$ contains a cycle.
Then we have $(W_O^{\perp x})_{\mathrm{fin}} \neq 1$ if and only if the following five conditions are satisfied:
\begin{enumerate}
\item \label{thm_item:fin_oddconn_cycl_bipyramid}
There exists a subset $K \subseteq O$ with $|K| \geq 4$ and $x_1,x_2 \in K$ such that $m_{x_1,x_2}=2$, $m_{x_1,y}=m_{x_2,y}=3$ for $y \in K^-=K \smallsetminus \{x_1,x_2\}$, and $m_{y,y'}=\infty$ for any distinct $y,y' \in K^-$.
\item \label{thm_item:fin_oddconn_cycl_trunk}
The graph $\Gamma_K^{\mathrm{odd}}$ contains all simple closed paths in $\Gamma_O^{\mathrm{odd}}$.
Hence $\Gamma_O^{\mathrm{odd}}$ admits a decomposition $\Gamma_O^{\mathrm{odd}}=\Gamma_K^{\mathrm{odd}} \cup \bigcup_{y \in K}T_y$ as in Remark \ref{rem:decompintotrees}.
\item \label{thm_item:fin_oddconn_cycl_almostallinfty}
For distinct $s,t \in O$, the order $m_{s,t}$ of $st$ is $2$, odd or $\infty$.
Moreover, if $m_{s,t}=2$, then one of the following three conditions is satisfied:
\begin{itemize}
\item $\{s,t\}=\{x_1,x_2\}$;
\item one of $s,t$ is $x_i$ with $i \in \{1,2\}$, and the other is an element of $T_y \smallsetminus y$ with $y \in K^-$;
\item $s,t$ are comparable elements of the same poset $T_y$ with $y \in K^-$.
\end{itemize}
\item \label{thm_item:fin_oddconn_cycl_orderideal}
For $i \in \{1,2\}$ and $y \in K^-$, the set
\[
T_y^{(i)}=y \cup \{s \in T_y \smallsetminus y \mid m_{x_i,s}=2\}
\]
is an order ideal of $T_y$. 
We have $T_y^{(1)}=T_y^{(2)}$, which we denote by $T_y^o$.
Moreover, we have $m_{y,s}=3$ if $s \in T_y^o$ is an atom of $T_y$.
\item \label{thm_item:fin_oddconn_cycl_sametree}
Suppose that $s,t \in T_y$ with $y \in K^-$, $s \prec_y t$ and $m_{s,t}=2$.
Then $t \in T_y^o \smallsetminus y$, and for distinct $s',t' \in T_y$, we have $m_{s',t'}=2$ if $s' \prec_y t' \preceq_y t$, $s' \preceq_y s$ and two vertices $s',t'$ of $T_y$ are not adjacent.
Moreover, if $t$ covers some $z$ and this $z$ covers $s$ in $T_y$, then $m_{s,z}=m_{z,t}=3$.
\end{enumerate}
Moreover, if these conditions are satisfied, then we have
\[
(\mathcal{E}_x^{(O)})_{\mathrm{fin}}=\{(x_1,x_2),(x_2,x_1)\},\ (x_1,x_2) \overset{1}{\sim}_O (x_2,x_1)
\]
and $(W_O^{\perp x})_{\mathrm{fin}}$ is generated by a single generator $r_x^{(O)}(x_1,x_2)=r_x^{(O)}(x_2,x_1)$.
\end{thm}
\begin{thm}
\label{thm:finitepart_oddconnected_acyclic_trankisfinitetype}
Put $O=S_{\sim x}^{\mathrm{odd}}$, and suppose that $\Gamma_O^{\mathrm{odd}}$ is acyclic and there exists a subset $K \subseteq O$ of type $A_n$ ($n \geq 3$), $D_n$ ($n \geq 4$), $E_6$, $E_7$, $E_8$, $H_3$, $H_4$ or $P(m)$ (recall from Example \ref{ex:WperpxforP} the definition of type $P(m)$) such that
\begin{equation}
\label{eq:Kistrunk}
\textrm{if } s,t \in O \textrm{ and } m_{s,t} \textrm{ is even, then } s,t \in K.
\end{equation}
Then for any $x' \in K$, we have $W_O^{\perp x'}=W_K^{\perp x'}$, $R_{x'}^{(O)}=R_{x'}^{(K)}$, $(\mathcal{E}_{x'}^{(O)})_{\mathrm{fin}}=\mathcal{E}_{x'}^{(O)}=\mathcal{E}_{x'}^{(K)}$ and $|W_K^{\perp x'}|<\infty$.
Hence $|W_O^{\perp x}|<\infty$ by Theorem \ref{thm:propertyofWperpx} (\ref{thm_item:conjugateinduceisomorphism}).
\end{thm}
\begin{thm}
\label{thm:finitepart_oddconnected_acyclic_trankisnonfinitetype}
Put $O=S_{\sim x}^{\mathrm{odd}}$, and suppose that $\Gamma_O^{\mathrm{odd}}$ is acyclic and there exists a subset $K=\{x_1,x_2,x_3\} \subseteq O$ satisfying the following four conditions:
\begin{enumerate}
\item \label{thm_item:fin_oddconn_acycl_nonftype_KisA3}
We have $m_{x_1,x_3}=2$ and $m_{x_1,x_2}=m_{x_2,x_3}=3$.
Hence $K$ induces a decomposition $\Gamma_O^{\mathrm{odd}}=\Gamma_K^{\mathrm{odd}} \cup \bigcup_{i=1}^{3}T_{x_i}$ as in Remark \ref{rem:decompintotrees}.
\item \label{thm_item:fin_oddconn_acycl_nonftype_almostallinfty}
For distinct $s,t \in O$, the order $m_{s,t}$ of $st$ is $2$, odd or $\infty$.
Moreover, if $m_{s,t}=2$, then one of the following three conditions is satisfied:
\begin{itemize}
\item $\{s,t\}=\{x_1,x_3\}$;
\item one of $s,t$ is $x_i$ with $i \in \{1,3\}$, and the other belongs to $T_{x_2} \smallsetminus x_2$;
\item $s,t$ are comparable elements of $T_{x_2}$.
\end{itemize}
\item \label{thm_item:fin_oddconn_acycl_nonftype_orderideal}
For $i \in \{1,3\}$, the set
\[
T_{x_2}^{(i)}=x_2 \cup \{s \in T_{x_2} \smallsetminus x_2 \mid m_{x_i,s}=2\}
\]
is an order ideal of $T_{x_2}$.
We have $T_{x_2}^{(1)}=T_{x_2}^{(3)}$, which we denote by $T_{x_2}^o$.
Moreover, we have $m_{x_2,s}=3$ if $s \in T_{x_2}^o$ is an atom of $T_{x_2}$.
\item \label{thm_item:fin_oddconn_acycl_nonftype_Tx2}
Suppose that $s,t \in T_{x_2}$, $s \prec_{x_2} t$ and $m_{s,t}=2$.
Then $t \in T_{x_2}^o \smallsetminus x_2$, and for distinct $s',t' \in T_{x_2}$, we have $m_{s',t'}=2$ if $s' \prec_{x_2} t' \preceq_{x_2} t$, $s' \preceq_{x_2} s$ and two vertices $s',t'$ of $T_{x_2}$ are not adjacent.
Moreover, if $t$ covers some $z$ and this $z$ covers $s$, then $m_{s,z}=m_{z,t}=3$.
\end{enumerate}
If $K'=K \cup V(T_{x_2}^o)$ is of type $A_3$ or $D_n$ ($n \geq 4$), then the hypothesis of Theorem \ref{thm:finitepart_oddconnected_acyclic_trankisfinitetype} is satisfied (with $K'$ playing the role of $K$ there).
Otherwise, we have $(\mathcal{E}_x^{(O)})_{\mathrm{fin}}=\{(x_1,x_3),(x_3,x_1)\}$, $(x_1,x_3) \overset{1}{\sim}_O (x_3,x_1)$ and $(W_O^{\perp x})_{\mathrm{fin}}$ is generated by a single generator $r_x^{(O)}(x_1,x_3)=r_x^{(O)}(x_3,x_1)$.
\end{thm}
\begin{thm}
\label{thm:finitepart_oddconnected_acyclic}
Put $O=S_{\sim x}^{\mathrm{odd}}$, and suppose that $\Gamma_O^{\mathrm{odd}}$ is acyclic.
Then we have $(W_O^{\perp x})_{\mathrm{fin}} \neq 1$ if and only if the hypothesis of one of Theorems \ref{thm:finitepart_oddconnected_acyclic_trankisfinitetype} and \ref{thm:finitepart_oddconnected_acyclic_trankisnonfinitetype} is satisfied.
\end{thm}
\subsection{Proof of Theorem \ref{thm:finitepart_oddconnected_cyclic}}
\label{sec:specialcase_proof_cyclic}

This subsection is devoted to the proof of Theorem \ref{thm:finitepart_oddconnected_cyclic}.
First, we prove the ``only if'' part.
Suppose that $(W_O^{\perp x})_{\mathrm{fin}} \neq 1$, and let $(x',x'') \in (\mathcal{E}_x^{(O)})_{\mathrm{fin}}$.
Then by Lemmas \ref{lem:cyclecontainsyonly}--\ref{lem:cyclecontainsyands}, any simple closed path in $\Gamma_O^{\mathrm{odd}}$ contains both $x'$ and $x''$, and is as in the conclusion of Lemma \ref{lem:cyclecontainsyands}.
This implies that the union $K$ of the vertex sets of all simple closed paths in $\Gamma_O^{\mathrm{odd}}$ satisfies Conditions \ref{thm_item:fin_oddconn_cycl_bipyramid} and \ref{thm_item:fin_oddconn_cycl_trunk} of Theorem \ref{thm:finitepart_oddconnected_cyclic}, where $(x_1,x_2)=(x',x'')$.
This condition also implies that $(x_1,x_2),(x_2,x_1) \in \left[x',x''\right]_O$, so $(x_1,x_2),(x_2,x_1) \in (\mathcal{E}_x^{(O)})_{\mathrm{fin}}$ (see Corollary \ref{cor:inthesamecomponentofWperpx}).

From now, we verify the remaining three conditions.

\noindent
\textbf{Step 1: For any $i \in \{1,2\}$, we have $m_{z,x_i}=\infty$ if $z$ is an atom of $T_{x_{3-i}}$.}

Note that $m_{z,x_i}$ is not odd.
Assume contrary that $m_{z,x_i}$ is even.
Then two generators $r_x^{(O)}(x_i,x_{3-i})$ and $r_x^{(O)}(x_i,z)$ of $W_O^{\perp x}$ do not commute, since their conjugates $r_{x_i}^{(O)}(x_i,x_{3-i})$ and $r_{x_i}^{(O)}(x_i,z)$ do not commute by Lemma \ref{lem:conditionfornonorthogonal} (note that $z$ is adjacent to $x_{3-i}$ in $\Gamma_O^{\mathrm{odd}}$); while $(x_i,z) \not\in (\mathcal{E}_x^{(O)})_{\mathrm{fin}}$ by Lemma \ref{lem:cyclecontainsyonly} and $(x_i,x_{3-i}) \in (\mathcal{E}_x^{(O)})_{\mathrm{fin}}$.
This is a contradiction, so Step 1 is concluded.

By Conditions \ref{thm_item:fin_oddconn_cycl_bipyramid} and \ref{thm_item:fin_oddconn_cycl_trunk}, any neighbor $z$ of $x_i$ in $\Gamma_O^{\mathrm{odd}}$ satisfies either $z \in K^-$, or $z \in T_{x_i}$ and $m_{z,x_{3-i}}=\infty$ (by Step 1).
Thus a direct computation shows that $\left[x_1,x_2\right]_O=\{(x_1,x_2),(x_2,x_1)\}$; and if $(z,z') \overset{2}{\sim}_O \left[x_1,x_2\right]_O$, then we have $z=x_i$ with $i \in \{1,2\}$, $z' \in O \smallsetminus K$ and $m_{x_{3-i},z'}=2$.
Moreover, Lemmas \ref{lem:cyclecontainsyonly} and \ref{lem:cyclenotcontainy} imply that $(\mathcal{E}_x^{(O)})_{\mathrm{fin}}=\left[x_1,x_2\right]_O$ (by Condition \ref{thm_item:fin_oddconn_cycl_trunk}).
Thus by Corollaries \ref{cor:necessaryconditionforfinitecomponent} and \ref{cor:inthesamecomponentofWperpx}, any $\sim_O$-equivalence class in $\mathcal{E}_x^{(O)}$ other than $\left[x_1,x_2\right]_O$ must be related to $\left[x_1,x_2\right]_O$ via the relation $\overset{2}{\sim}_O$, hence
\begin{equation}
\label{eq:finitepart_oddconnected_cyclic_element}
\textrm{it contains some } (x_i,z) \textrm{ with } i \in \{1,2\},\ z \in O \smallsetminus K \textrm{ and } m_{x_{3-i},z}=2.
\end{equation}
This yields the following observation:
\begin{eqnarray}
\nonumber
&&\textrm{if a subset } \mathcal{E}' \textrm{ of } \mathcal{E}_x^{(O)} \smallsetminus \left[x_1,x_2\right]_O \textrm{ is closed under the relation}\\
\label{eq:finitepart_oddconnected_cyclic_equivclass}
&&\overset{1}{\sim}_O \textrm{ and contains no element as in (\ref{eq:finitepart_oddconnected_cyclic_element}), then } \mathcal{E}'=\emptyset.
\end{eqnarray}

\noindent
\textbf{Step 2: For any distinct $s,t \in O$, the order $m_{s,t}$ is $2$, odd or $\infty$.}

Assume contrary that $m_{s,t}$ is even and $m_{s,t} \neq 2$.
Then, since $m_{x_1,x_2}=2$, one of $s$ and $t$ is not in $\{x_1,x_2\}$.
Say, $s \not\in \{x_1,x_2\}$.
Then $\left[s,t\right]_O=\{(s,t)\}$ by Example \ref{ex:misevenneq2}, so $\mathcal{E}'=\{(s,t)\}$ contradicts (\ref{eq:finitepart_oddconnected_cyclic_equivclass}).

\noindent
\textbf{Step 3: For any $i \in \{1,2\}$, we have $m_{z,x_i}=\infty$ if $z \in T_{x_{3-i}} \smallsetminus x_{3-i}$.}

By Step 1 (and Figure \ref{fig:localsubgraph}), it is straightforward to show that the set
\[
\mathcal{E}'=\{(t,x_i) \mid t \in T_{x_{3-i}} \smallsetminus x_{3-i} \textrm{ and } m_{t,x_i}<\infty\}
\]
satisfies the hypothesis of (\ref{eq:finitepart_oddconnected_cyclic_equivclass}), so $\mathcal{E}'=\emptyset$, proving the claim.

\noindent
\textbf{Step 4: For any $i \in \{1,2\}$, we have $m_{z,z'}=\infty$ if $z \in T_{x_i}$ and $z' \in T_{x_{3-i}} \smallsetminus x_{3-i}$.}

By Step 3, we may assume that $z \neq x_i$.
Then by Step 3 again, the set
\[
\mathcal{E}'=\{(t,z) \mid t \in T_{x_{3-i}} \smallsetminus x_{3-i} \textrm{ and } m_{t,z}<\infty\}
\]
satisfies the hypothesis of (\ref{eq:finitepart_oddconnected_cyclic_equivclass}), so $\mathcal{E}'=\emptyset$, proving the claim.

\noindent
\textbf{Step 5: For any $i \in \{1,2\}$ and $y \in K^-$, we have $m_{z,z'}=\infty$ if $z \in T_{x_i} \smallsetminus x_i$ and $z' \in T_y$.}

In the case $z'=y$, Step 1 and Condition \ref{thm_item:fin_oddconn_cycl_bipyramid} imply that the union $\mathcal{E}'$ of the two sets
\[
\{(t,y) \mid t \in T_{x_i} \smallsetminus x_i \textrm{ and } m_{t,y}<\infty\}
\]
and
\[
\{(t,t') \mid t \in T_y,\ t' \textrm{ is an atom of } T_{x_i} \textrm{ and } m_{t,t'}<\infty\}
\]
satisfies the hypothesis of (\ref{eq:finitepart_oddconnected_cyclic_equivclass}), so $\mathcal{E}'=\emptyset$, proving the claim.
Now the claim for general $z'$ follows in the same way as Step 4 following from Step 3.

\noindent
\textbf{Step 6: For any $i \in \{1,2\}$, we have $m_{z,z'}=\infty$ if $z,z'$ are distinct, non-adjacent vertices of the tree $T_{x_i}$.
(Hence by Steps 4--6, if $z,z' \in O$ and $m_{z,z'}=2$, then we have $z,z' \not\in T_{x_i} \smallsetminus x_i$.)}

Steps 3, 5 and Condition \ref{thm_item:fin_oddconn_cycl_bipyramid} imply that the set
\[
\mathcal{E}'=\{(t,t') \mid t,t' \in T_{x_i} \textrm{ and } m_{t,t'} \textrm{ is even}\}=\mathcal{E}_{x_i}^{(V(T_{x_i}))}
\]
satisfies the hypothesis of (\ref{eq:finitepart_oddconnected_cyclic_equivclass}); indeed, by Step 3, no element $(t,t')$ of $\mathcal{E}'$ satisfies that $t' \neq x_i$ and $m_{t',x_{3-i}}=2$.
Thus $\mathcal{E}'=\emptyset$, proving the claim.

\noindent
\textbf{Step 7: For any distinct $y,y' \in K^-$, we have $m_{z,z'}=\infty$ if $z \in T_y$ and $z' \in T_{y'}$.}

Since $m_{y,y'}=\infty$ (see Condition \ref{thm_item:fin_oddconn_cycl_bipyramid}), the set
\[
\mathcal{E}'=\{(t,y) \mid t \in T_{y'} \textrm{ and } m_{t,y}<\infty\}
\]
satisfies the hypothesis of (\ref{eq:finitepart_oddconnected_cyclic_equivclass}), so $\mathcal{E}'=\emptyset$, proving the claim in the case $z=y$.
For a general case, since $m_{y,z'}=\infty$ as above, the set
\[
\mathcal{E}''=\{(t,z') \mid t \in T_y \textrm{ and } m_{t,z'}<\infty\}
\]
satisfies the hypothesis of (\ref{eq:finitepart_oddconnected_cyclic_equivclass}), so $\mathcal{E}''=\emptyset$, proving the claim.

\noindent
\textbf{Step 8: For any $y \in K^-$, we have $m_{z,z'}=\infty$ if $z,z'$ are incomparable elements of the poset $T_y$.
(Hence by Steps 2, 6, 7 and 8, Condition \ref{thm_item:fin_oddconn_cycl_almostallinfty} of the theorem is satisfied.)}

Recall that the meet $z \wedge z'$ of $z$ and $z'$ in the poset $T_y$ exists.
Let
\[
z \wedge z'=s_0 \prec_y s_1 \prec_y \cdots \prec_y s_k=z \textrm{ and } z \wedge z'=t_0 \prec_y t_1 \prec_y \cdots \prec_y t_\ell=z' 
\]
be the unique saturated chains in $T_y$ from $z \wedge z'$ to $z$ and $z'$, respectively.
Note that $k,\ell \geq 1$ since $z$ and $z'$ are incomparable.
Now if $k=1$, then the union $\mathcal{E}'$ of the two sets
\[
\{(t,s_1) \mid t_1 \preceq_y t \in T_y \textrm{ and } m_{t,s_1}<\infty\}
\]
and
\[
\{(s,t_1) \mid s_1 \preceq_y s \in T_y \textrm{ and } m_{s,t_1}<\infty\}
\]
satisfies the hypothesis of (\ref{eq:finitepart_oddconnected_cyclic_equivclass}), so $\mathcal{E}'=\emptyset$, proving $m_{z',s_1}=\infty$ as desired.
Moreover, if $k \geq 2$, then $m_{s_1,z'}=\infty$ as above, so the set
\[
\mathcal{E}''=\{(s,z') \mid s_1 \prec_y s \in T_y \textrm{ and } m_{s,z'}<\infty\}
\]
satisfies the hypothesis of (\ref{eq:finitepart_oddconnected_cyclic_equivclass}); thus $\mathcal{E}''=\emptyset$, proving the claim.

\noindent
\textbf{Step 9: Condition \ref{thm_item:fin_oddconn_cycl_orderideal} of the theorem is satisfied.}

First, if $s \in T_y^{(i)} \smallsetminus y$, then we have $m_{x_{3-i},s}<\infty$ by Lemma \ref{lem:deniedtriangle} (since $(x_i,x_{3-i}) \in (\mathcal{E}_x^{(O)})_{\mathrm{fin}}$ and $m_{x_i,s}=2$), so $m_{x_{3-i},s}=2$ by Condition \ref{thm_item:fin_oddconn_cycl_almostallinfty}.
Thus we have $T_y^{(1)}=T_y^{(2)}$ by symmetry.

For the remaining claims of this step, let $y=y_0 \prec_y y_1 \prec_y \cdots \prec_y y_k=s$ be the unique saturated chain in $T_y$ from $y$ to $s \in T_y^o \smallsetminus y$, so $k \geq 1$.
Now if $k=1$ and $m_{y,y_1} \neq 3$, then the nonempty set
\[
\{(t,x_i) \mid y_1 \preceq_y t \in T_y \textrm{ and } m_{t,x_i}<\infty\}
\]
satisfies the hypothesis of (\ref{eq:finitepart_oddconnected_cyclic_equivclass}), a contradiction.
Thus we have $m_{y,s}=3$ whenever $k=1$ (namely $s$ is an atom of $T_y$).
Moreover, if $k \geq 2$ and $m_{x_i,y_{k-1}} \neq 2$, then the nonempty set
\[
\{(t,x_i) \mid y_k \preceq_y t \in T_y \textrm{ and } m_{t,x_i}<\infty\}
\]
satisfies the hypothesis of (\ref{eq:finitepart_oddconnected_cyclic_equivclass}), a contradiction.
Thus we have $m_{x_i,y_{k-1}}=2$ whenever $k \geq 2$; so we have $m_{x_i,y_j}=2$ inductively for $1 \leq j \leq k$.
Hence the claim follows.

\noindent
\textbf{Step 10: Condition \ref{thm_item:fin_oddconn_cycl_sametree} of the theorem is satisfied.}

Let $y=y_0 \prec_y y_1 \prec_y \cdots \prec_y y_k=t$ be the unique saturated chain in $T_y$ from $y$ to $t$, and $s=y_\ell$, so $k \geq \ell+2 \geq 2$.
First, we show that $t \in T_y^o \smallsetminus y$.
If this fails, then $m_{x_1,t}=m_{x_2,t}=\infty$ by Conditions \ref{thm_item:fin_oddconn_cycl_almostallinfty} and \ref{thm_item:fin_oddconn_cycl_orderideal}; so Conditions \ref{thm_item:fin_oddconn_cycl_bipyramid} and \ref{thm_item:fin_oddconn_cycl_almostallinfty} imply that the union of the two sets
\[
\{(y_i,y_k) \mid i \leq k-2 \textrm{ and } m_{y_i,y_k}=2\}
\]
and
\[
\{(z,y_{k-2}) \mid t \preceq_y z \in T_y \textrm{ and } m_{z,y_{k-2}}=2\}
\]
is nonempty and satisfies the hypothesis of (\ref{eq:finitepart_oddconnected_cyclic_equivclass}), a contradiction.
Thus $t \in T_y^o \smallsetminus y$ as desired.

Secondly, we show that $m_{s,z}=m_{z,t}=3$ if $\ell=k-2$, where $z=y_{k-1}$.
If this fails, then the nonempty set
\[
\{(z',y_{k-2}) \mid y_k \preceq_y z' \in T_y \textrm{ and } m_{z',y_{k-2}}=2\}
\]
satisfies the hypothesis of (\ref{eq:finitepart_oddconnected_cyclic_equivclass}), a contradiction.
Thus $m_{s,z}=m_{z,t}=3$ as desired.

Finally, we show that $m_{y_i,y_j}=2$ if $0 \leq i \leq \ell$, $j \leq k$ and $i \leq j-2$, concluding the proof of Step 10.
By induction, it suffices to consider the following two cases; $(i,j)=(\ell-1,k)$, or $\ell<k-2$ and $(i,j)=(\ell,k-1)$.
Now for the first case, if $m_{y_{\ell-1},y_k} \neq 2$, then Condition \ref{thm_item:fin_oddconn_cycl_almostallinfty} implies that the union of the two sets
\[
\{(y_p,y_k) \mid \ell \leq p \leq k-2 \textrm{ and } m_{y_p,y_k}=2\}
\]
and
\[
\{(z,y_{k-2}) \mid y_k \preceq_y z \in T_y \textrm{ and } m_{z,y_{k-2}}=2\}
\]
is nonempty and satisfies the hypothesis of (\ref{eq:finitepart_oddconnected_cyclic_equivclass}), a contradiction.
Thus $m_{y_{\ell-1},y_k}=2$ as desired.
On the other hand, for the second case, if $\ell<k-2$ and $m_{y_\ell,y_{k-1}} \neq 2$, then the nonempty set
\[
\{(z,y_\ell) \mid y_k \preceq_y z \in T_y \textrm{ and } m_{z,y_\ell}=2\}
\]
satisfies the hypothesis of (\ref{eq:finitepart_oddconnected_cyclic_equivclass}), a contradiction.
Thus $m_{y_\ell,y_{k-1}}=2$ as desired.
Hence Step 10 is concluded.

By Steps 8--10, the proof of the ``only if'' part is concluded.

From now, we prove the `if' part, and verify the description of $(\mathcal{E}_x^{(O)})_{\mathrm{fin}}$ given in the statement.
By Remark \ref{rem:Efinisindependentofx}, we may assume that $x=x_1$.
First, Conditions \ref{thm_item:fin_oddconn_cycl_bipyramid}--\ref{thm_item:fin_oddconn_cycl_almostallinfty} imply that
\[
\left[x_1,x_2\right]_O=\{(x_1,x_2),(x_2,x_1)\}=\mathcal{E}_x^{(K)},\ (\mathcal{E}_x^{(O)})_{\mathrm{fin}} \subseteq \mathcal{E}_x^{(K)}
\]
(see Lemmas \ref{lem:cyclecontainsyonly} and \ref{lem:cyclenotcontainy} for the latter inclusion), and the elements $(x_1,x_2)$ and $(x_2,x_1)$ of $\mathcal{E}_x^{(O)}$ satisfy the condition (\ref{eq:gammaisindependentofc}) in Proposition \ref{prop:gammaisindependentofc}.
Indeed, we have $\mathscr{F}_{x_i}^{(O)}=\mathscr{F}_{x_i}^{(K)}$ by the shape of $\Gamma_O^{\mathrm{odd}}$; while for any $y \in K^-$, the path $(x_{3-i},y,x_i)$ in $\Gamma_O^{\mathrm{odd}}$ is the image of the edge of $\mathcal{I}_x^{(O)}$ from $(x_i,x_{3-i})$ to $(x_{3-i},x_i)$, induced by the subgraph $\Gamma_I$ of $\Gamma_O$ (where $I=\{x_1,x_2,y\}$), by the map $\iota$ (see Figure \ref{fig:localsubgraph}).

Thus it suffices to show that $(x_1,x_2) \in (\mathcal{E}_x^{(O)})_{\mathrm{fin}}$, or $\left[s,t\right]_O \overset{2}{\sim}_O \left[x_1,x_2\right]_O$ for all $(s,t) \in \mathcal{E}_x^{(O)} \smallsetminus \left[x_1,x_2\right]_O$ (see Lemma \ref{lem:ifgammaisindependentofc} (\ref{lem_item:whengammaisindependentofc_relation})).
By Condition \ref{thm_item:fin_oddconn_cycl_almostallinfty}, we have $m_{s,t}=2$, and we can divide the proof into the following four cases:

\noindent
\textbf{Case 1: $s=x_i$ with $i \in \{1,2\}$, and $t \in T_y \smallsetminus y$ with $y \in K^-$.}

Now we have $m_{x_{3-i},t}=2$ by Condition \ref{thm_item:fin_oddconn_cycl_orderideal}, so the claim follows since $(x_i,t) \overset{2}{\sim}_O (x_i,x_{3-i}) \in \left[x_1,x_2\right]_O$.

\noindent
\textbf{Case 2: $t=x_i$ with $i \in \{1,2\}$, and $s \in T_y \smallsetminus y$ with $y \in K^-$.}

Let $y=y_0 \prec_y y_1 \prec_y \cdots \prec_y y_k=s$ be the unique saturated chain in $T_y$ from $y$ to $s$ (so $k \geq 1$).
Then by Conditions \ref{thm_item:fin_oddconn_cycl_bipyramid} and \ref{thm_item:fin_oddconn_cycl_orderideal}, we have $m_{x_1,y_j}=m_{x_2,y_j}=2$ for $1 \leq j \leq k$, and $m_{x_i,y}=m_{y,y_1}=3$.
This implies that
\[
(s,t)=(y_k,x_i) \overset{1}{\sim}_O (y_{k-1},x_i) \overset{1}{\sim}_O \cdots \overset{1}{\sim}_O (y_1,x_i) \overset{1}{\sim}_O (x_i,y_1) \overset{2}{\sim}_O (x_i,x_{3-i}),
\]
proving the claim.

\noindent
\textbf{Case 3: $s,t \in T_y$ with $y \in K^-$, and $s \prec_y t$.}

Let $y=y_0 \prec_y y_1 \prec_y \cdots \prec_y y_k=t$ be the unique saturated chain in $T_y$ from $y$ to $t$, and $s=y_\ell$ with $\ell \leq k-2$.
Then by Condition \ref{thm_item:fin_oddconn_cycl_sametree}, we have $m_{y_i,y_k}=2$ for $0 \leq i \leq \ell$, and $m_{x_1,y_k}=m_{x_2,y_k}=2$.
This implies that
\[
(s,t)=(y_\ell,y_k) \overset{1}{\sim}_O (y_{\ell-1},y_k) \overset{1}{\sim}_O \cdots \overset{1}{\sim}_O (y_0,y_k) \overset{1}{\sim}_O (x_1,y_k) \overset{2}{\sim}_O (x_1,x_2),
\]
proving the claim.

\noindent
\textbf{Case 4: $s,t \in T_y$ with $y \in K^-$, and $t \prec_y s$.}

Let $y=y_0 \prec_y y_1 \prec_y \cdots \prec_y y_k=s$ be the unique saturated chain in $T_y$ from $y$ to $s$, and $t=y_\ell$ with $\ell \leq k-2$.
Then by Condition \ref{thm_item:fin_oddconn_cycl_sametree}, we have $m_{y_\ell,y_i}=2$ for $\ell+2 \leq i \leq k$, $m_{y_\ell,y_{\ell+1}}=m_{y_{\ell+1},y_{\ell+2}}=3$, $m_{y_i,y_{\ell+2}}=2$ for $0 \leq i \leq \ell$, and $m_{x_1,y_{\ell+2}}=m_{x_2,y_{\ell+2}}=2$.
This implies that
\begin{eqnarray*}
(s,t)=(y_k,y_\ell) \overset{1}{\sim}_O (y_{k-1},y_\ell) \overset{1}{\sim}_O \cdots \overset{1}{\sim}_O (y_{\ell+2},y_\ell) \overset{1}{\sim}_O (y_\ell,y_{\ell+2})\\
\overset{1}{\sim}_O (y_{\ell-1},y_{\ell+2}) \overset{1}{\sim}_O \cdots \overset{1}{\sim}_O (y_0,y_{\ell+2}) \overset{1}{\sim}_O (x_1,y_{\ell+2}) \overset{2}{\sim}_O (x_1,x_2),
\end{eqnarray*}
proving the claim.

Hence the proof of Theorem \ref{thm:finitepart_oddconnected_cyclic} is concluded.
\subsection{Proofs of Theorems \ref{thm:finitepart_oddconnected_acyclic_trankisfinitetype}--\ref{thm:finitepart_oddconnected_acyclic}}
\label{sec:specialcase_proof_acyclic}

In this subsection, we suppose that $\Gamma_O^{\mathrm{odd}}$ is acyclic, so every element of $\mathcal{E}_x^{(O)}$ satisfies the condition (\ref{eq:gammaisindependentofc}) in Proposition \ref{prop:gammaisindependentofc} (where $I=O$).
\begin{proof}
[Proof of Theorem \ref{thm:finitepart_oddconnected_acyclic_trankisfinitetype}]
We may assume that $x'=x$.
Since $\Gamma_O^{\mathrm{odd}}$ and $\Gamma_K^{\mathrm{odd}}$ are trees, we have $r_x^{(O)}(y,s)=r_x^{(K)}(y,s)$ for $(y,s) \in \mathcal{E}_x^{(K)}$; so the condition (\ref{eq:Kistrunk}) implies that $\mathcal{E}_x^{(O)}=\mathcal{E}_x^{(K)}$, $R_x^{(O)}=R_x^{(K)}$ and so $W_O^{\perp x}=W_K^{\perp x}$.
Moreover, the hypothesis means that $W_K^{\perp x}$ is finite (see Theorem \ref{thm:Wperpxforoddpath} for type $P(m)$).
Hence the claim follows.
\end{proof}
\begin{proof}
[Proof of Theorem \ref{thm:finitepart_oddconnected_acyclic_trankisnonfinitetype}]
The proof is similar to the `if' part of Theorem \ref{thm:finitepart_oddconnected_cyclic}, but slightly more complicated, since now $\Gamma_O^{\mathrm{odd}}$ is acyclic and so Lemmas \ref{lem:cyclecontainsyonly} and \ref{lem:cyclenotcontainy} do not work.

First, we consider the case that $K'$ is of type $A_3$ or $D_n$.
Our aim is to verify the condition (\ref{eq:Kistrunk}), namely $s,t \in K'$ whenever $s,t \in O$ and $m_{s,t}$ is even.
This holds if $s,t \in K$; so suppose not.
Then Condition \ref{thm_item:fin_oddconn_acycl_nonftype_almostallinfty} of the theorem shows that $m_{s,t}=2$, and we have (up to symmetry) either $s=x_i$ with $i \in \{1,3\}$ and $t \in T_{x_2} \smallsetminus x_2$, or $s,t \in T_{x_2}$ and $s \prec_{x_2} t$.
In the first case, we have $t \in T_{x_2}^o \subseteq K'$ by definition (see Condition \ref{thm_item:fin_oddconn_acycl_nonftype_orderideal}); while in the second case, we have $t \in T_{x_2}^o$ (by Condition \ref{thm_item:fin_oddconn_acycl_nonftype_Tx2}) and $s \in T_{x_2}^o$ (by Condition \ref{thm_item:fin_oddconn_acycl_nonftype_orderideal}).
Hence the claim follows in any case.

So suppose that $K'$ is not of type $A_3$ or $D_n$.
As well as Theorem \ref{thm:finitepart_oddconnected_cyclic}, we have $\left[x_1,x_3\right]_O=\{(x_1,x_3),(x_3,x_1)\}$ by Conditions \ref{thm_item:fin_oddconn_acycl_nonftype_KisA3} and \ref{thm_item:fin_oddconn_acycl_nonftype_almostallinfty}.
Moreover, the same argument as the proof of the `if' part of Theorem \ref{thm:finitepart_oddconnected_cyclic} proves that, for any $(s,t) \in \mathcal{E}_x^{(O)} \smallsetminus \left[x_1,x_3\right]_O$, the set $\left[s,t\right]_O$ contains an element $(x_i,s')$ with $i \in \{1,3\}$ and $s' \in T_{x_2}^o \smallsetminus x_2$; so $\left[s,t\right]_O \overset{2}{\sim}_O \left[x_1,x_3\right]_O$ since $(x_i,s') \overset{2}{\sim}_O (x_i,x_{4-i})$.
Thus by Lemma \ref{lem:ifgammaisindependentofc}, the element $r_x^{(O)}(x_1,x_3)=r_x^{(O)}(x_3,x_1)$ generates a (finite) irreducible component of $W_O^{\perp x}$, proving that $\left[x_1,x_3\right]_O \subseteq (\mathcal{E}_x^{(O)})_{\mathrm{fin}}$.

Now our remaining task is to show that $(\mathcal{E}_x^{(O)})_{\mathrm{fin}} \subseteq \left[x_1,x_3\right]_O$.
By the above argument and Corollary \ref{cor:inthesamecomponentofWperpx}, it suffices to verify that $(x_i,s) \not\in (\mathcal{E}_x^{(O)})_{\mathrm{fin}}$ for any $i \in \{1,3\}$ and $s \in T_{x_2}^o \smallsetminus x_2$.
Assume contrary that $(x_i,s) \in (\mathcal{E}_x^{(O)})_{\mathrm{fin}}$.

First, we show that there is an element $t$ of $T_{x_2}^o$ maximal subject to the condition $s \preceq_{x_2} t$.
If $t$ is some element of $T_{x_2}^o$ with $s \preceq_{x_2} t$, and $s=z_0 \prec_{x_2} z_1 \prec_{x_2} \cdots \prec_{x_2} z_k=t$ is the saturated chain from $s$ to $t$, then $z_j \in T_{x_2}^o$ and $(x_i,z_{j-1}) \overset{m_j}{\sim}_O (x_i,z_j)$ (where $m_j=m_{z_{j-1},z_j} \neq 2$) for $1 \leq j \leq k$ (see Condition \ref{thm_item:fin_oddconn_acycl_nonftype_orderideal}).
This implies that all the \emph{distinct} generators $r_x^{(O)}(x_i,z_j)$ of $W_O^{\perp x}$ belong to the irreducible component of $W_O^{\perp x}$ containing $r_x^{(O)}(x_i,s)$.
Since the component is finite and independent of $k$, the number $k+1$ of those generators must be bounded above, yielding such a maximal element.

Moreover, the above argument also shows that $(x_i,t) \in (\mathcal{E}_x^{(O)})_{\mathrm{fin}}$.
Now Lemma \ref{lem:deniedtriangle} implies that $m_{t,t'}<\infty$ for all $t' \in T_{x_2}^o \smallsetminus x_2$, so $t' \preceq_{x_2} t$ by Condition \ref{thm_item:fin_oddconn_acycl_nonftype_almostallinfty} and maximality of $t$.
Thus $T_{x_2}^o=\wedge^t$ and it is the saturated chain from $x_2$ to $t$, say $x_2=z_0 \prec_{x_2} z_1 \prec_{x_2} \cdots \prec_{x_2} z_k=t$ with $k \geq 1$.
Put $m_{p,q}=m_{z_p,z_q}$.
Then $(x_i,z_j) \overset{m_{j,j+1}}{\sim}_O (x_i,z_{j+1})$ and $(x_i,z_j) \in (\mathcal{E}_x^{(O)})_{\mathrm{fin}}$ for $1 \leq j \leq k-1$ by a similar argument to the previous paragraph.

We show that $(z_0,z_1,\dots,z_k)$ is of type $A_{k+1}$ (see (\ref{eq:ruleforfinitepart}) for terminology).
Note that $m_{0,1}=3$ by Condition \ref{thm_item:fin_oddconn_acycl_nonftype_orderideal}; so the claim holds if $k=1$.
Secondly, if $k \geq 2$ and $m_{0,2} \neq 2$, then we have (by Condition \ref{thm_item:fin_oddconn_acycl_nonftype_almostallinfty}) $\left[x_1,z_2\right]_O=\{(x_1,z_2)\}$, $\left[x_3,z_2\right]_O=\{(x_3,z_2)\}$ and $\left[x_1,z_2\right]_O \overset{\infty}{\sim}_O \left[x_3,z_2\right]_O$, contradicting the fact $(x_i,z_2) \in (\mathcal{E}_x^{(O)})_{\mathrm{fin}}$.
Thus we have $m_{0,2}=2$ whenever $k=2$, so $m_{1,2}=3$ (by Condition \ref{thm_item:fin_oddconn_acycl_nonftype_Tx2}), proving the claim if $k=2$.
Moreover, if $k \geq 3$, then $(x_i,z_k) \in (\mathcal{E}_x^{(O)})_{\mathrm{fin}}$ and $m_{x_i,z_k}=m_{x_i,z_{k-2}}=2$, so Lemma \ref{lem:deniedtriangle} implies that $m_{k-2,k}<\infty$, therefore $m_{k-2,k}=2$ by Condition \ref{thm_item:fin_oddconn_acycl_nonftype_almostallinfty}.
Thus the claim follows from Condition \ref{thm_item:fin_oddconn_acycl_nonftype_Tx2}.

Hence we have shown that $T_{x_2}^o$ is of type $A_{k+1}$, which means that $K'$ is of type $D_{k+3}$, contradicting the above assumption on $K'$.
Thus the proof is concluded.
\end{proof}
The rest of this subsection is devoted to the proof of Theorem \ref{thm:finitepart_oddconnected_acyclic}.
The `if' part has been proved, so we show the ``only if'' part.
First we prepare the following two lemmas.
\begin{lem}
\label{lem:lemmaforalmostinfinite}
Under the hypothesis of Theorem \ref{thm:finitepart_oddconnected_acyclic}, let $I=\{x_1,x_2,\dots,x_n\}$ be a subset of $O$ such that $(x_1,x_2,\dots,x_n)$ is of type $A_n$, $H_3$, $H_4$ or $P(m)$ (see (\ref{eq:ruleforfinitepart}) and Example \ref{ex:WperpxforP} for terminology).
Suppose that $x \in I$ and $(\mathcal{E}_x^{(O)})_{\mathrm{fin}} \cap \mathcal{E}_x^{(I)} \neq \emptyset$ (so $n \geq 3$).
Let $x_0 \in O \smallsetminus I$ be a neighbor of $x_i$ with $2 \leq i \leq n-1$ in the tree $\Gamma_O^{\mathrm{odd}}$.
\begin{enumerate}
\item \label{lem_item:lemforalmostinf_twopossibilities}
If $(x_1,\dots,x_n)$ is of type $A_n$, then one of the following two conditions is satisfied:
\begin{itemize}
\item $m_{x_0,x_j}=2$ for $1 \leq j \leq n$ with $j \neq i$, and $m_{x_0,x_i}=3$;
\item $m_{x_0,x_j}=\infty$ for $1 \leq j \leq n$ with $j \neq i$.
\end{itemize}
Moreover, in the former case, we have either $i \in \{2,n-1\}$ (so $I \cup x_0$ is of type $D_{n+1}$), or $5 \leq n \leq 7$ and $i \in \{3,n-2\}$ (so $I \cup x_0$ is of type $E_{n+1}$).
\item \label{lem_item:lemforalmostinf_allinf}
If $(x_1,\dots,x_n)$ is of type $H_3$, $H_4$ or $P(m)$, then $m_{x_0,x_j}=\infty$ for $1 \leq j \leq n$ with $j \neq i$.
\end{enumerate}
\end{lem}
\begin{proof}
Put $J=I \cup x_0$, $m_{j,k}=m_{x_j,x_k}$ and $\xi_{j,k}=(x_j,x_k)$.
Note that $m_{0,j}$ is even or $\infty$ for $1 \leq j \leq n$ with $j \neq i$, since $\Gamma_O^{\mathrm{odd}}$ is acyclic.

First, we consider the case that $(x_1,\dots,x_n)$ is of type $A_n$ or $H_4$.
In this case, Theorem \ref{thm:Wperpxforoddpath} (\ref{thm_item:Wperpxforoddpath_AHP}) implies that $\mathcal{E}_x^{(I)} \subseteq (\mathcal{E}_x^{(O)})_{\mathrm{fin}}$.
Now if $m_{0,j}<\infty$ for all $j$, then $\Gamma_J^{\mathrm{odd}}$ is acyclic and $J$ is irreducible and $2$-spherical; so Corollary \ref{cor:no2sphericalacyclicsupset} implies that $J$ is of finite type, which means that $(x_1,\dots,x_n)$ is of type $A_n$ and $J$ is as in the statement of Claim \ref{lem_item:lemforalmostinf_twopossibilities}.
So suppose that $m_{0,j}=\infty$ for some $j \neq i$.
Assume that $j<i$, since the other case $j>i$ is similar.
Then for $i+1 \leq k \leq n$, since $\xi_{k,j} \in (\mathcal{E}_x^{(O)})_{\mathrm{fin}}$, Lemma \ref{lem:deniedtriangle} implies that $m_{k,0}$ is not even, so $m_{k,0}=\infty$.
Similarly, for $1 \leq k \leq i-1$, we have $m_{k,0}=\infty$ by Lemma \ref{lem:deniedtriangle} since $\xi_{k,i+1} \in (\mathcal{E}_x^{(O)})_{\mathrm{fin}}$ and $m_{0,i+1}=\infty$ as above.
Thus the claim follows in this case.

Secondly, we consider the case that $(x_1,\dots,x_n)$ is of type $H_3$, so $n=3$.
Then $\mathcal{E}_x^{(I)}=\{\xi_{1,3},\xi_{3,1}\}$ and so $\xi_{j,k} \in (\mathcal{E}_x^{(O)})_{\mathrm{fin}}$ where $\{j,k\}=\{1,3\}$.
Now if $m_{0,k}<\infty$, then a direct computation shows that
\[
\left[\xi_{j,k}\right]_J=\{\xi_{j,k}\} \overset{\infty}{\sim}_J \{\xi_{0,k},\xi_{k,0}\} \supseteq \left[\xi_{0,k}\right]_J \textrm{ and so } \left[\xi_{j,k}\right]_J \overset{\infty}{\sim}_J \left[\xi_{0,k}\right]_J,
\]
contradicting the fact $\xi_{j,k} \in (\mathcal{E}_x^{(J)})_{\mathrm{fin}}$ (see Lemma \ref{lem:restrictionoffinitepart}).
Thus we have $m_{0,k}=\infty$, so Lemma \ref{lem:deniedtriangle} implies that $m_{j,0}$ is not even, namely $m_{j,0}=\infty$.
Hence the claim holds.

Finally, we consider the case that $(x_1,\dots,x_n)$ is of type $P(m)$, so $n=4$.
Note that $\xi_{4,1} \in (\mathcal{E}_x^{(O)})_{\mathrm{fin}}$ (by Theorem \ref{thm:Wperpxforoddpath} (\ref{thm_item:Wperpxforoddpath_AHP})) and $\left[\xi_{4,1}\right]_J=\{\xi_{4,1}\}$ (since $m_{3,1}=\infty$).
Now if $1 \leq j \leq i-1$ and $m_{0,j}<\infty$, then we have
\[
\left[\xi_{j,0}\right]_J \subseteq \{\xi_{1,0},\dots,\xi_{i-1,0},\xi_{0,i-1}\} \overset{\infty}{\sim}_J \xi_{4,1} \textrm{ and so } \left[\xi_{4,1}\right]_J \overset{\infty}{\sim}_J \left[\xi_{j,0}\right]_J,
\]
a contradiction (see Lemma \ref{lem:restrictionoffinitepart}).
Thus $m_{0,j}=\infty$ for $1 \leq j \leq i-1$, so $m_{0,4}$ is not even by Lemma \ref{lem:deniedtriangle} (since $m_{0,1}=\infty$), therefore $m_{0,4}=\infty$.
Moreover, if $i=2$, then $m_{0,3}=\infty$; otherwise, we have
\[
\left[\xi_{0,3}\right]_J \subseteq \{\xi_{0,3},\xi_{3,0}\} \overset{\infty}{\sim}_J \xi_{4,1} \textrm{ and so } \left[\xi_{4,1}\right]_J \overset{\infty}{\sim}_J \left[\xi_{0,3}\right]_J,
\]
a contradiction.
Hence the proof is concluded.
\end{proof}
\begin{lem}
\label{lem:conditionfortrunk}
Under the hypothesis of Theorem \ref{thm:finitepart_oddconnected_acyclic}, let $K$ be a subset of $O$ connected in the tree $\Gamma_O^{\mathrm{odd}}$.
Suppose further that $x \in K$, $(\mathcal{E}_x^{(O)})_{\mathrm{fin}} \cap \mathcal{E}_x^{(K)} \neq \emptyset$ (so $K \neq \emptyset$) and
\begin{equation}
\label{eq:conditionfortrunk_adjacentinfty}
\textrm{if } y \in O \smallsetminus K \textrm{ is adjacent to } z \in K \textrm{ in } \Gamma_O^{\mathrm{odd}}, \textrm{ then } m_{y,s}=\infty \textrm{ for } s \in K \smallsetminus z.
\end{equation}
Then $K$ satisfies the condition (\ref{eq:Kistrunk}) in Theorem \ref{thm:finitepart_oddconnected_acyclic_trankisfinitetype}.
\end{lem}
\begin{proof}
Let $\Gamma_O^{\mathrm{odd}}=\Gamma_K^{\mathrm{odd}} \cup \bigcup_{y \in K}T_y$ be the decomposition of $\Gamma_O^{\mathrm{odd}}$ as in Remark \ref{rem:decompintotrees}, and $\xi \in (\mathcal{E}_x^{(O)})_{\mathrm{fin}} \cap \mathcal{E}_x^{(K)}$.
Then (\ref{eq:conditionfortrunk_adjacentinfty}) implies that the subset $\mathcal{E}_x^{(K)}$ of $\mathcal{E}_x^{(O)}$ is closed under the relation $\overset{1}{\sim}_O$, so $\left[\xi\right]_O \subseteq \mathcal{E}_x^{(K)}$ and (by Corollary \ref{cor:necessaryconditionforfinitecomponent})
\begin{equation}
\label{lem_eq:conditionfortrunk_relatedtoK}
\textrm{for any } (s,t) \in \mathcal{E}_x^{(O)} \smallsetminus \mathcal{E}_x^{(K)}, \textrm{ we have } \left[s,t\right]_O \overset{k}{\sim}_O \mathcal{E}_x^{(K)} \textrm{ for some } 2 \leq k<\infty.
\end{equation}

First, we show that $m_{s,t}=\infty$ for $s \in K$ and $t \in O \smallsetminus K$ unless $s$ and $t$ are adjacent in $\Gamma_O^{\mathrm{odd}}$.
Let $y \in K$ such that $t \in T_y \smallsetminus y$.
The claim follows from (\ref{eq:conditionfortrunk_adjacentinfty}) if $t$ is an atom of $T_y$; so suppose not.
Let $z_0 \prec_y z_1 \prec_y \cdots \prec_y z_k$ be the saturated chain from $z_0=y$ to $z_k=t$ (so $k \geq 2$).
If $k=2$, $s \neq y$ and $m_{s,t}<\infty$, then it follows from (\ref{eq:conditionfortrunk_adjacentinfty}) that
\[
\left[t,s\right]_O \subseteq \{(t',s) \mid t' \in \vee_t \textrm{ and } m_{s,t'}<\infty\} \overset{\infty}{\sim}_O \mathcal{E}_x^{(K)} \textrm{ and so } \left[t,s\right]_O \overset{\infty}{\sim}_O \mathcal{E}_x^{(K)},
\]
contradicting (\ref{lem_eq:conditionfortrunk_relatedtoK}).
Thus $m_{s,t}=\infty$ whenever $k=2$ and $s \neq y$.
Moreover, if $k=2$, $s=y$ and $m_{s,t}<\infty$, then $m_{t,z}=\infty$ for all $z \in K \smallsetminus y$; so we have
\[
\left[s,t\right]_O \subseteq \{(y,t)\} \cup \{(t',y) \mid t' \in \vee_t \textrm{ and } m_{y,t'}<\infty\} \overset{\infty}{\sim}_O \mathcal{E}_x^{(K)}
\]
and $\left[s,t\right]_O \overset{\infty}{\sim}_O \mathcal{E}_x^{(K)}$, contradicting (\ref{lem_eq:conditionfortrunk_relatedtoK}).
Thus $m_{s,t}=\infty$ whenever $k=2$.
Finally, if $k \geq 3$ and $m_{s,t}<\infty$, then (since $m_{s,z_2}=\infty$ as above) we have
\[
\left[t,s\right]_O \subseteq \{(t',s) \mid t' \in \vee_{z_3} \textrm{ and } m_{s,t'}<\infty\} \overset{\infty}{\sim}_O \mathcal{E}_x^{(K)} \textrm{ and so } \left[t,s\right]_O \overset{\infty}{\sim}_O \mathcal{E}_x^{(K)},
\]
contradicting (\ref{lem_eq:conditionfortrunk_relatedtoK}).
Thus the claim holds for any $s,t$.

Now our remaining task is to show that $m_{s,t}=\infty$ for any distinct $s,t \in O \smallsetminus K$ which are not adjacent in $\Gamma_O^{\mathrm{odd}}$.
However, the result in the previous paragraph shows that the subset $\mathcal{E}_x^{(K)}$ of $\mathcal{E}_x^{(O)}$ is closed under all the relations $\overset{k}{\sim}_O$ with $k<\infty$, so we have $\mathcal{E}_x^{(O)} \smallsetminus \mathcal{E}_x^{(K)}=\emptyset$ by (\ref{lem_eq:conditionfortrunk_relatedtoK}), proving the claim.
Hence the proof is concluded.
\end{proof}
We come back to the proof of Theorem \ref{thm:finitepart_oddconnected_acyclic}.
Fix a finite irreducible component $G$ of $W_O^{\perp x}$.
We divide the proof into two cases.

\noindent
\textbf{Case 1: There does not exist a path $(x_1,x_2,\dots,x_n)$ in $\Gamma_O^{\mathrm{odd}}$ of type $A_n$ such that $r_x^{(O)}(x_j,x_k) \in G$ for some $j,k \in \{1,2,\dots,n\}$.}

Take a generator $r_x^{(O)}(x',x'')$ of $G$, so $(x',x'') \in (\mathcal{E}_x^{(O)})_{\mathrm{fin}}$.
Then Theorem \ref{thm:Wperpxforoddpath} implies that the non-backtracking path $P$ in the tree $\Gamma_O^{\mathrm{odd}}$ between $x'$ and $x''$ is of type $A_n$ (with $n \geq 3$), $H_3$, $H_4$ or $P(m)$ (otherwise $(\mathcal{E}_{x'}^{(V(P))})_{\mathrm{fin}}=\emptyset$, while $(x',x'') \in (\mathcal{E}_{x'}^{(O)})_{\mathrm{fin}} \cap \mathcal{E}_{x'}^{(V(P))}$ by Remark \ref{rem:Efinisindependentofx}, contradicting Lemma \ref{lem:restrictionoffinitepart}).
However, the same theorem implies further that, unless $P$ is of type $H_3$, any subpath of $P$ of type $A_3$ contradicts the hypothesis of Case 1.
Thus $P$ is of type $H_3$, say $P=(x_1,x_2,x_3)$ with $m_{x_1,x_2}=3$, $m_{x_2,x_3}=5$, $m_{x_1,x_3}=2$ and $\{x',x''\}=\{x_1,x_3\}$.

We show that $K=V(P)$ satisfies the hypothesis of Theorem \ref{thm:finitepart_oddconnected_acyclic_trankisfinitetype}, proving Theorem \ref{thm:finitepart_oddconnected_acyclic} in this case.
By Lemma \ref{lem:conditionfortrunk}, it suffices to verify the condition (\ref{eq:conditionfortrunk_adjacentinfty}).
This follows immediately from Lemma \ref{lem:lemmaforalmostinfinite} (\ref{lem_item:lemforalmostinf_allinf}) in the case $z=x_2$.

Secondly, suppose that $z=x_1$, and we show that $m_{y,x_2}=m_{y,x_3}=\infty$.
If $m_{y,x_3}<\infty$, then since $(x',x'') \in (\mathcal{E}_x^{(O)})_{\mathrm{fin}}$, Theorem \ref{thm:Wperpxforoddpath} (applied to $K'=K \cup y$) implies that $K'$ is of type $H_4$ (note that $K'$ cannot be of type $A_4$ or $P(m)$).
However, now the same theorem also shows that $r_x^{(O)}(y,x_2) \in G$, so the path $(y,z,x_2)$ of type $A_3$ contradicts the hypothesis of Case 1.
Thus we have $m_{y,x_3}=\infty$.
Moreover, if $m_{y,x_2}<\infty$, then we have
\[
\left[y,x_2\right]_{K'} \subseteq \{(y,x_2),(x_2,y)\} \overset{\infty}{\sim}_{K'} \{(x',x'')\}=\left[x',x''\right]_{K'} 
\]
and $\left[x',x''\right]_{K'} \overset{\infty}{\sim}_{K'} \left[y,x_2\right]_{K'}$, contradicting the fact $(x',x'') \in (\mathcal{E}_x^{(O)})_{\mathrm{fin}}$.
Thus we have $m_{y,x_2}=\infty$ as desired.

Finally, the proof of the remaining case $z=x_3$ is similar to the previous paragraph; indeed, now Theorem \ref{thm:Wperpxforoddpath} implies immediately that $m_{x_1,y}=\infty$ (note that now $K'=K \cup y$ cannot be of type $A_4$, $H_4$ or $P(m)$), so $m_{x_2,y}=\infty$ (otherwise we have $\left[x',x''\right]_{K'} \overset{\infty}{\sim}_{K'} \left[x_2,y\right]_{K'}$, a contradiction).
Hence we have verified the condition (\ref{eq:conditionfortrunk_adjacentinfty}), proving the claim in this case.

\noindent
\textbf{Case 2: There exists a path $(x_1,x_2,\dots,x_n)$ in $\Gamma_O^{\mathrm{odd}}$ of type $A_n$ such that $r_x^{(O)}(x_j,x_k) \in G$ for some $j,k \in \{1,2,\dots,n\}$ (so $n \geq 3$).}

Now by Theorem \ref{thm:Wperpxforoddpath} (\ref{thm_item:Wperpxforoddpath_AHP}), all elements $r_x^{(O)}(x_j,x_k)$ with $j,k \in \{1,2,\dots,n\}$ belong to the finite $G$, so the total number $n-2$ of those generators (without repetitions) is bounded above.
We may assume without loss of generality that the quantity $n$ now attains the upper bound, so the path is maximal subject to the condition.
Owing to Theorem \ref{thm:propertyofWperpx} (\ref{thm_item:conjugateinduceisomorphism}), we suppose further that $x \in I=\{x_1,\dots,x_n\}$.
Note that $\mathcal{E}_x^{(I)} \subseteq (\mathcal{E}_x^{(O)})_{\mathrm{fin}}$ as above.

First we prepare the following lemma.
\begin{lem}
\label{lem:extendedfromterminal}
In this case, suppose further that $y \in O \smallsetminus I$ is adjacent to $x_n$ in $\Gamma_O^{\mathrm{odd}}$.
Then one of the following three conditions is satisfied:
\begin{enumerate}
\item \label{lem_item:extendterminal_H4}
$n=3$ and $(x_1,x_2,x_3,y)$ is of type $H_4$;
\item \label{lem_item:extendterminal_P}
$n=3$ and $(y,x_3,x_2,x_1)$ is of type $P(m)$;
\item \label{lem_item:extendterminal_inf}
$m_{x_i,y}=\infty$ for $1 \leq i \leq n-1$.
\end{enumerate}
\end{lem}
\begin{proof}
It suffices to verify Condition \ref{lem_item:extendterminal_inf} under the assumtion that $J=I \cup y$ is not of type $H_4$ or $P(m)$.
Now $J$ cannot be of type $A_{n+1}$ by maximality of $n$, while $(x_1,x_n) \in (\mathcal{E}_x^{(O)})_{\mathrm{fin}}$; so Theorem \ref{thm:Wperpxforoddpath} (applied to $J$) implies that $m_{x_1,y}=\infty$.
Moreover, if $3 \leq i \leq n-1$, then $(x_i,x_1) \in (\mathcal{E}_x^{(O)})_{\mathrm{fin}}$ and $m_{x_1,y}=\infty$, so we have $m_{x_i,y}=\infty$ by Lemma \ref{lem:deniedtriangle}.

Finally, assume contrary that $m_{x_2,y}<\infty$.
Then we have $\left[x_2,y\right]_J=\{(x_2,y)\}$ if $n \geq 4$, and $\left[x_2,y\right]_J \subseteq \{(x_2,y),(y,x_2)\}$ if $n=3$; while
\[
\left[x_1,x_3\right]_J=\{(x_1,x_3),(x_3,x_1),(x_4,x_1),\dots,(x_n,x_1)\} \overset{\infty}{\sim}_J \{(x_2,y),(y,x_2)\}.
\]
This implies that $\left[x_1,x_3\right]_J \overset{\infty}{\sim}_J \left[x_2,y\right]_J$, contradicting the fact $(x_1,x_3) \in (\mathcal{E}_x^{(O)})_{\mathrm{fin}}$.
Hence $m_{x_2,y}=\infty$, so the proof is concluded.
\end{proof}
Now we divide the remaining proof of Case 2 into the following four cases.

\noindent
\textbf{Case 2-1: $n \geq 4$.}

First, we consider the case that, for any $y \in O \smallsetminus I$ adjacent to $x_i$ with $2 \leq i \leq n-1$ in $\Gamma_O^{\mathrm{odd}}$, we have $m_{y,x_j}=\infty$ for $1 \leq j \leq n$ with $j \neq i$.
Now Lemma \ref{lem:extendedfromterminal} implies that the set $K=I$ satisfies the condition (\ref{eq:conditionfortrunk_adjacentinfty}), so the claim follows from Lemma \ref{lem:conditionfortrunk}.

So suppose that a vertex $x_0 \in O \smallsetminus I$ is adjacent to $x_i$ with $2 \leq i \leq n-1$ in $\Gamma_O^{\mathrm{odd}}$, and $m_{x_0,x_j}<\infty$ for some $1 \leq j \leq n$ with $j \neq i$.
In this case, Lemma \ref{lem:lemmaforalmostinfinite} (\ref{lem_item:lemforalmostinf_twopossibilities}) implies (up to symmetry) that $K=I \cup x_0$ is either of type $D_{n+1}$ (with $i=2$) or of type $E_{n+1}$ (with $5 \leq n \leq 7$ and $i=3$).
We verify the condition (\ref{eq:conditionfortrunk_adjacentinfty}) for this $K$, proving the claim by Lemma \ref{lem:conditionfortrunk}.
So take $z \in K$ and $y \in O \smallsetminus K$ adjacent in $\Gamma_O^{\mathrm{odd}}$ with each other.
Put $K'=K \cup y$.

\noindent
\textbf{Step 1: If $z=x_0$, then $m_{x_j,y}=\infty$ for $1 \leq j \leq n$.}

If $K$ is of type $D_{n+1}$, then the path $(x_n,\dots,x_3,x_2,x_0)$ (with $n$ vertices) also satisfies the hypothesis of Case 2, and is maximal subject to the condition as well as the original path $(x_1,\dots,x_n)$.
Thus Lemma \ref{lem:extendedfromterminal} can also be applied to the path, proving that $m_{x_j,y}=\infty$ for $2 \leq j \leq n$.
Moreover, since $(x_1,x_3) \in (\mathcal{E}_x^{(O)})_{\mathrm{fin}}$ and $m_{x_3,y}=\infty$, we have $m_{x_1,y}=\infty$ by Lemma \ref{lem:deniedtriangle}.
Thus the claim follows.

So suppose that $K$ is of type $E_{n+1}$ (with $5 \leq n \leq 7$ and $i=3$).
First, we show that either $m_{x_1,y}=\infty$ or $m_{x_n,y}=\infty$.
If both of them fail, then since $(x_1,x_3),(x_n,x_3) \in (\mathcal{E}_x^{(O)})_{\mathrm{fin}}$, Theorem \ref{thm:Wperpxforoddpath} (applied to the paths $(x_1,x_2,x_3,x_0,y)$ and $(x_n,x_{n-1},\dots,x_3,x_0,y)$ which cannot be of type $H_3$, $H_4$ or $P(m)$) implies that these two paths are of type $A_k$ for some $k$.
However, this means that the set $K'$ is irreducible, $2$-spherical and not of finite type, and $\Gamma_{K'}^{\mathrm{odd}}$ is acyclic; while $(\mathcal{E}_x^{(K)})_{\mathrm{fin}} \neq \emptyset$ (see Lemma \ref{lem:restrictionoffinitepart}).
This contradicts Corollary \ref{cor:no2sphericalacyclicsupset}, so we have either $m_{x_1,y}=\infty$ or $m_{x_n,y}=\infty$.

We may assume that $m_{x_1,y}=\infty$, since the other case is similar.
Now for $3 \leq j \leq n$, since $(x_j,x_1) \in (\mathcal{E}_x^{(O)})_{\mathrm{fin}}$, we have $m_{x_j,y}=\infty$ by Lemma \ref{lem:deniedtriangle}.
Moreover, since $m_{x_4,y}=\infty$ and $(x_2,x_4) \in (\mathcal{E}_x^{(O)})_{\mathrm{fin}}$, we have $m_{x_2,y}=\infty$ by Lemma \ref{lem:deniedtriangle}.
Thus the proof of Step 1 is concluded.

\noindent
\textbf{Step 2: If $z=x_j$ with $2 \leq j \leq n-1$, then $m_{x_k,y}=\infty$ for $0 \leq k \leq n$ with $k \neq j$.}

Assume contrary that the conclusion of Step 2 fails.
Then by Lemma \ref{lem:lemmaforalmostinfinite} (\ref{lem_item:lemforalmostinf_twopossibilities}) applied to $I$ and $y$, we have $m_{y,x_j}=3$, and $m_{y,x_k}=2$ for $1 \leq k \leq n$ with $k \neq j$.
Now if $m_{x_0,y}<\infty$, then the set $K'$ is irreducible, $2$-spherical and not of finite type, and $\Gamma_{K'}^{\mathrm{odd}}$ is acyclic.
This contradicts Corollary \ref{cor:no2sphericalacyclicsupset} (since $(\mathcal{E}_x^{(K)})_{\mathrm{fin}} \neq \emptyset$), so we have $m_{x_0,y}=\infty$.
Moreover, since $i \leq n-2$, $(x_n,x_i) \in (\mathcal{E}_x^{(O)})_{\mathrm{fin}}$ and $(x_n,x_i) \overset{3}{\sim}_O (x_n,x_0)$, we have $(x_n,x_0) \in (\mathcal{E}_x^{(O)})_{\mathrm{fin}}$ by Corollary \ref{cor:inthesamecomponentofWperpx}; while $m_{x_n,y}=2$ as above.
This contradicts Lemma \ref{lem:deniedtriangle}, so the proof of Step 2 is concluded.

\noindent
\textbf{Step 3: If $z=x_j$ with $j \in \{1,n\}$, then $m_{x_k,y}=\infty$ for $0 \leq k \leq n$ with $k \neq j$.}

By Lemma \ref{lem:extendedfromterminal} (applied to $I$ and $y$) and symmetry, we have $m_{x_k,y}=\infty$ for $1 \leq k \leq n$ with $k \neq j$.
This implies that
\[
\left[x_{i-1},x_{i+1}\right]_{K'}=\{(x_1,x_{i+1}),\dots,(x_{i-1},x_{i+1}),(x_{i+1},x_{i-1}),\dots,(x_n,x_{i-1})\}.
\]
Now if $m_{x_0,y}<\infty$, then we have $\left[x_0,y\right]_{K'}=\{(x_0,y)\}$ (since $m_{x_i,y}=\infty$ as above) and so $\left[x_0,y\right]_{K'} \overset{\infty}{\sim}_{K'} \left[x_{i-1},x_{i+1}\right]_{K'}$, contradicting the fact $(x_{i-1},x_{i+1}) \in (\mathcal{E}_x^{(O)})_{\mathrm{fin}}$.
Thus we have $m_{x_0,y}=\infty$, concluding the proof of Step 3.

Hence we have verified the condition (\ref{eq:conditionfortrunk_adjacentinfty}), concluding Case 2-1.

\noindent
\textbf{Case 2-2: $n=3$, and there exists an element $x_4 \in O \smallsetminus I$ such that $K=I \cup x_4$ is of type $H_4$.}

By symmetry, we may assume that $m_{x_3,x_4}=5$.
Since $(x_1,x_3) \in (\mathcal{E}_x^{(O)})_{\mathrm{fin}}$, we have $\mathcal{E}_x^{(K)} \subseteq (\mathcal{E}_x^{(O)})_{\mathrm{fin}}$ by Theorem \ref{thm:Wperpxforoddpath} (\ref{thm_item:Wperpxforoddpath_AHP}).
We verify the condition (\ref{eq:conditionfortrunk_adjacentinfty}) for $K$, which proves the claim similarly.
So take a vertex $y \in O \smallsetminus K$ adjacent to $x_i \in K$ in $\Gamma_O^{\mathrm{odd}}$.
First, Lemma \ref{lem:lemmaforalmostinfinite} (\ref{lem_item:lemforalmostinf_allinf}) implies that, if $i \in \{2,3\}$, then $m_{x_j,y}=\infty$ for $1 \leq j \leq 4$ with $j \neq i$.

Secondly, suppose that $i=1$.
Then Theorem \ref{thm:Wperpxforoddpath} (applied to the path $(y,x_1,\dots,x_4)$) implies that $m_{x_4,y}=\infty$, since this path cannot be of type $A_k$, $H_3$, $H_4$ or $P(m)$.
On the other hand, since $(x_2,x_4) \in (\mathcal{E}_x^{(O)})_{\mathrm{fin}}$, we have $m_{x_2,y}=\infty$ by Lemma \ref{lem:deniedtriangle}.
Moreover, we have $m_{x_3,y}=\infty$, since otherwise (by putting $K'=K \cup y$) we have
\[
\left[x_3,y\right]_{K'}=\{(x_3,y)\} \overset{\infty}{\sim}_{K'} \{(x_4,x_2)\}=\left[x_4,x_2\right]_{K'},
\]
contradicting the fact $(x_4,x_2) \in (\mathcal{E}_x^{(O)})_{\mathrm{fin}}$.

Finally, suppose that $i=4$.
Then Theorem \ref{thm:Wperpxforoddpath} (applied to the paths $(x_1,\dots,x_4,y)$ and $(x_2,x_3,x_4,y)$) implies that $m_{x_1,y}=m_{x_2,y}=\infty$, since $(x_2,x_4) \in (\mathcal{E}_x^{(O)})_{\mathrm{fin}}$ and any of these two paths cannot be of type $A_k$, $H_3$, $H_4$ or $P(m)$.
Moreover, we have $m_{x_3,y}=\infty$, since otherwise we have
\[
\left[x_3,y\right]_{K'}=\{(x_3,y)\} \overset{\infty}{\sim}_{K'} \{(x_4,x_2)\}=\left[x_4,x_2\right]_{K'}
\]
similarly, contradicting the fact $(x_4,x_2) \in (\mathcal{E}_x^{(O)})_{\mathrm{fin}}$.
Hence the condition (\ref{eq:conditionfortrunk_adjacentinfty}) has been verified.

\noindent
\textbf{Case 2-3: $n=3$, and there exists an element $x_0 \in O \smallsetminus I$ such that $K=I \cup x_0$ is of type $P(m)$.}

By symmetry, we may assume that $m_{x_0,x_1}=m$ and $m_{x_0,x_2}=\infty$.
The proof is similar to Case 2-2; we have $\mathcal{E}_x^{(K)} \subseteq (\mathcal{E}_x^{(O)})_{\mathrm{fin}}$, and it suffices to verify the condition (\ref{eq:conditionfortrunk_adjacentinfty}) for $K$.
So take $y \in O \smallsetminus K$ adjacent to $x_i \in K$ in $\Gamma_O^{\mathrm{odd}}$.
The proof follows immediately from Lemma \ref{lem:lemmaforalmostinfinite} (\ref{lem_item:lemforalmostinf_allinf}) if $i \in \{1,2\}$.

Secondly, suppose that $i=3$.
Then Theorem \ref{thm:Wperpxforoddpath} (applied to the path $(x_0,\dots,x_3,y)$) implies that $m_{x_0,y}=\infty$.
Moreover, we have $m_{x_j,y}=\infty$ for $j \in \{1,2\}$, since otherwise (by putting $K'=K \cup y$) we have
\[
\left[x_j,y\right]_{K'} \subseteq \{(x_1,y),(x_2,y),(y,x_2)\} \overset{\infty}{\sim}_{K'} \{(x_3,x_0)\}=\left[x_3,x_0\right]_{K'},
\]
contradicting the fact $(x_3,x_0) \in (\mathcal{E}_x^{(O)})_{\mathrm{fin}}$.

Finally, suppose that $i=0$.
Then Theorem \ref{thm:Wperpxforoddpath} (applied to $(y,x_0,\dots,x_3)$) implies that $m_{x_3,y}=\infty$.
Moreover, we have $m_{x_j,y}=\infty$ for $j \in \{1,2\}$, since otherwise we have
\[
\left[x_j,y\right]_{K'} \subseteq \{(y,x_1),(x_1,y),(x_2,y)\} \overset{\infty}{\sim}_{K'} \{(x_3,x_0)\}=\left[x_3,x_0\right]_{K'},
\]
contradicting the fact $(x_3,x_0) \in (\mathcal{E}_x^{(O)})_{\mathrm{fin}}$.
Hence (\ref{eq:conditionfortrunk_adjacentinfty}) has been verified.

\noindent
\textbf{Case 2-4: $n=3$, and the hypothesis of Case 2-2 or Case 2-3 is not satisfied.}

Our aim here is to show that $K=I$ satisfies the hypothesis of Theorem \ref{thm:finitepart_oddconnected_acyclic_trankisnonfinitetype}.
Our proof traces the proof of the ``only if'' part of Theorem \ref{thm:finitepart_oddconnected_cyclic} (given in Section \ref{sec:specialcase_proof_cyclic}), with slight modification.
For instance, the elements $x_1$ and $x_3$ here play a role of elements $x_1$ and $x_2$ in the proof of Theorem \ref{thm:finitepart_oddconnected_cyclic}.
On the other hand, the analogue (Step 7') of Step 7 in the proof of Theorem \ref{thm:finitepart_oddconnected_cyclic} is trivial (since the set $K^-=K \smallsetminus \{x_1,x_3\}$ now consists of only one element $x_2$), so Step 7' will be skipped in the following proof.

We start the proof.
Note that Condition \ref{thm_item:fin_oddconn_acycl_nonftype_KisA3} of Theorem \ref{thm:finitepart_oddconnected_acyclic_trankisnonfinitetype} holds by the choice of $K$.
Decompose $\Gamma_O^{\mathrm{odd}}$ as $\Gamma_K^{\mathrm{odd}} \cup \bigcup_{i=1}^{3}T_{x_i}$ (see Remark \ref{rem:decompintotrees}).
Note that $(x_1,x_3),(x_3,x_1) \in (\mathcal{E}_x^{(O)})_{\mathrm{fin}}$ since now $(x_1,x_3) \overset{1}{\sim}_O (x_3,x_1)$.

\noindent
\textbf{Step 1': For any $i \in \{1,3\}$, we have $m_{z,x_2}=m_{z,x_{4-i}}=\infty$ if $z$ is an atom of $T_{x_i}$.}

By the hypothesis of Case 2-4 and symmetry, the claim follows immediately from Lemma \ref{lem:extendedfromterminal}.

By Step 1', we have $\left[x_1,x_3\right]_O=\{(x_1,x_3),(x_3,x_1)\}$; and if $(z,z') \in \mathcal{E}_x^{(O)}$ and $(z,z') \overset{m}{\sim}_O \left[x_1,x_3\right]_O$ for some $2 \leq m<\infty$, then we have $z=x_i$ with $i \in \{1,3\}$, $z' \in O \smallsetminus K$ and $m_{x_{4-i},z'}<\infty$.
Thus (since $\left[x_1,x_3\right]_O \subseteq (\mathcal{E}_x^{(O)})_{\mathrm{fin}}$) Corollary \ref{cor:necessaryconditionforfinitecomponent} implies that any $\sim_O$-equivalence class in $\mathcal{E}_x^{(O)}$ other than $\left[x_1,x_3\right]_O$ must be related to $\left[x_1,x_3\right]_O$ via some $\overset{m}{\sim}_O$ with $2 \leq m<\infty$, hence
\begin{equation}
\label{eq:finitepart_oddconnected_acyclic_element}
\textrm{it contains some } (x_i,z) \textrm{ with } i \in \{1,3\},\ z \in O \smallsetminus K \textrm{ and } m_{x_{4-i},z}<\infty.
\end{equation}
This yields the following observation:
\begin{eqnarray}
\nonumber
&&\textrm{if a subset } \mathcal{E}' \textrm{ of } \mathcal{E}_x^{(O)} \smallsetminus \left[x_1,x_3\right]_O \textrm{ is closed under the relation}\\
\label{eq:finitepart_oddconnected_acyclic_equivclass}
&&\overset{1}{\sim}_O \textrm{ and contains no element as in (\ref{eq:finitepart_oddconnected_acyclic_element}), then } \mathcal{E}'=\emptyset.
\end{eqnarray}

\noindent
\textbf{Step 2': For any distinct $s,t \in O$, the order $m_{s,t}$ is $2$, odd or $\infty$.}

\noindent
\textbf{Step 3': For any $i \in \{1,3\}$, we have $m_{z,x_i}=\infty$ if $z \in T_{x_{4-i}} \smallsetminus x_{4-i}$.}

\noindent
\textbf{Step 4': For any $i \in \{1,3\}$, we have $m_{z,z'}=\infty$ if $z \in T_{x_i}$ and $z' \in T_{x_{4-i}} \smallsetminus x_{4-i}$.}

The proofs of Steps 2'--4' are the same as Steps 2--4, respectively, in the proof of Theorem \ref{thm:finitepart_oddconnected_cyclic}.

\noindent
\textbf{Step 5': For any $i \in \{1,3\}$, we have $m_{z,z'}=\infty$ if $z \in T_{x_i} \smallsetminus x_i$ and $z' \in T_{x_2}$.}

In the case $z'=x_2$, Step 1' implies that the set
\[
\mathcal{E}'=\{(t,x_2) \mid t \in T_{x_i} \smallsetminus x_i \textrm{ and } m_{t,x_2}<\infty\}
\]
satisfies the hypothesis of (\ref{eq:finitepart_oddconnected_acyclic_equivclass}), so $\mathcal{E}'=\emptyset$, proving the claim.
The remaining proof is the same as Step 5 in the proof of Theorem \ref{thm:finitepart_oddconnected_cyclic}.

\noindent
\textbf{Step 6': For any $i \in \{1,3\}$, we have $m_{z,z'}=\infty$ if $z,z'$ are distinct, non-adjacent vertices of the tree $T_{x_i}$.
(Hence by Steps 4'--6', if $z,z' \in O$ and $m_{z,z'}=2$, then we have $z,z' \not\in T_{x_i} \smallsetminus x_i$.)}

\noindent
\textbf{Step 8': We have $m_{z,z'}=\infty$ if $z,z'$ are incomparable elements of the poset $T_{x_2}$.
(Hence by Steps 2', 6' and 8', Condition \ref{thm_item:fin_oddconn_acycl_nonftype_almostallinfty} of Theorem \ref{thm:finitepart_oddconnected_acyclic_trankisnonfinitetype} is satisfied.)}

\noindent
\textbf{Step 9': Condition \ref{thm_item:fin_oddconn_acycl_nonftype_orderideal} of Theorem \ref{thm:finitepart_oddconnected_acyclic_trankisnonfinitetype} is satisfied.}

\noindent
\textbf{Step 10': Condition \ref{thm_item:fin_oddconn_acycl_nonftype_Tx2} of Theorem \ref{thm:finitepart_oddconnected_acyclic_trankisnonfinitetype} is satisfied.}

The proofs of Steps 6'--10' are the same as Steps 6--10, respectively, in the proof of Theorem \ref{thm:finitepart_oddconnected_cyclic}.

By Steps 8'--10', the hypothesis of Theorem \ref{thm:finitepart_oddconnected_acyclic_trankisnonfinitetype} is satisfied, concluding Case 2-4.

Hence Case 2, therefore the proof of Theorem \ref{thm:finitepart_oddconnected_acyclic}, is concluded.
\section{The structure of $(W^{\perp x})_{\mathrm{fin}}$}
\label{sec:finitepart}

This section completes our description of the finite part $(W^{\perp x})_{\mathrm{fin}}$ of the Coxeter group $W^{\perp x}$.
Throughout this section, we put
\[
O=S_{\sim x}^{\mathrm{odd}} \textrm{ and } E=\{y \in S \smallsetminus O \mid m_{y,z}<\infty \textrm{ for some } z \in O\}
\]
(see Definition \ref{defn:oddCoxetergraph} for notation).
For $I \subseteq S$ and $s \in S$, write
\[
I_k(s)=\{t \in I \mid m_{s,t}=k\} \textrm{ for } 2 \leq k \leq \infty
\]
and
\[
I_{k_1,k_2,\dots}(s)=I_{k_1}(s) \cup I_{k_2}(s) \cup \cdots \textrm{ and } I_{\mathrm{even}}(s)=\bigcup_{k \textrm{ even}}I_k(s).
\]
Note that $O=O_{\mathrm{even}}(s) \cup O_\infty(s) \neq O_\infty(s)$ for $s \in E$, while $O=O_\infty(s)$ for $s \in S \smallsetminus (O \cup E)$.
For simplicity, the symbol `$S$' in some notations will be omitted unless some ambiguity occurs.
For instance, the set $\mathcal{E}_x^{(O)}$ is a subset of $\mathcal{E}_x=\mathcal{E}_x^{(S)}$ closed under the relation $\overset{1}{\sim}\,=\,\overset{1}{\sim}_S$.
Note also that $s \in E$ if $(y,s) \in \mathcal{E}_x \smallsetminus \mathcal{E}_x^{(O)}$.
\subsection{Lemmas for main theorems}
\label{sec:finitepart_lemma}

Here we prepare some lemmas, summarizing preliminary observations.
\begin{lem}
\label{lem:relationsoutsideO}
Let $(y_1,s_1),(y_2,s_2) \in \mathcal{E}_x \smallsetminus \mathcal{E}_x^{(O)}$, so $s_1,s_2 \in E$.
\begin{enumerate}
\item \label{lem_item:reloutsideO_1simcomp}
If $m_{y_1,s_1} \geq 4$, then $\left[y_1,s_1\right]=\{(y_1,s_1)\}$.
If $m_{y_1,s_1}=2$, then $\left[y_1,s_1\right]=\{(z,s_1) \mid z \in I\}$ and $\iota(\pi_1(\mathcal{I}_x;(y_1,s_1)))=\mathcal{F}_{y_1}^{(I)}$ where $I=O_2(s_1)_{\sim y_1}^{\mathrm{odd}}$.
\item \label{lem_item:reloutsideO_vtx_sissame}
Suppose that $s_1=s_2$ and $y_1 \neq y_2$.
Then
\[
\begin{cases}
(y_1,s_1) \overset{1}{\sim} (y_2,s_2) & \textrm{if } m_{y_1,s_1}=m_{y_2,s_1}=2 \textrm{ and } m_{y_1,y_2} \textrm{ is odd};\\
(y_1,s_1) \overset{2}{\sim} (y_2,s_2) & \textrm{if } m_{y_1,y_2}=3 \textrm{ and } \{m_{y_1,s_1},m_{y_2,s_1}\}=\{2,4\};\\
(y_1,s_1) \overset{\infty}{\sim} (y_2,s_2) & \textrm{otherwise}.
\end{cases}
\]
\item \label{lem_item:reloutsideO_vtx_sisdistinct}
Suppose that $y_1=y_2$ and $s_1 \neq s_2$.
Then
\[
\begin{cases}
(y_1,s_1) \overset{m}{\sim} (y_2,s_2) & \textrm{if } m_{y_1,s_1}=m_{y_1,s_2}=2 \textrm{ and } m_{s_1,s_2}=m<\infty;\\
(y_1,s_1) \overset{2}{\sim} (y_2,s_2) & \textrm{if } m_{s_1,s_2}=2 \textrm{ and } m_{y_1,s_i}=2 \textrm{ for some } i \in \{1,2\};\\
(y_1,s_1) \overset{4}{\sim} (y_2,s_2) & \textrm{if } m_{s_1,s_2}=3 \textrm{ and } \{m_{y_1,s_1},m_{y_1,s_2}\}=\{2,4\};\\
(y_1,s_1) \overset{\infty}{\sim} (y_2,s_2) & \textrm{otherwise}.
\end{cases}
\]
\item \label{lem_item:reloutsideO_comp_sissame}
Suppose that $s_1=s_2$ and $\left[y_1,s_1\right] \neq \left[y_2,s_2\right]$.
Then
\[
\begin{cases}
\left[y_1,s_1\right] \overset{2}{\sim} \left[y_2,s_2\right] & \textrm{if } m_{y_i,s_1}=2,\ m_{y_{3-i},s_1}=4 \textrm{ and } m_{z,y_{3-i}}=3\\
& \qquad \textrm{ for some } i \in \{1,2\} \textrm{ and } z \in O_2(s_1)_{\sim y_i}^{\mathrm{odd}};\\
\left[y_1,s_1\right] \overset{\infty}{\sim} \left[y_2,s_2\right] & \textrm{otherwise}.
\end{cases}
\]
\item \label{lem_item:reloutsideO_comp_sisdistinct}
Suppose that $s_1 \neq s_2$.
Then $\left[y_1,s_1\right] \neq \left[y_2,s_2\right]$, and
\[
\begin{cases}
\left[y_1,s_1\right] \overset{m}{\sim} \left[y_2,s_2\right] & \textrm{if } m_{y_1,s_1}=m_{y_2,s_2}=2,\ m_{s_1,s_2}=m<\infty\\
& \qquad \textrm{ and } O_2(s_1)_{\sim y_1}^{\mathrm{odd}} \cap O_2(s_2)_{\sim y_2}^{\mathrm{odd}} \neq \emptyset;\\
\left[y_1,s_1\right] \overset{2}{\sim} \left[y_2,s_2\right] & \textrm{if } m_{s_1,s_2}=m_{y_i,s_i}=2,\ m_{y_{3-i},s_{3-i}} \neq 2\\
& \qquad \textrm{ and } y_{3-i} \in O_2(s_i)_{\sim y_i}^{\mathrm{odd}} \textrm{ for some } i \in \{1,2\};\\
\left[y_1,s_1\right] \overset{4}{\sim} \left[y_2,s_2\right] & \textrm{if } m_{s_1,s_2}=3,\ m_{y_i,s_i}=2,\ m_{y_{3-i},s_{3-i}}=4\\
& \qquad \textrm{ and } y_{3-i} \in O_2(s_i)_{\sim y_i}^{\mathrm{odd}} \textrm{ for some } i \in \{1,2\};\\
\left[y_1,s_1\right] \overset{\infty}{\sim} \left[y_2,s_2\right] & \textrm{otherwise}.
\end{cases}
\]
\end{enumerate}
\end{lem}
\begin{proof}
Claim \ref{lem_item:reloutsideO_1simcomp} follows from Examples \ref{ex:misevenneq2} and \ref{ex:mis2}, while Claims \ref{lem_item:reloutsideO_vtx_sissame} and \ref{lem_item:reloutsideO_vtx_sisdistinct} follow from definition (see Figure \ref{fig:localsubgraph}).
Claims \ref{lem_item:reloutsideO_comp_sissame} and \ref{lem_item:reloutsideO_comp_sisdistinct} are deduced from the previous claims by straightforward observations.
\end{proof}
\begin{lem}
\label{lem:irrsubsetofE}
Let $E'$ be an irreducible subset of $E$ with $|E'| \geq 2$.
Then the set $\{r_x(c;y,s) \in R_x \mid s \in E'\}$ is an irreducible subset of the generating set $R_x=R_x^{(S)}$ of the Coxeter group $W^{\perp x}$.
\end{lem}
\begin{proof}
By the hypothesis, it suffices to show that two generators $r_x(c;y,s)$ and $r_x(c';y',s')$ do not commute whenever $s,s' \in E'$ and $m_{s,s'} \geq 3$.
Now by Lemma \ref{lem:relationsoutsideO} (\ref{lem_item:reloutsideO_comp_sisdistinct}), we have $\left[y,s\right] \neq \left[y',s'\right]$, and the relation $\left[y,s\right] \overset{2}{\sim} \left[y',s'\right]$ does not hold.
Thus the claim follows from Theorem \ref{thm:presentationofWperpx}.
\end{proof}
\begin{lem}
\label{lem:ncsrconditionforfinoutsideO}
Let $(y,s) \in (\mathcal{E}_x)_{\mathrm{fin}} \smallsetminus \mathcal{E}_x^{(O)}$, so $s \in E$.
\begin{enumerate}
\item \label{lem_item:ncsrforfinoutsideO_mis2}
Suppose that $m_{y,s}=2$.
Then we have $O=O_{2,4,\infty}(s)$, and in the graph $\Gamma_O^{\mathrm{odd}}$, the set $O_2(s)$ is connected and contains all simple closed paths.
If $z \in O_4(s)$, then $m_{z,z'}=3$ for some $z' \in O_2(s)$.
Moreover, for any $t \in E \smallsetminus s$, we have $m_{s,t}<\infty$, every connected component of $\Gamma_{O_2(t)}^{\mathrm{odd}}$ intersects $O_2(s)$, and we have
\[
\begin{cases}
O=O_{2,\infty}(t) & \textrm{if } m_{s,t} \geq 4;\\
O=O_{2,4,\infty}(t) \textrm{ and } O_4(t) \subseteq O_2(s) & \textrm{if } m_{s,t}=3;\\
O_{\mathrm{even}}(t) \smallsetminus O_2(t) \subseteq O_2(s) & \textrm{if } m_{s,t}=2.
\end{cases}
\]
\item \label{lem_item:ncsrforfinoutsideO_mis4}
Suppose that $m_{y,s}=4$.
Then we have $O=O_{2,\infty}(s) \cup y$, the graph $\Gamma_O^{\mathrm{odd}}$ is acyclic, and every connected component of $\Gamma_{O_2(s)}^{\mathrm{odd}}$ contains a vertex $z$ with $m_{z,y}=3$.
Moreover, for any $t \in E \smallsetminus s$, we have $O=O_{2,\infty}(t)$, $m_{y,t}=2$, $m_{s,t} \in \{2,3\}$ and $\Gamma_{O_2(t)}^{\mathrm{odd}}$ is connected.
\item \label{lem_item:ncsrforfinoutsideO_other}
Suppose that $m_{y,s} \geq 6$.
Then $O=O_\infty(s) \cup y$ and $\Gamma_O^{\mathrm{odd}}$ is acyclic.
Moreover, for any $t \in E \smallsetminus s$, we have $O=O_{2,\infty}(t)$, $m_{y,t}=m_{s,t}=2$ and $\Gamma_{O_2(t)}^{\mathrm{odd}}$ is connected.
\end{enumerate}
\end{lem}
\begin{proof}
Here we use Corollary \ref{cor:necessaryconditionforfinitecomponent} several times, without references.

\noindent
\textbf{(\ref{lem_item:ncsrforfinoutsideO_mis2})}
First, Proposition \ref{prop:gammaisindependentofc} and Lemma \ref{lem:relationsoutsideO} (\ref{lem_item:reloutsideO_1simcomp}) imply that
\[
\mathscr{F}_y=\iota(\pi_1(\mathcal{I}_x;(y,s)))=\mathscr{F}_y^{(I)} \textrm{ where } I=O_2(s)_{\sim y}^{\mathrm{odd}},
\]
which means that $\Gamma_I^{\mathrm{odd}}=(\Gamma_{O_2(s)}^{\mathrm{odd}})_{\sim y}$ contains all simple closed paths in $\Gamma_O^{\mathrm{odd}}$.
Let $z \in O_{\mathrm{even}}(s)$.
Then, since $(y,s) \in (\mathcal{E}_x)_{\mathrm{fin}}$, we have $(z,s) \in \left[y,s\right]$ or $\left[z,s\right] \overset{m}{\sim} \left[y,s\right]$ for some $2 \leq m<\infty$.
Now if $m_{z,s}=2$, then we have $\left[z,s\right]=\left[y,s\right]$ and $z \in I$ (by Lemma \ref{lem:relationsoutsideO} (\ref{lem_item:reloutsideO_1simcomp}) and (\ref{lem_item:reloutsideO_comp_sissame})); so $O_2(s)=I$, therefore $\Gamma_{O_2(s)}^{\mathrm{odd}}$ is connected.
On the other hand, if $m_{z,s} \neq 2$, then (by Lemma \ref{lem:relationsoutsideO} (\ref{lem_item:reloutsideO_1simcomp}) and (\ref{lem_item:reloutsideO_comp_sissame})) we have $\left[z,s\right] \neq \left[y,s\right]$ and $\left[z,s\right] \overset{2}{\sim} \left[y,s\right]$; so $z \in O_4(s)$ and $m_{z,z'}=3$ for some $z' \in I=O_2(s)$, therefore $O=O_{2,4,\infty}(s)$.

Finally, suppose that $t \in E \smallsetminus s$, and let $z' \in O_{\mathrm{even}}(t)$.
We use Lemma \ref{lem:relationsoutsideO} (\ref{lem_item:reloutsideO_comp_sisdistinct}) several times.
Then, since $(y,s) \in (\mathcal{E}_x)_{\mathrm{fin}}$, we have $\left[y,s\right] \neq \left[z',t\right]$ and so $\left[z',t\right] \overset{m}{\sim} \left[y,s\right]$ for some $2 \leq m<\infty$.
This implies that $m_{s,t}<\infty$; and we have $O_2(s) \cap O_2(t)_{\sim z'}^{\mathrm{odd}} \neq \emptyset$ whenever $z' \in O_2(t)$, proving that every connected component of $\Gamma_{O_2(t)}^{\mathrm{odd}}$ intersects $O_2(s)$.
Moreover, if $m_{s,t} \geq 4$, then we have $z' \in O_2(t)$ and $\left[y,s\right] \overset{m_{s,t}}{\sim} \left[z',t\right]$; so $O=O_{2,\infty}(t)$.
If $m_{s,t}=3$ and $z' \not\in O_2(t)$, then we have $\left[y,s\right] \overset{4}{\sim} \left[z',t\right]$, $z' \in O_4(t)$ and $z' \in O_2(s)$; so $O=O_{2,4,\infty}(t)$ and $O_4(t) \subseteq O_2(s)$.
Finally, if $m_{s,t}=2$ and $m_{z',t} \neq 2$, then we have $\left[y,s\right] \overset{2}{\sim} \left[z',t\right]$ and $z' \in I=O_2(s)$; so $O_{\mathrm{even}}(t) \smallsetminus O_2(t) \subseteq O_2(s)$.
Hence the claim holds.

\noindent
\textbf{(\ref{lem_item:ncsrforfinoutsideO_mis4})}
First, since $(y,s) \in (\mathcal{E}_x)_{\mathrm{fin}}$, Corollary \ref{cor:miseven} implies that $\Gamma_O^{\mathrm{odd}}$ is acyclic.
Moreover, for any $y' \in O_{\mathrm{even}}(s) \smallsetminus y$, we have $\left[y,s\right] \neq \left[y',s\right]$ (by Lemma \ref{lem:relationsoutsideO} (\ref{lem_item:reloutsideO_1simcomp})) and so $\left[y',s\right] \overset{m}{\sim} \left[y,s\right]$ for some $2 \leq m<\infty$.
Thus Lemma \ref{lem:relationsoutsideO} (\ref{lem_item:reloutsideO_comp_sissame}) implies that $y' \in O_2(s)$ and $m_{z,y}=3$ for some $z \in O_2(s)_{\sim y'}^{\mathrm{odd}}$; therefore $O=O_{2,\infty} \cup y$.

Finally, suppose that $t \in E \smallsetminus s$, and let $z \in O_{\mathrm{even}}(t)$.
Then, since $(y,s) \in (\mathcal{E}_x)_{\mathrm{fin}}$, we have $\left[z,t\right] \overset{m}{\sim} \left[y,s\right]$ for some $2 \leq m<\infty$ (note that $\left[y,s\right] \neq \left[z,t\right]$).
Now by Lemma \ref{lem:relationsoutsideO} (\ref{lem_item:reloutsideO_comp_sisdistinct}), we have either $m_{s,t}=2$, $z \in O_2(t)$ and $y \in O_2(t)_{\sim z}^{\mathrm{odd}}$ (now $\left[y,s\right] \overset{2}{\sim} \left[z,t\right]$); or $m_{s,t}=3$, $z \in O_2(t)$ and $y \in O_2(t)_{\sim z}^{\mathrm{odd}}$ (now $\left[y,s\right] \overset{4}{\sim} \left[z,t\right]$).
Thus we have $O_{\mathrm{even}}(t)=O_2(t)$ (so $O=O_{2,\infty}(t)$), $m_{y,t}=2$, $m_{s,t} \in \{2,3\}$, and every connected component of $\Gamma_{O_2(t)}^{\mathrm{odd}}$ contains the common vertex $y$ (so $\Gamma_{O_2(t)}^{\mathrm{odd}}$ is connected).
Hence the claim holds.

\noindent
\textbf{(\ref{lem_item:ncsrforfinoutsideO_other})}
The graph $\Gamma_O^{\mathrm{odd}}$ is acyclic as well as Claim \ref{lem_item:ncsrforfinoutsideO_mis4}.
Now if $y' \in O_{\mathrm{even}}(s) \smallsetminus y$, then (since $(y,s) \in (\mathcal{E}_x)_{\mathrm{fin}}$) we have either $\left[y,s\right]=\left[y',s\right]$ or $\left[y',s\right] \overset{m}{\sim} \left[y,s\right]$ for some $2 \leq m<\infty$.
However, Lemma \ref{lem:relationsoutsideO} (\ref{lem_item:reloutsideO_1simcomp}) and (\ref{lem_item:reloutsideO_comp_sissame}) show that this is impossible.
Thus we have $O_{\mathrm{even}}(s)=\{y\}$ and so $O=O_\infty(s) \cup y$.

Finally, suppose that $t \in E \smallsetminus s$, and let $z \in O_{\mathrm{even}}(t)$.
Then we have $\left[z,t\right]=\left[y,s\right]$ or $\left[z,t\right] \overset{m}{\sim} \left[y,s\right]$ for some $2 \leq m<\infty$; so Lemma \ref{lem:relationsoutsideO} (\ref{lem_item:reloutsideO_comp_sisdistinct}) implies that $z \in O_2(t)$, $m_{s,t}=2$ and $y \in O_2(t)_{\sim z}^{\mathrm{odd}}$.
Thus we have $O=O_{2,\infty}(t)$, $y \in O_2(t)$ and every connected component of $\Gamma_{O_2(t)}^{\mathrm{odd}}$ contains $y$, which means that $\Gamma_{O_2(t)}^{\mathrm{odd}}$ is connected.
Hence the claim holds.
\end{proof}
\subsection{Main theorems}
\label{sec:finitepart_maintheorem}

Our description of $(W^{\perp x})_{\mathrm{fin}}$ is divided into some parts.
First, we consider the case that $(W^{\perp x})_{\mathrm{fin}} \cap W_O^{\perp x} \neq 1$, or equivalently, $(\mathcal{E}_x)_{\mathrm{fin}} \cap \mathcal{E}_x^{(O)} \neq \emptyset$.
\begin{thm}
\label{thm:finpart_O_typeB}
Suppose that $S$ contains a sequence $(x_1,x_2,\dots,x_n,s_0)$ of type $B_{n+1}$ (see (\ref{eq:ruleforfinitepart}) for terminology) with $n \geq 3$, satisfying the following three conditions:
\begin{enumerate}
\item \label{thm_item:finpart_O_typeB_trunk}
The graph $\Gamma_O^{\mathrm{odd}}$ is acyclic, $s_0 \in E$, and $K=\{x_1,\dots,x_n\}$ is a subset of $O$ (of type $A_n$) satisfying the condition (\ref{eq:Kistrunk}).
\item \label{thm_item:finpart_O_typeB_s0}
The set $\{s_0\}$ is an irreducible component of $E$, $O=O_{2,\infty}(s_0) \cup x_n$, and $\Gamma_{O_2(s_0)}^{\mathrm{odd}}$ is connected.
\item \label{thm_item:finpart_O_typeB_nots0}
For any $s \in E \smallsetminus s_0$, we have $O=O_{2,\infty}(s)$, and the set $O_2(s)$ is connected in $\Gamma_O^{\mathrm{odd}}$ and contains $K$.
\end{enumerate}
Then $W^{\perp x_i}=W_{K \cup s_0}^{\perp x_i} \times W_{E \smallsetminus s_0}$ for any $x_i \in K$.
Hence $(W^{\perp x_i})_{\mathrm{fin}}=W_{K \cup s_0}^{\perp x_i} \times (W_{E \smallsetminus s_0})_{\mathrm{fin}}$ and $(W^{\perp x_i})_{\mathrm{fin}} \cap W_O^{\perp x_i}=W_K^{\perp x_i} \neq 1$.
\end{thm}
\begin{proof}
Condition \ref{thm_item:finpart_O_typeB_trunk} implies that any element of $\mathcal{E}_{x_i}$ satisfies (\ref{eq:gammaisindependentofc}), and $\mathcal{E}_{x_i}^{(O)}=\mathcal{E}_{x_i}^{(K)}=\bigsqcup_{j=3}^{n}\left[x_1,x_j\right]$ (see Example \ref{ex:WperpxforAn}).
On the other hand, owing to Lemma \ref{lem:relationsoutsideO}, Conditions \ref{thm_item:finpart_O_typeB_s0} and \ref{thm_item:finpart_O_typeB_nots0} imply that
\[
\mathcal{E}_{x_i} \smallsetminus \mathcal{E}_{x_i}^{(O)}=\left[x_1,s_0\right] \sqcup \left[x_n,s_0\right] \sqcup \bigsqcup_{s \in E \smallsetminus s_0}\left[x_1,s\right].
\]
Now by Lemma \ref{lem:relationsoutsideO}, for $s \in E \smallsetminus s_0$, we have $(x_1,x_j) \overset{2}{\sim} (x_1,s)$ for $3 \leq j \leq n$ (since $K \subseteq O_2(s)$ by Condition \ref{thm_item:finpart_O_typeB_nots0}), $(x_1,s_0) \overset{2}{\sim} (x_1,s)$ (since $m_{s_0,s}=2$ by Condition \ref{thm_item:finpart_O_typeB_s0}) and $(x_n,s_0) \overset{2}{\sim} (x_n,s) \in \left[x_1,s\right]$ (since $x_n \in O_2(s)_{\sim x_1}^{\mathrm{odd}}$ by Condition \ref{thm_item:finpart_O_typeB_nots0}).
This means that $W^{\perp x_i}$ is the direct product of two subgroups $W_{E \smallsetminus s_0}$ and $W_{K \cup s_0}^{\perp x_i}$, generated by
\[
\{r_{x_i}(x_1,s) \mid s \in E \smallsetminus s_0\}=E \smallsetminus s_0
\]
(note that $K \subseteq O_2(s)$ for any $s \in E \smallsetminus s_0$) and
\[
\{r_{x_i}(x_1,x_j) \mid 3 \leq j \leq n\} \cup \{r_{x_i}(x_1,s_0),r_{x_i}(x_n,s_0)\}=R_{x_i}^{(K \cup s_0)},
\]
respectively.
Hence the claim holds.
\end{proof}
\begin{thm}
\label{thm:finpart_O_trunk_acyclic}
Suppose that the hypothesis of Theorem \ref{thm:finitepart_oddconnected_acyclic_trankisfinitetype} is satisfied, and for any $s \in E$, we have $O=O_{2,\infty}(s)$, and $O_2(s)$ is connected in $\Gamma_O^{\mathrm{odd}}$ and contains $K$.
Then for any $x' \in K$, we have $W^{\perp x'}=W_K^{\perp x'} \times W_E$.
Hence $(W^{\perp x'})_{\mathrm{fin}}=W_K^{\perp x'} \times (W_E)_{\mathrm{fin}}$ and $(W^{\perp x'})_{\mathrm{fin}} \cap W_O^{\perp x'}=W_K^{\perp x'} \neq 1$.
\end{thm}
\begin{proof}
In this case, the graph $\Gamma_O^{\mathrm{odd}}$ is acyclic, while $|W_K^{\perp x'}|<\infty$ by Theorem \ref{thm:finitepart_oddconnected_acyclic_trankisfinitetype}.
Now by Lemma \ref{lem:relationsoutsideO}, we have $\mathcal{E}_{x'}^{(O)}=\mathcal{E}_{x'}^{(K)}$ (see Theorem \ref{thm:finitepart_oddconnected_acyclic_trankisfinitetype}), $\mathcal{E}_{x'} \smallsetminus \mathcal{E}_{x'}^{(O)}=\bigsqcup_{s \in E}\left[x',s\right]$ (since $\Gamma_{O_2(s)}^{\mathrm{odd}}$ is connected and $K \subseteq O_2(s)$), and $\left[x',s\right] \overset{2}{\sim} \left[y,y'\right]$ for any $s \in E$ and $(y,y') \in \mathcal{E}_{x'}^{(K)}$ (since $(y,s) \in \left[x',s\right]$ and $y,y' \in O_2(s)$).
Hence the claim follows, since $\{r_{x'}(x',s) \mid s \in E\}=E$.
\end{proof}
\begin{thm}
\label{thm:finpart_O_trunk_cyclic}
Assume the conditions in one of Theorems \ref{thm:finitepart_oddconnected_cyclic} or \ref{thm:finitepart_oddconnected_acyclic_trankisnonfinitetype}.
Put $(x'_1,x'_2)=(x_1,x_2)$ in the case of Theorem \ref{thm:finitepart_oddconnected_cyclic}, and $(x'_1,x'_2)=(x_1,x_3)$ in the case of Theorem \ref{thm:finitepart_oddconnected_acyclic_trankisnonfinitetype}, respectively.
Suppose further that, for any $s \in E$, we have $O=O_{2,\infty}(s)$, $x'_1,x'_2 \in O_2(s)$ and every connected component of $\Gamma_{O_2(s)}^{\mathrm{odd}}$ contains $x'_1$ or $x'_2$ (so $\Gamma_{O_2(s)}^{\mathrm{odd}}$ consists of at most two connected components).
Put $x'=x'_1$, $K^-=K \smallsetminus \{x'_1,x'_2\}$ and $K'=K \cup \bigcup_{y \in K^-}V(T_y^o)$.
\begin{enumerate}
\item \label{thm:finpart_O_trunk_cyclic_reducedtoabove}
Suppose that the hypothesis of Theorem \ref{thm:finitepart_oddconnected_acyclic_trankisnonfinitetype} is satisfied, $K'$ is of type $A_3$ or $D_n$ (with $n \geq 4$), and for any $s \in E$, the set $O_2(s)$ is connected in $\Gamma_O^{\mathrm{odd}}$ and contains $K'$.
Then the hypothesis of Theorem \ref{thm:finpart_O_trunk_acyclic} is satisfied, with the set $K'$ playing the role of $K$.
\item \label{thm:finpart_O_trunk_cyclic_notreducedtoabove}
Otherwise, we have $(W^{\perp x'})_{\mathrm{fin}}=\langle x'_2 \rangle \times W_{E'}$ and $(W^{\perp x'})_{\mathrm{fin}} \cap W_O^{\perp x'}=\langle x'_2 \rangle \neq 1$, where $E'$ is the union of the irreducible components $I$ of $E$ of finite type such that $K' \subseteq O_2(s)$ for all $s \in I$.
\end{enumerate}
\end{thm}
\begin{proof}
Claim \ref{thm:finpart_O_trunk_cyclic_reducedtoabove} follows from Theorem \ref{thm:finitepart_oddconnected_acyclic_trankisnonfinitetype} and direct verification.
For Claim \ref{thm:finpart_O_trunk_cyclic_notreducedtoabove}, first we show the following lemma.
\begin{lem}
\label{thm_lem:finpart_O_trunk_cyclic}
Under the hypothesis of Theorem \ref{thm:finpart_O_trunk_cyclic}, suppose further that $(\mathcal{E}_{x'})_{\mathrm{fin}} \cap \mathcal{E}_{x'}^{(O)} \not\subseteq \{(x'_1,x'_2),(x'_2,x'_1)\}$.
Then Condition \ref{thm:finpart_O_trunk_cyclic_reducedtoabove} of the theorem is satisfied.
\end{lem}
\begin{proof}
[Proof of Lemma \ref{thm_lem:finpart_O_trunk_cyclic}.]
Now we have $(\mathcal{E}_{x'}^{(O)})_{\mathrm{fin}} \not\subseteq \{(x'_1,x'_2),(x'_2,x'_1)\}$ by Lemma \ref{lem:restrictionoffinitepart}.
Thus by Theorems \ref{thm:finitepart_oddconnected_cyclic} and \ref{thm:finitepart_oddconnected_acyclic_trankisnonfinitetype}, the condition in Theorem \ref{thm:finitepart_oddconnected_cyclic} is not satisfied, and the hypothesis of Theorem \ref{thm:finitepart_oddconnected_acyclic_trankisnonfinitetype} is satisfied, with $K'$ of type $A_3$ or $D_n$ (where $n \geq 4$).
So it follows that the hypothesis of Theorem \ref{thm:finitepart_oddconnected_acyclic_trankisfinitetype} is satisfied (with $K'$ playing the role of $K$).
Note that $K'$ is not of type $A_3$ (since otherwise $(\mathcal{E}_{x'}^{(O)})_{\mathrm{fin}}=\{(x'_1,x'_2),(x'_2,x'_1)\}$, a contradiction), so the order ideal $T_{x_2}^o$ of $T_{x_2}$ is a finite chain $y_1 \prec_{x_2} y_2 \prec_{x_2} \cdots \prec_{x_2} y_k$, where $k=n-2 \geq 2$, $y_1=x_2$, and $(y_1,\dots,y_k)$ is of type $A_k$.

Our remaining task is to show that, for any $s \in E$, the set $O_2(s)$ is connected in the tree $\Gamma_O^{\mathrm{odd}}$ and contains $K'$.
Now it suffices to verify that $y_k \in O_2(s)$.
Indeed, if this holds, then by the hypothesis and the shape of $\Gamma_O^{\mathrm{odd}}$, the connected component $(\Gamma_{O_2(s)}^{\mathrm{odd}})_{\sim y_k}$ containing either $x'_1$ or $x'_2$ must contain the whole $K'$ (note that $x'_1,x'_2 \in O_2(s)$).
Thus every connected component of $\Gamma_{O_2(s)}^{\mathrm{odd}}$ also contains $K'$, so $\Gamma_{O_2(s)}^{\mathrm{odd}}$ is connected, as desired.

Now let $s \in E$, and take an element $(z,z') \in (\mathcal{E}_{x'})_{\mathrm{fin}} \cap \mathcal{E}_{x'}^{(O)}$ with $(z,z') \not\in \{(x'_1,x'_2),(x'_2,x'_1)\}$.
Then $z,z' \in K' \smallsetminus x'_{3-i}$ for some $i \in \{1,2\}$; so we have $(x'_i,y_j) \in \left[z,z'\right]$ for some $2 \leq j \leq k$ (see Example \ref{ex:WperpxforAn}), while $\left[x'_i,y_j\right] \overset{3}{\sim} \left[x'_i,y_{j+1}\right] \overset{3}{\sim} \cdots \overset{3}{\sim} \left[x'_i,y_k\right]$, therefore $(x'_i,y_k) \in (\mathcal{E}_{x'})_{\mathrm{fin}}$ by Corollary \ref{cor:inthesamecomponentofWperpx}.
Since $x'_i \in O_2(s)$ by the hypothesis, we have $m_{y_k,s}<\infty$ by Lemma \ref{lem:deniedtriangle}, so $y_k \in O_2(s)$ by the hypothesis $O=O_{2,\infty}(s)$.
Hence the claim holds.
\end{proof}
We come back to the proof of Claim \ref{thm:finpart_O_trunk_cyclic_notreducedtoabove}.
Now we have $(\mathcal{E}_{x'})_{\mathrm{fin}} \cap \mathcal{E}_{x'}^{(O)} \subseteq \{(x'_1,x'_2),(x'_2,x'_1)\}$ by Lemma \ref{thm_lem:finpart_O_trunk_cyclic}.
On the other hand, by Theorems \ref{thm:finitepart_oddconnected_cyclic} and \ref{thm:finitepart_oddconnected_acyclic_trankisnonfinitetype}, we have $\left[x'_1,x'_2\right]=\{(x'_1,x'_2),(x'_2,x'_1)\}$, and $(x'_1,x'_2)$ satisfies (\ref{eq:gammaisindependentofc}).
Moreover, it is also shown in the proofs of Theorems \ref{thm:finitepart_oddconnected_cyclic} and \ref{thm:finitepart_oddconnected_acyclic_trankisnonfinitetype} that
\begin{eqnarray}
\nonumber
\textrm{for any } (y,y') \in \mathcal{E}_{x'}^{(O)} \smallsetminus \left[x'_1,x'_2\right], \textrm{ we have } (x'_i,s) \in \left[y,y'\right]\\
\label{thm_eq:finpart_O_trunk_cyclic_equivclass}
\textrm{for some } i \in \{1,2\} \textrm{ and } s \in K' \smallsetminus K \textrm{ (so $m_{s,x'_{3-i}}=2$)}.
\end{eqnarray}

We divide the proof into two steps.

\noindent
\textbf{Step 1: $(W^{\perp x'})_{\mathrm{fin}} \subseteq \langle x'_2 \rangle \times W_{E'}$.}

By Proposition \ref{prop:gammaisindependentofc}, any generator of $(W^{\perp x'})_{\mathrm{fin}}$ is of the form $r_{x'}(y_1,y_2)$ with $(y_1,y_2) \in (\mathcal{E}_{x'})_{\mathrm{fin}}$ satisfying (\ref{eq:gammaisindependentofc}).
Since $(\mathcal{E}_{x'})_{\mathrm{fin}} \cap \mathcal{E}_{x'}^{(O)} \subseteq \left[x'_1,x'_2\right]$, we have $r_{x'}(y_1,y_2)=r_{x'}(x'_1,x'_2)=x'_2$ if $(y_1,y_2) \in \mathcal{E}_{x'}^{(O)}$.

We show that $r_{x'}(y_1,y_2) \in E'$ if $(y_1,y_2) \not\in \mathcal{E}_{x'}^{(O)}$, concluding Step 1.
By the hypothesis, we have $O=O_{2,\infty}(y_2)$, so $m_{y_1,y_2}=2$.
Now Lemma \ref{lem:ncsrconditionforfinoutsideO} (\ref{lem_item:ncsrforfinoutsideO_mis2}) shows that the graph $\Gamma_{O_2(y_2)}^{\mathrm{odd}}$ is connected and contains all simple closed paths in $\Gamma_O^{\mathrm{odd}}$.
In the case of Theorem \ref{thm:finitepart_oddconnected_cyclic}, this implies that $K \subseteq O_2(y_2)$, since $K$ is the union of vertex sets of simple closed paths in $\Gamma_O^{\mathrm{odd}}$.
In the case of Theorem \ref{thm:finitepart_oddconnected_acyclic_trankisnonfinitetype}, the connectedness implies that $O_2(y_2)$ contains the unique vertex $x_2$ between two vertices $x'_1$ and $x'_2$ in $O_2(y_2)$ (see the hypothesis), so $K \subseteq O_2(y_2)$.
In any case, we have $K \subseteq O_2(y_2)$, so $(x'_1,y_2) \in \left[y_1,y_2\right]$ (see Lemma \ref{lem:relationsoutsideO} (\ref{lem_item:reloutsideO_1simcomp})), therefore $(x'_1,y_2) \in (\mathcal{E}_{x'})_{\mathrm{fin}}$ and $r_{x'}(y_1,y_2)=r_{x'}(x'_1,y_2)=y_2$.
Now Lemma \ref{lem:deniedtriangle} implies that $m_{z,y_2}<\infty$ for $z \in K' \smallsetminus K$ (since $m_{x'_1,z}=2$), so $K' \smallsetminus K \subseteq O_2(y_2)$ (since $O=O_{2,\infty}(y_2)$) and $K' \subseteq O_2(y_2)$.

Let $I$ be the irreducible component of $E$ containing $y_2$.
If $I=\{y_2\}$, then $I \subseteq E'$ by the above result, proving the claim.
So suppose that $|I| \geq 2$.
Then for $s \in I \smallsetminus y_2$ and $z \in O_{\mathrm{even}}(y_2)$ (note that $O_{\mathrm{even}}(y_2) \neq \emptyset$ since $z \in E$), Lemma \ref{lem:irrsubsetofE} shows that the generator $r_{x'}(z,s)$ belongs to the same (finite) irreducible component of $W^{\perp x'}$ as $r_{x'}(y_1,y_2)$, therefore $(z,s) \in (\mathcal{E}_{x'})_{\mathrm{fin}}$.
Thus we have $K' \subseteq O_2(s)$ by applying the argument in the previous paragraph to $(z,s)$.
Finally, by replacing $z$ with $x'_1$, the above argument implies that all generators $r_{x'}(x'_1,s)=s$ (with $s \in I$) belong to the same finite irreducible component of $W^{\perp x'}$, proving that $I$ is of finite type.
Hence we have $I \subseteq E'$, so Step 1 is concluded.

\noindent
\textbf{Step 2: $\langle x'_2 \rangle \times W_{E'} \subseteq (W^{\perp x'})_{\mathrm{fin}}$.}

Recall that $(x'_1,x'_2)$ satisfies (\ref{eq:gammaisindependentofc}).
First, we show that $r_{x'}(x'_1,x'_2)=x'_2$ commutes with all generators of $W^{\perp x'}$, proving that $\langle x'_2 \rangle$ is an irreducible component of $(W^{\perp x'})_{\mathrm{fin}}$.
It suffices to show that $\left[x'_1,x'_2\right] \overset{2}{\sim} \left[y,y'\right]$ for all $(y,y') \in \mathcal{E}_{x'} \smallsetminus \left[x'_1,x'_2\right]$.
This follows from (\ref{thm_eq:finpart_O_trunk_cyclic_equivclass}) if $(y,y') \in \mathcal{E}_{x'}^{(O)}$, since now $(x'_i,s) \overset{2}{\sim} (x'_i,x'_{3-i})$.
So suppose that $(y,y') \not\in \mathcal{E}_{x'}^{(O)}$ (note that $y \in O_2(y')$ since $O=O_{2,\infty}(y')$).
Then $O_2(y')_{\sim y}^{\mathrm{odd}}$ contains $x'_i$ for some $i \in \{1,2\}$, and $x'_{3-i} \in O_2(y')$, by the hypothesis.
Thus we have $(x'_i,x'_{3-i}) \overset{2}{\sim} (x'_i,y') \in \left[y,y'\right]$ (see Lemma \ref{lem:relationsoutsideO}), as desired.

On the other hand, let $I \subseteq E'$ be an irreducible component of $E$ (of finite type) as in the statement of Claim \ref{thm:finpart_O_trunk_cyclic_notreducedtoabove}.
We show that $W_I$ is an irreducible component of $W^{\perp x'}$, implying that $W_I \subseteq (W^{\perp x'})_{\mathrm{fin}}$ and concluding the proof.
Now for $s \in I$, every connected component of $\Gamma_{O_2(s)}^{\mathrm{odd}}$ contains either $x'_1$ or $x'_2$ (by the hypothesis), so it contains the whole $K'$ in common (since $K' \subseteq O_2(s)$ by the choice of $I$).
This means that $\Gamma_{O_2(s)}^{\mathrm{odd}}$ is connected and contains $\Gamma_{K'}^{\mathrm{odd}}$, so $(x'_1,s) \in \mathcal{E}_{x'}$ satisfies the condition (\ref{eq:gammaisindependentofc}) by the shape of $\Gamma_O^{\mathrm{odd}}$.
Thus we have $\{r_{x'}(x'_1,s) \mid s \in I\}=I$; so it suffices to show that, for any $(y,y') \in \mathcal{E}_{x'}$, we have either $\left[y,y'\right]=\left[x'_1,s\right]$ for some $s \in I$, or $\left[y,y'\right] \overset{2}{\sim} \left[x'_1,s\right]$ for all $s \in I$.

If $\left[y,y'\right]=\left[x'_1,x'_2\right]$, then $(x'_1,x'_2) \overset{2}{\sim} (x'_1,s)$ for all $s \in I$ (since $x'_2 \in O_2(s)$), so the claim follows.
If $(y,y') \in \mathcal{E}_{x'}^{(O)} \smallsetminus \left[x'_1,x'_2\right]$, then (\ref{thm_eq:finpart_O_trunk_cyclic_equivclass}) shows that $(x'_i,t) \in \left[y,y'\right]$ for some $i \in \{1,2\}$ and $t \in K' \smallsetminus K$, while $(x'_i,s) \in \left[x'_1,s\right]$ and $(x'_i,s) \overset{2}{\sim} (x'_i,t)$ for all $s \in I$ (since $K' \subseteq O_2(s)$).
Thus the claim follows in this case.
Finally, suppose that $(y,y') \in \mathcal{E}_{x'} \smallsetminus \mathcal{E}_{x'}^{(O)}$.
If $y' \in I$, then $y \in O_2(y')_{\sim x'_1}^{\mathrm{odd}}$ as above, so $\left[y,y'\right]=\left[x'_1,y'\right]$.
On the other hand, if $y' \in E \smallsetminus I$, then $O_2(y')_{\sim y}^{\mathrm{odd}}$ contains $x'_i$ for some $i \in \{1,2\}$.
Now for $s \in I$, we have $m_{y',s}=2$ since $I$ is an irreducible component of $E$, while $(x'_i,y') \in \left[y,y'\right]$ and $(x'_i,y') \overset{2}{\sim} (x'_i,s) \in \left[x_1,s\right]$ (since $K' \subseteq O_2(s)$).
Thus the claim holds for any $(y,y')$, concluding Step 2.

Hence the proof of Theorem \ref{thm:finpart_O_trunk_cyclic} is concluded.
\end{proof}
The next theorem is our first main result in this paper.
The proof will be given in Section \ref{sec:finitepart_maintheorem_O}.
\begin{thm}
\label{thm:finpart_O}
We have $(\mathcal{E}_x)_{\mathrm{fin}} \cap \mathcal{E}_x^{(O)} \neq \emptyset$ if and only if the hypothesis of one of Theorems \ref{thm:finpart_O_typeB}, \ref{thm:finpart_O_trunk_acyclic} and \ref{thm:finpart_O_trunk_cyclic} is satisfied.
\end{thm}

From now, we consider the case that $(\mathcal{E}_x)_{\mathrm{fin}} \cap \mathcal{E}_x^{(O)}=\emptyset$.
\begin{thm}
\label{thm:finpart_nonO_Otrunk_F4}
Suppose the following conditions:
\begin{enumerate}
\item \label{thm_item:finpart_nonO_Otrunk_F4_trunk}
The graph $\Gamma_O^{\mathrm{odd}}$ is acyclic, and $m_{y,z}$ is odd or $\infty$ for any $y,z \in O$.
\item \label{thm_item:finpart_nonO_Otrunk_F4_F4exist}
There is a sequence $(y',y,s_0,s_1)$ in $S$ of type $F_4$ (see (\ref{eq:ruleforfinitepart}) for terminology) such that $y,y' \in O$, $\{s_0,s_1\}$ is an irreducible component of $E$, $O=O_{2,\infty}(s_0) \cup y$ and $\Gamma_{O_2(s_0)}^{\mathrm{odd}}$ is connected.
\item \label{thm_item:finpart_nonO_Otrunk_F4_E}
For any $s \in E \smallsetminus s_0$, we have $O=O_{2,\infty}(s)$, and $\Gamma_{O_2(s)}^{\mathrm{odd}}$ is connected and contains $y$ and $y'$.
\end{enumerate}
Then we have $W^{\perp y}=W_J^{\perp y} \times W_{E \smallsetminus J}$ where $J=\{y',y,s_0,s_1\}$.
Hence we have $(W^{\perp y})_{\mathrm{fin}}=W_J^{\perp y} \times (W_{E \smallsetminus J})_{\mathrm{fin}}$, so $m_{y,s_0}=4$ and $(y,s_0) \in (\mathcal{E}_{y})_{\mathrm{fin}}$.
\end{thm}
\begin{proof}
By Condition \ref{thm_item:finpart_nonO_Otrunk_F4_trunk}, we have $\mathcal{E}_y^{(O)}=\emptyset$, and any element of $\mathcal{E}_y$ satisfies (\ref{eq:gammaisindependentofc}).
Thus by Conditions \ref{thm_item:finpart_nonO_Otrunk_F4_F4exist}, \ref{thm_item:finpart_nonO_Otrunk_F4_E} and Lemma \ref{lem:relationsoutsideO} (\ref{lem_item:reloutsideO_1simcomp}), we have
\[
\mathcal{E}_y=\left[y,s_0\right] \sqcup \left[y',s_0\right] \sqcup \left[y,s_1\right] \sqcup \bigsqcup_{s \in E \smallsetminus J}\left[y,s\right].
\]
Moreover, three generators $r_y(y,s_0)$, $r_y(y',s_0)$ and $r_y(y,s_1)$ generate $W_J^{\perp y}$, and $r_y(y,s)=s$ for any $s \in E \smallsetminus J$.
Thus our remaining task is to show that every $s \in E \smallsetminus J$ commutes with the first three generators.
Now we have $m_{s_0,s}=m_{s_1,s}=2$ since $\{s_0,s_1\}$ is an irreducible component of $E$ by Condition \ref{thm_item:finpart_nonO_Otrunk_F4_F4exist}; while $y,y' \in O_2(s)$ by Condition \ref{thm_item:finpart_nonO_Otrunk_F4_E}.
Thus we have $J=J_2(s)$, proving the claim.
Hence the proof is concluded.
\end{proof}
\begin{thm}
\label{thm:finpart_nonO_Otrunk_B}
Suppose the following conditions:
\begin{enumerate}
\item \label{thm_item:finpart_nonO_Otrunk_B_trunk}
The graph $\Gamma_O^{\mathrm{odd}}$ is acyclic, and $m_{y,z}$ is odd or $\infty$ for any $y,z \in O$.
\item \label{thm_item:finpart_nonO_Otrunk_B_Bexist}
There is either a sequence $(s_n,s_{n-1},\dots,s_0,y)$ in $S$ of type $B_{n+2}$, or a sequence $(s_0,y)$ in $S$ of type $I_2(m)$ with $m$ even and $m \neq 2$ where we put $n=0$ and $I_2(4)=B_2$ (see (\ref{eq:ruleforfinitepart}) for terminology), such that $y \in O$, $\{s_0,\dots,s_n\}$ is an irreducible component of $E$ and $O=O_\infty(s_0) \cup y$.
\item \label{thm_item:finpart_nonO_Otrunk_B_E}
For any $s \in E \smallsetminus s_0$, we have $O=O_{2,\infty}(s)$, and $\Gamma_{O_2(s)}^{\mathrm{odd}}$ is connected and contains $y$.
\end{enumerate}
Then we have $W^{\perp y}=W_J^{\perp y} \times W_{E \smallsetminus J}$ where $J=\{y,s_0,\dots,s_n\}$.
Hence $(W^{\perp y})_{\mathrm{fin}}=W_J^{\perp y} \times (W_{E \smallsetminus J})_{\mathrm{fin}}$, so $m_{y,s_0} \neq 2$ and $(y,s_0) \in (\mathcal{E}_y)_{\mathrm{fin}}$.
\end{thm}
\begin{proof}
The proof is similar to Theorem \ref{thm:finpart_nonO_Otrunk_F4}.
By Condition \ref{thm_item:finpart_nonO_Otrunk_B_trunk}, we have $\mathcal{E}_y^{(O)}=\emptyset$, and any element of $\mathcal{E}_y$ satisfies (\ref{eq:gammaisindependentofc}).
Thus we have $\mathcal{E}_y=\bigsqcup_{s \in E}\left[y,s\right]$ by Conditions \ref{thm_item:finpart_nonO_Otrunk_B_Bexist}, \ref{thm_item:finpart_nonO_Otrunk_B_E} and Lemma \ref{lem:relationsoutsideO} (\ref{lem_item:reloutsideO_1simcomp}).
Moreover, the generators $r_y(y,s_i)$ with $0 \leq i \leq n$ generate $W_J^{\perp y}$, and $r_y(y,s)=s$ for any $s \in E \smallsetminus J$.
Thus our remaining task is to show that every $s \in E \smallsetminus J$ commutes with the former generators.
Now we have $m_{s_i,s}=2$ for $0 \leq i \leq n$ by Condition \ref{thm_item:finpart_nonO_Otrunk_B_Bexist}, while $y \in O_2(s)$ by Condition \ref{thm_item:finpart_nonO_Otrunk_B_E}.
Thus we have $J=J_2(s)$, proving the claim.
Hence the proof is concluded.
\end{proof}
\begin{thm}
\label{thm:finpart_nonO_Onontrunk_previouscase}
Suppose the following conditions:
\begin{enumerate}
\item \label{thm_item:finpart_nonO_Onontrunk_prev_Oistree}
The graph $\Gamma_O^{\mathrm{odd}}$ is acyclic.
\item \label{thm_item:finpart_nonO_Onontrunk_prev_B2exist}
There is a sequence $(y,s_0)$ in $S$ of type $B_2$ (see (\ref{eq:ruleforfinitepart}) for terminology) such that $y \in O$, $\{s_0\}$ is an irreducible component of $E$, $O=O_{2,\infty}(s_0) \cup y$, $O_2(s_0) \neq \emptyset$ and the set $T=O_2(s_0) \cup y$ is connected in the graph $\Gamma_O^{\mathrm{odd}}$.
Thus $\Gamma_T^{\mathrm{odd}}$ is a tree, so $T$ admits the tree order $\preceq$ with root vertex $y$.
\item \label{thm_item:finpart_nonO_Onontrunk_prev_Tistrunk}
For any $y_1,y_2 \in O$, the order $m_{y_1,y_2}$ is $2$, odd or $\infty$.
Moreover, if $m_{y_1,y_2}=2$, then $y_1,y_2$ are comparable elements of the poset $T$.
\item \label{thm_item:finpart_nonO_Onontrunk_prev_atomofT}
We have $m_{y,y'}=3$ for any atom $y'$ of $T$.
\item \label{thm_item:finpart_nonO_Onontrunk_prev_orderideal}
Suppose that $y_1,y_2 \in T$, $m_{y_1,y_2}=2$ and $y_1 \prec y_2$.
Then we have $m_{z_1,z_2}=2$ if two vertices $z_1,z_2$ of $T$ are not adjacent, $z_1 \prec z_2 \preceq y_2$ and $z_1 \preceq y_1$.
Moreover, if $y_2$ covers some $z$ and this $z$ covers $y_1$ in $T$, then $m_{y_1,z}=m_{z,y_2}=3$.
\item \label{thm_item:finpart_nonO_Onontrunk_prev_E}
For $s \in E \smallsetminus s_0$, we have $O=O_{2,\infty}(s)$, and $\Gamma_{O_2(s)}^{\mathrm{odd}}$ is connected and contains $y$.
\end{enumerate}
Put
\[
T^o=\{y' \in O \mid m_{y,y'}<\infty\} \textrm{ (so $y \in T^o$)}.
\]
Suppose further that $T^o$ is a finite chain $y=y_0 \prec y_1 \prec \cdots \prec y_n$ with $n \geq 2$, $m_{y_{n-2},y_n}=2$, and $y_n \in O_2(s)$ for all $s \in E \smallsetminus s_0$.
Then the hypothesis of Theorem \ref{thm:finpart_O_typeB} is satisfied, where $(y_n,\dots,y_1,y,s_0)$ plays the role of $(x_1,\dots,x_n,s_0)$.
\end{thm}
\begin{proof}
By definition, all atoms of $T$ are contained in $T^o$.
Conditions \ref{thm_item:finpart_nonO_Onontrunk_prev_Tistrunk} and \ref{thm_item:finpart_nonO_Onontrunk_prev_orderideal} imply that $T^o$ is an order ideal of $T$; namely, for any $z \in T^o \smallsetminus y$ which is not an atom of $T$, we have $y \prec z$, and $m_{y,z}=2$ by Condition \ref{thm_item:finpart_nonO_Onontrunk_prev_Tistrunk}, so Condition \ref{thm_item:finpart_nonO_Onontrunk_prev_orderideal} shows that $m_{y,z'}=2$ for any $z' \in \wedge^z \smallsetminus y$ which is not an atom of $T$.
Moreover, Condition \ref{thm_item:finpart_nonO_Onontrunk_prev_Tistrunk} can be strengthened by replacing $T$ with $T^o$; indeed, if $y_1,y_2 \in T$, $m_{y_1,y_2}=2$ and $y_1 \prec y_2$, then Condition \ref{thm_item:finpart_nonO_Onontrunk_prev_orderideal} implies that $m_{y,y_2}=2$ since $y \prec y_2$, so $y_2 \in T^o$, therefore $y_1 \in T^o$ since $T^o$ is an order ideal.
In particular, the set $K=T^o$ satisfies the condition (\ref{eq:Kistrunk}).

Now the chain $T^o$ is saturated since it is an order ideal of $T$.
Since $m_{y_{n-2},y_n}=2$, Conditions \ref{thm_item:finpart_nonO_Onontrunk_prev_atomofT} and \ref{thm_item:finpart_nonO_Onontrunk_prev_orderideal} imply that $(y_0,\dots,y_n)$ is of type $A_{n+1}$, while $m_{y_0,s_0}=4$ and $T^o \smallsetminus y_0 \subseteq T \smallsetminus y_0=O_2(s_0)$ by Condition \ref{thm_item:finpart_nonO_Onontrunk_prev_B2exist}.
Thus $(y_n,\dots,y_1,y,s)$ is of type $B_{n+2}$.
Put $K=T^o$, which satisfies (\ref{eq:Kistrunk}) as above.
We verify the three conditions in Theorem \ref{thm:finpart_O_typeB}, with $(y_n,\dots,y_1,y,s_0)$ playing the role of $(x_1,\dots,x_n,s_0)$.
The condition \ref{thm_item:finpart_O_typeB_trunk} in Theorem \ref{thm:finpart_O_typeB} is already verified.
Secondly, now $T$ possesses a unique atom $y_1$ since $T^o$ does; this means that $O_2(s_0)=T \smallsetminus y$ is connected in $\Gamma_O^{\mathrm{odd}}$, so Condition \ref{thm_item:finpart_nonO_Onontrunk_prev_B2exist} here verifies the condition \ref{thm_item:finpart_O_typeB_s0} in Theorem \ref{thm:finpart_O_typeB}.
Finally, the condition \ref{thm_item:finpart_O_typeB_nots0} in Theorem \ref{thm:finpart_O_typeB} is verified by Conditions \ref{thm_item:finpart_nonO_Onontrunk_prev_Oistree} and \ref{thm_item:finpart_nonO_Onontrunk_prev_E} here; namely, for any $s \in E \smallsetminus s_0$, the connected subgraph $\Gamma_{O_2(s)}^{\mathrm{odd}}$ in the tree $\Gamma_O^{\mathrm{odd}}$ containing $y=y_0$ and $y_n$ also contains all the vertices $y_1,\dots,y_{n-1}$ between them.
Hence the claim holds.
\end{proof}
\begin{thm}
\label{thm:finpart_nonO_Onontrunk}
Suppose that the conditions \ref{thm_item:finpart_nonO_Onontrunk_prev_Oistree}--\ref{thm_item:finpart_nonO_Onontrunk_prev_E} in Theorem \ref{thm:finpart_nonO_Onontrunk_previouscase} are satisfied.
\begin{enumerate}
\item \label{thm_item:finpart_nonO_Onontrunk_emptytrunk}
Suppose further that $m_{y_1,y_2}$ is odd or $\infty$ for any $y_1,y_2 \in O$, $T$ has a unique atom $y'$ and $y' \in O_2(s)$ for all $s \in E \smallsetminus s_0$.
Then $J=\{y,y',s_0\}$ is of type $B_3$, and we have $W^{\perp y}=W_J^{\perp y} \times W_{E \smallsetminus s_0}$.
Hence $(W^{\perp y})_{\mathrm{fin}}=W_J^{\perp y} \times (W_{E \smallsetminus s_0})_{\mathrm{fin}}$, so $m_{y,s_0}=4$ and $(y,s_0) \in (\mathcal{E}_y)_{\mathrm{fin}}$.
\item \label{thm_item:finpart_nonO_Onontrunk_other}
Suppose further that the hypothesis of Theorem \ref{thm:finpart_nonO_Onontrunk_previouscase} or Claim \ref{thm_item:finpart_nonO_Onontrunk_emptytrunk} is not satisfied.
Then $(W^{\perp y})_{\mathrm{fin}}=\langle r_y(y,s_0) \rangle \times W_{E'}$, so $m_{y,s_0}=4$ and $(y,s_0) \in (\mathcal{E}_y)_{\mathrm{fin}}$, where $E'$ is the union of the irreducible components $K$ of $E \smallsetminus s_0$ of finite type such that $T^o \subseteq O_2(s)$ for any $s \in K$.
\end{enumerate}
\end{thm}
\begin{proof}
First we prepare two lemmas.
\begin{lem}
\label{thm_lem:finpart_nonO_Onontrunk_previoustheorem}
Suppose that the conditions \ref{thm_item:finpart_nonO_Onontrunk_prev_Oistree}--\ref{thm_item:finpart_nonO_Onontrunk_prev_E} in Theorem \ref{thm:finpart_nonO_Onontrunk_previouscase} are satisfied, and $(\mathcal{E}_x)_{\mathrm{fin}} \cap \mathcal{E}_x^{(O)} \neq \emptyset$.
Then the hypothesis of Theorem \ref{thm:finpart_nonO_Onontrunk_previouscase} is satisfied.
\end{lem}
\begin{proof}
[Proof of Lemma \ref{thm_lem:finpart_nonO_Onontrunk_previoustheorem}.]
By Theorem \ref{thm:finpart_O}, the hypothesis of one of Theorems \ref{thm:finpart_O_typeB}, \ref{thm:finpart_O_trunk_acyclic} and \ref{thm:finpart_O_trunk_cyclic} is satisfied.
Moreover, now $O \neq O_{2,\infty}(s_0)$ by Condition \ref{thm_item:finpart_nonO_Onontrunk_prev_B2exist}; so the latter two possibilities are denied, and the hypothesis of Theorem \ref{thm:finpart_O_typeB} is satisfied.
The element $s_0$ here plays the role of $s_0$ in Theorem \ref{thm:finpart_O_typeB}, since this is the unique element $s$ of $E$ with $O \neq O_{2,\infty}(s)$, by Condition \ref{thm_item:finpart_nonO_Onontrunk_prev_E} here and the condition \ref{thm_item:finpart_O_typeB_nots0} in Theorem \ref{thm:finpart_O_typeB}.
Similarly, Condition \ref{thm_item:finpart_nonO_Onontrunk_prev_B2exist} here and the condition \ref{thm_item:finpart_O_typeB_s0} imply that $x_n=y$.
The condition \ref{thm_item:finpart_O_typeB_trunk} in Theorem \ref{thm:finpart_O_typeB} shows that $K$ is a saturated chain $y=x_n \prec x_{n-1} \prec \cdots \prec x_1$, $K \subseteq T^o$ and $m_{x_3,x_1}=2$; while $x_1 \in O_2(s)$ for all $s \in E \smallsetminus s_0$ by the condition \ref{thm_item:finpart_O_typeB_nots0} in Theorem \ref{thm:finpart_O_typeB}.

Finally, we show that $T^o \subseteq K$, concluding the proof.
By the condition \ref{thm_item:finpart_O_typeB_s0} in Theorem \ref{thm:finpart_O_typeB}, the set $O_2(s_0)=T \smallsetminus y$ is connected in $\Gamma_O^{\mathrm{odd}}$ (see Condition \ref{thm_item:finpart_nonO_Onontrunk_prev_B2exist} here), which means that $T$ possesses a unique atom, namely $z_{n-1} \in K$.
Moreover, if $z \in T^o \smallsetminus y$ is not an atom, then $m_{y,z}=2$ and so $z \in K$, since $K$ satisfies (\ref{eq:Kistrunk}) by the condition \ref{thm_item:finpart_O_typeB_trunk} in Theorem \ref{thm:finpart_O_typeB}.
Hence the claim holds.
\end{proof}
\begin{lem}
\label{thm_lem:finpart_nonO_Onontrunk_previousclaim}
Suppose that the conditions \ref{thm_item:finpart_nonO_Onontrunk_prev_Oistree}--\ref{thm_item:finpart_nonO_Onontrunk_prev_E} in Theorem \ref{thm:finpart_nonO_Onontrunk_previouscase} are satisfied, we have $(\mathcal{E}_x)_{\mathrm{fin}} \cap \mathcal{E}_x^{(O)}=\emptyset$, and $(z,s_0) \in (\mathcal{E}_x)_{\mathrm{fin}}$ for some $z \in O_2(s_0)$.
Then the hypothesis of Theorem \ref{thm:finpart_nonO_Onontrunk} (\ref{thm_item:finpart_nonO_Onontrunk_emptytrunk}) is satisfied.
\end{lem}
\begin{proof}
[Proof of Lemma \ref{thm_lem:finpart_nonO_Onontrunk_previousclaim}.]
By Lemma \ref{lem:ncsrconditionforfinoutsideO} (\ref{lem_item:ncsrforfinoutsideO_mis2}), the nonempty set $O_2(s_0)=T \smallsetminus y$ is connected in $\Gamma_O^{\mathrm{odd}}$ (see Condition \ref{thm_item:finpart_nonO_Onontrunk_prev_B2exist}).
Thus $T$ possesses at most one atom, while $T \neq \{y\}$; therefore $T$ has a unique atom, say $y'$.
Since $\Gamma_{O_2(s_0)}^{\mathrm{odd}}$ is connected, Lemma \ref{lem:relationsoutsideO} (\ref{lem_item:reloutsideO_1simcomp}) implies that $(y',s_0) \in \left[z,s_0\right]$, so $(y',s_0) \in (\mathcal{E}_x)_{\mathrm{fin}}$.
Now for $s \in E \smallsetminus s_0$, some $z' \in O_2(s)$ exists by Condition \ref{thm_item:finpart_nonO_Onontrunk_other}, so we have $\left[y',s_0\right]=\left[z',s\right]$ or $\left[y',s_0\right] \overset{m}{\sim} \left[z',s\right]$ for some $2 \leq m<\infty$.
Thus Lemma \ref{lem:relationsoutsideO} (\ref{lem_item:reloutsideO_comp_sisdistinct}) implies that some $z'' \in O_2(s) \cap O_2(s_0)$ exists.
By Condition \ref{thm_item:finpart_nonO_Onontrunk_prev_E}, the connected subgraph $\Gamma_{O_2(s)}^{\mathrm{odd}}$ of $\Gamma_O^{\mathrm{odd}}$ contains all vertices between $z'' \in O_2(s_0)=T \smallsetminus y$ and $y$, so an atom of $T$, which is $y'$ as above.

Thus our remaining task is to show that $m_{y_1,y_2}$ is odd or $\infty$ for any $y_1,y_2 \in O$.
Assume contrary that $m_{y_1,y_2}=2$, $y_1,y_2 \in T$ and $y_1 \prec y_2$ (see Condition \ref{thm_item:finpart_nonO_Onontrunk_prev_Tistrunk}).
Then by the first paragraph of the proof of Theorem \ref{thm:finpart_nonO_Onontrunk_previouscase}, we have $y_1,y_2 \in T^o$ and $T^o$ is an order ideal of $T$.
Since $y_2$ cannot be the unique atom $y'$ of $T$, some element $y_3$ of $T^o$ covers $y'$ in $T$.
Now if $y_3 \in O_2(s_0)$, then we have $(y_3,y) \overset{4}{\sim} (y_3,s_0) \sim (y',s_0)$ and so $(y_3,y) \in (\mathcal{E}_x)_{\mathrm{fin}}$; however, this contradicts the hypothesis $(\mathcal{E}_x)_{\mathrm{fin}} \cap \mathcal{E}_x^{(O)}=\emptyset$.
On the other hand, if $y_3 \not\in O_2(s_0)$, then Lemma \ref{lem:relationsoutsideO} (\ref{lem_item:reloutsideO_1simcomp}) implies that
\[
\left[y',s_0\right] \subseteq \{(t,s_0) \mid t \in \vee_{y'} \smallsetminus \vee_{y_3} \textrm{ and } m_{t,s_0}=2\},
\]
while Condition \ref{thm_item:finpart_nonO_Onontrunk_prev_Tistrunk} implies that
\[
\left[y,y_3\right] \subseteq \{(y,y_3)\} \cup \{(t,y) \mid t \in \vee_{y_3} \textrm{ and } m_{t,y}=2\}.
\]
Thus we have $\left[y',s_0\right] \overset{\infty}{\sim} \left[y,y_3\right]$, contradicting the fact $(y',s_0) \in (\mathcal{E}_x)_{\mathrm{fin}}$.
In any case, we have a contradiction, concluding the proof.
\end{proof}
Now we come back to the proof of Theorem \ref{thm:finpart_nonO_Onontrunk}.
First, Condition \ref{thm_item:finpart_nonO_Onontrunk_prev_Oistree} implies that any element of $\mathcal{E}_y$ satisfies (\ref{eq:gammaisindependentofc}).

\noindent
\textbf{(\ref{thm_item:finpart_nonO_Onontrunk_emptytrunk})}
Since $T$ has a unique atom $y'$ by the hypothesis, Condition \ref{thm_item:finpart_nonO_Onontrunk_prev_B2exist} implies that the set $O_2(s_0)=T \smallsetminus y$ contains $y'$ and is connected in $\Gamma_O^{\mathrm{odd}}$.
Thus by the hypothesis, Conditions \ref{thm_item:finpart_nonO_Onontrunk_prev_B2exist} and \ref{thm_item:finpart_nonO_Onontrunk_prev_E} and Lemma \ref{lem:relationsoutsideO}, we have $\mathcal{E}_x^{(O)}=\emptyset$,
\[
\mathcal{E}_x=\left[y,s_0\right] \sqcup \left[y',s_0\right] \sqcup \bigsqcup_{s \in E \smallsetminus s_0}\left[y,s\right],
\]
and for any $s \in E \smallsetminus s_0$, we have $y,y' \in O_2(s)$, $(y',s) \in \left[y,s\right]$ and $m_{s,s_0}=2$.
Moreover, Condition \ref{thm_item:finpart_nonO_Onontrunk_prev_atomofT} shows that $m_{y,y'}=3$, so $J$ is of type $B_3$.

Now two generators $r_y(y,s_0)$ and $r_y(y',s_0)$ generate $W_J^{\perp y}$, while $r_y(y,s)=s$ for all $s \in E \smallsetminus s_0$.
Moreover, the above properties show that $(y,s_0) \overset{2}{\sim} (y,s)$ and $(y',s_0) \overset{2}{\sim} (y',s) \in \left[y,s\right]$ for all $s \in E \smallsetminus s_0$.
Thus for all $s \in E \smallsetminus s_0$, the generator $r_y(y,s)=s$ commutes with $r_y(y,s_0)$ and $r_y(y',s_0)$.
Hence the claim follows.

\noindent
\textbf{(\ref{thm_item:finpart_nonO_Onontrunk_other})}
In this case, Lemmas \ref{thm_lem:finpart_nonO_Onontrunk_previoustheorem} and \ref{thm_lem:finpart_nonO_Onontrunk_previousclaim} implies that $(\mathcal{E}_y)_{\mathrm{fin}} \cap \mathcal{E}_y^{(O)}=\emptyset$ and $(z,s_0) \not\in (\mathcal{E}_y)_{\mathrm{fin}}$ for any $z \in O_2(s_0)$ (see Remark \ref{rem:Efinisindependentofx}).
We divide the proof into two steps.

\noindent
\textbf{Step 1: $(W^{\perp y})_{\mathrm{fin}} \subseteq \langle r_y(y,s_0) \rangle \times W_{E'}$.}

By the above remark and Conditions \ref{thm_item:finpart_nonO_Onontrunk_prev_B2exist} and \ref{thm_item:finpart_nonO_Onontrunk_prev_E}, Lemma \ref{lem:relationsoutsideO} (\ref{lem_item:reloutsideO_1simcomp}) implies that any element of $(\mathcal{E}_y)_{\mathrm{fin}}$ belongs to either $\left[y,s_0\right]$ or $\left[y,s\right]$ for some $s \in E \smallsetminus s_0$.
Thus it suffices to show that $r_y(y,s) \in E'$ for all $s \in E \smallsetminus s_0$ such that $(y,s) \in (\mathcal{E}_y)_{\mathrm{fin}}$.
Note that $r_y(y,t)=t$ for all $t \in E \smallsetminus s_0$ by Condition \ref{thm_item:finpart_nonO_Onontrunk_prev_E}.

So suppose that $s \in E \smallsetminus s_0$ and $(y,s) \in (\mathcal{E}_y)_{\mathrm{fin}}$, and let $K$ be the irreducible component of $E \smallsetminus s_0$ containing $s$.
We show that $K \subseteq E'$, which yields $r_y(y,s)=s \in E'$ and concluding Step 1.
Now by the choice of $K$, all generators $r_y(y,t)=t$ with $t \in K$ belong to the same (finite) irreducible component of $W^{\perp y}$ as $r_y(y,s)=s$; therefore $K$ is also of finite type.
On the other hand, let $t \in K$, so $(y,t) \in (\mathcal{E}_y)_{\mathrm{fin}}$ as above.
Then we have $z \in O_2(t)$ if $z \in T^o \smallsetminus y$ and $z$ is not an atom of $T$; indeed, now $m_{y,z}=2$, $(y,t) \in (\mathcal{E}_y)_{\mathrm{fin}}$ and $O=O_{2,\infty}(t)$ by Condition \ref{thm_item:finpart_nonO_Onontrunk_prev_E}, so $z \in O_2(t)$ by Lemma \ref{lem:deniedtriangle}.
Moreover, if $z$ is an atom of $T$ and $z \not\in O_2(t)$, then $y \not\in O_2(s_0)$ and $z \in T \smallsetminus y=O_2(s_0)$ by Condition \ref{thm_item:finpart_nonO_Onontrunk_prev_B2exist}, so we have $O_2(t)_{\sim y}^{\mathrm{odd}} \subseteq T \smallsetminus \vee_z$ and $O_2(s_0)_{\sim z}^{\mathrm{odd}} \subseteq \vee_z$.
Now Lemma \ref{lem:relationsoutsideO} implies that $\left[y,t\right] \neq \left[z,s_0\right]$ and $\left[y,t\right] \overset{\infty}{\sim} \left[z,s_0\right]$, contradicting the fact $(y,t) \in (\mathcal{E}_y)_{\mathrm{fin}}$.
Thus $O_2(t)$ contains all atoms of $T$, so $T^o \subseteq O_2(t)$.
Hence we have $K \subseteq E'$ as desired.

\noindent
\textbf{Step 2: $\langle r_y(y,s_0) \rangle \times W_{E'} \subseteq (W^{\perp y})_{\mathrm{fin}}$.}

Note that, by Condition \ref{thm_item:finpart_nonO_Onontrunk_prev_B2exist}, we have $O_2(s_0)=T \smallsetminus y$ and every connected component of $\Gamma_{O_2(s_0)}^{\mathrm{odd}}$ contains an atom of $T$.
Let $\Lambda$ be the set of the atoms of $T$, so $\Lambda \subseteq O_2(s_0)$.
Then by Conditions \ref{thm_item:finpart_nonO_Onontrunk_prev_B2exist}, \ref{thm_item:finpart_nonO_Onontrunk_prev_Tistrunk} and \ref{thm_item:finpart_nonO_Onontrunk_prev_E}, and Lemma \ref{lem:relationsoutsideO} (\ref{lem_item:reloutsideO_1simcomp}), any generator of $W^{\perp y}$ is of one of the following five forms: $r_y(y,s_0)$; $r_y(y',s_0)$ with $y' \in \Lambda$; $r_y(y,s)=s$ with $s \in E \smallsetminus s_0$; $r_y(y_1,y_2)$ with $y_1,y_2 \in T$, $y_1 \prec y_2$ and $m_{y_1,y_2}=2$; and $r_y(y_2,y_1)$ with $y_1,y_2 \in T$, $y_1 \prec y_2$ and $m_{y_1,y_2}=2$.

We show that any generator of one of the last two forms coincides with some $r_y(y,y')$ such that $y' \in T^o$.
Suppose that $y_1,y_2 \in T$, $y_1 \prec y_2$ and $m_{y_1,y_2}=2$.
Then, as shown in the first paragraph of the proof of Theorem \ref{thm:finpart_nonO_Onontrunk_previouscase}, we have $y_1,y_2 \in T^o$, and $T^o$ is an order ideal of $T$.
Take the saturated chain $z_0 \prec z_1 \prec \cdots \prec z_k$ in $T^o$ from $z_o=y$ to $z_k=y_2$, and let $y_1=z_\ell$ with $\ell \leq k-2$.
Now Condition \ref{thm_item:finpart_nonO_Onontrunk_prev_orderideal} implies that $m_{z_i,z_k}=2$ for all $0 \leq i \leq \ell$, so we have $(z_0,z_k) \in \left[y_1,y_2\right]$ and $r_y(y_1,y_2)=r_y(y,z_k)$.
On the other hand, Condition \ref{thm_item:finpart_nonO_Onontrunk_prev_orderideal} also shows that $m_{z_\ell,z_i}=2$ for all $\ell+2 \leq i \leq k$, $m_{z_\ell,z_{\ell+1}}=m_{z_{\ell+1},z_{\ell+2}}=3$, and $m_{z_i,z_{\ell+2}}=2$ for all $0 \leq i \leq \ell$.
Thus we have $(z_0,z_{\ell+2}) \in \left[z_k,z_\ell\right]$ and $r_y(y_2,y_1)=r_y(y,z_{\ell+2})$, as desired.

By using Lemma \ref{lem:relationsoutsideO}, we show that $r_y(y,s_0)$ commutes with all generators of $W^{\perp y}$, proving that $r_y(y,s_0) \in (W^{\perp y})_{\mathrm{fin}}$.
First, for $y' \in \Lambda$, we have $m_{y,y'}=3$ by Condition \ref{thm_item:finpart_nonO_Onontrunk_prev_atomofT}, $m_{y,s_0}=4$ and $m_{y',s_0}=2$; therefore $\left[y,s_0\right] \overset{2}{\sim} \left[y',s_0\right]$.
Secondly, Condition \ref{thm_item:finpart_nonO_Onontrunk_prev_B2exist} shows that $\{s_0\}$ is an irreducible component of $E$, so $m_{s,s_0}=2$ and $\left[y,s_0\right] \overset{2}{\sim} \left[y,s\right]$ for any $s \in E \smallsetminus s_0$.
Moreover, if $y' \in T^o$ and $(y,y') \in \mathcal{E}_y$, then $y' \in T \smallsetminus y=O_2(s_0)$, while $y' \not\in \Lambda$ and so $m_{y,y'}=2$.
Thus we have $\left[y,s_0\right] \overset{2}{\sim} \left[y,y'\right]$, so the claim holds.

Let $K$ be an irreducible component of $E \smallsetminus s_0$ such that $K \subseteq E'$; so $K$ is of finite type and $T^o \subseteq O_2(s)$ for any $s \in K$.
Then the generators $r_y(y,s)=s$ with $s \in K$ generate a finite group.
We show that, for any $s \in K$, the generator $r_y(y,s)$ commutes with all generators of $W^{\perp y}$ other than $r_y(y,t)$ with $t \in K$.
Once this is shown, the generators $r_y(y,s)=s$ with $s \in K$ generate a finite irreducible component of $W^{\perp y}$, so $K \subseteq (W^{\perp y})_{\mathrm{fin}}$ and Step 2 is concluded.
The generator $r_y(y,s)$ commutes with $r_y(y,s_0)$ as shown above.
We have $(y',s_0) \overset{2}{\sim} (y',s) \in \left[y,s\right]$ for any $y' \in \Lambda$, since $y' \in T^o \subseteq O_2(s)$ as above and $m_{s,s_0}=2$ by Condition \ref{thm_item:finpart_nonO_Onontrunk_prev_B2exist}.
For any $t \in E \smallsetminus (K \cup s_0)$, we have $m_{s,t}=2$ by the choice of $K$, so $\left[y,s\right] \overset{2}{\sim} \left[y,t\right]$.
Moreover, if $y' \in T^o$ and $(y,y') \in \mathcal{E}_y$, then we have $m_{y,y'}=2$ and $y' \in O_2(s)$, so $\left[y,s\right] \overset{2}{\sim} \left[y,y'\right]$.
Hence the claim follows.
\end{proof}
Our second main result is the following, proved in Section \ref{sec:finitepart_maintheorem_nonO}.
\begin{thm}
\label{thm:finpart_nonO}
Suppose that $(\mathcal{E}_x)_{\mathrm{fin}} \cap \mathcal{E}_x^{(O)}=\emptyset$.
Then we have $m_{y,s_0} \neq 2$ for some $(y,s_0) \in (\mathcal{E}_x)_{\mathrm{fin}}$ if and only if the hypothesis of one of Theorems \ref{thm:finpart_nonO_Otrunk_F4}, \ref{thm:finpart_nonO_Otrunk_B} and \ref{thm:finpart_nonO_Onontrunk} is satisfied.
\end{thm}
Finally, the next theorem concludes our description of $(W^{\perp x})_{\mathrm{fin}}$.
The proof will be given in Section \ref{sec:finitepart_maintheorem_generalcase}.
\begin{thm}
\label{thm:finpart_generalcase}
Suppose that the hypothesis of Theorem \ref{thm:finpart_O_typeB}, \ref{thm:finpart_O_trunk_acyclic}, \ref{thm:finpart_O_trunk_cyclic}, \ref{thm:finpart_nonO_Otrunk_F4}, \ref{thm:finpart_nonO_Otrunk_B} or \ref{thm:finpart_nonO_Onontrunk} is not satisfied.
Let $\mathcal{K}_1$ be the family of the irreducible components $K$ of $E$ of finite type satisfying the following conditions:
\begin{enumerate}
\item \label{thm_item:finpart_general_family1_O2containtrunk}
For any $s \in K$, we have $O=O_{2,\infty}(s)$, the graph $\Gamma_{O_2(s)}^{\mathrm{odd}}$ is connected and contains all simple closed paths in $\Gamma_O^{\mathrm{odd}}$, and
\begin{eqnarray}
\label{thm_eq:finpart_general_1_trunk_trunk}
\textrm{if } y_1,y_2 \in O \textrm{ and } m_{y_1,y_2} \textrm{ is even, then } y_1,y_2 \in O_2(s).
\end{eqnarray}
\item \label{thm_item:finpart_general_family1_allO2intersect}
We have $\bigcap_{s \in K}O_2(s) \neq \emptyset$.
\item \label{thm_item:finpart_general_family1_O2containOeven}
If $s \in K$ and $t \in E \smallsetminus K$, then
\begin{eqnarray}
\nonumber
O_{\mathrm{even}}(t) \smallsetminus O_2(t) \subseteq O_2(s), \textrm{ and every connected}\\
\label{thm_eq:finpart_general_1_Oeven_even}
\textrm{component of } \Gamma_{O_2(t)}^{\mathrm{odd}} \textrm{ intersects the set } O_2(s).
\end{eqnarray}
\end{enumerate}
Let $\mathcal{K}_2$ be the family of the irreducible components of $E$ of the form $\{s\}$ satisfying the following conditions:
\begin{enumerate}
\item \label{thm_item:finpart_general_family2_O2containtrunk}
The graph $\Gamma_{O_2(s)}^{\mathrm{odd}}$ is connected and contains all simple closed paths in $\Gamma_O^{\mathrm{odd}}$, and $s$ satisfies the condition (\ref{thm_eq:finpart_general_1_trunk_trunk}).
\item \label{thm_item:finpart_general_family2_Oeven}
We have $O=O_{2,4,\infty}(s)$, $O_4(s) \neq \emptyset$, and for any $y \in O_4(s)$, we have $m_{y,y'}=3$ for some $y' \in O_2(s)$.
(Hence $O_2(s) \neq \emptyset$.)
\item \label{thm_item:finpart_general_family2_O2containOeven}
The condition (\ref{thm_eq:finpart_general_1_Oeven_even}) is satisfied for $s$ and any $t \in E \smallsetminus s$.
\end{enumerate}
Fix an element $y_K \in \bigcap_{s \in K}O_2(s)$ for each $K \in \mathcal{K}_1 \cup \mathcal{K}_2$, and let $G_K$ be the subgroup of $W^{\perp x}$ generated by $\{r_x(y_K,s) \mid s \in K\}$.
Then $(W^{\perp x})_{\mathrm{fin}}$ is the (restricted) direct product of all $G_K$ with $K \in \mathcal{K}_1 \cup \mathcal{K}_2$, and each $G_K$ is isomorphic to $W_K$.
\end{thm}
\subsection{Proof of Theorem \ref{thm:finpart_O}}
\label{sec:finitepart_maintheorem_O}

This subsection is devoted to the proof of Theorem \ref{thm:finpart_O}.
Indeed, here we only prove the ``only if'' part, since the `if' part has been verified in Theorems \ref{thm:finpart_O_typeB}, \ref{thm:finpart_O_trunk_acyclic} and \ref{thm:finpart_O_trunk_cyclic}.
Before starting the proof, we give some preliminary observations.
\begin{lem}
\label{lem:proofoffinpartO_supportofequivclass}
Let $(y,y') \in (\mathcal{E}_x)_{\mathrm{fin}} \cap \mathcal{E}_x^{(O)}$ and $I \subseteq O$, and suppose that $z,z' \in I$ for all $(z,z') \in \left[y,y'\right]$.
Let $s \in E$.
\begin{enumerate}
\item \label{lem_item:proofoffinpartO_supp_even}
If $O \neq O_{2,\infty}(s)$, then $\emptyset \neq O_{\mathrm{even}}(s) \smallsetminus O_2(s) \subseteq I$ and $I_2(s) \neq \emptyset$.
\item \label{lem_item:proofoffinpartO_supp_noeven}
If $O=O_{2,\infty}(s)$, then $|I_2(s)| \geq 2$, and every connected component of $\Gamma_{O_2(s)}^{\mathrm{odd}}$ intersects $I$.
\end{enumerate}
\end{lem}
\begin{proof}
\textbf{(\ref{lem_item:proofoffinpartO_supp_even})}
The hypothesis means that $O_{\mathrm{even}}(s) \smallsetminus O_2(s) \neq \emptyset$.
Let $y'' \in O_{\mathrm{even}}(s) \smallsetminus O_2(s)$.
Then we have $\left[y'',s\right] \overset{m}{\sim} \left[y,y'\right]$ for some $2 \leq m<\infty$ (by Corollary \ref{cor:necessaryconditionforfinitecomponent}), while $\left[y'',s\right]=\{(y'',s)\}$ (see Lemma \ref{lem:relationsoutsideO} (\ref{lem_item:reloutsideO_1simcomp})).
Thus $(y'',s) \overset{m}{\sim} (z,z')$ for some $(z,z') \in \left[y,y'\right]$; so $z=y''$ and $z' \in O_2(s)$ since $m_{y'',s} \neq 2$ (see Figure \ref{fig:localsubgraph}), while $z,z' \in I$ by the hypothesis.
Hence $y'' \in I$ and $z' \in I_2(s)$, proving the claim.

\noindent
\textbf{(\ref{lem_item:proofoffinpartO_supp_noeven})}
Note that $O_2(s) \neq \emptyset$ since $O \neq O_\infty(s)$.
Now for $y'' \in O_2(s)$, a similar argument to Claim \ref{lem_item:proofoffinpartO_supp_even} shows that $(z,z') \overset{m}{\sim} \left[y'',s\right]$ for some $2 \leq m<\infty$ and $(z,z') \in \left[y,y'\right]$; so $z,z' \in I$, $z \in O_2(s)_{\sim y''}^{\mathrm{odd}}$ and $(z,z') \overset{m}{\sim} (z,s)$ (see Lemma \ref{lem:relationsoutsideO} (\ref{lem_item:reloutsideO_1simcomp})), therefore $m_{z',s}<\infty$ and so $m_{z',s}=2$ (since $O=O_{2,\infty}(s)$).
Thus $z,z' \in I_2(s)$, and the connected component of $\Gamma_{O_2(s)}^{\mathrm{odd}}$ containing $y''$ intersects $I$ (at $z$).
hence the claim holds.
\end{proof}
\begin{lem}
\label{lem:proofoffinpartO_trunkandeven}
Suppose that $(\mathcal{E}_x)_{\mathrm{fin}} \cap \mathcal{E}_x^{(O)} \neq \emptyset$, and let $K$ be a subset of $O$ satisfying (\ref{eq:Kistrunk}).
Then $|K_{\mathrm{even}}(s)| \geq 2$ for all $s \in E$.
\end{lem}
\begin{proof}
Take an element $(y,y') \in (\mathcal{E}_x)_{\mathrm{fin}} \cap \mathcal{E}_x^{(O)}$.
Then the condition (\ref{eq:Kistrunk}) implies that $\mathcal{E}_x^{(O)}=\mathcal{E}_x^{(K)}$, so $(y,y')$ and $I=K$ satisfy the hypothesis of Lemma \ref{lem:proofoffinpartO_supportofequivclass}.
Now the claim follows from Lemma \ref{lem:proofoffinpartO_supportofequivclass}.
\end{proof}
\begin{lem}
\label{lem:proofoffinpartO_singleton}
Suppose that there is an element $(y,y') \in (\mathcal{E}_x)_{\mathrm{fin}} \cap \mathcal{E}_x^{(O)}$ such that $\left[y,y'\right]=\{(y,y')\}$.
Then for $s \in E$, we have $O=O_{2,\infty}(s)$, and the set $O_2(s)$ is connected in $\Gamma_O^{\mathrm{odd}}$ and contains both $y$ and $y'$.
\end{lem}
\begin{proof}
Note that $(y,y')$ and $I=\{y,y'\}$ satisfy the hypothesis of Lemma \ref{lem:proofoffinpartO_supportofequivclass}.
Let $s \in E$.
First, assume contrary that $O \neq O_{2,\infty}(s)$.
Then by Lemma \ref{lem:proofoffinpartO_supportofequivclass} (\ref{lem_item:proofoffinpartO_supp_noeven}), we have either $y \in O_{\mathrm{even}}(s) \smallsetminus O_2(s)$ and $y' \in O_2(s)$, or $y \in O_2(s)$ and $y' \in O_{\mathrm{even}}(s) \smallsetminus O_2(s)$.
However, in any case, Lemma \ref{lem:relationsoutsideO} (\ref{lem_item:reloutsideO_1simcomp}) implies that $(y,s) \not\in \left[y',s\right]$, so we have $\left[y,y'\right]=\{(y,y')\} \overset{\infty}{\sim} (y',s)$ (see Figure \ref{fig:localsubgraph}), contradicting the hypothesis $(y,y') \in (\mathcal{E}_x)_{\mathrm{fin}}$.

Thus we have shown that $O=O_{2,\infty}(s)$, so Lemma \ref{lem:proofoffinpartO_supportofequivclass} (\ref{lem_item:proofoffinpartO_supp_noeven}) implies that $I \subseteq O_2(s)$ (since $|I|=2$), and every connected component of $\Gamma_{O_2(s)}^{\mathrm{odd}}$ contains $y$ or $y'$.
Moreover, that component containing $y'$ also contains $y$; otherwise, we have $(y,s) \not\in \left[y',s\right]$ (see Lemma \ref{lem:relationsoutsideO} (\ref{lem_item:reloutsideO_1simcomp})), so $\left[y',s\right] \overset{\infty}{\sim} \{(y,y')\}=\left[y,y'\right]$, contradicting the hypothesis $(y,y') \in (\mathcal{E}_x)_{\mathrm{fin}}$.
This means that every connected component of $\Gamma_{O_2(s)}^{\mathrm{odd}}$ contains $y$, namely $\Gamma_{O_2(s)}^{\mathrm{odd}}$ is connected.
Hence the claim holds.
\end{proof}
Here we prepare a definition.
The \emph{complement} $\overline{G}$ of a graph $G$ is the simple graph with vertex set $V(G)$ in which two distinct vertices are joined if and only if these are not adjacent in $G$.
\begin{lem}
\label{lem:proofoffinpartO_complement}
Let $I$ be a subset of $O$ such that $(y,y') \in (\mathcal{E}_x)_{\mathrm{fin}}$ whenever $(y,y') \in \mathcal{E}_x$ and $y,y' \in I$.
Then for $s \in E$ and $y \in I_{\mathrm{even}}(s)$, the set $I_{\mathrm{even}}(s)$ contains all vertices of the connected component of $\overline{\Gamma_I}$ containing $y$.
\end{lem}
\begin{proof}
It suffices to show that $y' \in I_{\mathrm{even}}(s)$ if $y'$ is a neighbor of $y$ in $\overline{\Gamma_I}$, namely if $y' \in I$ and $m_{y,y'}=2$.
Now we have $(y,y') \in (\mathcal{E}_x)_{\mathrm{fin}}$ by the hypothesis, while $m_{y,s}$ is even; so the claim follows from Lemma \ref{lem:deniedtriangle}.
\end{proof}

Now we start the proof of Theorem \ref{thm:finpart_O}.
Suppose that $(x',x'') \in (\mathcal{E}_x)_{\mathrm{fin}} \cap \mathcal{E}_x^{(O)}$.
Then we have $(x',x'') \in (\mathcal{E}_x^{(O)})_{\mathrm{fin}}$ by Lemma \ref{lem:restrictionoffinitepart}, so it follows from Theorems \ref{thm:finitepart_oddconnected_cyclic} and \ref{thm:finitepart_oddconnected_acyclic} that the conditions in one of Theorems \ref{thm:finitepart_oddconnected_cyclic}, \ref{thm:finitepart_oddconnected_acyclic_trankisfinitetype} and \ref{thm:finitepart_oddconnected_acyclic_trankisnonfinitetype} are now satisfied.
Note that the case of Theorem \ref{thm:finitepart_oddconnected_acyclic_trankisfinitetype} with $K$ of type $A_3$ or $D_n$ is actually included in the case of Theorem \ref{thm:finitepart_oddconnected_acyclic_trankisnonfinitetype} (with the set $K$ suitably chosen).
By Remark \ref{rem:Efinisindependentofx}, we may assume that $x$ belongs to the set $K$.
Now the situation is divided into the following five cases.

\noindent
\textbf{Case 1: The conditions in Theorem \ref{thm:finitepart_oddconnected_acyclic_trankisfinitetype} are satisfied, with $K$ of type $E_6$, $E_7$, $E_8$ or $H_4$.}

We verify the hypothesis of Theorem \ref{thm:finpart_O_trunk_acyclic}.
Now Theorem \ref{thm:finitepart_oddconnected_acyclic_trankisfinitetype} implies that $\mathcal{E}_x^{(O)}=\mathcal{E}_x^{(K)}$, so $(x',x'') \in (\mathcal{E}_x)_{\mathrm{fin}} \cap \mathcal{E}_x^{(K)}$; while a direct computation shows that $W_K^{\perp x}$ is irreducible (see Theorem \ref{thm:Wperpxforoddpath} for the case of type $H_4$).
Thus we have $\mathcal{E}_x^{(K)} \subseteq (\mathcal{E}_x)_{\mathrm{fin}}$.

Now $\overline{\Gamma_K}$ is connected, so for $s \in E$, the combination of Lemmas \ref{lem:proofoffinpartO_trunkandeven} and \ref{lem:proofoffinpartO_complement} implies that $K \subseteq O_{\mathrm{even}}(s)$.
Moreover, we have $K \subseteq O_2(s)$, since otherwise $K \cup s$ is a finite, irreducible, $2$-spherical subset of $O$ with $\Gamma_{K \cup s}^{\mathrm{odd}}$ acyclic, contradicting Corollary \ref{cor:no2sphericalacyclicsupset}.
Thus Lemma \ref{lem:proofoffinpartO_supportofequivclass} (\ref{lem_item:proofoffinpartO_supp_even}) implies that $O=O_{2,\infty}(s)$, so every connected component of $\Gamma_{O_2(s)}^{\mathrm{odd}}$, which intersects $K$ by Lemma \ref{lem:proofoffinpartO_supportofequivclass} (\ref{lem_item:proofoffinpartO_supp_noeven}), must contain the whole $K$ in common (since $\Gamma_K^{\mathrm{odd}}$ is connected).
This means that $\Gamma_{O_2(s)}^{\mathrm{odd}}$ is connected.
Hence the hypothesis of Theorem \ref{thm:finpart_O_trunk_acyclic} is satisfied, concluding Case 1.

\noindent
\textbf{Case 2: The conditions in Theorem \ref{thm:finitepart_oddconnected_acyclic_trankisfinitetype} are satisfied, with $K$ of type $H_3$ or $P(m)$.}

We verify the hypothesis of Theorem \ref{thm:finpart_O_trunk_acyclic}.
Put $K=\{x_1,\dots,x_n\}$ (with $n=|K|$), where the path $(x_1,\dots,x_n)$ is of type $H_3$ or $P(m)$ (see (\ref{eq:ruleforfinitepart}) and Example \ref{ex:WperpxforP} for terminology).
Now $\mathcal{E}_x^{(O)}=\mathcal{E}_x^{(K)}$ by Theorem \ref{thm:finitepart_oddconnected_acyclic_trankisfinitetype}, so $(\mathcal{E}_x)_{\mathrm{fin}} \cap \mathcal{E}_x^{(K)} \neq \emptyset$.

In the case of type $H_3$, we have $\mathcal{E}_x^{(K)}=\{(x_1,x_3),(x_3,x_1)\}$, so $(x_i,x_j) \in (\mathcal{E}_x)_{\mathrm{fin}}$ for some $(i,j) \in \{(1,3),(3,1)\}$.
Now we have $\left[x_i,x_j\right]=\left[x_i,x_j\right]_K=\{(x_i,x_j)\}$, so for $s \in E$, Lemma \ref{lem:proofoffinpartO_singleton} implies that $O=O_{2,\infty}(s)$, and $\Gamma_{O_2(s)}^{\mathrm{odd}}$ is connected and contains both $x_i$ and $x_j$; therefore it also contains the vertex $x_2$ of $\Gamma_O^{\mathrm{odd}}$ between $x_i$ and $x_j$.
Hence the claim holds.

The case of type $P(m)$ is similar; we have $\mathcal{E}_x^{(K)} \subseteq (\mathcal{E}_x)_{\mathrm{fin}}$ by Theorem \ref{thm:Wperpxforoddpath}, so $(x_4,x_1) \in (\mathcal{E}_x)_{\mathrm{fin}}$, while $\left[x_4,x_1\right]=\{(x_4,x_1)\}$.
Thus the claim is deduced by the same argument, since $x_2$ and $x_3$ are all the vertices of $\Gamma_O^{\mathrm{odd}}$ between $x_1$ and $x_4$.
Hence Case 2 is concluded.

\noindent
\textbf{Case 3: The conditions in Theorem \ref{thm:finitepart_oddconnected_acyclic_trankisfinitetype} are satisfied, with $K$ of type $A_n$ and $n \geq 4$.}

Put $K=\{x_1,\dots,x_n\}$, where $(x_1,\dots,x_n)$ is of type $A_n$ (see (\ref{eq:ruleforfinitepart}) for terminology).
Note that $\mathcal{E}_x^{(O)}=\mathcal{E}_x^{(K)}$ by Theorem \ref{thm:finitepart_oddconnected_acyclic_trankisfinitetype}, so $(x',x'')$ and $K$ satisfy the hypothesis of Lemma \ref{lem:proofoffinpartO_supportofequivclass}.
Now since $W_K^{\perp x}$ is irreducible (see Theorem \ref{thm:Wperpxforoddpath}) and $\overline{\Gamma_K}$ is connected, we have $\mathcal{E}_x^{(K)} \subseteq (\mathcal{E}_x)_{\mathrm{fin}}$ and $K \subseteq O_{\mathrm{even}}(s)$ for all $s \in E$ by the same argument as Case 1.
Moreover, the same argument as Case 1 shows further that, if $K \subseteq O_2(s)$, then $O=O_{2,\infty}(s)$ and the graph $\Gamma_{O_2(s)}^{\mathrm{odd}}$ is connected and contains $K$.
Thus the hypothesis of Theorem \ref{thm:finpart_O_trunk_acyclic} is satisfied if $K \subseteq O_2(s)$ for all $s \in E$.

So suppose that $K \not\subseteq O_2(s_0)$ for some $s_0 \in E$.
We verify the hypothesis of Theorem \ref{thm:finpart_O_typeB}; the condition \ref{thm_item:finpart_O_typeB_trunk} is already satisfied.
Now since $K \subseteq O_{\mathrm{even}}(s_0)$ and $(\mathcal{E}_x^{(K)})_{\mathrm{fin}} \neq \emptyset$, we have (up to symmetry) $K \smallsetminus x_n \subseteq O_2(s_0)$ and $m_{x_n,s_0}=4$ (so $(x_1,\dots,x_n,s_0)$ is of type $B_{n+1}$); otherwise $K \cup s_0$ is finite, irreducible, $2$-spherical and not of finite type, and $\Gamma_{K \cup s_0}^{\mathrm{odd}}$ is acyclic, contradicting Corollary \ref{cor:no2sphericalacyclicsupset}.
On the other hand, Lemma \ref{lem:proofoffinpartO_supportofequivclass} (\ref{lem_item:proofoffinpartO_supp_even}) implies that $O_{\mathrm{even}}(s_0) \smallsetminus O_2(s_0) \subseteq K$; thus we have $O=O_{2,\infty}(s_0) \cup x_n$.
Moreover, we have $(x_1,s_0) \in (\mathcal{E}_x)_{\mathrm{fin}}$ since $(x_1,s_0) \overset{4}{\sim} (x_1,x_n) \in (\mathcal{E}_x)_{\mathrm{fin}}$; so Lemma \ref{lem:ncsrconditionforfinoutsideO} (\ref{lem_item:ncsrforfinoutsideO_mis2}) shows that $\Gamma_{O_2(s_0)}^{\mathrm{odd}}$ is connected.

Finally, let $s \in E \smallsetminus s_0$.
Since $(x_1,s_0) \in (\mathcal{E}_x)_{\mathrm{fin}}$ and $K \subseteq O_{\mathrm{even}}(s)$, we have $m_{s,s_0}<\infty$ by Lemma \ref{lem:deniedtriangle}.
Moreover, we have $K \subseteq O_2(s)$ and $m_{s,s_0}=2$, since otherwise the $2$-spherical set $K \cup \{s,s_0\}$ is not of finite type and so contradicts Corollary \ref{cor:no2sphericalacyclicsupset}.
Thus $\{s_0\}$ is an irreducible component of $E$, so the condition \ref{thm_item:finpart_O_typeB_s0} is satisfied.
On the other hand, the condition \ref{thm_item:finpart_O_typeB_nots0} is also satisfied by the argument in the first paragraph of Case 3.
Hence the hypothesis of Theorem \ref{thm:finpart_O_typeB} has been verified, concluding Case 3.

\noindent
\textbf{Case 4: The conditions in one of Theorems \ref{thm:finitepart_oddconnected_cyclic} and \ref{thm:finitepart_oddconnected_acyclic_trankisnonfinitetype} are satisfied, and $(\mathcal{E}_x)_{\mathrm{fin}} \cap \mathcal{E}_x^{(K)} \neq \emptyset$.}

Put $(x'_1,x'_2)=(x_1,x_2)$ in the case of Theorem \ref{thm:finitepart_oddconnected_cyclic}, and $(x'_1,x'_2)=(x_1,x_3)$ in the case of Theorem \ref{thm:finitepart_oddconnected_acyclic_trankisnonfinitetype}.
Then we have $\mathcal{E}_x^{(K)} \subseteq (\mathcal{E}_x)_{\mathrm{fin}}$, since the set $\mathcal{E}_x^{(K)}$ forms a single $\sim$-equivalence class $\left[x'_1,x'_2\right]$ (see Theorems \ref{thm:finitepart_oddconnected_cyclic} and \ref{thm:finitepart_oddconnected_acyclic_trankisnonfinitetype}).
Now $(x'_1,x'_2)$ and $I=\{x'_1,x'_2\}$ satisfy the hypothesis of Lemma \ref{lem:proofoffinpartO_supportofequivclass}.
Thus for $s \in E$, Lemma \ref{lem:proofoffinpartO_supportofequivclass} shows that $O_{\mathrm{even}}(s) \smallsetminus O_2(s)=\{x'_i\}$ (so $O=O_{2,\infty}(s) \cup x'_i$) and $x'_{3-i} \in O_2(s)$ for some $i \in \{1,2\}$ if $O \neq O_{2,\infty}(s)$; while $x'_1,x'_2 \in O_2(s)$ and every connected component of $\Gamma_{O_2(s)}^{\mathrm{odd}}$ contains either $x'_1$ or $x'_2$ if $O=O_{2,\infty}(s)$.
This means that the hypothesis of Theorem \ref{thm:finpart_O_trunk_cyclic} is satisfied if $O=O_{2,\infty}(s)$ for all $s \in E$.

So suppose that $O \neq O_{2,\infty}(s_0)$ for some $s_0 \in E$.
By symmetry, we may assume that $x'_1 \in O_2(s_0)$ and $x'_2 \in O_{\mathrm{even}}(s_0) \smallsetminus O_2(s_0)$.
This implies that $(x'_1,s_0) \overset{m}{\sim} (x'_1,x'_2) \in (\mathcal{E}_x)_{\mathrm{fin}}$ (where $m=m_{x'_2,s_0} \neq 2$), so $(x'_1,s_0) \in (\mathcal{E}_x)_{\mathrm{fin}}$.
Now in the case of Theorem \ref{thm:finitepart_oddconnected_cyclic}, the graph $\Gamma_O^{\mathrm{odd}}$ has a simple closed path containing $x'_2 \not\in O_2(s_0)$; this contradicts Lemma \ref{lem:ncsrconditionforfinoutsideO} (\ref{lem_item:ncsrforfinoutsideO_mis2}) applied to $(x'_1,s_0)$.
Thus this is the case of Theorem \ref{thm:finitepart_oddconnected_acyclic_trankisnonfinitetype}, so $(x'_1,x'_2)=(x_1,x_3)$.

By applying Lemma \ref{lem:ncsrconditionforfinoutsideO} (\ref{lem_item:ncsrforfinoutsideO_mis2}) to $(x_1,s_0)$, we have $O=O_{2,4,\infty}(s_0)$, and $\Gamma_{O_2(s_0)}^{\mathrm{odd}}$ is connected.
The former property implies that $x_3 \in O_4(s_0)$, so we have $m_{z,x_3}=3$ for some $z \in O_2(s_0)$ by Lemma \ref{lem:ncsrconditionforfinoutsideO} (\ref{lem_item:ncsrforfinoutsideO_mis2}) again.
Since $\Gamma_O^{\mathrm{odd}}$ is a tree, the connectedness of $\Gamma_{O_2(s_0)}^{\mathrm{odd}}$ implies that $x_3$ does not exist between $x_1$ and $z$; therefore $z=x_2$.
Thus $(x_1,x_2,x_3,s_0)$ is of type $B_4$; while $O=O_{2,\infty}(s_0) \cup x_3$ as shown in the first paragraph of Case 4.

We show that $T_{x_2}^o=\{x_2\}$.
If this fails, then by Condition \ref{thm_item:fin_oddconn_acycl_nonftype_orderideal} in Theorem \ref{thm:finitepart_oddconnected_acyclic_trankisnonfinitetype}, some $y \in T_{x_2}^o$ is an atom of $T_{x_2}$.
Since $(x_1,s_0) \in (\mathcal{E}_x)_{\mathrm{fin}}$ and $m_{x_1,y}=2$, Lemma \ref{lem:deniedtriangle} shows that $m_{y,s_0}<\infty$.
However, now the set $K \cup \{s_0,y\}$ is irreducible, $2$-spherical and not of finite type, contradicting Corollary \ref{cor:no2sphericalacyclicsupset}.
Thus we have $T_{x_2}^o=\{x_2\}$, so $K'=K$; this implies (by Theorem \ref{thm:finitepart_oddconnected_acyclic_trankisnonfinitetype}) that $K$ satisfies the condition (\ref{eq:Kistrunk}).

We show that the hypothesis of Theorem \ref{thm:finpart_O_typeB} is satisfied.
The condition \ref{thm_item:finpart_O_typeB_trunk} has been verified.
For $s \in E \smallsetminus s_0$, we have $x_1 \in O_{\mathrm{even}}(s)$ as shown in the first paragraph of Case 4, so Lemma \ref{lem:deniedtriangle} implies that $m_{s,s_0}<\infty$ since $(x_1,s_0) \in (\mathcal{E}_x)_{\mathrm{fin}}$.
Now if $O \neq O_{2,\infty}(s)$, then the above argument applied to $s_0$ also works for $s$, proving that $K \cup s$ is of type $B_4$.
However, now $K \cup \{s,s_0\}$ is irreducible, $2$-spherical and not of finite type, contradicting \ref{cor:no2sphericalacyclicsupset}.
Thus we have $O=O_{2,\infty}(s)$, so as shown in the first paragraph of Case 4, we have $x_1,x_3 \in O_2(s)$ and every connected component of $\Gamma_{O_2(s)}^{\mathrm{odd}}$ contains $x_1$ or $x_3$.
Now if $x_2 \not\in O_2(s)$, then since $x_3 \not\in O_2(s_0)$, two sets $O_2(s_0)_{\sim x_1}^{\mathrm{odd}}$ and $O_2(s)_{\sim x_3}^{\mathrm{odd}}$ do not intersect; so $\left[x_1,s_0\right] \overset{\infty}{\sim} \left[x_3,s\right]$ by Lemma \ref{lem:relationsoutsideO} (\ref{lem_item:reloutsideO_comp_sisdistinct}), contradicting the fact $(x_1,s_0) \in (\mathcal{E}_x)_{\mathrm{fin}}$.
Thus we have $x_2 \in O_2(s)$, so $\Gamma_{O_2(s)}^{\mathrm{odd}}$ is connected and $K \subseteq O_2(s)$; therefore the condition \ref{thm_item:finpart_O_typeB_nots0} is satisfied.
Finally, if $m_{s,s_0} \neq 2$, then $K \cup \{s,s_0\}$ is irreducible, $2$-spherical and not of finite type, contradicting \ref{cor:no2sphericalacyclicsupset}.
Thus we have $m_{s,s_0}=2$, so $\{s_0\}$ is an irreducible component of $E$; therefore the condition \ref{thm_item:finpart_O_typeB_s0} is satisfied.
Hence the hypothesis of Theorem \ref{thm:finpart_O_typeB} is satisfied, concluding Case 4.

\noindent
\textbf{Case 5: The conditions in one of Theorems \ref{thm:finitepart_oddconnected_cyclic} and \ref{thm:finitepart_oddconnected_acyclic_trankisnonfinitetype} are satisfied, and $(\mathcal{E}_x)_{\mathrm{fin}} \cap \mathcal{E}_x^{(K)}=\emptyset$.}

In this case, we have $(\mathcal{E}_x)_{\mathrm{fin}} \cap (\mathcal{E}_x^{(O)} \smallsetminus \mathcal{E}_x^{(K)}) \neq \emptyset$, so $(\mathcal{E}_x^{(O)})_{\mathrm{fin}} \not\subseteq \mathcal{E}_x^{(K)}$ by Lemma \ref{lem:restrictionoffinitepart}.
By Theorems \ref{thm:finitepart_oddconnected_cyclic} and \ref{thm:finitepart_oddconnected_acyclic_trankisnonfinitetype}, this happens only if the hypothesis of Theorem \ref{thm:finitepart_oddconnected_acyclic_trankisnonfinitetype} is satisfied and $K'$ is of type $A_3$ or $D_n$ (with $n \geq 4$); so the hypothesis of Theorem \ref{thm:finitepart_oddconnected_acyclic_trankisfinitetype} is satisfied, where $K'$ plays the role of $K$ in Theorem \ref{thm:finitepart_oddconnected_acyclic_trankisfinitetype}.
Note that $\mathcal{E}_x^{(O)}=\mathcal{E}_x^{(K')}$.
The possibility of type $A_3$ is also denied, since $K'=K$ and so $\mathcal{E}_x^{(O)}=\mathcal{E}_x^{(K)}$ in the case of type $A_3$.
Thus $K'$ is of type $D_n$; so $T_{x_2}^o$ is a chain, say $y_0 \prec_{x_2} y_1 \prec_{x_2} \cdots \prec_{x_2} y_{n-3}$, such that $(y_0,\dots,y_{n-3})$ is of type $A_{n-2}$.
Moreover, if $n=4$, then $\mathcal{E}_x^{(O)}=\left[x_1,x_3\right] \sqcup \left[x_1,y_1\right] \sqcup \left[x_3,y_1\right]$ and so $(x_i,y_1) \in (\mathcal{E}_x)_{\mathrm{fin}}$ for some $i \in \{1,3\}$.
Now by the first remark of the proof, the hypothesis of Theorem \ref{thm:finitepart_oddconnected_acyclic_trankisnonfinitetype} is satisfied, where the set $\{x_i,x_2,y_1\}$ plays the role of $K$; therefore the proof is reduced to Case 4.
Thus we may assume that $n \geq 5$.

By Remark \ref{rem:Efinisindependentofx}, we may assume that $x=x_2$.
Then by putting $I'=K' \smallsetminus x_3$ and $I''=K' \smallsetminus x_1$, we have $\mathcal{E}_x^{(O)}=\mathcal{E}_x^{(K')}=\mathcal{E}_x^{(K)} \cup \mathcal{E}_x^{(I')} \cup \mathcal{E}_x^{(I'')}$, while $\mathcal{E}_x^{(I')} \cap \mathcal{E}_x^{(I'')} \neq \emptyset$.
Now we show that $\mathcal{E}_x^{(I')} \cup \mathcal{E}_x^{(I'')} \subseteq (\mathcal{E}_x)_{\mathrm{fin}}$.
Recall that $(\mathcal{E}_x)_{\mathrm{fin}} \cap (\mathcal{E}_x^{(O)} \smallsetminus \mathcal{E}_x^{(K)}) \neq \emptyset$, so we may assume by symmetry that $(\mathcal{E}_x)_{\mathrm{fin}} \cap \mathcal{E}_x^{(I')} \neq \emptyset$.
Then we have $\mathcal{E}_x^{(I')} \subseteq (\mathcal{E}_x)_{\mathrm{fin}}$ by Theorem \ref{thm:Wperpxforoddpath} (\ref{thm_item:Wperpxforoddpath_AHP}) applied to $I'$, so $(\mathcal{E}_x)_{\mathrm{fin}} \cap \mathcal{E}_x^{(I'')} \neq \emptyset$, therefore $\mathcal{E}_x^{(I'')} \subseteq (\mathcal{E}_x)_{\mathrm{fin}}$ similarly.

We show that the hypothesis of Theorem \ref{thm:finpart_O_trunk_acyclic} is satisfied.
Let $s \in E$.
Then Lemma \ref{lem:proofoffinpartO_trunkandeven} implies that $(K')_{\mathrm{even}}(s) \neq \emptyset$; by symmetry, we may assume that $(I')_{\mathrm{even}}(s) \neq \emptyset$ (note that $K'=I' \cup I''$).
Now both $\overline{\Gamma_{I'}}$ and $\overline{\Gamma_{I''}}$ are connected, so Lemma \ref{lem:proofoffinpartO_complement} implies that $I' \subseteq O_{\mathrm{even}}(s)$, and further that $I'' \subseteq O_{\mathrm{even}}(s)$ (since $I' \cap I'' \neq \emptyset$), therefore $K' \subseteq O_{\mathrm{even}}(s)$.
Moreover, we have $K' \subseteq O_2(s)$, since otherwise $K' \cup s$ is irreducible, $2$-spherical and not of finite type, contradicting Corollary \ref{cor:no2sphericalacyclicsupset}.
Finally, Lemma \ref{lem:proofoffinpartO_supportofequivclass} (applied to $(x',x'')$ and $K'$) denies the possibility $O \neq O_{2,\infty}(s)$ and implies that every connected component of $\Gamma_{O_2(s)}^{\mathrm{odd}}$ intersects $K'$.
Since $K' \subseteq O_2(s)$ and $\Gamma_{K'}^{\mathrm{odd}}$ is connected, these components contain $K'$ in common, forcing $\Gamma_{O_2(s)}^{\mathrm{odd}}$ to be connected.
Hence the hypothesis of Theorem \ref{thm:finpart_O_trunk_acyclic} is satisfied, concluding Case 5.

Thus the proof of Theorem \ref{thm:finpart_O} is concluded.
\subsection{Proof of Theorem \ref{thm:finpart_nonO}}
\label{sec:finitepart_maintheorem_nonO}

This subsection is devoted to the proof of Theorem \ref{thm:finpart_nonO}.
The `if' part has been verified in Theorems \ref{thm:finpart_nonO_Otrunk_F4}, \ref{thm:finpart_nonO_Otrunk_B} and \ref{thm:finpart_nonO_Onontrunk}; so we show the ``only if'' part.
Moreover, to verify the hypothesis of Theorem \ref{thm:finpart_nonO_Onontrunk}, it suffices to check the conditions \ref{thm_item:finpart_nonO_Onontrunk_prev_Oistree}--\ref{thm_item:finpart_nonO_Onontrunk_prev_E} in Theorem \ref{thm:finpart_nonO_Onontrunk_previouscase} only.
Indeed, now the hypothesis of one of Theorems \ref{thm:finpart_nonO_Onontrunk_previouscase} and \ref{thm:finpart_nonO_Onontrunk} is satisfied; while the hypothesis $(\mathcal{E}_x)_{\mathrm{fin}} \cap \mathcal{E}_x^{(O)}=\emptyset$ denies the former possibility, since Theorem \ref{thm:finpart_O} shows that the hypothesis of Theorem \ref{thm:finpart_O_typeB} cannot be satisfied.

We divide the proof into four cases.
Note that $\left[y,s_0\right]=\{(y,s_0)\}$.

\noindent
\textbf{Case 1: $O=O_{\infty}(s_0) \cup y$ and $\{s_0\}$ is an irreducible component of $E$.}

We verify the three conditions in Theorem \ref{thm:finpart_nonO_Otrunk_B}, with $n=0$ and $m=m_{y,s_0}$.
By applying Lemma \ref{lem:ncsrconditionforfinoutsideO} (\ref{lem_item:ncsrforfinoutsideO_mis4}) or (\ref{lem_item:ncsrforfinoutsideO_other}) to $(y,s_0) \in (\mathcal{E}_x)_{\mathrm{fin}}$, it follows that $\Gamma_O^{\mathrm{odd}}$ is acyclic, and for all $s \in E \smallsetminus s_0$, we have $O=O_{2,\infty}(s)$ and the graph $\Gamma_{O_2(s)}^{\mathrm{odd}}$ is connected and contains $y$.
Thus Conditions \ref{thm_item:finpart_nonO_Otrunk_B_Bexist} and \ref{thm_item:finpart_nonO_Otrunk_B_E} have been verified.
Now our remaining task is to show that $m_{z,z'}$ is odd or $\infty$ for any $z,z' \in O$.
We have $(t,t') \overset{\infty}{\sim} (y,s_0)$ for any $(t,t') \in \mathcal{E}_x^{(O)}$; otherwise either $m_{s_0,t}$ or $m_{s_0,t'}$ must be finite, but this is impossible since $O=O_\infty(s_0) \cup y$.
Thus if $m_{z,z'}$ is even, then $\left[z,z'\right] \subseteq \mathcal{E}_x^{(O)}$ and so $\left[y,s_0\right]=\{(y,s_0)\} \overset{\infty}{\sim} \left[z,z'\right]$, contradicting the hypothesis $(y,s_0) \in (\mathcal{E}_x)_{\mathrm{fin}}$.
Hence the claim follows, so Case 1 is concluded.

\noindent
\textbf{Case 2: $O \neq O_{\infty}(s_0) \cup y$ and $\{s_0\}$ is an irreducible component of $E$.}

In this case, Lemma \ref{lem:ncsrconditionforfinoutsideO} (\ref{lem_item:ncsrforfinoutsideO_other}) denies the possibility $m_{y,s_0} \geq 6$, so $m_{y,s_0}=4$.
We verify the conditions \ref{thm_item:finpart_nonO_Onontrunk_prev_Oistree}--\ref{thm_item:finpart_nonO_Onontrunk_prev_E} in Theorem \ref{thm:finpart_nonO_Onontrunk_previouscase}.
Conditions \ref{thm_item:finpart_nonO_Onontrunk_prev_Oistree} and \ref{thm_item:finpart_nonO_Onontrunk_prev_E} follow from Lemma \ref{lem:ncsrconditionforfinoutsideO} (\ref{lem_item:ncsrforfinoutsideO_mis4}) applied to $(y,s_0)$.
Now $\Gamma_O^{\mathrm{odd}}$ is a tree; we denote by $\preceq$ the tree order on $O$ with root vertex $y$.

Lemma \ref{lem:ncsrconditionforfinoutsideO} (\ref{lem_item:ncsrforfinoutsideO_mis4}) also implies that $O=O_{2,\infty}(s_0) \cup y$ (so $O_2(s_0) \neq \emptyset$ by the hypothesis of Case 2), and every connected component of $\Gamma_{O_2(s_0)}^{\mathrm{odd}}$ contains a vertex $z$ with $m_{z,y}$.
This $z$ is a neighbor of $y$ in $\Gamma_O^{\mathrm{odd}}$; so it follows that the set $T=O_2(s_0) \cup y$ is connected in $\Gamma_O^{\mathrm{odd}}$.
Thus Condition \ref{thm_item:finpart_nonO_Onontrunk_prev_B2exist} has been verified.
Note that the poset $T$ is now an order ideal of $O$, by connectedness of $\Gamma_T^{\mathrm{odd}}$.
Moreover, if this connected component of $\Gamma_{O_2(s_0)}^{\mathrm{odd}}$ contains an atom $y'$ of $T$, then this $z$ must be $y'$; thus Condition \ref{thm_item:finpart_nonO_Onontrunk_prev_atomofT} has been verified.

From now, we verify Conditions \ref{thm_item:finpart_nonO_Onontrunk_prev_Tistrunk} and \ref{thm_item:finpart_nonO_Onontrunk_prev_orderideal}.
We prepare some preliminary observations.
If $(z_1,z_2) \in \mathcal{E}_x^{(O)}$, then $\left[z_1,z_2\right] \neq \left[y,s_0\right]$ and so we have $\left[z_1,z_2\right] \overset{m}{\sim} \left[y,s_0\right]=\{(y,s_0)\}$ for some $2 \leq m<\infty$, since $(y,s_0) \in (\mathcal{E}_x)_{\mathrm{fin}}$.
This implies that $(z'_1,z'_2) \overset{m}{\sim} (y,s_0)$ for some $(z'_1,z'_2) \in \left[z_1,z_2\right]$, which means that $z'_1=y$ and $m_{z'_2,s_0}<\infty$, therefore $z'_2 \in O_2(s_0) \smallsetminus y=T$ since $O=O_{2,\infty}(s_0) \cup y$.
Thus we have the following:
\begin{equation}
\label{thm_eq:finpart_nonO_equivclass}
\textrm{every $\sim$-equivalence class in } \mathcal{E}_x^{(O)} \textrm{ contains some } (y,t) \textrm{ with } t \in T.
\end{equation}

First, we show that $m_{y_1,y_2}$ is $2$, odd or $\infty$ for any $y_1,y_2 \in O$.
If $m_{y_1,y_2}$ is even and $m_{y_1,y_2} \neq 2$, then $y_i \neq y$ for some $i \in \{1,2\}$, so the set $\left[y_i,y_{3-i}\right]=\{(y_i,y_{3-i})\}$ contradicts (\ref{thm_eq:finpart_nonO_equivclass}).
Thus the claim follows.

We show that, if $y_1,y_2 \in O$ and $m_{y_1,y_2}=2$, then $y_1$ and $y_2$ are comparable in $O$.
Assume contrary that $y_1$ and $y_2$ are incomparable, and let $z_0 \prec z_1 \prec \cdots \prec z_k$ and $z'_0 \prec z'_1 \prec \cdots \prec z'_\ell$ be the saturated chains in $O$ from $y_1 \wedge y_2$ to $y_1=z_k$ and $y_2=z'_\ell$, respectively, with $k,\ell \neq 0$.
If $k=1$, then we have
\[
\left[y_2,y_1\right] \subseteq \{(t,z_1) \mid t \in \vee_{z'_1}\} \cup \{(t,z'_1) \mid t \in \vee_{z_1}\},
\]
so $\left[y_2,y_1\right]$ contradicts (\ref{thm_eq:finpart_nonO_equivclass}).
If $k \geq 2$, then the above argument shows that $m_{z_1,y_2} \neq 2$, so we have
\[
\left[y_1,y_2\right] \subseteq \{(t,y_2) \mid t \in \vee_{z_2}\},
\]
therefore $\left[y_1,y_2\right]$ contradicts (\ref{thm_eq:finpart_nonO_equivclass}).
Thus the claim follows.

Take any $y_1,y_2 \in O$ such that $m_{y_1,y_2}=2$ and $y_1 \prec y_2$.
Let $z_0 \prec z_1 \prec \cdots \prec z_k$ be the saturated chain in $O$ from $z_0=y$ to $z_k=y_2$, and $y_1=z_\ell$ with $0 \leq \ell \leq k-2$.
We show that $m_{z_i,z_k}=2$ for all $0 \leq i<\ell$, and $z_k \in T$.
If $m_{z_i,z_k} \neq 2$ for some $0 \leq i<\ell$, then we have
\[
\left[z_\ell,z_k\right] \subseteq \{(t,z_k) \mid t \in \vee_{z_{i+1}} \smallsetminus \vee_{z_{k-1}}\} \cup \{(t,z_{k-2}) \mid t \in \vee_{z_k}\},
\]
so $\left[z_\ell,z_k\right]$ contradicts (\ref{thm_eq:finpart_nonO_equivclass}).
Thus the former claim follows.
Moreover, since $m_{y,z_k}=2$ and $(y,s_0) \in (\mathcal{E}_x)_{\mathrm{fin}}$, Lemma \ref{lem:deniedtriangle} implies that $m_{z_k,s_0}<\infty$, so $z_k \in O_2(s_0) \subseteq T$ since $O=O_{2,\infty}(s_0) \cup y$, as desired.
Note that $z_i \in T$ for all $0 \leq i \leq k$, since $T$ is an order ideal.
Thus we have $y_1,y_2 \in T$, so Condition \ref{thm_item:finpart_nonO_Onontrunk_prev_Tistrunk} has been verified.

We show that $m_{z_{k-2},z_{k-1}}=m_{z_{k-1},z_k}=3$ whenever $\ell=k-2$.
If this fails, then we have
\[
\left[z_k,z_{k-2}\right] \subseteq \{(t,z_{k-2}) \mid t \in \vee_{z_k}\},
\]
so $\left[z_k,z_{k-2}\right]$ contradicts (\ref{thm_eq:finpart_nonO_equivclass}).
Thus the claim follows.

Finally, we show that $m_{z_i,z_j}=2$ whenever $0 \leq i \leq \ell$ and $i+2 \leq j \leq k$.
Note that $m_{z_i,z_k}=2$ as above, so we may assume that $j<k$.
Now if $m_{z_i,z_j} \neq 2$, then we have
\[
\left[z_k,z_i\right] \subseteq \{(t,z_i) \mid t \in \vee_{z_{j+1}}\},
\]
so $\left[z_k,z_i\right]$ contradicts (\ref{thm_eq:finpart_nonO_equivclass}).
Thus $m_{z_i,z_j}=2$ as desired.

Hence Condition \ref{thm_item:finpart_nonO_Onontrunk_prev_orderideal} has been verified, so Case 2 is concluded.

\noindent
\textbf{Case 3: $O=O_{\infty}(s_0) \cup y$ and $\{s_0\}$ is not an irreducible component of $E$.}

We verify the hypothesis of Theorem \ref{thm:finpart_nonO_Otrunk_B} with $n \geq 1$.
By Remark \ref{rem:Efinisindependentofx}, we may assume that $x=y$.
Let $K$ be the irreducible component of $E$ containing $s_0$, so $K \neq \{s_0\}$.
Now $m_{s,s_0} \neq 2$ for some $s \in K \smallsetminus s_0$; so Lemma \ref{lem:ncsrconditionforfinoutsideO} (\ref{lem_item:ncsrforfinoutsideO_other}) denies the possibility $m_{y,s_0} \geq 6$, therefore $m_{y,s_0}=4$.
Condition \ref{thm_item:finpart_nonO_Otrunk_B_E} follows from Lemma \ref{lem:ncsrconditionforfinoutsideO} (\ref{lem_item:ncsrforfinoutsideO_mis4}) applied to $(y,s_0)$.
Lemma \ref{lem:ncsrconditionforfinoutsideO} (\ref{lem_item:ncsrforfinoutsideO_mis4}) also shows that $\Gamma_O^{\mathrm{odd}}$ is acyclic.
Moreover, since $O=O_{\infty}(s_0) \cup y$, the same argument as Case 1 implies that Condition \ref{thm_item:finpart_nonO_Otrunk_B_trunk} is satisfied.

Finally, we verify Condition \ref{thm_item:finpart_nonO_Otrunk_B_Bexist}, concluding Case 3.
Note that $r_y(y,s)=s$ for any $s \in E \smallsetminus s_0$ by Condition \ref{thm_item:finpart_nonO_Otrunk_B_E}.
By Lemma \ref{lem:irrsubsetofE}, all the generators $r_y(y,s)=s$ with $s \in K \smallsetminus s_0$ belong to the same (finite) irreducible component of $W^{\perp y}$ as $r_y(y,s_0)$; so $K \smallsetminus s_0$ is of finite type.
On the other hand, we have $m_{s,s_0} \in \{2,3\}$ for all $s \in E \smallsetminus s_0$, by Lemma \ref{lem:ncsrconditionforfinoutsideO} (\ref{lem_item:ncsrforfinoutsideO_mis4}).
Now if $s_0$ has two neighbors $s,s'$ in $\Gamma_K$, then $m_{s_0,s}=m_{s_0,s'}=3$ and $m_{y,s}=m_{y,s'}=2$, so we have $(y,s) \overset{4}{\sim} (y,s_0) \overset{4}{\sim} (y,s')$.
This means that the three generators $r_y(y,s)$, $r_y(y,s_0)$ and $r_y(y,s')$, which belong to the same finite irreducible component of $W^{\perp y}$ as above, generate an infinite group.
This contradiction implies that $s_0$ has a unique neighbor in $\Gamma_K$, say $s_1$.
Since $K \smallsetminus s_0$ is of finite type, the graph $\Gamma_{K \smallsetminus s_0}$ is acyclic, and $K \smallsetminus s_0$ is finite and $2$-spherical, so $\Gamma_K$ is also acyclic and $K$ is also finite and $2$-spherical.
Now by Condition \ref{thm_item:finpart_nonO_Otrunk_B_E}, the set $K \cup y$ is finite, irreducible and $2$-spherical, and $\Gamma_{K \cup y}^{\mathrm{odd}}$ is acyclic.
Thus Corollary \ref{cor:no2sphericalacyclicsupset} implies that $K \cup y$ must be of finite type, since $(y,s_0) \in (\mathcal{E}_y)_{\mathrm{fin}}$.
Since $m_{y,s_0}=4$ and $m_{y,s}=2$ for all $s \in K \smallsetminus s_0$, this is possible only if $K \cup y$ is of type $B_n$ with $3 \leq n<\infty$.
Hence Condition \ref{thm_item:finpart_nonO_Otrunk_B_Bexist} has been verified, so Case 3 is concluded.

\noindent
\textbf{Case 4: $O \neq O_{\infty}(s_0) \cup y$ and $\{s_0\}$ is not an irreducible component of $E$.}

We verify the hypothesis of Theorem \ref{thm:finpart_nonO_Otrunk_F4}.
We may assume that $x=y$.
Let $K$ be the irreducible component of $E$ containing $s_0$.
Since $K \neq \{s_0\}$, the same argument as the first paragraph of Case 3 implies that $m_{y,s_0}=4$, the graph $\Gamma_O^{\mathrm{odd}}$ is acyclic, and for any $s \in E \smallsetminus s_0$, we have $O=O_{2,\infty}(s)$ and the graph $\Gamma_{O_2(s)}^{\mathrm{odd}}$ is connected and contains $y$.
Moreover, the same argument as the second paragraph of Case 3 implies that the set $K \cup y$ is of type $B_{n+1}$ with $2 \leq n=|K|<\infty$.
Write $K=\{s_0,s_1,\dots,s_n\}$, where the sequence $(s_n,\dots,s_1,s_0,y)$ is of type $B_{n+1}$ (see (\ref{eq:ruleforfinitepart}) for terminology).
Note that $(y,s_i) \in (\mathcal{E}_y)_{\mathrm{fin}}$ for $1 \leq i \leq n$, since $r_y(y,s_i)$ and $r_y(y,s_0)$ belong to the same finite irreducible component of $W^{\perp y}$, by Lemma \ref{lem:irrsubsetofE}.

Lemma \ref{lem:ncsrconditionforfinoutsideO} (\ref{lem_item:ncsrforfinoutsideO_mis4}) implies that $O=O_{2,\infty}(s_0) \cup y$, so $O_2(s_0) \neq \emptyset$ by the hypothesis of Case 4.
This lemma also shows that a connected component of $\Gamma_{O_2(s_0)}^{\mathrm{odd}}$ contains a neighbor $y'$ of $y$ in $\Gamma_O^{\mathrm{odd}}$.
Now since $y \not\in O_2(s_0)$ and $\Gamma_O^{\mathrm{odd}}$ is acyclic, if $y' \not\in O_2(s_i)$ for some $1 \leq i \leq n$, then $O_2(s_0)_{\sim y'}^{\mathrm{odd}}$ and $O_2(s_i)_{\sim y}^{\mathrm{odd}}$ cannot intersect with each other, so $\left[y',s_0\right] \overset{\infty}{\sim} \left[y,s_i\right]$ by Lemma \ref{lem:relationsoutsideO} (\ref{lem_item:reloutsideO_comp_sisdistinct}).
This contradicts the fact $(y,s_i) \in (\mathcal{E}_y)_{\mathrm{fin}}$, so we have $y' \in O_2(s_i)$ for all $1 \leq i \leq n$.
In particular, the set $K \cup \{y,y'\}$ is finite, irreducible and $2$-spherical, and $\Gamma_{K \cup \{y,y'\}}^{\mathrm{odd}}$ is acyclic; so Corollary \ref{cor:no2sphericalacyclicsupset} implies that $K \cup \{y,y'\}$ is of finite type, therefore $n=2$ and $(y',y,s_0,s_1)$ is of type $F_4$ (see (\ref{eq:ruleforfinitepart}) for terminology).
Note that $(y',s_0),(y',s_1) \in (\mathcal{E}_y)_{\mathrm{fin}}$, since Lemma \ref{lem:irrsubsetofE} implies that the two generators $r_y(y',s_0)$ and $r_y(y',s_1)$ belong to the same finite irreducible component of $W^{\perp y}$.

To verify Condition \ref{thm_item:finpart_nonO_Otrunk_F4_F4exist}, the remaining task is to show that $\Gamma_{O_2(s_0)}^{\mathrm{odd}}$ is connected.
Since $y' \in O_2(s_0)$ and $(y',s_0) \in (\mathcal{E}_y)_{\mathrm{fin}}$ as above, this follows from Lemma \ref{lem:ncsrconditionforfinoutsideO} (\ref{lem_item:ncsrforfinoutsideO_mis2}) applied to $(y',s_0)$.
Thus Condition \ref{thm_item:finpart_nonO_Otrunk_F4_F4exist} is satisfied.

To verify Condition \ref{thm_item:finpart_nonO_Otrunk_F4_E}, the remaining task is to show that $y' \in O_2(s)$ for any $s \in E \smallsetminus \{s_0,s_1\}$.
If this fails, then $O_2(s)_{\sim y}^{\mathrm{odd}}$ cannot intersect with $O_2(s_0)_{\sim y'}^{\mathrm{odd}}$, since $y \not\in O_2(s_0)$.
Thus $\left[y,s\right] \overset{\infty}{\sim} \left[y',s_0\right]$ by Lemma \ref{lem:relationsoutsideO} (\ref{lem_item:reloutsideO_comp_sisdistinct}), contradicting the fact $(y',s_0) \in (\mathcal{E}_y)_{\mathrm{fin}}$.
Hence Condition \ref{thm_item:finpart_nonO_Otrunk_F4_E} is satisfied.

Finally, to verify Condition \ref{thm_item:finpart_nonO_Otrunk_F4_trunk}, the remaining task is to show that $m_{z,z'}$ is odd or $\infty$ for any $z,z' \in O$.
Assume contrary that $m_{z,z'}$ is even.
Since $(y,s_0) \in (\mathcal{E}_y)_{\mathrm{fin}}$ and $\left[y,s_0\right]=\{(y,s_0)\}$, we have $(y,s_0) \overset{m}{\sim} \left[z,z'\right]$ for some $2 \leq m<\infty$, so $(y,s_0) \overset{m}{\sim} (t,t')$ for some $(t,t') \in \left[z,z'\right] \subseteq \mathcal{E}_y^{(O)}$.
Since $m_{y,s_0}=4$, $m_{t,t'}$ is even and $m_{s_0,t'}$ is not odd, this implies that $t=y$ and $m_{y,t'}=m_{s_0,t'}=2$.
Since $\Gamma_{O_2(s_0)}^{\mathrm{odd}}$ is connected and $y',t' \in O_2(s_0)$, we have $(t',s_0) \in \left[y',s_0\right]$ and $(t',s_0) \overset{4}{\sim} (t',y)$, so $(t',y) \in (\mathcal{E}_y)_{\mathrm{fin}} \cap \mathcal{E}_y^{(O)}$ since $(y',s_0) \in (\mathcal{E}_y)_{\mathrm{fin}}$.
However, this contradicts the hypothesis of Theorem \ref{thm:finpart_nonO}.
Hence Condition \ref{thm_item:finpart_nonO_Otrunk_F4_trunk} is satisfied, so Case 4 is concluded.

Thus the proof of Theorem \ref{thm:finpart_nonO} is concluded.
\subsection{Proof of Theorem \ref{thm:finpart_generalcase}}
\label{sec:finitepart_maintheorem_generalcase}

This subsection is devoted to the proof of Theorem \ref{thm:finpart_generalcase}.
By Theorems \ref{thm:finpart_O} and \ref{thm:finpart_nonO}, and the hypothesis of Theorem \ref{thm:finpart_generalcase}, we have $(\mathcal{E}_x)_{\mathrm{fin}} \cap \mathcal{E}_x^{(O)}=\emptyset$ and $m_{y,s}=2$ for any $(y,s) \in (\mathcal{E}_x)_{\mathrm{fin}}$.
First, we show that $G_K$ is isomorphic to $W_K$ for every $K \in \mathcal{K}_1 \cup \mathcal{K}_2$.
Indeed, by definition, all the generators $r_x(y_K,s)$ of $G_K$ are mapped to $r_{y_K}(y_K,s)$, respectively, by taking the conjugation by a common element.
Since $y_K \in O_2(s)$, we have $r_{y_K}(y_K,s)=s$, so this map is the desired isomorphism from $G_K$ to $W_K$.

We divide the remaining proof into two steps.

\noindent
\textbf{Step 1: Every generator of $(W^{\perp x})_{\mathrm{fin}}$ belongs to some $G_K$ with $K \in \mathcal{K}_1 \cup \mathcal{K}_2$.}

Let $r_x(c;y_0,s_0)$ be a generator of $(W^{\perp x})_{\mathrm{fin}}$.
Then $(y_0,s_0) \in (\mathcal{E}_x)_{\mathrm{fin}}$, so $(y_0,s_0) \not\in \mathcal{E}_x^{(O)}$ and $m_{y_0,s_0}=2$ as above.
Now Lemma \ref{lem:ncsrconditionforfinoutsideO} (\ref{lem_item:ncsrforfinoutsideO_mis2}) shows that $\Gamma_{O_2(s_0)}^{\mathrm{odd}}$ is connected and contains all simple closed paths in $\Gamma_O^{\mathrm{odd}}$.
This implies that $(y_0,s_0)$ satisfies the condition (\ref{eq:gammaisindependentofc}), so we have $r_x(c;y_0,s_0)=r_x(y_0,s_0)$.
Let $K$ be the irreducible component of $E$ containing $s_0$.
Now it suffices to show that $K \in \mathcal{K}_1 \cup \mathcal{K}_2$; indeed, if this holds, then Lemma \ref{lem:relationsoutsideO} (\ref{lem_item:reloutsideO_1simcomp}) shows that $(y_K,s_0) \in \left[y_0,s_0\right]$, so $r_x(y_0,s_0)=r_x(y_K,s_0) \in G_K$.
We divide the proof into two cases.

\noindent
\textbf{Case 1-1: $O \neq O_{2,\infty}(s'_0)$ for some $s'_0 \in K$.}

First, we show that $K=\{s_0\}$, so $O \neq O_{2,\infty}(s_0)$.
Assume contrary that $K \neq \{s_0\}$, and take $y \in O_{\mathrm{even}}(s'_0) \smallsetminus O_2(s'_0)$.
Then Lemma \ref{lem:irrsubsetofE} implies that two generators $r_x(y,s'_0)$ and $r_x(y_0,s_0)$ belong to the same irreducible component of $W^{\perp x}$, so $r_x(y,s'_0) \in (W^{\perp x})_{\mathrm{fin}}$ since $r_x(y_0,s_0) \in (W^{\perp x})_{\mathrm{fin}}$, therefore $(y,s'_0) \in (\mathcal{E}_x)_{\mathrm{fin}}$.
However, this contradicts the remark in the first paragraph of the proof.
Thus we have $K=\{s_0\}$.

We show that $K \in \mathcal{K}_2$.
Note that $m_{s_0,t}=2$ for all $t \in E \smallsetminus s_0$, and it will follow that $O_4(s_0) \neq \emptyset$ once we have $O=O_{2,4,\infty}(s_0)$, since $O \neq O_{2,\infty}(s_0)$.
Thus Conditions \ref{thm_item:finpart_general_family2_Oeven} and \ref{thm_item:finpart_general_family2_O2containOeven} in the definition of $\mathcal{K}_2$ follow immediately from Lemma \ref{lem:ncsrconditionforfinoutsideO} (\ref{lem_item:ncsrforfinoutsideO_mis2}) applied to $(y_0,s_0)$.
Now the remaining task is to show that $s_0$ satisfies the condition (\ref{thm_eq:finpart_general_1_trunk_trunk}).
Note that, since $\Gamma_{O_2(s_0)}^{\mathrm{odd}}$ is connected and contains all simple closed paths in $\Gamma_O^{\mathrm{odd}}$ as shown above, the graph $\Gamma_O^{\mathrm{odd}}$ admits a decomposition $\Gamma_{O_2(s_0)}^{\mathrm{odd}} \cup \bigcup_{y \in O_2(s_0)}T_y$ as in Remark \ref{rem:decompintotrees}.

First, assume contrary that $m_{z,z'}$ is even for some $z \in O_2(s_0)$ and $z' \in O \smallsetminus O_2(s_0)=O_{4,\infty}(s_0)$.
Then we have $(z,s_0) \in \left[y_0,s_0\right]$ by Lemma \ref{lem:relationsoutsideO} (\ref{lem_item:reloutsideO_1simcomp}), so $(z,s_0) \in (\mathcal{E}_x)_{\mathrm{fin}}$, therefore $m_{z',s_0}<\infty$ by Lemma \ref{lem:deniedtriangle}, and $z' \in O_4(s_0)$.
Moreover, we have $m_{z,z'}=2$, since otherwise the set $\{z,z',s_0\}$ is irreducible, $2$-spherical and not of finite type, and $\Gamma_{\{z,z',s_0\}}^{\mathrm{odd}}$ is acyclic, contradicting Corollary \ref{cor:no2sphericalacyclicsupset}.
Thus we have $(z,s_0) \overset{4}{\sim} (z,z')$, so $(z,z') \in (\mathcal{E}_x)_{\mathrm{fin}} \cap \mathcal{E}_x^{(O)}$.
This contradicts the remark in the first paragraph of the proof.

Secondly, assume contrary that $m_{z,z'}$ is even for some $z \in T_y \smallsetminus y$ and $z' \in T_{y'} \smallsetminus y'$ such that $y,y' \in O_2(s_0)$ and $y \neq y'$.
Note that $m_{y,z'}$ is not odd, so $m_{y,z'}=\infty$ by the previous paragraph.
Now we have
\[
\left[z,z'\right] \subseteq \{(t,z') \mid t \in T_y \smallsetminus y\} \textrm{ and } \left[y_0,s_0\right] \subseteq \{(t',s_0) \mid t' \in O_2(s_0)\},
\]
so $\left[z,z'\right] \overset{\infty}{\sim} \left[y_0,s_0\right]$ since $(T_y \smallsetminus y) \cap O_2(s_0)=\emptyset$.
This contradicts the fact $(y_0,s_0) \in (\mathcal{E}_x)_{\mathrm{fin}}$.

Finally, assume contrary that $m_{z,z'}$ is even for some $z,z' \in T_y \smallsetminus y$ with $y \in O_2(s_0)$.
Then $m_{y,t}$ is odd or $\infty$ for any $t \in T_y \smallsetminus y$ as shown above, so we have
\[
\left[z,z'\right] \subseteq \{(t,t') \mid t,t' \in T_y \smallsetminus y\}.
\]
Thus we have $\left[z,z'\right] \overset{\infty}{\sim} \left[y_0,s_0\right]$ by a similar argument to the previous paragraph.
This contradicts the fact $(y_0,s_0) \in (\mathcal{E}_x)_{\mathrm{fin}}$.

Hence the condition (\ref{thm_eq:finpart_general_1_trunk_trunk}) for $s_0$ has been verified, concluding Case 1-1.

\noindent
\textbf{Case 1-2: $O=O_{2,\infty}(s)$ for all $s \in K$.}

We show that $K \in \mathcal{K}_1$.
Choose an element $y_s \in O_2(s)$ for each $s \in K$, where $y_{s_0}=y_0$.
First, we show that $|K|<\infty$ and $(y_s,s) \in (\mathcal{E}_x)_{\mathrm{fin}}$ for any $s \in K$.
This is obvious if $K=\{s_0\}$; so suppose that $K \neq \{s_0\}$.
Then, since $(y_0,s_0) \in (\mathcal{E}_x)_{\mathrm{fin}}$, Lemmas \ref{lem:irrsubsetofE} and \ref{lem:relationsoutsideO} (\ref{lem_item:reloutsideO_comp_sisdistinct}) imply that all generators $r_x(y_s,s)$ with $s \in K$ are distinct and belong to the same finite irreducible component of $W^{\perp x}$.
Thus the claim follows.

By applying Lemma \ref{lem:ncsrconditionforfinoutsideO} (\ref{lem_item:ncsrforfinoutsideO_mis2}) to $(y_s,s)$ for each $s \in K$, it follows that $\Gamma_{O_2(s)}^{\mathrm{odd}}$ is connected and contains all simple closed paths in $\Gamma_O^{\mathrm{odd}}$, every connected component of $\Gamma_{O_2(t)}^{\mathrm{odd}}$ intersects $O_2(s)$ for any $t \in E \smallsetminus s$, and $O_{\mathrm{even}}(t) \smallsetminus O_2(t) \subseteq O_2(s)$ for any $t \in E \smallsetminus K$ (since $m_{s,t}=2$ by the choice of $K$).
In particular, Condition \ref{thm_item:finpart_general_family1_O2containOeven} in the definition of $\mathcal{K}_1$ is satisfied, and we have $O_2(s) \cap O_2(t) \neq \emptyset$ for any $s,t \in K$.

We show that Condition \ref{thm_item:finpart_general_family1_allO2intersect} is satisfied.
If $\Gamma_O^{\mathrm{odd}}$ has a simple closed path, then it is contained in every $\Gamma_{O_2(s)}^{\mathrm{odd}}$ with $s \in K$ as shown above, so all the $\Gamma_{O_2(s)}^{\mathrm{odd}}$ have a common vertex.
On the other hand, if $\Gamma_O^{\mathrm{odd}}$ is acyclic, then all $\Gamma_{O_2(s)}^{\mathrm{odd}}$ with $s \in K$ are nonempty subtrees of the tree $\Gamma_O^{\mathrm{odd}}$, any two of which intersect with each other.
Thus all the $\Gamma_{O_2(s)}^{\mathrm{odd}}$ have a common vertex by Lemma \ref{lem:intersectionoftree}, since $|K|<\infty$.
Hence Condition \ref{thm_item:finpart_general_family1_allO2intersect} has been verified.
Moreover, now Lemma \ref{lem:irrsubsetofE} implies that $G_K$ is contained in the finite irreducible component of $W^{\perp x}$ containing $r_x(y_0,s_0)=r_x(y_K,s_0)$.
Since $G_K \simeq W_K$ as above, this means that $K$ is of finite type.

Finally, we verify Condition \ref{thm_item:finpart_general_family1_O2containtrunk}.
The remaining task is to show that any $s \in K$ satisfies (\ref{thm_eq:finpart_general_1_trunk_trunk}).
Now since $\Gamma_{O_2(s)}^{\mathrm{odd}}$ is connected and contains all simple closed paths, and $(y_s,s) \in (\mathcal{E}_x)_{\mathrm{fin}}$, the same argument as that verifying the condition (\ref{thm_eq:finpart_general_1_trunk_trunk}) for $s_0$ in Case 1-1 also implies that any $s \in K$ satisfies (\ref{thm_eq:finpart_general_1_trunk_trunk}).
Thus Condition \ref{thm_item:finpart_general_family1_O2containtrunk} has been verified, so Case 1-2 is concluded.

Hence Step 1 is concluded.

\noindent
\textbf{Step 2: For $K \in \mathcal{K}_1 \cup \mathcal{K}_2$, the subgroup $G_K$ is a finite irreducible component of $W^{\perp x}$.}

Recall that $G_K \simeq W_K$ as above, while $W_K$ is of finite type by definition of $\mathcal{K}_1$ and $\mathcal{K}_2$.
Thus the finiteness of $G_K$ follows.
Now it suffices to show that any generator $r_x(y_K,s)$ with $s \in K$ commutes with all generators of $W^{\perp x}$ other than those of $G_K$.

First, Condition \ref{thm_item:finpart_general_family1_O2containtrunk} in definition of $\mathcal{K}_1$ and Condition \ref{thm_item:finpart_general_family2_O2containtrunk} in definition of $\mathcal{K}_2$ imply that, for any $s \in K$, the graph $\Gamma_{O_2(s)}^{\mathrm{odd}}$ is connected and contains all simple closed paths in $\Gamma_O^{\mathrm{odd}}$.
Thus (\ref{eq:gammaisindependentofc}) is satisfied for any element of $\mathcal{E}_x$ of the form $(y,s)$ with $s \in K$ and $y \in O_2(s)$.
Now the two conditions also imply that, for every $s \in K$, the generator $r_x(y_K,s)$ commutes with any generator $r_x(c;z,z')$ with $z,z' \in O$.
Indeed, the condition (\ref{thm_eq:finpart_general_1_trunk_trunk}) shows that $z,z' \in O_2(s)$, so $(z,z') \overset{2}{\sim} (z,s) \in \left[y_K,s\right]$ by Lemma \ref{lem:relationsoutsideO} (\ref{lem_item:reloutsideO_1simcomp}) since $\Gamma_{O_2(s)}^{\mathrm{odd}}$ is connected, therefore the claim follows from Lemma \ref{lem:ifgammaisindependentofc}.

We show that $r_x(y_K,s)$ commutes with any generator of the form $r_x(c;z,t)$ such that $t \in E \smallsetminus K$ and $z \in O_{\mathrm{even}}(t)$.
Note that $m_{s,t}=2$.
By condition (\ref{thm_eq:finpart_general_1_Oeven_even}), if $z \in O_{\mathrm{even}}(t) \smallsetminus O_2(t)$, then we have $z \in O_2(s)$ and so $(z,t) \overset{2}{\sim} (z,s) \in \left[y_K,s\right]$, since $\Gamma_{O_2(s)}^{\mathrm{odd}}$ is connected.
Thus $r_x(y_K,s)$ commutes with $r_x(c;z,t)$ by Lemma \ref{lem:ifgammaisindependentofc}.
On the other hand, if $z \in O_2(t)$, then we have $O_2(t)_{\sim z}^{\mathrm{odd}} \cap O_2(s) \neq \emptyset$, so $\left[z,t\right] \overset{2}{\sim} \left[y_K,s\right]$ by Lemma \ref{lem:relationsoutsideO} (\ref{lem_item:reloutsideO_comp_sisdistinct}), therefore $r_x(y_K,s)$ commutes with $r_x(c;z,t)$ by Lemma \ref{lem:ifgammaisindependentofc}.

Now if $K \in \mathcal{K}_1$, then by Condition \ref{thm_item:finpart_general_family1_O2containtrunk}, the generators of $W^{\perp x}$ are either those already considered, or those of $G_K$.
Thus $G_K$ is an irreducible component of $W^{\perp x}$, as desired.
Similarly, if $K \in \mathcal{K}_2$, then the generators of $W^{\perp x}$, which are not already considered above and not those of $G_K$, are of the form $r_x(y,s)$ with $y \in O_4(s)$ (see Condition \ref{thm_item:finpart_general_family2_Oeven} in definition of $\mathcal{K}_2$).
Now we have $m_{y,y'}=3$ for some $y' \in O_2(s)$ by Condition \ref{thm_item:finpart_general_family2_Oeven}, so Lemma \ref{lem:relationsoutsideO} (\ref{lem_item:reloutsideO_comp_sissame}) implies that $\left[y,s\right] \overset{2}{\sim} \left[y_K,s\right]$, therefore $r_x(y_K,s)$ commutes with $r_x(y,s)$.
Thus $G_K$ is an irreducible component of $W^{\perp x}$, so Step 2 is concluded.

Hence the proof of Theorem \ref{thm:finpart_generalcase} is concluded.
\section{On Reflection-independent Coxeter groups}
\label{sec:resultsonrefindep}

In a paper \cite{Bah}, Bahls introduced the notion of reflection independence of Coxeter groups.
Namely, a Coxeter group $W$ is \emph{reflection independent} if the set $S^W$ of reflections in $W$ is independent of the choice of the generating set $S$ of $W$.
This condition is equivalent to that $f(S^W)=S'{}^{W'}$ for any Coxeter system $(W',S')$ and any group isomorphism $f:W \overset{\sim}{\to} W'$.
Not all Coxeter groups are reflection independent; it is well known that the symmetric group of degree $6$, which is a Coxeter group of type $A_5$ with the five adjacent transpositions as generators, admits another generating set whose member is not a transposition.
Reflection independence for Coxeter groups in several classes have been examined.

In a preceding paper \cite{Nui_ext}, the author gave a new sufficient condition for a Coxeter group (possibly of infinite rank) to be reflection independent, in terms of the structure of the finite parts $(W^{\perp x})_{\mathrm{fin}}$ of the Coxeter groups $W^{\perp x}$ studied in the previous sections.
Owing to this result, we give some new classes of reflection-independent Coxeter groups.
Let $(W,S)$ denote a Coxeter system throughout this section.

\subsection{Key lemmas}
\label{lem:refindep_keylemma}

This subsection summarizes some preliminary observations.
The following fact is fundamental in a study of reflection independence of Coxeter groups.
\begin{prop}
[{\cite[Lemma 3.7]{Bra-McC-Muh-Neu}}]
\label{prop:relationbetweenreflections}
Let $(W',S')$ be another Coxeter system, and $f:W \overset{\sim}{\to} W'$ a group isomorphism.
If $f(S) \subseteq S'{}^{W'}$, then $f(S^W)=S'{}^{W'}$.

Hence $W$ is reflection independent if and only if, for any group isomorphism $f$ from $W$ to another Coxeter group $W'$, the image $f(x)$ of every generator $x \in S$ is a reflection in $W'$.
\end{prop}
Note that, in the statement, it suffices to check the condition for only one $x$ in every conjugacy class of $S$, namely in every connected component of the odd Coxeter graph $\Gamma^{\mathrm{odd}}$ of $W$ (see Proposition \ref{prop:oddCoxetergraph}).

In preceding papers, several observations based on properties of maximal finite subgroups of Coxeter groups have been given.
The following facts are important in the strategy.
\begin{lem}
[{\cite[Corollary 7]{Fra-How_rank3}}]
\label{lem:intersectionofparabolic}
The intersection of two (hence finitely many) parabolic subgroups of $W$ is also parabolic.
\end{lem}
\begin{lem}
[{\cite[Lemma 8]{Fra-How_rank3}}]
\label{lem:maxfinstandardismaxfin}
If $W_I$ is a maximal finite standard parabolic subgroup of $W$, then it is a maximal finite subgroup of $W$.
\end{lem}
On the other hand, Theorem \ref{thm:Tits_finitesubgroup} yields the following immediately.
\begin{lem}
\label{lem:maxfinisparabolic}
Any maximal finite subgroup of $W$ is a maximal finite parabolic subgroup.
\end{lem}
The following theorem of Roger W.\ Richardson is also required.
\begin{thm}
[{\cite[Theorem A]{Ric}}]
\label{thm:richardson_involution}
Let $W$ be a Coxeter group.
\begin{enumerate}
\item \label{thm_item:ric_inv_exist}
Any involution in $W$ is conjugate to the longest element $w_0(I)$ of a finite standard parabolic subgroup $W_I$ of $W$ such that $w_0(I) \in Z(W_I)$.
\item \label{thm_item:ric_inv_unique}
Let $W_I$ and $W_J$ be two finite standard parabolic subgroups of $W$, and suppose that $w_0(I) \in Z(W_I)$ and $w_0(J) \in Z(W_J)$.
Then $w_0(I)$ and $w_0(J)$ are conjugate in $W$ if and only if $I$ and $J$ are conjugate in $W$.
\end{enumerate}
\end{thm}
Now we obtain the following key fact, which improves {\cite[Lemma 1.6]{Fra}}.
\begin{thm}
\label{thm:oldconditionforrefindep}
Suppose that a subset $I \subseteq S$ consists of at most three mutually commuting generators which are conjugate in $W$, and $W_I$ is the intersection of finitely many maximal finite subgroups of $W$.
Then for any Coxeter system $(W',S')$ and any group isomorphism $f:W \overset{\sim}{\to} W'$, there is an element $w \in W'$ such that $f(I) \subseteq wS'w^{-1}$.
In particular, we have $f(I) \subseteq S'{}^{W'}$.
\end{thm}
\begin{proof}
Put $I=\{x_i \mid 1 \leq i \leq n\}$ with $n=|I| \leq 3$.
Since $W_I$ is the intersection of finitely many maximal finite subgroups, its image $f(W_I)$ has the same property in $W'$.
Lemma \ref{lem:maxfinisparabolic} implies further that $f(W_I)$ is the intersection of finitely many parabolic subgroups of $W'$, so $f(W_I)$ is also parabolic by Lemma \ref{lem:intersectionofparabolic}.
We may assume that $f(W_I)$ is a standard parabolic subgroup $W'_J$ of $W'$.
Moreover, since $W_I$ is an elementary abelian $2$-group with $n$ generators, its image $W'_J$ has the same property, so $J$ also consists of $n$ mutually commuting generators, say $y_i \in S'$ ($1 \leq i \leq n$).

By the hypothesis, all $x_i$ are conjugate, so all $f(x_i)$ are also conjugate.
Let $K_i$ be the subset of $J$ such that $f(x_i)=\prod_{y_j \in K_i}y_j$.
Then $f(x_i)$ is the longest element of $W'_{K_i}$, so Theorem \ref{thm:richardson_involution} (\ref{thm_item:ric_inv_unique}) shows that all $K_i$ are conjugate in $W'$, therefore their cardinality $k$ are equal.
Regarding $W'_J$ as an $n$-dimensional vector space over $\mathbb{F}_2$ with basis $y_1,\dots,y_n$, the property that all the $f(x_i)$ generate $W'_J$ is interpreted as that the corresponding $n$ vectors spans $W'_J$; so we obtain a nonsingular $n \times n$ matrix over $\mathbb{F}_2$ each of whose columns contains precisely $k$ 1's.
Now since $n \leq 3$, a straightforward observation shows that this is possible only if $k=1$, namely $f(x_i) \in J$.
Hence the claim holds.
\end{proof}

On the other hand, the author gave some observations for the reflection independence in a preceding paper \cite{Nui_ext}, by a different approach.
Here we summarize them with slight modification.
We start with a wider setting; let $(W,S)$ and $(W',S')$ be Coxeter systems, $f:W \overset{\sim}{\to} W'$ a group isomorphism, and $I$ a subset of $S$ of finite type such that $w_0(I) \in Z(W_I)$ and the Coxeter graph $\Gamma_I$ of $W_I$ admits no nontrivial automorphisms.
Let $W^{\perp I}$ be the subgroup of $W$ generated by the reflections which fix the subset $\Pi_I \subseteq \Pi$ pointwise; so $W^{\perp I}=W^{\perp x}$ in the above notations whenever $I=\{x\}$.
A general theorem of Deodhar \cite{Deo_refsub} or Dyer \cite{Dye} proves that $W^{\perp I}$ is a Coxeter group, so it admits the finite part $(W^{\perp I})_{\mathrm{fin}}$ (see Definition \ref{defn:finitepart}).
Now $f(w_0(I))$ is an involution as well as $w_0(I)$, so Theorem \ref{thm:richardson_involution} (\ref{thm_item:ric_inv_exist}) shows that it is conjugate to the longest element $w_0(J)$ of some finite standard parabolic subgroup $W'_J$ of $W'$ such that $w_0(J) \in Z(W'_J)$.
The subgroup $W'{}^{\perp J}$ of $W'$, which is also a Coxeter group, is similarly defined.

In the setting, the argument in {\cite[Section 3.3]{Nui_ext}} yields the following:
\begin{prop}
[{\cite[Equation (3.6)]{Nui_ext}}]
\label{prop:Nui_ext_correspondence}
There is a finite subgroup $G_{\rho'}$ of $W'$ such that the image of the subgroup $W_I \times (W^{\perp I})_{\mathrm{fin}}$ of $W$ under $f$ is conjugate to $W'_J \rtimes ((W'{}^{\perp J})_{\mathrm{fin}} \rtimes G_{\rho'})$.
\end{prop}
The result in another papar \cite{Nui_centra} of the author describes the structure of the group $W^{\perp I}$.
This proposition suggests that, even if $|I| \neq 1$, studying semidirect product decompositions of the Coxeter groups $W_I \times (W^{\perp I})_{\mathrm{fin}}$ would yield some information on the set $J$, namely on the image $f(w_0(I))$.

In particular, for the case that $|I|=1$ in which we are interested, Proposition \ref{prop:Nui_ext_correspondence} yields the following observation.
\begin{thm}
[{\cite[Theorem 3.7]{Nui_ext}}]
\label{thm:Nui_ext_scforrefindep}
Let $x \in S$, and suppose that the finite part $(W^{\perp x})_{\mathrm{fin}}$ of $W^{\perp x}$ is either trivial, or generated by a single element which is conjugate to $x$ in $W$.
Then $f(x)$ is a reflection in $W'$ for any Coxeter system $(W',S')$ and any group isomorphism $f:W \overset{\sim}{\to} W'$.
\end{thm}
Note that by Proposition \ref{prop:relationbetweenreflections}, a Coxeter group $W$ is reflection independent if every $x \in S$ satisfies the condition in this theorem.
Since the structure of $(W^{\perp x})_{\mathrm{fin}}$ has been completely determined in the previous sections, this theorem actually gives rise to a sufficient condition for $W$ to be reflection independent, which is efficiently verifiable.
Moreover, in some sense, our description of $(W^{\perp x})_{\mathrm{fin}}$ would suggest that the group $(W^{\perp x})_{\mathrm{fin}}$ is trivial in a `generic' case, so $W$ is reflection independent in that case.
\subsection{Some reflection-independent Coxeter groups}
\label{sec:someclassesofrefindep}

The aim of this subsection is to examine the reflection independence of Coxeter groups in certain classes.
Some further properties of isomorphisms between Coxeter groups in those classes will be investigated.
\subsubsection{$2$-spherical Coxeter groups}
\label{sec:refindep_2spherical}

Recall that $W$ is called $2$-spherical if $m_{s,t}<\infty$ for any $s,t \in S$.
Then the description of the groups $(W^{\perp x})_{\mathrm{fin}}$ given in Sections \ref{sec:mainresult_specialcase} and \ref{sec:finitepart} yield the following, generalizing Theorem \ref{thm:specialcase_2spherical_acyclic}.
\begin{prop}
\label{prop:Wperpxfinfor2spherical}
Let $W$ be an infinite, irreducible, $2$-spherical Coxeter group, and $x \in S$.
Then we have $(W^{\perp x})_{\mathrm{fin}} \neq 1$ if and only if one of the conditions \ref{thm_item:specialcase_2sp_ac_B} and \ref{thm_item:specialcase_2sp_ac_D} in Theorem \ref{thm:specialcase_2spherical_acyclic} is satisfied.
Moreover, in this case, the group $(W^{\perp x})_{\mathrm{fin}}$ is generated by a single element conjugate to $x$.
\end{prop}
\begin{proof}
The `if' part and the latter claim are deduced by Theorem \ref{thm:specialcase_2spherical_acyclic} and Proposition \ref{prop:oddCoxetergraph}.
For the ``only if'' part, put $O$ and $E$ as in Section \ref{sec:finitepart}.
Note that $O \neq \emptyset$ and $S=O \cup E$ since $W$ is $2$-spherical.
First, we give some preliminary observations.
\begin{lem}
\label{prop_lem:Wperpx2sp_facts}
In this case, we have the followings:
\begin{enumerate}
\item \label{prop_lem_item:Wperpx2sp_facts_compofE}
There is no irreducible component $K$ of $E$ such that $O=O_{2,\infty}(s)$ for all $s \in K$.
\item \label{prop_lem_item:Wperpx2sp_facts_trunk}
If $K \subseteq O$ satisfies (\ref{eq:Kistrunk}), and $\Gamma_K^{\mathrm{odd}}$ is nonempty, connected and contains all simple closed paths in $\Gamma_O^{\mathrm{odd}}$, then either $K=O$, or $|K|=1$ and $|O|=2$.
\item \label{prop_lem_item:Wperpx2sp_facts_Oissmall}
If $\Gamma_O^{\mathrm{odd}}$ is acyclic, and $m_{y,z}$ is odd or $\infty$ for any $y,z \in O$, then $|O| \leq 2$.
\end{enumerate}
\end{lem}
\begin{proof}
[Proof of Lemma \ref{prop_lem:Wperpx2sp_facts}.]
\textbf{(\ref{prop_lem_item:Wperpx2sp_facts_compofE})}
Since $W$ is $2$-spherical, such $K$ satisfies that $O=O_2(s)$ for all $s \in K$, so $K$ is an irreducible component of $S=O \cup E$, contradicting the irreducibleness of $S$.

\noindent
\textbf{(\ref{prop_lem_item:Wperpx2sp_facts_trunk})}
The condition (\ref{eq:Kistrunk}) now implies that $m_{y,z}$ is odd for any $y,z \in O$ with $\{y,z\} \not\subseteq K$.
Now since $\Gamma_K^{\mathrm{odd}}$ is connected, if $K \neq O$ and $|K| \geq 2$, then two adjacent vertices of $\Gamma_K^{\mathrm{odd}}$ and a vertex of $O \smallsetminus K$ form a simple closed path, contradicting the assumption on $K$.
Similarly, if $|O \smallsetminus K| \geq 2$, then two vertices of $O \smallsetminus K$ and a vertex of $K$ form a simple closed path, a contradiction.
Thus the claim follows.

\noindent
\textbf{(\ref{prop_lem_item:Wperpx2sp_facts_Oissmall})}
This follows from Claim \ref{prop_lem_item:Wperpx2sp_facts_trunk}, since now $K=\{x\} \subseteq O$ satisfies the hypothesis of Claim \ref{prop_lem_item:Wperpx2sp_facts_trunk}.
\end{proof}
Now we come back to the proof of the proposition.
In this case, the hypothesis of Theorem \ref{thm:finpart_O_typeB} cannot be satisfied.
Indeed, if this is satisfied, then Lemma \ref{prop_lem:Wperpx2sp_facts} (\ref{prop_lem_item:Wperpx2sp_facts_compofE}) and Conditions \ref{thm_item:finpart_O_typeB_s0} and \ref{thm_item:finpart_O_typeB_nots0} imply that $E=\{s_0\}$, while Lemma \ref{prop_lem:Wperpx2sp_facts} (\ref{prop_lem_item:Wperpx2sp_facts_trunk}) and Condition \ref{thm_item:finpart_O_typeB_trunk} show that $O=K$ since $|K| \geq 3$.
Thus $S=O \cup E$ is of type $B_{n+1}$, contradicting the infiniteness of $W$.

Similarly, in the case of Theorem \ref{thm:finpart_O_trunk_acyclic}, Lemma \ref{prop_lem:Wperpx2sp_facts} implies that $K=O$ and $E=\emptyset$, while $K$ now cannot be of type $P(m)$, therefore $S=K$ is of finite type, a contradiction.
In the case of Theorem \ref{thm:finpart_nonO_Otrunk_F4}, Lemma \ref{prop_lem:Wperpx2sp_facts} implies that $|O| \leq 2$, so $O=\{y,y'\}$, and $E=\{s_0,s_1\}$, therefore $S=J$ is of type $F_4$, a contradiction.
In the case of Theorem \ref{thm:finpart_nonO_Otrunk_B}, we have $O=O_\infty(s_0) \cup y=\{y\}$ by Condition \ref{thm_item:finpart_nonO_Otrunk_B_Bexist}, and $E=\{s_0,\dots,s_n\}$ by Lemma \ref{prop_lem:Wperpx2sp_facts}, so $S=J$ is of type $B_{n+2}$ or $I_2(m)$, a contradiction.

Suppose that the hypothesis of Theorem \ref{thm:finpart_generalcase} is satisfied.
Now Lemma \ref{prop_lem:Wperpx2sp_facts} and Condition \ref{thm_item:finpart_general_family1_O2containtrunk} in definition of $\mathcal{K}_1$ imply that $\mathcal{K}_1=\emptyset$.
We show that $\mathcal{K}_2=\emptyset$.
If $\{s\} \in \mathcal{K}_2$, then $\emptyset \neq O_2(s) \neq O$ by Condition \ref{thm_item:finpart_general_family2_Oeven}, so $|O_2(s)|=1$ and $|O|=2$ by Condition \ref{thm_item:finpart_general_family2_O2containtrunk} and Lemma \ref{prop_lem:Wperpx2sp_facts} (\ref{prop_lem_item:Wperpx2sp_facts_trunk}) applied to $O_2(s) \subseteq O$.
Put $O_2(s)=\{y\}$ and $O=\{y,y'\}$, so $y' \in O_4(s)$ and $m_{y,y'}=3$ by Condition \ref{thm_item:finpart_general_family2_Oeven}.
Now for any $t \in E \smallsetminus s$, Condition \ref{thm_item:finpart_general_family2_O2containOeven} implies that $O_{\mathrm{even}}(t) \smallsetminus O_2(t)=\emptyset$; otherwise we have $O_{\mathrm{even}}(t) \smallsetminus O_2(t)=\{y\}$, so $O_2(t)=\{y'\}$ since $O=O_{\mathrm{even}}(t)$, therefore $\Gamma_{O_2(t)}^{\mathrm{odd}}$ does not intersect $O_2(t)$, a contradiction.
Thus $E=\{s\}$ by Lemma \ref{prop_lem:Wperpx2sp_facts} (\ref{prop_lem_item:Wperpx2sp_facts_compofE}), so $S=\{y,y',s\}$ which is of type $B_3$ as above, a contradiction.
Hence we have shown that $\mathcal{K}_1 \cup \mathcal{K}_2=\emptyset$, therefore $(W^{\perp x})_{\mathrm{fin}}=1$.

The remaining cases are now two cases; Theorems \ref{thm:finpart_O_trunk_cyclic} and \ref{thm:finpart_nonO_Onontrunk}.
In the case of Theorem \ref{thm:finpart_O_trunk_cyclic}, we have $E=\emptyset$ by Lemma \ref{prop_lem:Wperpx2sp_facts} (\ref{prop_lem_item:Wperpx2sp_facts_compofE}).
Now Condition \ref{thm_item:fin_oddconn_cycl_bipyramid} in Theorem \ref{thm:finitepart_oddconnected_cyclic} cannot be satisfied; so this is the case of Theorem \ref{thm:finitepart_oddconnected_acyclic_trankisnonfinitetype}.
We have $T_{x_1}=\{x_1\}$, since otherwise Condition \ref{thm_item:fin_oddconn_acycl_nonftype_almostallinfty} in Theorem \ref{thm:finitepart_oddconnected_acyclic_trankisnonfinitetype} implies that $m_{z,x_3}=\infty$ for any $z \in T_{x_1} \smallsetminus x_1$, a contradiction.
We also have $T_{x_3}=\{x_3\}$ by symmetry.
Condition \ref{thm_item:fin_oddconn_acycl_nonftype_almostallinfty} also implies that any two elements of the poset $T_{x_2}$ are comparable; so $T_{x_2}$ is a (possibly infinite) chain $x_2=y_0 \prec_{x_2} y_1 \prec_{x_2} \cdots$.
It also follows that $T_{x_2}^o=T_{x_2}$, $m_{y_0,y_1}=3$ by Condition \ref{thm_item:fin_oddconn_acycl_nonftype_orderideal}, and $m_{y_i,y_j}=2$ for any $i,j \geq 0$ with $i \leq j-2$; so Condition \ref{thm_item:fin_oddconn_acycl_nonftype_Tx2} implies that $m_{y_{i-1},y_i}$ for $i \geq 1$.
Thus we have shown that $S=O=K'$ is of type $A_3$, $D_n$ with $4 \leq n<\infty$, or $D_\infty$; the former two cases are denied by infiniteness of $W$.
Hence now Condition \ref{thm_item:specialcase_2sp_ac_D} in Theorem \ref{thm:specialcase_2spherical_acyclic} is satisfied.

Finally, suppose that the hypothesis of Theorem \ref{thm:finpart_nonO_Onontrunk} is satisfied.
Lemma \ref{prop_lem:Wperpx2sp_facts} (\ref{prop_lem_item:Wperpx2sp_facts_compofE}) implies that $E=\{s_0\}$ by Conditions \ref{thm_item:finpart_nonO_Onontrunk_prev_B2exist} and \ref{thm_item:finpart_nonO_Onontrunk_prev_E}.
By Conditions \ref{thm_item:finpart_nonO_Onontrunk_prev_B2exist} and \ref{thm_item:finpart_nonO_Onontrunk_prev_Tistrunk}, we have $O=O_2(s_0) \cup y=T$ and any two elements of $T$ are comparable, so $O$ is a chain $y=y_0 \prec y_1 \prec \cdots$.
Now as well as the previous paragraph, Conditions \ref{thm_item:finpart_nonO_Onontrunk_prev_atomofT} and \ref{thm_item:finpart_nonO_Onontrunk_prev_orderideal} imply that $O$ is of type $A_n$ with $n<\infty$ or type $A_\infty$.
Thus by Condition \ref{thm_item:finpart_nonO_Onontrunk_prev_B2exist}, $S=O \cup s_0$ is of type $B_{n+1}$ or $B_\infty$; the former possibility is denied since $|W|=\infty$.
Hence now Condition \ref{thm_item:specialcase_2sp_ac_B} in Theorem \ref{thm:specialcase_2spherical_acyclic} is satisfied, concluding the proof.
\end{proof}
Hence by Theorem \ref{thm:Nui_ext_scforrefindep} and Proposition \ref{prop:relationbetweenreflections}, we have the following result.
\begin{thm}
\label{thm:refindep_2spherical}
Let $W$ be an infinite, irreducible and $2$-spherical Coxeter group, possibly of infinite rank.
Then $W$ is reflection independent.
\end{thm}
\subsubsection{Odd-connected Coxeter groups}
\label{sec:refindep_oddconnected}

In this paper, we say that a Coxeter group $W$ is \emph{odd-connected} if the odd Coxeter graph $\Gamma^{\mathrm{odd}}$ of $W$ is connected.
Owing to Proposition \ref{prop:oddCoxetergraph}, this condition is equivalent to that all generators $x \in S$, hence all reflections in $W$, are conjugate with each other.
Moreover, an observation shows further that, this is also equivalent to that $W$ possesses a unique subgroup of index two.
Thus it follows that the odd-connectedness of $W$ is actually independent of the choice of the set $S$.
For such a Coxeter group $W$, we have the following results.
Note that $O=S_{\sim x}^{\mathrm{odd}}=S$ in this case, so the result in Section \ref{sec:specialcase_oddconnected} now works.
\begin{prop}
\label{prop:Wperpxfinforoddconnected}
Let $W$ be an infinite odd-connected Coxeter group, and $x \in S$.
Suppose that $(W^{\perp x})_{\mathrm{fin}}$ is neither trivial nor generated by a single reflection.
Then the hypothesis of Theorem \ref{thm:finitepart_oddconnected_acyclic_trankisfinitetype} is satisfied, with $O=S$.
Moreover, we have $K \neq S$ unless $K$ is of type $P(m)$.
\end{prop}
\begin{proof}
If $(W^{\perp x})_{\mathrm{fin}} \neq 1$, then by Theorems \ref{thm:finitepart_oddconnected_cyclic} and \ref{thm:finitepart_oddconnected_acyclic}, the condition in one of Theorems \ref{thm:finitepart_oddconnected_cyclic}, \ref{thm:finitepart_oddconnected_acyclic_trankisfinitetype} and \ref{thm:finitepart_oddconnected_acyclic_trankisnonfinitetype} is satisfied.
Now in the case of Theorem \ref{thm:finitepart_oddconnected_cyclic} or \ref{thm:finitepart_oddconnected_acyclic_trankisnonfinitetype}, it follows that $(W^{\perp x})_{\mathrm{fin}}$ is generated by a single generator, which is a reflection.
Thus the hypothesis of Theorem \ref{thm:finitepart_oddconnected_acyclic_trankisfinitetype} is now satisfied.
Now the last claim in the statement follows from the infiniteness of $W$.
\end{proof}
\begin{thm}
\label{thm:refindep_oddconnected}
Let $W$ be an infinite odd-connected Coxeter group, possibly of infinite rank.
Then $W$ is reflection independent.
\end{thm}
\begin{proof}
Let $x \in S$.
Note that now all reflections in $W$ are conjugate; so by Proposition \ref{prop:relationbetweenreflections}, it suffices to show that $f(x) \in S'{}^{W'}$ for any Coxeter system $(W',S')$ and any group isomorphism $f:W \overset{\sim}{\to} W'$.
By Theorem \ref{thm:Nui_ext_scforrefindep} and Proposition \ref{prop:Wperpxfinforoddconnected}, it suffices to consider the case that the hypothesis of Theorem \ref{thm:finitepart_oddconnected_acyclic_trankisfinitetype} is satisfied, and $K \neq S$ unless $K$ is of type $P(m)$.
Now the condition (\ref{eq:Kistrunk}) implies that
\begin{equation}
\label{thm_eq:refindep_oddconn_infty}
\text{if } s \in S \smallsetminus K, \textrm{ then } m_{s,t}<\infty \textrm{ for at most one } t \in K.
\end{equation}
Indeed, since $\Gamma^{\mathrm{odd}}$ is acyclic and $\Gamma_K^{\mathrm{odd}}$ is connected, this $s$ is adjacent in $\Gamma^{\mathrm{odd}}$ to at most one vertex in $K$.

First, suppose that $K$ is not of type $P(m)$, and take some $y \in K$ and $z \in S \smallsetminus K$ which are adjacent in the connected graph $\Gamma^{\mathrm{odd}}$.
Put $I=\{y,z\}$.
We show that $W_K$ and $W_I$ are maximal finite standard parabolic subgroups of $W$.
By (\ref{thm_eq:refindep_oddconn_infty}), this holds for $W_K$ since $|K| \geq 2$.
Let $s \in S \smallsetminus I$.
Then (\ref{thm_eq:refindep_oddconn_infty}) implies that $m_{s,z}=\infty$ if $s \in K$.
Moreover, if $s \not\in K$, then at most one of $m_{s,z}$ and $m_{s,t}$ is odd, since $\Gamma^{\mathrm{odd}}$ is acyclic.
Thus (\ref{eq:Kistrunk}) implies that one of them is $\infty$.
Hence the claim holds for $W_I$.
Now Lemma \ref{lem:maxfinstandardismaxfin} shows that $W_K$ and $W_I$ are maximal finite subgroups of $W$, so the set $I \cap K=\{y\}$ satisfies the hypothesis of Theorem \ref{thm:oldconditionforrefindep}.
Thus we have $f(y) \in S'{}^{W'}$, so $f(x) \in S'{}^{W'}$ since $x$ is conjugate to $y$, as desired.

Secondly, suppose that $K$ is of type $P(m)$.
Write $K=\{x_1,x_2,x_3,x_4\}$ such that $m_{x_1,x_2}=m$, $m_{x_2,x_3}=m_{x_3,x_4}=3$, $m_{x_1,x_3}=\infty$ and $m_{x_1,x_4}=m_{x_2,x_4}=2$.
Put $I=K \smallsetminus x_1$ and $J=K \smallsetminus x_3$.
We show that $W_I$ and $W_J$ are maximal finite standard parabolic subgroups of $W$.
Now if $s \in S \smallsetminus K$, then (\ref{thm_eq:refindep_oddconn_infty}) implies that $W_{I \cup s}$ and $W_{J \cup s}$ are infinite, since $|I| \geq 2$ and $|J| \geq 2$.
Thus the claim holds, since $|K \smallsetminus I|=|K \smallsetminus J|=1$ and $|W_K|=\infty$.
Now $I \cap J$ consists of two commuting generators $x_2$ and $x_4$ which are conjugate, so Lemma \ref{lem:maxfinstandardismaxfin} and Theorem \ref{thm:oldconditionforrefindep} prove that $f(I \cap J) \subseteq S'{}^{W'}$, therefore $f(x) \in S'{}^{W'}$ similarly.
Hence the proof is concluded.
\end{proof}

\textbf{Koji Nuida}\\
Graduate School of Mathematical Sciences, University of Tokyo\\
Supported by JSPS Research Fellowship (No.\ 16-10825)\\
E-mail: \textit{nuida@ms.u-tokyo.ac.jp}
\end{document}